\documentclass[12pt]{amsart}
\usepackage{mathpazo} 
\usepackage{avant}           
\usepackage{tensor}
\usepackage{todonotes}
\usepackage[backend=biber, style=alphabetic]{biblatex}

\usepackage{graphicx}
\usepackage{amssymb}
\makeindex
\usepackage{tikz-cd}

\newtheorem{thm}{Theorem}[section]
\newtheorem{cor}[thm]{Corollary}

\newtheorem{theorem}{Theorem}[section]

\newtheorem{lem}[thm]{Lemma}
\newtheorem{prop}[thm]{Proposition}
\theoremstyle{definition}
\newtheorem{example}[thm]{Example}
\newtheorem{defn}[thm]{Definition}
\theoremstyle{remark}
\newtheorem{rem}[thm]{Remark}
\numberwithin{equation}{section}

\newcommand{\del}{\partial}

\renewcommand{\del}{\partial}
\newcommand{\delb}{\overline{\partial}}

\newcommand{\CM}{{\mathcal M}}

\newcommand{\CC}{{\mathsf C}}
\newcommand{\CD}{{\mathsf D}}
\newcommand{\CA}{{\mathcal A}}
\newcommand{\CG}{{\mathcal G}}
\newcommand{\CP}{{\mathbb P}}
\newcommand{\CF}{{\mathcal F}}

\newcommand{\BE}{{\mathbb E}}
\newcommand{\CN}{{\mathcal N}}
\newcommand{\cm}{{\mathcal{M}}}
\newcommand{\BR}{{\mathbf R}}
\newcommand{\CW}{{\mathcal{W}}}

\newcommand{\CS}{{\mathcal S}}
\newcommand{\CO}{{\mathcal O}}
\newcommand{\D}{{\mathcal D}}

\newcommand{\cinf}{{\mathcal C}^{\infty}}

\newcommand{\A}{{A}}

\newcommand{\Z}{{\mathbb{Z}}}

\newcommand{\R}{{\mathbb{R}}}
\newcommand{\C}{{\mathbb{C}}}
\newcommand{\N}{{\mathbb{N}}}
\newcommand{\1}{{\mathbb{1}}}
\newcommand{\E}{{\mathcal{E}}}
\newcommand{\EE}{{\mathbb{E}}}

\newcommand{\e}{\mathcal{E}}
\newcommand{\y}{\mathsf{y}}
\newcommand{\p}{\mathsf{p}}

\newcommand{\f}{\mathcal{F}}
\newcommand{\po}{\mathcal{P}}
\newcommand{\g}{\mathfrak{g}}
\newcommand{\gl}{\mathfrak{gl}}

\newcommand{\h}{\mathfrak{h}}

\newcommand{\Diff}{\operatorname{Diff}}
\newcommand{\diff}{\operatorname{diff}}

\newcommand{\ob}{\operatorname{Obj}}
\newcommand{\op}{\operatorname{op}}
\newcommand{\CE}{\operatorname{CE}}
\newcommand{\At}{\operatorname{At}}
\newcommand{\GL}{\operatorname{GL}}
\newcommand{\SL}{\operatorname{SL}}
\newcommand{\inv}{\operatorname{Inv}}
\newcommand{\gr}{\operatorname{gr}}

\newcommand{\set}{\mathsf{Set}}
\newcommand{\dm}{\mathsf{dM}}
\newcommand{\dmfd}{\mathsf{dMfd}}
\newcommand{\dan}{\mathsf{dAnSp}}
\newcommand{\ban}{\mathsf{Ban}}
\newcommand{\dban}{\mathsf{dBan}}
\newcommand{\mfd}{\mathsf{Mfd}}

\newcommand{\sch}{\mathsf{Sch}}
\newcommand{\Mod}{\mathsf{Mod}}
\newcommand{\Modi}{\mathtt{Mod}}
\newcommand{\ho}{\mathsf{Ho}}
\newcommand{\cohmod}{\mathsf{Mod}^{Coh}}
\newcommand{\comgrp}{\mathsf{ComGrp}}
\newcommand{\common}{\mathsf{ComMon}}
\newcommand{\cohmodi}{\mathtt{Mod}^{Coh}}
\newcommand{\perf}{\mathsf{Perf}}

\newcommand{\dcoh}{\mathsf{D^b_{Coh}}}
\newcommand{\dperf}{\mathsf{D_{Perf}}}
\newcommand{\dr}{\mathsf{dR}}
\newcommand{\sing}{\operatorname{Sing}}

\newcommand{\gp}{\mathsf{Grp}}
\newcommand{\gpd}{\mathsf{Grpd}}
\newcommand{\aff}{\mathsf{Aff}}\newcommand{\alg}{\mathsf{Alg}}
\newcommand{\algd}{\mathsf{Algd}}
\newcommand{\ring}{\mathsf{Ring}}
\newcommand{\cring}{\mathsf{CRing}}
\newcommand{\calg}{\mathsf{CAlg}}
\newcommand{\dgcalg}{\mathsf{dgCAlg}}
\newcommand{\dcalg}{\mathsf{dg^{\ge 0}CAlg}}
\newcommand{\fol}{{\operatorname{Fol}}}
\newcommand{\sfol}{{\operatorname{SFol}}}
\newcommand{\Flat}{{\operatorname{flat}}}
\newcommand{\sh}{\mathsf{Sh}}
\newcommand{\psh}{\mathsf{PSh}}
\newcommand{\cov}{\mathsf{cov}}
\newcommand{\cstar}{\mathsf{C^*Alg}}
\newcommand{\Top}{\mathsf{Top}}
\newcommand{\Ch}{\mathsf{Ch}}
\newcommand{\dgcat}{\mathsf{dgCat}}
\newcommand{\rep}{\mathsf{Rep}}
\newcommand{\Loc}{\mathsf{Loc}}
\newcommand{\QCoh}{\mathsf{QCoh}}
\newcommand{\Lx}{{\mathbb{L}}}
\newcommand{\M}{{\mathfrak{M}}}
\newcommand{\Sp}{{\mathcal{S}}}

\newcommand{\cat}{\mathsf{Cat}}
\newcommand{\cocycle}{\mathsf{Cocycle}}
\newcommand{\wcocycle}{\mathsf{wCocycle}}
\newcommand{\repi}{\mathtt{Rep}}
\newcommand{\loci}{\mathtt{Loc}}
\newcommand{\mfdi}{\mathtt{Mfd}}
\newcommand{\RH}{\operatorname{RH}}
\newcommand{\n}{\operatorname{N}}

\newcommand{\ku}{\mathsf{ku}}

\newcommand{\holim}{\operatorname{Holim}}
\newcommand{\hocolim}{\operatorname{Hocolim}}
\newcommand{\igpd}{\mathsf{Grpd_{\infty}}}
\newcommand{\ilgpd}{\mathsf{Lie_{\infty}Grpd}}

\newcommand{\lalgd}{\mathsf{L_{\infty}Algd}^{\operatorname{dg}}}
\newcommand{\lalgdi}{\mathtt{L_{\infty}Algd}}
\newcommand{\open}{\operatorname{Open}}
\newcommand{\Char}{\operatorname{Char}}
\newcommand{\shi}{\mathsf{Sh}_{\infty}}
\newcommand{\pshi}{\mathsf{PSh_{\infty}}}
\newcommand{\vect}{\mathsf{Vect}}
\newcommand{\nerv}{\mathsf{N}}
\newcommand{\bt}{\bullet}
\newcommand{\io}{\operatorname{(\infty, 1)}}
\newcommand{\algi}{\mathtt{Alg}}

\newcommand{\free}{{\operatorname{Free}}}
\newcommand{\dg}{{\operatorname{dg}}}
\newcommand{\coh}{{\operatorname{coh}}}
\newcommand{\qcoh}{{\operatorname{qcoh}}}
\newcommand{\st}{{\operatorname{St}}}
\newcommand{\dgmod}{\mathsf{Mod}^{\operatorname{dg}}}
\newcommand{\cialg}{\mathcal{C}^{\infty}\mathsf{Alg}}
\newcommand{\ci}{\mathcal{C}^{\infty}}

\newcommand{\Hom}{\operatorname{Hom}}
\newcommand{\HOM}{\underline{\operatorname{Hom}}}
\newcommand{\rhom}{\operatorname{\mathbf{R}Hom}}
\newcommand{\coker}{{\operatorname{Coker}}}

\newcommand{\map}{\operatorname{Map}}
\newcommand{\kan}{\operatorname{Kan}}
\newcommand{\End}{\operatorname{End}}

\newcommand{\ext}{\operatorname{Ext}}
\newcommand{\fun}{\operatorname{Fun}}

\newcommand{\funi}{\operatorname{Fun_{\infty}}}

\newcommand{\id}{\operatorname{Id}}
\newcommand{\tr}{\operatorname{Tr}}
\newcommand{\im}{\operatorname{Im}}

\newcommand{\ann}{\operatorname{Ann}}
\newcommand{\codim}{\operatorname{codim}}
\newcommand{\LLP}{\operatorname{LLP}}

\newcommand{\Span}{\operatorname{Span}}
\newcommand{\spec}{\operatorname{Spec}}
\newcommand{\germ}{\operatorname{Germ}}
\newcommand{\aut}{\operatorname{Aut}}

\newcommand{\der}{\operatorname{Der}}
\newcommand{\mon}{\operatorname{Mon}}
\newcommand{\hol}{\operatorname{Hol}}
\newcommand{\holtrans}{\operatorname{HolTrans}}
\newcommand{\cech}{\operatorname{\check{C}ech}}
\newcommand{\rk}{\operatorname{rk}}
\newcommand{\unsh}{\operatorname{UnSh}}
\newcommand{\ider}{\operatorname{InnDer}}

\newcommand{\its}{\textit}

\newcommand{\ev}{\operatorname{ev}}
\newcommand{\sym}{\operatorname{Sym}}
\newcommand{\ber}{\operatorname{Ber}}
\newcommand{\td}{\operatorname{Td}}

\newcommand{\pr}{{\operatorname{pr}}}
\newcommand{\li}{L_{\infty}}

\newcommand{\sep}{{\operatorname{sep}}}

\newcommand{\str}{\operatorname{Str}}
\newcommand{\pont}{\operatorname{Pont}}
\newcommand{\ch}{\operatorname{ch}}

\makeatletter
\newcommand{\colim@}[2]{%
	\vtop{\m@th\ialign{##\cr
			\hfil$#1\operator@font colim$\hfil\cr
			\noalign{\nointerlineskip\kern1.5\ex@}#2\cr
			\noalign{\nointerlineskip\kern-\ex@}\cr}}%
}
\newcommand{\colim}{%
	\mathop{\mathpalette\colim@{\rightarrowfill@\textstyle}}\nmlimits@
}

\newcommand{\extp}{\@ifnextchar^\@extp{\@extp^{\,}}}
\def\@extp^#1{\mathop{\bigwedge\nolimits^{\!#1}}}

\makeatother

\begin{document}

\title{Derived Lie $\infty$-groupoids and algebroids in higher differential geometry}

\author{Qingyun Zeng}
\address{Department of Mathematics,
	University of Pennsylvania, Philadelphia, PA 19104}
\curraddr{Department of Mathematics, University of Pennsylvania, Philadelphia, PA 19104}
\email{qze@math.upenn.edu}
\thanks{The first author was supported in part by NSF Grant \#000000.}


\subjclass[2000]{Primary 54C40, 14E20; Secondary 46E25, 20C20}

\date{January 1, 2001 and, in revised form, June 22, 2001.}


\keywords{Differential geometry, algebraic geometry}

\begin{abstract}
We study various problems arising in higher geometry using derived Lie $\infty$-groupoids and algebroids. We construct homotopical algebras for derived Lie $\infty$-groupoids and algebroids and study their homotopy-coherent representations. Then we apply these tools in studying singular foliations and their characteristic classes. Finally, we prove an $A_{\infty}$ de Rham theorem and higher Riemann-Hilbert correspondence for foliated manifolds.
\end{abstract}

\maketitle
\tableofcontents	
\part{Introduction}
\section{Introduction}
This thesis studies higher categorical and homotopical methods in differential geometry, which is also called {\it higher differential geometry}\index{higher differential geometry}. Recall that in traditional differential geometry, the geometric object we study are usually differentiable manifolds with some additional structures, like complex structures, symplectic structures, Calabi-Yau structures etc. Traditional manifolds theory does not permit singularities. Though there are some tools like stratifies spaces, orbifolds etc. which allow us to study singular manifolds, these tools were built to solve special problems rather than general ones. On the other hand, higher category theory and homotopy theory have been grown rapidly, and are adapted to algebraic geometry and algebraic topology. Hence, higher differential geometry is an adaption of (higher) categorical and homotopical methods in traditional differential geometry.

In order to motivate the necessity of higher categories and homotopy theory, let us first look at the following problems arising in algebraic and differential geometry.

\begin{example}[Moduli problems]
    Let $\M_g$ be the moduli space of curves of genus $g$, that is, a functor sending $\spec(A)$ to the classes of  curves over $\spec(A)$ for some $A\in \calg$. The differential geometric analogue of $\M_g$ is that for each base space $S$ we have a category such that its objects are fiber bundles $X\to S$ fibered in Riemann surfaces endowed with a fiberwise smoothly varying complex structure. 

Usually, we want to put geometric structures on moduli spaces, like manifold,  varieties, and schemes. However, $\M_g$ is not a sheaf on $\sch_k$, since two algebraic curves could be isomorphic under a base extension. That implies that we cannot represent $\M_g$ by schemes and even algebraic spaces. Note that the Yoneda embedding gives a functor $\y: \sch_k \to   \sh(\aff)$ by sending $X \to \hom (-, X)$, where $\hom(-,X)$ is a functor $\aff^{op}\to \set$. In order to solve our representability problem, we want to construct a functor similar to $\hom (-, X)$, but we want to categorify the codomain $\set$ replacing by it to the category of groupoids $\gpd$, which is equivalent to the 1-homotopy type. Then we recover the functor $\M_g:\aff^{\op} \to \gpd$ which is the moduli stack over algebraic curves of genus $g$. 
\end{example}

This illustrates the need for study stacks. The adaption of stacks in differential geometry, called {\it differentiable stacks}\index{differentiable stack} or {\it smooth stacks}\index{smooth stack} has been studied in \cite{Met03}\cite{BX06}\cite{Car11} etc. On the other hand, differential stacks can be presented by {\it Lie groupoids}\index{Lie groupoid}, which has been studied widely in operator algebras and non-commutative geometry.

We can define higher stacks by switching the codomain of $\aff^{op} \to \gpd$ to higher homotopy types. 
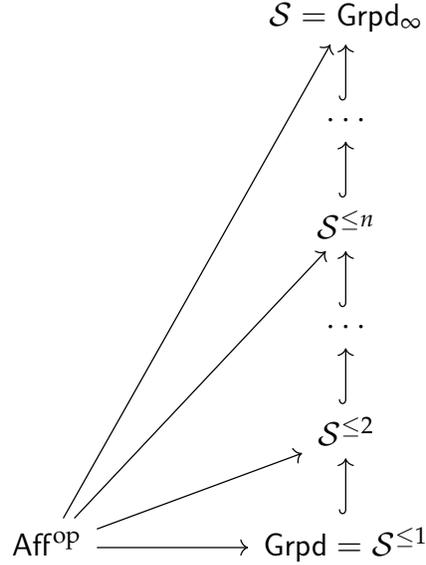
\begin{center}
	\begin{tikzcd}
&  & \mathcal{S} = \igpd                \\
&  & \cdots \arrow[u, hook] \\
&  & \mathcal{S}^{\le n} \arrow[u, hook] \\
&  & \cdots \arrow[u, hook] \\
&  & \Sp^{\le 2} \arrow[u, hook] \\
\aff^{\op} \arrow[rr] \arrow[rru]  \arrow[rruuu]  \arrow[rruuuuu] &  & \gpd = \Sp^{\le 1} \arrow[u, hook]
\end{tikzcd}
\end{center}

Here $\mathcal{S}$ denote the category of spaces, which is considered to be the $\infty$-homotopy types.

The need for  enhancing the codomain of $\aff^{\op} \to \gpd$ is shown in the following example.

\begin{example}[Pontryagin-Thom construction] 
Let $\mfd$ be the category of smooth manifolds. 	Let $X\in \mfd$ be a compact manifold and $\Omega$ the unoriented cobordism ring. $X$ represent a class $[X] \in \Omega$. The Pontryagin-Thom construction tells us that $[X]$ is classified by a homotopy class of maps $S^n$ to the Thom spectrum $MO$ for $n$ large enough. We can always pick a representative $f$ from this class such that $f$ is smooth (away from the base point) and meets the zero section $B\subset MO$ transversely. Then we have that $f^{-1}(B)$ is a manifold which is cobordant to $X$, i.e. $[f^{-1}(B)]=[X] \in \Omega$. We have the following pullback diagram
		\begin{center}
		\begin{tikzcd}
		f^{-1} B \arrow[d] \arrow[r]\arrow[dr, phantom, "\ulcorner", very near start] & B \arrow[d,""] \\
		S^n \arrow[r]           & MO         
		\end{tikzcd}
	\end{center}
\end{example}
The transversality is essential in the above construction. We first represent a class in $\Omega$ by a homotopy class of maps, which has a dense collection of smooth maps. Once we perturb the map to be transversal to the zero section, then we can get an actual manifold rather than just a class in $\Omega$. Suppose that the transversality is not required, we would have that a correspondence between smooth maps $S^n\to Mo$ and the zero loci of them.

However, transversality is essential in $\mfd$. For example, let $f:X\to Z$, $g:Y\to Z$ be arbitrary smooth maps, then the fiber product of $f$ and $g$ does not exist in $\mfd$. If we restrict to the case that $f$ and $g$ are transversal to each other, i.e. $f_{*}T_xX +g_{*}T_yY=T_zZ$ for $f(x)=g(y)=z$, then the fiber product $X\times_Z Y$ exists in $\mfd$. In particular, if either $f$ or $g$ is a submersion, then $X\times_Z Y$ exists.
\begin{rem}
	There may exist a pullback when $f$ and $g$ are not transversal to each other.
\end{rem}

This example illustrates that we want to enlarge our category of manifolds to have sufficient limits. To solve this problem, we want to pass the domain, for example, the category of commutative algebras, to differential graded commutative algebras or simplicial commutative algebras (c.f. \cite{Lur09a} \cite{TV02}\cite{TV08}). Hence, a derived $\infty$-stack should be modeled by higher categorifying both domain and codomain, i.e. as an $\infty$-functor
\begin{equation*}
    \dgcalg \to \igpd
\end{equation*}

Finally, let us look at modules over geometric objects. In classical algebraic geometry, (quasi)coherent modules play an important role. People study geometric properties using derived categories of (quasi)coherent modules.
\begin{example}[Derived categories]
    Consider $X\in \sch_k$ be a scheme over some field $k$, we can consider the category of quasi-coherent sheaves $\QCoh$ on $X$. Recall a quasi-coherent sheaf   $\mathcal{F}$ of $X$ is a sheaf of modules over the structure sheaf $\CO_X$ that is locally presentable, i.e. locally we have a following exact sequence
\begin{equation*}
\CO_X^{I_\alpha}|_{U_\alpha} \to \CO_X^{J_\alpha}|_{U_\alpha} \to \mathcal{F}|_{U_\alpha}\to 0,
\end{equation*}
where $\{{U_\alpha}\}_{\alpha}$ is a cover of $X$. The (unbounded) derived $D(X)=D\big(\QCoh(x)\big)$ is defined to be the homotopy category of a Quillen model structure on the category of unbounded chain complexes over $\QCoh(X)$. This is a powerful invariant of schemes, especially when $X$ is not smooth since it contains the cotangent complex $\Lx_X$ of $X$ and dualizing complex $\omega_X$ of $X$. Note that, if $X$ is not smooth, then $\Lx_X$ and $\omega_X$ may not be bounded.

One problem with the classical derived categories is that it does not behave well under gluing, i.e. in general we have $D(X) \not= \lim_{\{U_{\alpha}\}} D(U_{\alpha})$ where $\{U_{\alpha} \}_{\alpha}$ is a Zariski cover of $X$. An easy example is taking $X=\CP^1$ covered by two principal open sets $U_0$ and $U_1$, it's easy to verify that
\begin{equation*}
D(\CP^1) \to D(U_0) \times_{D(U_0\cap U_1)} D(U_1)
\end{equation*} 
is not faithful. In order to solve this problem, we want to pass to the $\infty$-derived category of $X$, denoted by $L_{\QCoh}(X)$ which behaves well under gluing by taking the homotopy fiber product (homotopy pullback). In particular, for our previous example $X=\CP^1$, we have
\begin{equation*}
L_{\QCoh}(X) = L_{\QCoh}(U_0)\times_{L_{\QCoh}(U_0 \cap U_1)} L_{\QCoh}(U_1)
\end{equation*}
$\infty$-derived categories are the main example of {\it stable $\infty$-categories}\index{stable $\infty$-category} introduced by Lurie \cite{Lur06, Lur17}, which is also an enhancement of {\it dg-categories}(c.f. \cite{Kel06}).
\end{example}
For differential geometry, (quasi)coherent sheaves do not really make sense since given a manifold $M$, its structure sheaf $\CO_M$, i.e. the sheaf of $\cinf$-functions over $M$, is not coherent. A substitute is to consider {\it perfect complexes} or {\it pseudo-coherent sheaves} introduced by SGA 6\cite{Ber06}. In {\cite{Blo05}}, Block constructed a dg-enhancement, called {\it cohesive modules}, for the derived category of $\CO_X$-modules over a complex manifold with coherent cohomology, which can be easily adapted to the case of various geometric structures on differentiable manifolds. We will take this idea in many of our construction and study the dg-categories related to them.

\section{Motivations}
Lie groups and Lie algebras are important tools in studying geometry, for example their actions on manifolds etc.

    \begin{thm}[Lie's 3rd theorem]
        There exists a (simply connected) Lie group $G$ corresponding to every finite dimensional Lie algebra $\g$.
    \end{thm}
    
    We can understand this as an integration functor $\int$ between Lie algebras and Lie groups
    \begin{equation*}
        \int : \mathsf{LieAlg} \to \mathsf{LieGrp}
    \end{equation*}
    
    We can also study the relation between Lie groups representations and Lie algebras representations.
    \begin{thm}
        If $G$ is simply connected, then every representation the Lie algebra ${\mathfrak {g}}$ of 
$G$ comes from a representation 
$G$ itself.
    \end{thm}
    
    We can understand this also as an integration functor
    \begin{equation*}
        \int: \rep(\g) \to \rep(G)
    \end{equation*}
Hence, under suitable assumptions, we have the following commutative squares
    \begin{center}
    \begin{equation*}
        \begin{tikzcd}[ampersand replacement=\&]
{\g} \arrow[r, "\int"] \arrow[d, "\rep"] \& {G} \arrow[d, "\rep"] \\
{\rep(\g)} \arrow[r, "\int"]                \& {\rep(G)}               
\end{tikzcd}
    \end{equation*}
    \end{center}
    We will mainly study the generalizations of these diagrams.
    
    Roughly speaking, {\it higher geometry} uses higher homotopical and categorical method in studying higher structures that traditional differential geometric method cannot handle.
	
	“Higher” usually means two directions of enhancing the classical structures (i.e. smooth/complex/symplectic manifolds, noncommutative space, etc.) which usually presents by algebras (i.e. commutative algebras, associative algebras, $C^*$-algebras):
	\begin{enumerate}
	    \item "Stacky" direction: positive grading in a dga
	    \item “Derived” direction: negative grading in a dga
	\end{enumerate}
	
	We will first generalize the objects in the previous diagram,
    \begin{itemize}
        \item Lie algebras $\	\Rightarrow$ $L_{\infty}$-algebroids. (Nuiten, Lavau etc.)
        \item Lie groups $\Rightarrow$ Lie ${\infty}$-groupoids. (Zhu, Pridham, Behrend-Getzler)
        \item Representations $\Rightarrow$ $\infty$-representations. (Abad-Crainic and Block).
    \end{itemize}
    
    The overall theme of this thesis is the study of (derived) Lie $\infty$-groupoids and (derived) $L_{\infty}$-algebroids and their representations, with applications focusing on (singular) foliations.

\section{Summary of results}
Our model of derived $\infty$-stacks for differential geometry are derived Lie $\infty$-groupoids, by which we mean Lie $\infty$-groupoid objects in some ($\infty$)-category of derived spaces, including derived manifolds (in the sense of \cite{Nui18}), derived $k$-analytic spaces (in the sense of \cite{Pri20b}), derived Banach manifolds, and derived non-commutative spaces (in the sense of \cite{Pri20c}). For the first two categories, or more generally {\it homotopy descent categories}, we construct {\it category of fibrant objects}(CFO) structures on them:
\begin{theorem}
   Let $(\dm, \mathcal{T})$ be a category with pretopology, then the category of derived Lie $\infty$-groupoids in $(\dm, \mathcal{T})$, $\ilgpd_{\dm}$,  carries a category of fibrant object structure, where fibrations are Kan fibrations, and weak equivalences are stalkwise weak equivalences.
\end{theorem}
This generalizes \cite{BG17} which considers descent categories (for example, the category of $\cinf$-schemes) and \cite{RZ20} which consider the category of Banach manifolds. We also develop a parallel theory which do not assume the underlying category has sufficient (homotopy) limits, which generalizes the result in \cite{RZ20}.

\begin{theorem}
     Given an incomplete category with locally stalkwise pretopology $(\dm, \mathcal{T})$, then the category of derived Lie $\infty$-groupoids in $(\dm, \mathcal{T})$ carries an incomplete category of fibrant object structure (iCFO), where fibrations are Kan fibrations, and weak equivalences are stalkwise weak equivalences.
\end{theorem}

These CFO or iCFO structures allow us to perform homotopical algebras explicitly, and in particular they present the $\infty$-categories associated to derived Lie $\infty$-groupoids.

The reason we want to use Lie $\infty$-groupoids rather than sheaf-theoretic $\infty$-stacks are coming from the 1-truncated case, where differential gasometers have already used Lie groupoids in studying various geometric problems, like foliations, non-commutative geometry, index theory etc. We hope that those analytic tools, like groupoid $C^*$-algebras, pseudodifferential calculus, $K$-theory, index theory etc, developed in Lie 1-groupoids can be adapted to Lie $\infty$-groupoids.

The infinitesimal counterpart of (derived) Lie $\infty$-groupoids are {\it (derived) $L_{\infty}$-algebroids}, which is a generalization of both {\it $L_{\infty}$-algebras} and {\it Lie algebroids}. The homotopy theory of derived $L_{\infty}$-algebroids is developed in \cite{Nui18}, which endows the category of derived $L_{\infty}$-algebroids over a derived manifold a {\it semi-model structure}. The Semi-model structure was introduced first in \cite{Hov98}, which is a  weaker notion than the usual model structure. This result elaborates the fact that derived $L_{\infty}$-algebroids do not have fibrant replacements in general. We study modules and representation of derived $L_{\infty}$-algebroids, and establish the equivalence between the $\infty$-representations of derived $L_{\infty}$-algebroids and the quasi-cohesive modules (c.f. \cite{Blo05}\cite{BD10}) over the Chevalley-Eilenberg algebra associated to derived $L_{\infty}$-algebroids. Note that in the Chevalley-Eilenberg algebra of a derived $L_{\infty}$-algebroid is actually a {\it stacky cdga} introduced by \cite{Pri17}. Following these, we develop Chern-Weil theory and characteristic classes for perfect $A^0$-modules with $\Z$-connections over (derived) $L_{\infty}$-algebroids. 

The main application of (derived) $\infty$-groupoids and $\li$-algebroids are (singular) foliations. Foliation studies partitions of a manifold into submanifold, which is an important tool in differential geometry and topology. The original idea of foliation can be traced back to Cartan's study on integration of PDEs, which lead to the notion of {\it exterior differential system}. 

{\it  Regular foliations}, i.e., foliations associated to integrable distributions, have been studies widely using traditional analytic tools as well as groupoids and algebroids. Singular foliations (in the sense of Stefen \cite{Ste74} and Sussmann \cite{Sus73}), are much more complicated. For example, it took many years for people, for example \cite{Pra85} \cite{Deb00}\cite{Deb01},  to try to construct holonomy groupoids of singular foliations.
It was until in \cite{AS06}, Androulidakis and Skandalis constructed holonomy groupoids for all singular foliations. Though their notion is good enough to do many constructions like $C^*$-algebras and pseudo-differential calculus, there are still many drawbacks. For example, the topology of holonomy groupoids can be pretty bad, hence the arrow spaces will not be manifolds in general. This issue reminds our principle of higher geometry: singular objects are truncation of higher homotopical objects. Hence, a natural question is, given a singular foliation, can we find a (derived) Lie $\infty$-groupoid $G_{\bt}$, such that the truncation of $G_{\bt}$ is equivalent to the holonomy groupoid?

\cite{LLS20} studied a special class of singular foliations, which admits resolution by  vector bundles. They construct $L_{\infty}$-algebroids structure on those singular foliations, and proved this construction is universal in a sense which is similar to universality in category theory. In some sense, they are looking at singular foliations (regard as an $\CO_M$-module) which admits resolution by finitely generated projective $\CO_M$-modules, and lift the dg-$\CO_M$-module structure to $L_{\infty}$-algebroid structure, i.e. the free functor
\begin{equation*}
    \free: \Mod_{\CO_M}^{\dg} \to \lalgd_{\CO_M}
\end{equation*}

However, their method does not work for many cases. For example, holomorphic singular foliations over a complex manifold $X$ only admits local resolution by finitely generated projective $\CO_X$-modules due to the finiteness property of coherent sheaves. By a result of \cite{Blo05}, we can construct a {\it cohesive module} resolving a coherent sheaf. Hence, using the tool of cohesive modules, we have
\begin{theorem}
	Given a holomorphic singular foliation $\CF$ on a compact complex manifold $\CF$ , there exist an $L_{\infty}$-algebroid $\g$ over $A$, where the linear part of $\g$ corresponds to the cohesive module $E^{\bullet}$ associated to $\CF^{\infty} = \CF\otimes_{\CO_X} \cinf(X)$.
\end{theorem}

Inspired by this result, we define {\it perfect singular foliations} to be singular foliations which are  perfect $\CO_M$-modules, i.e. foliations with local resolutions by finitely generated projective $\CO_M$-modules. With the similar method as holomorphic singular foliations, we can construct $L_{\infty}$-algebroids out of perfect singular foliations.

Next we turn to a specific class of foliations, which is called {\it elliptic involutive structures} \cite{Tre09}\cite{BCH14}\cite{Kor14}. This involutive structure is a  combinations of complex structure and foliation structure, which turns to be equivalent to transversely holomorphic foliations. We study modules over these foliations, and proves an extended version of Oka's theorem. This allows us to define $V$-coherent sheaves for an elliptic involutive structures. We also define $V$-coherent analytic sheaves for elliptic involutive structures, which generalizes Pali's $\Bar{\del}$-coherent analytic sheaves for complex manifolds in \cite{Pal03}. Using similar techniques as \cite{Blo05} on coherent sheaves,  we get 

\begin{theorem}
    Let $(X, V)$ be a compact manifold $X$ with an elliptic involutive structure $V$, then there exists an equivalence of categories between $\dcoh(X)$, the bounded derived category of
complexes of sheaves of $\CO_V$-modules with coherent $V$-analytic cohomology, and $\ho \mathcal{P}_{\CA^{\bullet}}$, the homotopy category of the dg-category of cohesive modules over $\CA^{\bullet} = \sym V^{\vee}[-1]$, i.e.   
\begin{equation*}
    \dcoh(X, \CO_V) \simeq \ho \mathcal{P}_{\CA^{\bullet}}
\end{equation*}
\end{theorem}

We then  study the homotopical structure on the category of singular foliated manifolds. \cite{GZ19} introduces the notion of {\it Hausdorff Morita equivalence} between singular foliated manifolds. We utilize their result and construct fibrations and path objects for singular foliated manifolds, and get 
\begin{theorem}
There exists an incomplete category of fibrant objects structure on the category of singular foliated manifolds $\mfd^{\sfol}$.
\end{theorem}
Following \cite{Bun18}, we also construct algebraic $K$-theory of singular foliations.

\cite{BS14} introduces the notion of {\it $\infty$-local systems} on smooth manifolds, which can be regarded as homotopical coherent representations of the fundamental $\infty$-groupoids $\Pi^{\infty}(M)$. This inspires us to define {\it $\infty$-representations} of derived Lie $\infty$-groupoids, which generalizes both \cite{BS14} for fundamental $\infty$-groupoids and \cite{AC11} for simplicial sets and Lie groupoids, and show the equivalence between $\infty$-local systems and $\infty$-representations with value in $\Mod^{\dg}_A$. We then prove an $A_{\infty}$ de Rham theorem for foliations:
\begin{theorem}
    	Let $(M, \CF)$ be a foliated manifold, there exists an $A_{\infty}$-quasi-isomorphism between $\big(\Omega^{\bullet}(\CF), -d, \wedge \big)$ and $\big(C^{\bullet}(\CF), \delta, \cup \big)$. 
\end{theorem}
and at the module level, we prove a Riemann-Hilbert correspondence for foliated $\infty$-local system,  
\begin{theorem}
    The $\infty$-category $\loci^{\infty}_{\Ch_k}\CF$ is equivalent to the $\infty$-category $\Modi_{\A}^{\coh}$, for $\A = \CE(\CF)$.
\end{theorem}
We can interpret this result as an equivalence between the  $\infty$-representations of the $L_{\infty}$-algebroid $\CF$ and the $\infty$-representations of the Monodromy $\infty$-groupoids $\mon^{\infty}(\CF)$ of $\CF$. Note that, regarding $\CF$ as an $L_{\infty}$-algebroid, its integration $\int \CF$ is equivalent to $\mon^{\infty}(\CF)$, where
\begin{equation*}
    \int: \lalgd_{\cinf{M}} \to \ilgpd_M
\end{equation*}
is the {\it Lie integration functor}\index{Lie integration functor} (c.f. \cite{Hen08}\cite{SS19}\cite{RZ20}), which is a generalization of Sullivan's integration functor for simply-connected groups \cite{Sul77}. Hence, generalizing the Riemann-Hilbert functor, we get an integration functions between $\infty$-representations of $L_{\infty}$-algebroids and Lie $\infty$-groupoids over $\cinf$-manifolds.
\begin{equation*}
    \int_{\rep}: \rep^{\infty}_A(\g) \to \rep^{\infty}_A\big(\int \g \big)
\end{equation*}
It's a natural question to ask when $\int_{\rep}$ will  be an ($\infty$-)equivalence. We won't address this problem in this thesis, and it will be one of the future topics.

\section{Organization of the paper}
{\bf Chapter 2} gives a brief introduction to the algebraic and homotopical language that we will use throughout this paper.

{\bf Chapter 3} studies homotopy theory of derived Lie $\infty$-groupoids in various derived spaces. We will construct either categories of fibrant objects (CFO) or incomplete categories of fibrant objects (iCFO) depending on the property of the underlying (homotopical) categories.

{\bf Chapter 4} studies derived $L_{\infty}$ algebroids and its representations, and we connect these to the theory of cohesive modules.  {\bf Chapter 5} studies Characteristic classes related to cohesive modules and $L_{\infty}$-algebroids.

In {\bf Chapter 6}, we study singular foliations in various categories and then use $L_{\infty}$-algebroids to study singular foliations. 

{\bf Chapter 7} studies higher monodromy and holonomy of regular and singular foliations. We study foliations on higher stacks, and gives an explicit presentation of foliations on tangent $\infty$-groupoids.

Finally, in {\bf Chapter 8}, we develop the notion of foliated $\infty$-local system, which is equivalent to the $\infty$-representation of the monodromy $\infty$-groupoid of a foliation. Then we prove the $A_{\infty}$ de Rham theorem and Higher Riemann-Hilbert correspondence for  foliated $\infty$-local system.

\part{Preliminaries}
In this section, we will recall some basic algebraic and homotopy theoretical language that we will use though out the thesis. First, we will talk about homotopical algebras, including model categories and simplicial sets. Then we will talk about dg-algebras and dg-categories and their homotopy generalizations. These will be the main tools to model derived spaces. Next, we will give a brief introduction to higher categories, which will be the main language of this thesis. Though sometimes we don't need all the generality of the language of $\infty$-categories, we still keep that direction in mind for future studies. Finally, we give an overview of {\it derived differential topology} which was developed comprehensively in \cite{Nui18} (see also \cite{Spi08}\cite{Lur09a}\cite{Pri20a}). This will be the foundation of this thesis.
\section{Homotopical algebras}

\subsection{Model categories}

Model category\index{model category} is one of the major tools in modern homotopy theory, which is originally introduced in \cite{Qui67}. Later we will see model categories are major sources of $\infty$-categories, which involves more 'higher' homotopical properties, and we will also see several variants or weaken notions of model categories, including semi-model categories, pseudo-model categories, and category of fibrant objects etc. First, let's introduce the most standard version of model categories, for which we will follow the definition in \cite{Hov07}.

\begin{defn}
    A {\it model structure}\index{model structure} on a category $\CC$ consists of three subcategories of $\CC$ called weak equivalences\index{weak equivalence!in a model category} $\mathcal{W}$, cofibrations\index{cofibration!in a model category} $\mathcal{C}$, and fibrations\index{fibration!in a model category} $\mathcal{F}$ satisfying the following properties:
    \begin{enumerate}
        \item $\CC$ is contains all finite limits and finite colimits. In particular, $\CC$ is both initial and terminal. We shall denote the initial object bye $\empty$ and the terminal object by $*$.
        \item (2-out-of-3) Let $f$ and $g$ be morphisms in $\CC$ such that $g\circ f$ is defined. If two of $f$, $g$, or $fg$ are weak equivalences, then so is the third.
        \item (Retracts) Weak equivalences, fibrations, and cofibrations are closed under retracts. Recall that a map $f:X\to Y$ is called a {\it retract}\index{retract} of $g: X'\to Y'$ if there is a commutative diagram
        \begin{center}
                \begin{tikzcd}
                {X} \arrow[r] \arrow[d, "f"] & {X'} \arrow[r] \arrow[d,"g"] & {X} \arrow[d,"f"] \\
                {Y} \arrow[r]           & {Y'} \arrow[r]           & {Y}    
                \end{tikzcd}
        \end{center}
        where the top row and the bottom row compose to $\id_X$ and $\id_Y$ respectively.
        \item (Lifting criterion) Consider the diagram
        \begin{center} 
            \begin{tikzcd}
                A \arrow[d, "f", hook] \arrow[r]    & X      \arrow[d, "g",  two heads]          \\
                B \arrow[r] \arrow[ru, "h", dotted] & Y 
                \end{tikzcd}
        \end{center}
        where $f\in \mathcal{C}$ and $g\in \mathcal{F}$. If one of the $f$ or $g$ is also a weak equivalence, then there exist a lift $g: B\to X$ such that the whole diagram commutes. We call the morphisms in $\mathcal{F}\cap \mathcal{W}$ {\it acyclic fibrations} \index{acyclic fibration! in a model category}, and the morphisms in $\mathcal{C}\cap \mathcal{W}$ {\it acyclic cofibrations}\index{acyclic cofibration! in a model category}.
        \item Every morphism 
        $f : X \to Y$ in $\CC$ can be factored as a composition $X\stackrel{\sim}{\hookrightarrow}A \twoheadrightarrow Y$ of
an acyclic cofibration followed by a fibration, and as a composition $X\hookrightarrow B \stackrel{\sim}{\twoheadrightarrow} Y$ of a
cofibration followed by an acyclic fibration.
    \end{enumerate}
    If $\CC$ have a model structure, we say $\CC$ is a model category.
\end{defn}

We call a model category {\it bicomplete}\index{bicomplete} if it contains all small limits and small colimits. We call a model category {\it factorizable} if the factorization in axiom (5) is functorial. Note that \cite{Hov07} require a model category to be bicomplete. 

An object $A$ in a model category $\CC$ is said to be {\it fibrant}\index{fibrant} if the unique map $A\to *$ is a fibration. Similarly, an object $B$ in a model category $\CC$ is said to be {\it cofibrant}\index{cofibrant} if the unique map $\empty \to B$ is a cofibration. If all objects in a model category $\CC$ is fibrant(cofibrant), then we say $\CC$ is fibrant(cofibrant).
\begin{example}[Quillen model structure on Spaces]
    Let $\Top$ be the category of topological spaces and continuous maps between them. The {\it Quillen model structure}\index{model structure!Quillen} on $\Top$ consists of the following data:
    \begin{enumerate}
        \item Weak equivalences are weak homotopy equivalences.
        \item Fibrations are Serre fibrations.
        \item Cofibrations are $\LLP(\mathcal{W}\cap \mathcal{F})$.
    \end{enumerate}
    Quillen model structure turns out to be fibrant, and the cofibrant objects are exactly CW-complexes.
\end{example}

\begin{example}[Simplicial sets]
    Consider the category of simplicial sets $s\set$. We can equip it with a model structure as follows
    \begin{enumerate}
        \item Cofibrations are monomorphisms, i.e. a map $f:X_{\bt}\to Y_{\bt}$ such that at each level $n$ we have an injection $f_n:X_n\to Y_n$.
        \item Fibrations are Kan fibrations.
        \item Weak equivalences are weak homotopy equivalences, i.e. morphisms whose geometric realization is a weak homotopy equivalence of topological spaces.
    \end{enumerate}
    
    Note that all objects are cofibrant in this model structure, and fibrant objects are call {\it Kan complexes}\index{Kan complex} or $\infty$-groupoids\index{$\infty$-groupoid}.
\end{example}
\begin{rem}
    
\end{rem}
\begin{example}[Chain complexes]
    Let $A$ be a unital associative ring. Consider $\Ch^{\ge 0}_A$ the category of non-negatively graded chain complexes of $A$-modules. We can put a model structure on $\Ch^{\ge 0}_A$ by the following data:
    \begin{itemize}
        \item The fibrations consists of all maps $f:X_{\bt} \to Y_{\bt}$ of complexes which are degreewise surjection of $A$-modules for $n > 0$.
        \item The weak equivalences are quasi-isomorphisms.
        \item Cofibrations are $\LLP(\mathcal{W}\cap\mathcal{F})$.
    \end{itemize}
    It turns out that the cofibrations in this model structure are exactly degreewise injection and $\coker(f_n)$ is a projective $A$-module for $n\ge 0$. This is called the {\it projective model structure}\index{model structure!for chain complexes} for chain complexes. This model structure is fibrant, and cofibrant objects are $X_{\bt}$ such that all components $X_i$ are projective $A$-modules.
    
    There is also a dual model structure on $\Ch^{\ge 0}_A$ called the {\it injective model structure}.
\end{example}
\subsection{Simplicial sets}
    We denote the category of simplicial sets by $s\set$. In this paper, a {\it simplicial category}\index{simplicial category} will always mean a category enriched in $s\set$, i.e. let $\CC$ be a simplicial category, for any object $x, y\in \CC$, there is a simplicial set $\HOM_{\CC}(x,y)$. At the same time, the underlying category has morphisms $\HOM_{\CC}(x,y)_0$, and the homotopy category $\ho(\CC)$ has morphisms $\pi_0 \HOM_{\CC}(x,y)$.
    
    A {\it simplicial structure}\index{simplicial structure} on a category $\CC$ is given by operations $\otimes: s\set \times \CC \to \CC$ or $(-)^-: s\set^{\op}\times \CC \to \CC$. If a category $\CC$ is equipped with a simplicial structure, then we have
    \begin{equation*}
        \Hom_{\CC}(x\otimes K,y) = \Hom_{s\set}(K,\HOM_{\CC}(x,y^K)) = \Hom_{\CC}(x,y)
    \end{equation*}
    for any $x,y \in \CC$ and $K\in s\set$.
    
    Given a category $\CC$, we write $s\CC$ for the category of simplicial objects in $\CC$.
    
    \begin{defn}
        Let $X_{\bt}$ be a simplicial object in a complete category $\CC$, we write $\Hom_{s\set}(-, x): s\set^{\op} \to \CC$ the right Kan extension of $x$ with respect to the Yoneda embedding $\y: \Delta^{\op} \to s\set$.
    \end{defn}
    Explicitly, $\Hom_{s\set}(-, x)$ can be constructed as the unique limit-preserving functor determined by $\Hom_{s\set}(\Delta^n, X_{\bt}) = X_n$ which is functorial for face and boundary maps.
    \begin{defn}
        Let $X_{\bt}$ be a simplicial object in a complete category $\CC$, and $K\in s\set$. Define the $K$-matching object in $\CC$ for $X_{\bt}$ by $M_K(X_{\bt}) = \Hom_{s\set}(K, X_{\bt})$.
    \end{defn}

If $\CC$ is also a model category, then we can equip $s\CC$ a model structure, called the {\it Reedy model structure}\index{model structure!Reedy}, as follows:
\begin{enumerate}
    \item A morphism $f:X_{\bt} \to Y_{\bt}$ is a {\it Reedy fibration}\index{fibration!Reedy} if
    \begin{equation*}
        X_n \to M_{\del \Delta^n}(X) \times_{M_{\del \Delta^n}(Y)} Y_n
    \end{equation*}
    are fibrations in $\CC$ for all $n$.
    \item Weak equivalences\index{weak equivalence!Reedy} are levelwise weak equivalences in $\CC$, i.e. $f$ is a weak equivalence if and only if each $f_n$ is a weak equivalence in $\CC$.
    \item Cofibrations are defined through the lifting properties.
\end{enumerate}

\begin{defn}
    Let $\CC$ be a model category, we write $\rhom_{s\set}(-,x):s\set^{\op} \to \CC$ the homotopy right Kan extension of $x$. We define the homotopy $K$-matching object $M_K^h(X_{\bt}) = \rhom_{s\set}(K, X)$. 
\end{defn}
Explicitly, we can realize $\rhom_{s\set}(-,x)$ by $\Hom_{s\set}(-, \BR x)$, where $\BR x$ is a fibrant replacement of $x$ is the Reedy model structure of $\CC$.

\section{Homotopy algebras}
     \subsection{dg algebras and dg categories}
    In this thesis, the main source of 'derived' and 'stacky'(or 'higher') parts of the geometry is presented by some differential graded algebras. 
    \begin{defn}
    A (cochain) {\it differential graded algebra}\index{differential graded algebra} $A= ( A^{\bt}, d)$ is a graded $k$-algebra $A^{\bt}$ with a differential $d: A^{\bt}\to A^{\bt}[1]$ satisfying
    \begin{enumerate}
        \item (Leibniz rule) $d$ is an (odd) derivation, i.e.
        \begin{equation*}
            d(ab) = (da)b + (-1)^{|a|}a(db)
        \end{equation*}
        for all $a, b \in A^{\bt}$.
        \item (Flatness) $d^2=0$.
    \end{enumerate}
    \end{defn}
    In this paper, non-negatively graded dga's will often be used as models for 'stacky' or 'higher' geometric objects. 
    Similarly, we can define chain dga's which concentrate on non-positive degrees, which are often used to model 'derived' geometric objects. We won't consider $\Z$-graded dga's and we will consider a substitute called {\it stacky dga}\index{stacky dga} in later chapters. Also, we shall consider all dga's to be unital unless otherwise mentioned explicitly.
    
    Morphisms between dga's  are degreewise morphisms which commute with differentials. A dga $A$ is called {\it (graded) commutative} if $ ab = (-1)^{|ab|}$ for any $a, b \in A$. Let $\dgcalg_k$ denote the category of commutative dga's (or cdga), and  $\dgcalg_k^{\ge 0}$ denote the category of non-negatively graded cdga's.
    
    We can equip  $\dgcalg_k^{\ge 0}$ a model structure by:
    \begin{enumerate}
        \item Fibrations are degreewise surjections.
        \item Weak equivalences are quasi-isomorphisms.
        \item Cofibrations are $\LLP(\mathcal{W}\cap \mathcal{F})$.
    \end{enumerate}
    Similar to the case of chain complexes, this model structure is called projective model structure\index{projective model structure} for cdga's. Again, all objects in this model structure are fibrant. 
    
    Next, we will look at modules over dga's.
    \begin{defn}
        Let $A=(A^{\bt}, d_A)$ be a dga. A (right) {\it dg-$A$-module}\index{dg-$A$-module}  $M=(M_{\bt}, d_M)$ over $A$ is a graded $A$-module with a differential $d_M$ such that
        \begin{enumerate}
            \item (Leibniz rule) For $a\in A, m \in M$
        \begin{equation*}
            d(m\cdot a ) = (d_M m)\cdot a + (-1)^{|m|}m\cdot (d_{A}a)
        \end{equation*}
        for all $a, b \in A^{\bt}$.
        \item (Flatness) $d^2=0$.
        \end{enumerate}
    \end{defn}
    We denote $\Mod_{A^{\bt}}^{\dg}$ the category of (unbounded) chain complexes of dg-$A$-modules. We can endow it a projective model structure similar to above. We also denote $\Mod_{A^{\bt}}^{\dg,\ge0}$ the category of non-negatively graded dg-$A$-modules, which also carries a similar model structure.

    \begin{defn}
           A {\it differential graded category}\index{differential graded category} (dg-category) $\CC$ is a category enriched over the category of $\Z$-graded cochain complexes of $k$-modules, i.e. $\CC$ consists of the following data:
        \begin{enumerate}
            \item A set of objects $\ob(\CC)$.
            \item For all $x, y \in \ob(\CC)$, a complex of morphisms $\CC(x,y)$. Write $(\CC(x,y),d)$ for this complex.
            \item The composition of morphisms is a morphism of complexes and factors through the tensor product of complexes
            \begin{equation*}
             \CC(y, z)\otimes \CC(x, y)\to \CC(x, z).
            \end{equation*}

satisfying the usual associativity condition
            
        \end{enumerate}
    \end{defn}
    For more details about dg-categories, see \cite{Kel06}. 
    \begin{example}
        A dga $A$ can be regard as a dg-category over a single object $*$.
    \end{example}
    \begin{example}
        Fix a dga $A$, let's consider $\Mod^{\dg}_{A}$. We can equip it a dg-category structure by enlarging it hom: define the morphism complex of $E, F \in \Mod^{\dg}_{A}$ to be 
        \begin{equation*}
            \HOM_{\Mod^{\dg}_{A}}^{\bt}({E,F}) = \bigoplus_{n\in \Z} \HOM_{\Mod^{\dg}_{A}}^{n}({E,F})
        \end{equation*}
        where 
        \begin{equation*}
            \HOM_{\Mod^{\dg}_{A}}^{n}({E,F}) = \prod_{i\in \Z} \Hom_{\Mod^{\dg}_{A}}(E^i,F^{i+n})
        \end{equation*}
        which are the degree $n$ morphisms of graded $A$-modules. We define a differential $d_{\HOM}:\HOM_{\Mod^{\dg}_{A}}^{\bt}({E,F})\to \HOM_{\Mod^{\dg}_{A}}^{\bt+1}({E,F})$ by
        \begin{equation*}
            d_{\HOM}(f) = [d,f] = d_F \circ f - (-1)^{|f|}f \circ d_E
        \end{equation*}
        It is easy to check that $d_{\HOM}^2 = 0$ by the fact that $d_E^2 = d_F^2 = 0$.
    \end{example}
\subsection{$L_{\infty}$-algebras}

We consider $L_{\infty}$ algebra on a graded vector space $V$. 

\begin{defn}
	Let $V=(V_{-i})$ be a graded vector space equipped with degree 1 graded symmetric bracket $\{\cdots\}_k: V\times V\times \cdots \times V \to V$ for all $k\ge 1$ such that the general Jacobi identity
	\begin{equation*}
	\sum_{i=1}^n \sum_{\sigma\in Un(i,n-i) }\epsilon(\sigma) \big\{\{ x_{\sigma(1)},\cdots, x_{\sigma(i)} \}_i, x_{\sigma(i+1)}, \cdots, x_{\sigma(n)} \big\}=0
	\end{equation*} holds, where $\epsilon$ is the sign function for graded symmetric permutations. 
\end{defn}Here we use the convention of graded symmetric bracket, which simplifies computations. In fact, $V$ is a $L_{\infty}[1]$-algebra in the usual graded antisymmetric bracket notation, where $L_{\infty}[1]$-algebra means $L_{\infty}$-algebra structure on $E[1]=\bigoplus_i E_i[1]=\bigoplus_i E_{i+1}$. We denote the 1-bracket $\{-\}_1$ by $d$, then the general Jacobi identity reads $\{-\}_1\circ \{-\}_1=d\circ d=0$, which implies $E$ is also a chain complex. Next consider $n=2$, we have 
\begin{equation*}
d\{x,y\}+\{dx,y\}+(-1)^{|x|}\{x,dy\}=0
\end{equation*} 
In short, we denote this as  $[d,m_2]=0$. For $n=3$, we have $\{\}_2\circ \{\}_2+[d,\{\}_3]=0$ by previous convention. Note that we also take the permutation into account. This equation says that the classical Jacobi identity holds up to homotopy of 3-bracket. For $n\ge 3$, we have a sequence of higher Jacobi identities:
\begin{equation*}
\sum_{i=1}^k\{\}_{k-i}\circ \{\}_i=0
\end{equation*}  
\begin{rem}
	Note that each $n$-ary bracket is a  multilinear and symmetric map, hence is determined uniquely by its values on even elements. Let $\xi \in V^{even}$. We can consider the following odd vector field
	\begin{equation*}
	Q=Q^i(\xi)\frac{\del}{\del \xi_i}=\sum_{n \ge 0} \frac{1}{n!}\{\xi, \cdots, \xi \}_n
	\end{equation*}
	where we identify $\xi=\xi^ie_i$ as $\xi^i\frac{\del}{\del \xi_i}$ as s constant vector field. Putting $\xi$ into the generalized Jacobi identity, we can define the $n$-th Jacobiator 
	\begin{equation*}
	J^n(\xi, \cdots, \xi)=\sum_{i=0}^n \frac{n!}{i!(n-i)!}\big\{ \{\xi,\cdots, \xi \}_{n-i}, \xi,\cdots ,\xi  \big\}_{i}
	\end{equation*}
\end{rem}

Define $J=\sum_{n \ge 0} \frac{1}{n!} J^n(\xi)$ which encounters all the general Jacobi identities. Observe that
\begin{align*}
Q^2=&\bigg( \sum_{n \ge 0} \frac{1}{n!}\{\xi, \cdots, \xi \}_n\bigg)\circ \bigg( \sum_{k \ge 0} \frac{1}{k!}\{\xi, \cdots, \xi \}_k\bigg)\\
=&\sum_{j\ge 0} \sum_{n+k=j}\frac{1}{n!}\frac{1}{(j-n)!}\big\{ \{\xi,\cdots, \xi \}_{j-n}, \xi,\cdots ,\xi  \big\}_{n} \\
=J
\end{align*}
Hence $Q$ is homological if and only if all Jacobiators vanish. 
\subsection{Derived algebras}
\subsubsection{Derived $\cinf$-rings}

\begin{defn}
	A $\cinf$-ring is a set $A$ such that for every $C^{\infty}$ function $\phi: \R^n\to \R^m$, there is an operation $\phi_*: A^{\times n}\to A^{\times m}$, and if we have another $C^{\infty}$ function $\psi: \R^m \to \R^k$, the following diagram commutes
	\begin{center}
		\begin{tikzcd}
		A^{\times n} \arrow[r, "\phi_*"] \arrow[rd, "(\psi\circ \phi)_*"'] & A^{\times m} \arrow[d, "\psi_*"] \\
		& A^{\times k}          
		\end{tikzcd}
	\end{center}
\end{defn}

In synthetic differential geometry, we define {\it affine $\cinf$-schemes} to be the opposite category of $\cinf$-rings, and then by gluing we get $\cinf$-schemes. $\mfd$ is a full subcategory of the category of $\cinf$-schemes $\cinf\sch$.

Let $X\in \sch_k$, we have a canonical functor $\y_X=\Hom(-,X): \aff\sch_k^{op}\to \set$. Recall, at the beginning of this note, we talked about how to categorify the codomain of this functor to get (higher) stacks. In derived algebraic geometry, we also want to pass the domain  $\aff\sch_k^{op}\simeq \calg$ to its (higher homotopical) derived version. Usually we replace commutative algebras by simplicial commutative algebras $s\calg$ or differential graded commutative algebras $\dgcalg$ (for $\Char(k)=0$ ). We will apply these constructions to $\cinf$-ring, and we will model derived $\cinf$-rings by (connective) dg-$\cinf$-rings.

\begin{defn}[\cite{CR12}]
	A {\it dg-$\cinf$-ring}\index{dg-$\cinf$-ring} is a non-negatively graded commutative dg-algebra over $\R$ such that $A_0$ has a structure of $\cinf$-ring. Denote the category of dg-$\cinf$-rings by $\cinf\alg^{\dg}$.
\end{defn}
\begin{example}[derived critical locus]
	Let $X\in\mfd$, and $\{f_i\}_{i=1}^n$ is a collection of $C^{\infty}$ functions on $X$. Consider a dg $\cinf$-ring defined by $A=C^{\infty}(M)[\eta_1,\cdots, \eta_n]$ which is the polynomial algebra generated by $\eta_1,\cdots, \eta_n$ in degree 1 over $ C^{\infty}(M)$, and satisfying $\del \eta_i =f_i$ for any $i$. $A$ models the {\it derived critical locus} of a function $f=(f_1,\cdots,f_n): M\to \R^n$ on $M$. We have $\pi_0(A)=C^{\infty}(M)/(f_1,\cdots,f_n)$. Note that if $0$ is a regular value of $f$, then $A$ is quasi-isomorphic to $C^{\infty}\big(f^{-1}(0) \big)$. 
\end{example}

\begin{prop}[\cite{CR12}]
There is a tractable model structure on $\cinf\alg^{dg}$, where weak equivalences are quasi-isomorphism (as dga), and fibrations are surjections on all non-zero degrees. 
\end{prop}
Denote the associated $\infty$-category of $\cinf\alg^{\dg}$ by $\cinf\algi$.

\subsubsection{Derived Banach manifolds}
In this section, we shall briefly construct another generalization of ordinary manifolds. First we denote by $\ban$ the category of {\it Banach manifolds}\index{Banach manifold}, i.e the objects are manifolds locally model on Banach spaces instead of $\R^n$, and maps are $\cinf$-maps (or just $C^n$-maps). For more details about Banach manifolds and geometry, see \cite{Lan95}. Note that we have a fully faithful embedding $\mfd \hookrightarrow \ban$.

We define a {\it submersion}\index{submersion} between two Banach manifolds $f:X\to Y$ to be a morphism such that for any $x\in X$, there exists an open neighborhood $U_x$ of $x$, and an open neighborhood $V_{f(x)}$ of $f(x)$, and a local section $\sigma: V_{f(x)} \to U_x$, i.e., $f\circ \sigma = \id$ and $\sigma(f(x)) = x$. Note that we will always take $U_x$ to be the connected component of $f^{-1}(V_{f(x)})$ containing $x$.

\begin{defn}
    We define {\it derived Banach manifold}\index{derived Banach manifold} to be a space locally modelled on a dga $A^{\bt}$, where  $A^0$ is of the form $\cinf{(M)}$ for some Banach manifold $M$. We denote the category of derived Banach manifolds by $\dban$. 
\end{defn}
Consider the subcategory $\dban^{\sep} \subset \dban$ whose objects consist of separable Banach manifolds which are locally modelled on separable Banach spaces B which admit ‘smooth bump functions’. The objects of $\dban^{\sep}$ carry natural affine $\cinf$-scheme structures \cite[Example 4.44]{Joy19}, hence there exists a fully faithful embedding $\dban^{\sep}\hookrightarrow \dmfd$.

Recall a morphism $f:X\to Y$ between Banach manifolds is said to be a {\it submersion}, if given any $x\in X$, there exists neighborhoods $U_x$ of $x$ and $V_{f(x)}$ of $f(x)$, such that there exists a local section $\sigma: V_{f(x)}\to U_x$.

\begin{defn}
    Let $f:X_{\bt} \to Y_{\bt}$ be a map between derived Banach manifolds, 
    \begin{enumerate}
        \item $f$ is a {\it submersion}\index{submersion}, if given any $x\in X$, there exists neighborhoods $U_x$ of $x$ and $V_{f(x)}$ of $f(x)$, such that there exists a local section $\sigma: V_{f(x)}\to U_x$, i.e. 
        $$f^*\sigma^*:\CO_{X_{\bt}}(U_x)\to \CO_{Y_{\bt}}(V_{f(x)}) \to \CO_{X_{\bt}}(U_x)\simeq {\id}
        $$
        	\item {\it $\acute{e}$tale} if the underlying map between topological spaces is local homeomorphism and the map $f^{-1}\CO_Y\to \CO_X$ is an equivalence of sheaves.
    \end{enumerate}
\end{defn}

\subsubsection{Derived EFC-algebras}
We consider derived EFC-algebras in the sense of \cite{Pri20c}.
\subsubsection{Derived non-commutative space}
We consider derived non-commutative space in the sense of \cite{Pri20b}.

\section{Higher categories and $\infty$-categories}

\subsection{Higher categories}
The basic idea of higher categories is that we don't consider only the morphisms between objects, but also want to keep track of higher morphisms, i.e., morphisms between morphisms, morphisms between morphisms between morphisms etc. 

\begin{example}
	Consider $\cat$ the category consisting of all small categories. The object of $\cat$ are just small categories, with morphisms as functors between categories. Note that we also have a notion of morphisms between morphisms here, which are just natural transformations between functors. Hence, $\cat$ is naturally a {\it 2-category}\index{2-category} with objects as small categories, 1-morhisms as functors, and 2-morhphisms as natural transformations between functors. 
\end{example}

Another 2-category which comes from geometry is that of stacks over a base scheme $S$.

Notice that in $\cat$, all morphisms between small categories in fact forms a category $\fun(\cat)$ with natural transformations between functors as morphisms. Hence, we can also think of $\cat$ as a category enriched in 1-categories. This leads to the definition of (strict) $n$-categories:

\begin{defn}
	A (strict) {\it $n$-category}\index{$n$-category} is a category enriched in (strict) $(n-1)$-categories.
\end{defn}

Unfortunately, many higher categories in geometry and topology are not strict, for example, higher structures like associativity holds only up to isomorphisms with some coherence relations. This leads to the definition of weak $n$-categories. Weak 2-categories are well-understood, but even for just weak 3-categories, the coherence conditions are very complicated and hard to work with. Hence, we would like to search for a better notion of (weak) higher categories and even $\infty$-categories.

First of all, we still want the weak $n$-categories to be enriched in weak $(n-1)$-categories. Next, we would like the weak $n$-groupoids to model the homotopy $n$-type of spaces. The latter is called the (strong) homotopy hypothesis. Followed these two principles, we have

\begin{defn}
	A (weak) {\it $\infty$-groupoid}\index{$\infty$-groupoid!weak} is a topological space.
\end{defn}
Note that the category of topological spaces clearly corresponds to the homotopy $\infty$-type. 

\begin{example}[Fundamental $\infty$-groupoid of a topological space]
	To see why the above definition is reasonable, we consider any $X\in \Top$ and construct its fundamental $\infty$-groupoid $\Pi_{\infty} X$. Define $\ob (\Pi_{\infty}(X))$ to be points in $X$, and 1-morphisms to be path in $X$. Note that path in $X$ is not strictly associative, and hence not strictly invertible as well. Define the 2 morphisms to be homotopies between paths. Observe that 1-morphisms are invertible up to homotopies, i.e. 2-morphisms. Then continuing this fashion, we can define $n$-morphisms to be homotopies between $(n-1)$-morphisms, and $(n-1)$-morphisms are then invertible up to $n$-morphisms.
\end{example}

It is still hard to see how to see what the structure of a weak $\infty$-category should be. In order to simplify our construction, we want to consider $\infty$-categories which have all morphisms invertible at some level.

\begin{defn}
	An {\it $(\infty,n)$-category}\index{$(\infty,n)$-category} is a weak $\infty$-category such that all $k$-morphisms are (weakly) invertible for $k>n$. 
\end{defn}

\begin{rem}
	Since a weak $\infty$-groupoid has all morphisms (weakly) invertible, it corresponds to an $(\infty,0)$-category. In principle, we still want the $(\infty, n)$-categories to be enriched in $(\infty, n-1)$ category.
\end{rem}
\begin{defn}
	An $(\infty,1)$-category\index{$(\infty,1)$-category} is a category enriched in topological spaces.
\end{defn}
This is one of the model for $\infty$-categories, namely the topological enriched categories. In the following sections we shall see variations of it.
\subsection{Categorical motivations of $\infty$-categories}

Recall that simplicial sets are designed to model spaces. For each category $\CC$, we can built a simplicial set related to $\CC$ by taking its nerve $\CN \CC$, where
\begin{equation*}
(\CN \CC)_n = \hom_{\cat}([n]. \CC)
\end{equation*}
\begin{example}
	Let $G$ be a group, which is considered as a category with one object, then canonically $|\CN G| \simeq BG$. Here $|-|$ denotes the geometric realization and $BG$ is the classifying space of $G$.
\end{example}

If we know information about the categories, then clearly we know information about their nerves. In fact, we have
\begin{prop}
	If $f:\CC \to \CD$ is an equivalence of categories, then $\CN(f): \CN \CC \to \CN \CD$ is a weak equivalence.
\end{prop}
This is not surprising. It is then natural to think whether the converse is true. It seems like the nerve captures all the information of the original category.

\begin{example}
	Let $[0]$ be the category $\bullet$ with one object and no non-identity morphisms, $I$ be $\bullet \leftrightarrow \bullet$. Both nerves are contractible. Consider $[1]$ being $\bullet \to \bullet$, then $\CN [1]$ is also contractible, but clearly $[1]$ is not equivalent to $I$ or $[0]$.
\end{example}

What is the problem here? Note that the weak equivalence between simplicial sets comes after taking geometric realization, where we lose the information of the directions of arrows, for example, we can not distinguish whether a 1-simplex comes an isomorphism or not. Nevertheless, the converse will hold if both $\CC$ and $\CD$ are groupoids.
\begin{prop}
	$f:\CC \to \CD$ is an equivalence of groupoids if and only if $\CN(f): \CN \CC \to \CN \CD$ is a weak equivalence.
\end{prop}

This tells us that in order to think of simplicial sets as spaces, there is a closer relation to groupoids than general categories.

In order to distinguish nerves from non-equivalent category, we have two possible constructions, and each leads to a model of $\infty$-category:
\begin{enumerate}
	\item We change the definition of weak equivalence so that non-equivalent categories will not have weakly equivalent nerves.
	\item We refine the nerve construction, which can distinguish isomorphisms from other morphisms.
\end{enumerate}
\subsection{Quasi-categories}

First, let us recall the definition of Kan complexes:
\begin{defn}
	A {\it Kan complex}\index{Kan complex} $X_{\bullet} \in s\set$ is a simplicial set such that the canonical map $X_{\bullet} \to \ast$ is a Kan fibration, i.e. for any $n\ge 0$, $0\le k \le n$, we have a lift
	\begin{center}
		\begin{tikzcd}
		\Lambda^k[n] \arrow[rr] \arrow[d, hook] &  & X_{\bullet} \\
		\Delta[n] \arrow[rru, dotted]        &  &  
		\end{tikzcd}
	\end{center}
\end{defn}

Let's look at lower dimensional case. Consider $n=2$, we have 
\begin{equation*}
\partial \Delta[2] = \begin{tikzcd}
& v_1 \arrow[rd] &   \\
v_0 \arrow[ru] \arrow[rr] &              & v_2
\end{tikzcd}
\end{equation*}
\begin{equation*}
\Lambda^0[2] = \begin{tikzcd}
& v_1  &   \\
v_0 \arrow[ru] \arrow[rr] &              & v_2
\end{tikzcd}
\end{equation*}

\begin{equation*}
\Lambda^1[2] = \begin{tikzcd}
& v_1 \arrow[rd] &   \\
v_0 \arrow[ru] &              & v_2
\end{tikzcd}
\end{equation*}
\begin{equation*}
\Lambda^2[2] = \begin{tikzcd}
& v_1 \arrow[rd] &   \\
v_0  \arrow[rr] &              & v_2
\end{tikzcd}
\end{equation*}

For example, consider the horn $i:\Lambda^0[2]\to X_{\bullet}$. This horn specifies two arrows in $X_{\bullet}$, call them $f: i(v_0)\to i(v_1)$ and $g:i(v_1)\to i(v_2)$. The horn filling property requires the extension of this horn to a 2-simplex by an arrow $h: i(v_0)\to i(v_2)$ together with a homotopy between $g\circ f$ and $h$.

\begin{example}
    Show that for any $X\in \top$, $\sing X$ is a Kan complex. Here $\sing: \Top \to s\set$ is the Singular complex functor which takes singular complexes for a given topological space. Note that $\sing$ is right adjoint to the geometric realization
	\begin{equation*}
	({\vert- \vert} \dashv \sing) \colon \Top \stackrel{\overset{{|-|}}{\longleftarrow}}{\underset{\sing}{\longrightarrow}} s\set
	\end{equation*}
	above is actually a Quillen adjunction.
\end{example}

We also have another large class of Kan complexes:
\begin{prop}
	The nerve of a groupoid is a Kan complex.
\end{prop}
\begin{proof}
	For example, for \begin{center}
		\begin{tikzcd}
	\Lambda^2[0] \arrow[rr] \arrow[d, hook] &  & \CN X_{\bullet} \\
	\Delta[2] \arrow[rru, dotted]        &  &  
	\end{tikzcd}
	\end{center}
the lift exists since we can invert $i(v_0)\to i(v_1)$ ( since $X_{\bullet}$ is a groupoid). 
\end{proof}
Since composition in a category is unique, if a simplicial set $X_{\bullet}$ is the nerve of a groupoid, all the lifts are unique. Hence,
\begin{prop}
	A Kan complex is the nerve of a groupoid iff all the lifts
		\begin{center}
		\begin{tikzcd}
		\Lambda^k[n] \arrow[rr] \arrow[d, hook] &  & X_{\bullet} \\
		\Delta[n] \arrow[rru, dotted]        &  &  
		\end{tikzcd}
	\end{center}
are unique, where $0\le n, 0\le k \le n$.
\end{prop}
This result tells us that Kan complexes are indeed groupoids 'up to homotopy'. Since Kan complexes are just fibrant objects in  $s\set$, we know that a fibrant replacement of a simplicial set behaves like the nerve of a groupoid. Now we might wonder what is the notion of categories 'up to homotopy'?

\begin{defn}
		A {\it quasi-category}\index{quasi-category} $X_{\bullet} \in s\set$ is a simplicial set such that for any $n\ge 0$, $1\le k \le n-1$, we have a lift
	\begin{center}
		\begin{tikzcd}
		\Lambda^k[n] \arrow[rr] \arrow[d, hook] &  & X_{\bullet} \\
		\Delta[n] \arrow[rru, dotted]        &  &  
		\end{tikzcd}
	\end{center}
Note that these lifts corresponding exactly the filling of 'inner horns', hence we also call a quasi-category to be a weak Kan complex.
\end{defn}

Similar to the proof of the nerve of groupoids, we have
\begin{prop}
A quasi-category is the nerve of a category iff all above lifts are unique.
\end{prop}
We can build a model structure on $s\set$ where all fibrant objects are quasi-categories, and then we would expect that there are less weak equivalence. This is the Joyal model structure on $s\set$.

\subsection{Simplicial localizations}
Let $(\cm,\CW)$ be a category with weak equivalences (homotopical category), we can get a localization $\cm[\CW^{-1}]$. For example, if $\cm$ is a model category, then $\cm[\CW^{-1}]$ corresponds to the homotopy category of $\cm$. The problem with this localization process is that it does not preserve limits and colimits.

\begin{example}
	Note that
	\begin{center}
		\begin{tikzcd}
		S^0 \arrow[d] \arrow[r] \arrow[dr, phantom, "\ulcorner", very near start]& S^1 \arrow[d,"\times 2"] \\
		\ast \arrow[r]           & S^1          
		\end{tikzcd}
	\end{center}
\end{example}
is a pullback diagram, but if we take the map from $D^2$ and localize, we get
\begin{center}
	\begin{tikzcd}
	\Hom_{\ho(\Top)}(D^2, S^0) \arrow[d] \arrow[r] & \Hom_{\ho(\Top)}(D^2, S^1) \arrow[d] \\
	\Hom_{\ho(\Top)}(D^2, \ast) \arrow[r]           & \Hom_{\ho(\Top)}(D^2, S^1)          
	\end{tikzcd} 
\end{center}
which is not a pull back since  $\Hom_{\ho(\Top)}(D^2, S^0)$ consists of two points but all the others consist of just one point. To fix this, we should take the mapping spaces and then

\begin{center}
	\begin{tikzcd}
	\map_{(\Top)}(D^2, S^0) \arrow[d] \arrow[r] \arrow[dr, phantom, "\ulcorner", very near start]& \map_{(\Top)}(D^2, S^1) \arrow[d] \\
	\map_{(\Top)}(D^2, \ast) \arrow[r]           & \map_{(\Top)}(D^2, S^1)          
	\end{tikzcd} 
\end{center}
is a homotopy pullback. In general, how should we define mapping spaces for $\cm$?

\subsection{Simplicial categories}

First, let's recall the definition of simplicial categories.

\begin{defn}
	Let $\CC$ be a category. $\CC$ is called a {\it simplicial category}\index{simplicial category} if it is enriched in simplicial sets. In particular, 
	\begin{enumerate}
		\item For any $X,Y\in \ob \CC$, we have a simplicial set $\map_{}(X,Y)$, called the mapping space between $X$ and $Y$.
		\item For any $X, Y, Z \in \ob \CC$, there is a composition map
		\begin{equation*}
		\map(X,Y) \times \map (Y,Z) \to \map (X,Z)
		\end{equation*}
		\item For any $X\in \ob \CC$, the canonical map $\Delta[0] \to \map(X,X)$ specifies the identity map.
		\item  For any $X,Y\in \ob \CC$, we have
		\begin{equation*}
		 \map(X,Y)_0 \simeq \Hom(X,Y)
		\end{equation*}
		which is compatible with compositions.
	\end{enumerate}
\end{defn}

\begin{rem}
	Note that simplicial categories could also mean simplicial objects in $\cat$, and what we presented before is simplicially enriched categories. These two notions are not equivalent. We will always mean simplicial categories to be simplicially enriched categories.
\end{rem}
\begin{rem}
	Since simplicial sets are designed to model spaces, simplicial categories provide another model for $(\infty,1)$-categories.
\end{rem}

Suppose we have a model category $\CM$ which is also a simplicial category, then we have a notion of simplicial model categories if these two notions are compatible.

\begin{defn}
	A {\it simplicial model category}\index{simplicial model category} $\CM$ is a model category as well as a simplicial category and satisfies:
	\begin{enumerate}
		\item For any $X, Y \in \ob(\CM)$ and $K\in s\set$, there exist an object $X\otimes K$ and $Y^K$ such that
		\begin{equation*}
		\map(X\otimes K, Y)\simeq \map (K, \map (X,Y))\simeq \map (X, Y^K)
		\end{equation*} 
		which is natural in $X, Y, K$.
		\item  For any $i: A\to B$ a cofibration, and $p:X\to Y$ a fibration,
		\begin{equation*}
		 \map(B,X) \to \map (A,X) \times_{\map(A,Y)} \map(B, Y)
		\end{equation*}
		is a fibration, and is a weak equivalence if either $i$ or $p$ is.
	\end{enumerate}
\end{defn}

\begin{example}
	$s\set$ is naturally a simplicial model category with $K\otimes L = K \times L$, and $\map(K,L)=L^K$ is given by
	\begin{equation*}
	\map(K,L)_n = \Hom_{s\set}(K\times \Delta[n], L)
	\end{equation*}
\end{example}
\subsection{Simplicial localizations}
Let $\CC$ be a category. The {\it free category on $\CC$}\index{free category} is a category $F\CC$ with the same objects as $\CC$ and morphisms which are freely generated by non-identity morphisms in $\CC$. There are two natural functors $\phi: F\CC \to \CC$ which takes any generating morphisms $Fc$ to the morphism $c\in \CC$, and $\psi: F\CC \to F^2\CC$ which takes the generating morphisms $Fc$ of $F\CC$ to the generating morphisms $F(F\CC)$.

\begin{defn}
	The {\it standard simplicial resolution}\index{standard simplicial resolution} of $\CC$ is a simplicial category $F_{\bullet}\CC$ which has $F^{k+1}\CC$ in degree $k$ with face map $d_i: F^{k+1} \CC \to F^k \CC$ given by $F^i\phi F^{k-i}$ and degeneracy map given by $F^i\psi F^{k-i}$. 
\end{defn}
Note that here $F_{\bullet}\CC$ is actually a simplicial object in $\cat$, but the free functor does not change objects, it could be easily shown that $F_{\bullet}\CC$ is actually a simplicially enriched category.

Now we have all the machinery to define the homotopical version of localizations with respect to weak equivalences.

\begin{defn}
	Let $(\cm, \CW)$ be a category with weak equivalences, the {\it simplicial localization}\index{simplicial localization} of $\cm$ with respect to $\CW$ is $(F_{\bullet}\CW)^{-1}(F_{\bullet}\CM)$, which is constructed by levelwise localizations. We denote $(F_{\bullet}\CW)^{-1}(F_{\bullet}\CM)$ by $L(\cm, \CW)$ or simply $L\cm$.
\end{defn}

For any simplicial categories, we can recover original categories by taking components. In fact, let $\CC$ be a simplicial category, we define its category of components $\pi_0\CC$ to be a category with $\ob(\pi_0 \CC)=\ob \CC$ and $\Hom_{\pi_0 \CC}(X,Y)=\pi_0 \map_{\CC}(X,Y)$. The following theorem tells us that the simplicial localization is indeed a higher homotopical version of homotopy categories.
\begin{thm}
	Let $(\cm, \CW)$ be a category with weak equivalences, then
	\begin{equation*}
	\pi_0 L(\cm, \CW) \simeq  \cm[\CW^{-1}]
	\end{equation*}
\end{thm}

The problem with the standard simplicial localization is that we might get just a category with proper classes of morphisms between fixed objects. This is what also happening in the ordinary localizations. Another way of producing simplicial categories is the Hammock localization.

\begin{defn}
	Let $(\cm, \CW)$ be a category with weak equivalences, the {\it hammock localization}\index{hammock localization} of $\cm$ with respect to $\CW$ is a simplicial category $L^H(\cm, \CW)$ such that
		\begin{enumerate}
			\item $\ob (L^H(\cm, \CW))=\ob(\cm)$.
			\item For any $x, y \in \cm$, $\map_{L^H\cm} (X,Y)$ has $k$-simplices the reduced hammock of width $k$ and arbitrary length $n$
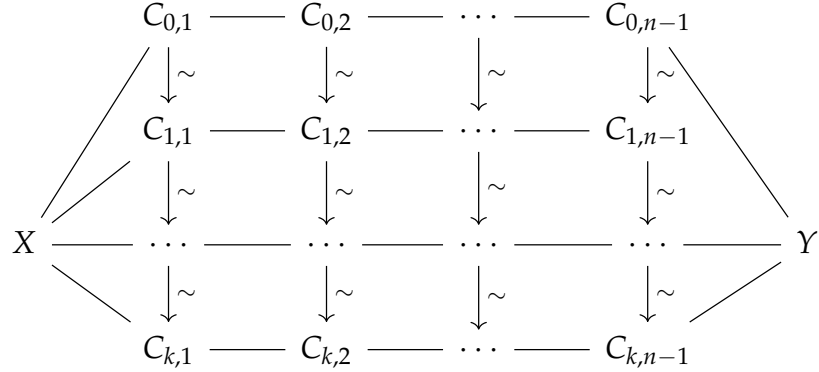
		\begin{center}
			\begin{tikzcd}
			& C_{0,1} \arrow[r, no head] \arrow[d,"\sim"] & C_{0,2} \arrow[r, no head] \arrow[d,"\sim"] & \cdots \arrow[r, no head] \arrow[d,"\sim"] & C_{0,n-1} \arrow[rdd, no head] \arrow[d,"\sim"] &   \\
			& C_{1,1} \arrow[r, no head] \arrow[d,"\sim"] & C_{1,2} \arrow[r, no head] \arrow[d,"\sim"] & \cdots \arrow[r, no head] \arrow[d,"\sim"] & C_{1,n-1} \arrow[d,"\sim"]                      &   \\
			X \arrow[ruu, no head] \arrow[ru, no head] \arrow[r, no head] \arrow[rd, no head] & \cdots \arrow[r, no head] \arrow[d,"\sim"] & \cdots \arrow[r, no head] \arrow[d,"\sim"] & \cdots \arrow[r, no head] \arrow[d,"\sim"] & \cdots \arrow[r, no head] \arrow[d,"\sim"]   & Y \\
			& C_{k,1} \arrow[r, no head]           & C_{k,2} \arrow[r, no head]           & \cdots \arrow[r, no head]           & C_{k,n-1} \arrow[ru, no head]            &  
			\end{tikzcd}
			\end{center}
			such that
			\begin{itemize}
				\item All vertical maps are in $\CW$.
				\item Horizontal maps are zig-zags, i.e. $\cdots \leftarrow \bullet\rightarrow \bullet \leftarrow \bullet \cdots$, and the arrows going to the left are in $\CW$.
				\item No column contains only identities.
				\item In each column, the horizontal arrows go in the same direction. 
			\end{itemize}
		
		\end{enumerate}
	In the case of model category, the description of hammocks localization can be greatly simplified, and it suffices to consider hammocks of length 3. For simplicity, we also denote the hammock localization by $L^H\cm$. Then a natural question is that are $L^H\cm$ and $L\cm$ the same? Or at least in some sense. 
	
	\begin{defn}
		Let $F:\CC \to \CD$ be a simplicial functor between simplicial categories, then $F$ is called a {\it Dwyer-Kan equivalence}\index{Dwyer-Kan equivalence} if 
	    \begin{enumerate}
	    	\item For any $X, Y \in \ob \CC$, $\map_{\CC}(X,Y) \to \map_{\CD}(FX,FX)$ is a weak equivalence.
	    	\item The induced functor $\pi_0 F: \pi_0 \CC \to \pi_0 \CD$ is an equivalence of categories.
	    \end{enumerate}
	\end{defn} 
\end{defn}

\begin{thm}
	Let $\cm$ be a model category, then $L^H\cm$ and $L\cm$ are Dwyer-Kan equivalent.
\end{thm}

In fact, up to Dwyer-Kan equivalence, any simplicial categories can be obtained as simplicial localizations from some categories with weak equivalences.

\begin{rem}
	We view a category with weak equivalences as a model for a homotopy theory, which determines a simplicial category by simplicial localization. Hence, the simplicial categories together with Dwyer-Kan equivalences actually form 'homotopy theory of homotopy theories'.
\end{rem}

\subsection{Homotopy mapping spaces}
Let $\M$ be a simplicial model category. First, as a consequence of the axioms for a simplicial model category, we have 
\begin{prop}
	Let $A, B, X \in \ob \cm$, $A\to B$ be a cofibration, and $X$ is a fibrant object, then 
	\begin{equation*}
	\map(B,X) \to \map (A,X)
	\end{equation*}
	is a fibration.
\end{prop}

If we have a weak equivalence $X\simeq X'$, in general $\map(X,Y)$ may not be weakly equivalent to $\map (X',Y)$, and similarly for $\map(Y,X)$ and $\map (Y,X')$. In order to get a homotopy invariant mapping space, we need to take the cofibrant/fibrant replacements.

\begin{defn}
	We define the {\it homotopy mapping space}\index{homotopy mapping space} to be $\map^h_{\cm} (X,Y) =\map_{\cm}(X^c, Y^f)$. Here $X^c$ and $Y^f$ denote the cofibrant replacement and fibrant replacement of $X$ and $Y$ respectively.
\end{defn}
Note that we can also define homotopy mapping space for a model category which is not simplicial by taking simplicial/cosimplicial resolution or $L\cm$. Let's go back to the case when $\cm$ is a model category. The following proposition justifies that the notion of homotopy mapping space is indeed a higher homotopical version of the ordinary hom set in $\ho(\cm)$.
\begin{prop}
	In a model category, $\map^h_{\cm}$ is fibrant, and we have
	\begin{equation*}
	\pi_0 \map^h_{\cm}(X,Y) \simeq \Hom_{\ho(\cm)}(X,Y)
	\end{equation*}
\end{prop} 
Now we can verify that the homotopy mapping spaces do solve problems about preserving limits at the beginning of this section.
\begin{prop}Let $\cm$ be a model category and $\CC$ be a small category.
	\begin{enumerate}
		\item Let $X\in \ob\cm$ be cofibrant, and $Y: \CC \to \cm$ a diagram of fibrant objects, then we have a weak equivalence
		\begin{equation*}
		\map_{\cm}(X, \holim_{\CC}Y_{\alpha} ) \simeq \holim_{\CC} \map_{\cm}(X, Y_{\alpha})
		\end{equation*}
		\item Let $Y\in \ob\cm$ be fibrant, and $X: \CC \to \cm$ a digram of cofibrant objects, then we have a weak equivalence
		\begin{equation*}
		\map_{\cm}(\hocolim_{\CC}X_{\alpha}, Y ) \simeq \holim_{\CC} \map_{\cm}(X_{\alpha}, Y)
		\end{equation*}
	\end{enumerate}
\end{prop}

Note that by our assumption, we can take the ordinary limits(colimits) for homotopy limits(colimits).
\section{Derived differential topology}

In this section, we will briefly introduce {\it derived differential topology}. Roughly speaking, derived differential topology is the $\cinf$ counterpart of derived algebraic geometry(DAG), where 'derived' is in the sense of Lurie and T\"{o}en-Vezzosi. Derived algebraic geometry is older and more developed. In general, derived geometry studies 'derived' spaces, which capture higher homotopical data of the classical spaces. The $\infty$-category of derived manifolds $\dm$ contains the ordinary smooth manifolds, but also many highly singular objects. People are using derived differential topology in studying moduli spaces, intersection theory, derived cobordisms etc. In order to do so, we need to apply the theory of $\infty$-categories heavily, especially Lurie's 'Structured space'. Below is a brief outline of the development of the theory of derived differential topology:

\begin{enumerate}
	\item Spivak \cite{Spi08} first defined the $\infty$-category of derived manifolds using homotopy sheaves of homotopy rings, which were introduced to study intersection theory and derived cobordisms.
	
	\item Lurie \cite{Lur09a} also gave a brief mentioning of derived differential topology in {\it DAG V: Structured space}, which will be further developed in {\it Spectral algebraic geometry} \cite{Lur18}
	
	\item Borisov-Noel\cite{BN11} gave an equivalent definition of derived manifolds using simplicial $\cinf$ rings.

	.  
	\item Joyce \cite{Joy12} introduced $\mathcal{D}$-manifolds, which form a strict 2-category. He also introduced $\mathcal{D}$-orbifolds. The main purpose of Joyce's work is to study moduli spaces arising in differential and symplectic geometry, including those used to define {\it Donaldson}, {\it Donaldson-Thomas}, {\it Gromov-Witten} and {\it Seiberg-Witten invariants}, {\it Floer theories}, and {\it Fukaya categories}.
	
	\item Nuiten \cite{Nui18} gave a comprehensive study of derived differential topology which is modeled on dg-$\cinf$-rings, based on the work of \cite{CR12}. \cite{Pri20a} took a similar approach, but restricts to simpler cases where derived manifolds are modeled on semi-free negatively graded dgas.
	
	\item During the writing process of this paper, Behrend, Liao, and Xu \cite{BLX21} develops a theory of derived manifolds modeled by bundles of curved $L_{\infty}[1]$-algebras, which is similar to \cite{Pri20a}. They prove that their derived manifolds form a {\it category of fibrant objects}, which gives an explicit presentation of its $\infty$-category.
\end{enumerate}

The idea of derived differential topology (geometry) is that we want to correct certain limits that exist in $\mfd$ but do not have the correct cohomological properties. In particular, we can form fiber products from non-transversal maps.

\subsection{Structured spaces}
Let $X\in \Top$, then we usually equip $X$ with some additional geometry structure on $X$ by associating $X$ with a sheaf $\CF$ on it.

	\begin{enumerate}
		\item Let $|X|$ be the underling topological space of a scheme $X$, then $\CF = \CO_X$ is the structure sheaf of $X$ with value in the category of commutative rings $\cring$.
		\item Again let $|X|$ be the underling topological space of a scheme $X$, we let $\CF$ be a quasi-coherent sheaf of $\CO_X$-modules on $X$.
		\item $|X|$ same as before. Let $\CF$ be an object of the derived category of quasi-coherent sheaves $D(\QCoh(X))$. This sheaf can be identified as a sheaf taking values in some $\infty$-category of module spectra.
		\item Let $X\in \mfd$, and $\CF$ be the sheaf of $C^{\infty}$ functions on $X$. This sheaf takes value in $\cring$ as well. Note that any smooth map $f:\R \to \R$ induce a morphism $\CF \to \CF$. In fact, it is easy to see that $C^{\infty}(X)$ has more delicate structure than simply being an $\R$-algebra.  
	\end{enumerate}

We want to define the {\it structured spaces}\index{structured spaces} introduced by Lurie which generalizes all the above examples and allow us to build the foundation of derived differential topology.

First, recall we say that a category is {\it locally presentable} if it is cocomplete and contains a small set $S$ of small objects such that every object in the category is a nice colimit over objects in $S$. We have a natural extension of this definition to $\infty$-categories:
\begin{defn}
	Let $\mathcal{D}$ be an $(\infty,1)$-category. We say $\mathcal{D}$ is {\it locally presentable}\index{locally presentable} if there is a small set $S$ of small objects such that every object of $\mathcal{D}$ can be presented by $(\infty,1)$-colimit over objects in $S$.
\end{defn}

\begin{rem}
	Suppose our $\infty$-categories are modeled by simplicial categories, and we assume mapping spaces are Kan complexes. We have the {\it homotopy coherent nerve}\index{homotopy coherent nerve} functor $N:s\set\cat \to s\set$ sending simplicial categories to quasi-categories. Then the $(\infty,1)$-(co)limits in quasi-categories correspond exactly homotopy (co)limits in simplicial categories. 
\end{rem}

Let $\mathcal{D}$ be a locally presentable $\infty$-category and $\mathcal{C}$ be a small $\infty$-category with finite limits. We put a Grothendieck topology on $\mathcal{C}$ generated by covers $\{U_i\to U \}$.

\begin{defn}
	A $\mathcal{D}$-valued sheaf on $\mathcal{C}$ is a functor $F: \mathcal{C}^{op}\to \mathcal{D}$ such that 
	\begin{equation*}
	F(U) \to \prod_i F(U_i) \rightrightarrows \prod_{j,k} F(U_j\times_{U}U_k)\substack{\rightarrow\\[-1em] \rightarrow \\[-1em] \rightarrow} \cdots
	\end{equation*}
	is a limit digram. Denote the category of $\mathcal{D}$-valued sheaves on $\mathcal{C}$ by $\sh(\mathcal{C};\mathcal{D})$.
\end{defn}
For example, let $X\in \Top$ and $\open(X)$ be the poset generated by open subspaces of $X$. Then $\sh(\open(X),\set)$ recovers the classical notion of sheaves.

Let $X, Y\in \Top$,  and $f: X\to Y$ be a morphism in $\Top$, i.e. a continuous function. We have an adjunction
\begin{equation*}
f^{-1}: \sh(Y,\mathcal{D}) \substack{\longrightarrow\\[-1em]  \longleftarrow}\sh(X,\mathcal{D}): f_*
\end{equation*}
where $f_*$ and $f^{-1}$ are the direct image functor and inverse image respectively. Consider the functor $\sh(-;\mathcal{D})^{op}: \Top\to \cat_{\infty}$ from topological spaces to $\infty$-categories, which sends continuous functions $f$ to direct image functors $f_{*}$ between the opposite categories of $\mathcal{D}$-valued sheaves.

In general. we can describe a functor $\mathcal{D} \to \cat_{\infty}$ equivalently by a locally cocartesian fibration $\mathcal{C}\to \mathcal{D}$. 
\begin{defn}
	Let $\pi: \mathcal{C}\to \mathcal{D}$ be a functor between $\infty$-categories. Let $\alpha: x \to y$ be a morphism in $\mathcal{D}$, we call a morphism $\tilde{\alpha}: a\to b$ in $\mathcal{D}$ {\it locally cocartesian lift}\index{locally cocartesian lift} if $\pi(\tilde{\alpha})= \alpha$, and precomposing  $\tilde{\alpha}$ induces an equivalence
	\begin{equation*}
	\tilde{\alpha}^*: \map_{\mathcal{C}_y}(b,c) \stackrel{-\circ \tilde{\alpha}}{\longrightarrow} \map_{\mathcal{C}}(a,c) \times_{\map_{\mathcal{C}}(x,y)} \{\alpha \}
	\end{equation*}
	where $\map_{\mathcal{C}_y}(b,c)$ is the mapping space in the fiber $\mathcal{C}_y$ over $y$. We called $\pi$ a {\it locally cocartesian fibration}\index{locally cocartesian fibration} if for any $\alpha: x\to y$ in $\mathcal{D}$ and $a\in \mathcal{C}_x$, we can find a locally cocartesian lift of $\alpha$. If all locally cocartesian arrows are closed under composition, we say $\pi$ is a cocartesian fibration.
\end{defn}

\begin{example}
	Let consider a simple case of cocartesian fibration. Consider the categories of modules $\Mod_A$ over some ring $A$. Consider a ring homomorphism $\phi: A\to B$, then naturally we have an induced map on modules $\phi_!: \Mod_A \to \Mod_B$ by extension of scalars, i.e. for any $M\in \Mod_A$, $\phi_!(M)= M\otimes_A B$. If we consider a category $\Mod$ of modules over all rings with objects $(A,M)$ where $M$ is a module over $A$, and morphisms have the form $(A,M) \to (B,N)$ where is a combination of ring homomorphism $A \to B$ and an $A$-linear map $M\to N$. It is easy to verify that this is a well-defined category. 
	
	Now consider a functor $\pi: \Mod\to \ring$ by mapping $(A,M)$ to $A$. Let $\phi: A\to B$ and $M$ an $A$-module, then we have a canonical map $\tilde{\phi}:(A,M)\to (B, \phi_!M)=(B, M\otimes_A B)$ induced a bijection by precomposition:
	
	\begin{equation*}
	\{(B,\phi_! M) \to (B,N)\ \text{in}\ \pi^{-1}(B)  \}\stackrel{\simeq}{\longrightarrow} \{(A,M)\stackrel{\psi}{\longrightarrow} (B,N)            \ \text{s.t}\ \pi(\psi)=\phi \}
	\end{equation*}
\end{example}

Now given a locally cocartesian fibration $\pi: \mathcal{C}\to \mathcal{D}$, let $\alpha: x\to y$ be a morphism in $\mathcal{D}$, then we have an induced functor $\alpha_! : \mathcal{C}_x \to \mathcal{C}_y$ between fiber of $x$ and $y$ respectively. In fact, let $a\in \mathcal{C}_x$, then $\alpha_!(a)=b$ for a locally cocartesian lift $a\to b$ of $\alpha$. In order to get $\alpha_!\beta_!= (\alpha\beta)_!$, we need the locally cocartesian arrows to be composable, which is ok if $\pi$ is a cocartesian fibration.

\begin{defn}[$\mathcal{D}$-structured spaces]
	Let $X\in \top$ and $\mathcal{D}$ be a locally presentable $\infty$-category, then we say $(X,\CO_X)$ is a $\mathcal{D}$-structured space if $\CO_X$ is a $\mathcal{D}$-valued sheaf. A map between two $\mathcal{D}$-structured space is a pair $(f, \tilde{f})$ where $f:X\to Y$ is a morphism in $\Top$ and $\tilde{f}: \CO_Y \to f_{*}\CO_X$ is a sheaf morphism. 
\end{defn}

Denote the $\infty$-category of $\mathcal{D}$-structured spaces by $\Top_{\mathcal{D}}$. The functor $\sh(-;\mathcal{D})^{op}: \Top\to \cat_{\infty}$ classifies a cocartesian fibration $\pi : \Top_{\mathcal{D}}\to \Top$. Denote the terminal object in $\Top$ by $*$. Consider the inclusion $i: \sh(*, \mathcal{D})\to \Top_{\mathcal{D}}$. Since $\pi: \Top_{\mathcal{D}}\to \Top$ is a cocartesian fibration, this inclusion functor has a left joint $\Gamma$ such that

\begin{equation*} 
\Gamma: \Top_{\mathcal{D}} \substack{\longrightarrow\\[-1em]  \longleftarrow}\sh(*, \mathcal{D})\simeq \mathcal{D}^{op} : i
\end{equation*}
which sends $(X, \CO_X)$ to its global sections $\CO_X(X)$.
\subsection{Construction of the $\infty$-category derived manifolds}

As we observed before, $\mfd$ does not have fiber products. In algebraic geometry, we have the category of schemes $\sch_k$ has fiber product since we have $\aff\sch^{op}_k\simeq \calg$ and we just need to compute the tensor product of commutative rings locally. Here we want to mimic the construction in algebraic geometry to extend the category of manifolds by looking at the algebraic structure on it. This method is developed in the context of {\it synthetic differential geometry}\index{synthetic differential geometry}.

As in the beginning of this section, any $X\in \mfd$ has an associated sheaf of rings of smooth function $\CO_X=C^{\infty}(X)$ on $X$. We can regard $X$ as a $\R$-scheme modeled on $\R^{\dim X}$ where the structure sheaf $\CO_X$ is a sheaf of local $\R$-algebras. Under this point of view, we can reinterpret many fundamental concepts in geometry and topology with more intrinsic constructions, for example
\begin{enumerate}
	\item The cotangent space at $x\in X$ is isomorphic to $I_p/I_p^2$, where $I_p$ is the unique maximal ideal of the stalk of $\CO_X$ at $x$.
	\item Consider the diagonal map $\Delta: X\to X\times X$. Let $\mathcal{I}$ be the sheaf of germs of smooth functions on $X\times X$ which vanish on the diagonal. Then consider the pullback of $\mathcal{I}/\mathcal{I}^2$ to $X$, denoted by $ \Delta^*( \mathcal{I}/\mathcal{I}^2)$. This construction yields a locally free sheaf called the {\it cotangent sheaf}\index{cotangent sheaf}. It is easy to verify that $ \Delta^*( \mathcal{I}/\mathcal{I}^2)$ corresponds to the cotangent bundle $T^*X$.
	\item We can also construct Taylor series (jets) similarly.
\end{enumerate}

However, a shortage of this method is that we lost the $C^{\infty}$ structure of manifolds. For example, $C^{\infty}(X)$ has much richer structures than simply being an $\R$-algebra. In order to solve this issue, we want to enlarge the category of manifolds to $\cinf$-schemes by constructions from $\cinf$-rings.

Consider the category of $\mathcal{D}=\cinf\alg_{\infty}$ structured spaces, called $\cinf$-ringed spaces, and we denote it $\Top_{\cinf}$.

\begin{defn}[Locally $\cinf$-ringed spaces]
	Define the category of {\it Locally $\cinf$-ringed spaces}\index{Locally $\cinf$-ringed spaces} $\Top_{\cinf}^{loc} \subset \Top_{\cinf}$ by
	\begin{enumerate}
		\item the objects of $\Top_{\cinf}^{loc}$ are structured spaces $(X,\CO_X)$ such that each stalk of the zeroth homotopy sheaf $\pi_0(\CO_X)_x$ is a local (discrete) $\cinf$-rings with residual field $\R$.
		\item  morphisms are morphisms $(X, \CO_X) \to (Y, \CO_Y)$ such that the map of stalks $\pi_0(\CO_{X,x})\to \pi_0 (\CO_{Y,f(x)})$ is a map of local rings. 
	\end{enumerate}
\end{defn}

\begin{prop}
	The global section functor $\Gamma$ fits into an adjunction with a right adjoint $\spec$ \begin{equation*}
	\Gamma: \Top_{\cinf}^{loc}\  \substack{\longrightarrow\\[-1em]  \longleftarrow}\   \cinf\alg :\spec
	\end{equation*}
\end{prop}

Now we define the essential image of the functor $\spec$ to be the ($\infty$-) category of {affine derived manifolds}\index{derived manifold!affine}, denoted by $\dm^{\aff}$. We call a locally $\cinf$-ringed space $(X,\CO_X)$ a {\it derived manifold}\index{derived manifold} if there exists an open cover $\{U_i\}_i$ of $X$ such that each $(U_i, \CO_X|_{U_i})\in \dm^{\aff}$. Denote the ($\infty$-) category of derived manifolds by $\dm$.

Clearly, $\mfd$ is a full subcategory of $\dm$, since for $M\in\mfd$, $M\simeq \spec\big(C^{\infty}(M) \big)$. In particular, we see that all smooth manifolds as derived manifolds are affine. 

\begin{example}
	The derived critical locus introduced before is a derived manifold. We have seen that derived critical locus is a derived enhancement of the classical critical locus. 
\end{example} 

\begin{example}
	Another large class of derived manifolds are given by {\it differential graded manifolds}. 
	A {\it graded manifold}\index{graded manifold} is defined to be a locally ringed space $\mathcal{M}=(M,\CO_{\mathcal{M}})$ where $M$ is a smooth manifold, and the structure sheaf $\CO_{\mathcal{M}}$ of $\mathcal{M}$ is locally isomorphic to $\CO(U)\otimes \sym (V^*)$ for an open set $U\subset M$ and $V$ a vector space. Here $\CO$ denotes the structure sheaf of $M$ as a smooth manifold and $\sym$ denotes the supercommutative tensor product. By a result of Batcher, any graded manifold $\mathcal{M}$ can be realized by a graded vector bundle $E\to M$ such that $\CO_{\mathcal{M}} \simeq \Gamma (\sym E^*)$. We say a graded manifold is a {\it differential graded manifold}\index{differential graded manifold} if it is equipped with a degree $+1$ vector field $Q$ with $Q^2=0$.
\end{example}

\begin{example}
	Joyce showed that many constructions in producing moduli spaces, for example, moduli spaces of $J$-holomorphic curves, yields derived manifolds. 
\end{example}

Let's go back to our motivating example of Pontryagin-Thom construction. We want to see whether $\dm$ solves the transversality problem in $\mfd$. 
\begin{prop}[\cite{Spi08}]
	The $\infty$-category $\dm$ has the following properties:
	\begin{enumerate}
		\item Let $X\in \mfd$ and $A,B$ be submanifolds of $X$, then the homotopy pull back $A\times_X^h B\in \dm$. We call $A\times_X^h B\in \dm$ the {\bf derived intersection}\index{derived intersection} of $A$ and $B$ in $X$.
		\item There exist an equivalence relation on the compact objects of $\dm$ which extend cobordism relation in $\mfd$, i.e. for any $X\in \mfd$, there is a ring $\Omega^{der}$, which is called the {\bf derived cobordism ring}\index{derived cobordism ring} over $X$, and a functor $i: \mfd \to \dm$ which induces a homomorphism $i_*: \Omega(T) \to \Omega^{der}(T)$.
		\item we have a derived cup product formula. Let $A,B$ be compact submanifolds of $X$, then we have
		\begin{equation*}
		[A]\smile[B] =[A\cap B]
		\end{equation*}
		in $\Omega^{der}(X)$.
	\end{enumerate}
\end{prop}

\begin{defn}
    Let $f:X\to Y$ be a morphism in $\dmfd$. We say $f$ is 
\begin{enumerate}
	\item a {\it closed (open) immersion}\index{immersion} if the underlying map between topological spaces is a closed (open) embedding, 
	and $\pi_0$ component of the morphism of sheaves $f^{-1}\CO_Y\to \CO_X$ is a surjection(equivalence).
	
	\item {\it $\acute{e}$tale}\index{$\acute{e}$tale morphism} if the underlying map between topological spaces is a local homeomorphism and the map $f^{-1}\CO_Y\to \CO_X$ is an equivalence of sheaves.
	
	\item {\it smooth}\index{smooth morphism} if for any $x\in X$, there are affine open neighborhood $U\ni x$, $V\ni f(x)$ such that the restricted map $f:U\to V$ is equivalent to a projection $V\times \R^n \to \R^n$. Note that this corresponds to the submersion in the classical differential geometry. In fact, $f:X\to Y$ in $\mfd$ is smooth as morphism in $\dm$ iff $f$ is a submersion.
	\item {\it locally finitely presented} if, for any point $x\in X$, there are affine open neighborhood $U\ni x$, $V\ni f(x)$ such that the restricted map $f:U\to V$ belongs to the smallest subcategory of $\aff/V$ containing $V\times \R \to V$ and is closed under finite limits.
\end{enumerate}
\end{defn}

\begin{lem} \label{lem: maps}
Let $P$ be one of the properties of maps above. then
	\begin{enumerate}
		\item The compositions of maps with property $P$ also has property $P$.
		\item Let  $f:X\to Y$ have property $P$, then the base change of $f$ under any morphism still has property $P$.
		\item Let $f:X\to Y$ be a morphism in $\dm$, and $\{U_i\to Y\}$ be an open cover of $Y$. Suppose that each base change $U_i\times_Y X\to U_i$ has property $P$, then $f$ has property $P$.
		\item Let $f:X\to Y$ be a smooth($\acute{e}tale $) surjection, and $g$ is any morphism. If $g\circ f$ is locally finitely presented or smooth($\acute{e}tale$), then $g$ is also locally finitely presented or smooth($\acute{e}tale$). 
	\end{enumerate}	
\end{lem}

    

\section{Differential geometric $L_{\infty}$ algebroids}
\subsection{$L_{\infty}$ algebroids}
Let $M$ be a smooth manifold and $E=(E_{-i})_{0\le i\le \infty}$ be a graded vector bundle over $M$. Let $\CO_M$ be the sheaf of $C^{\infty}$ functions on $M$. 
\begin{defn}
    An {\it $L_{\infty}$-algebroid}\index{$L_{\infty}$-algebroid!differential geometric} structure on $E$ is a sheaf of $L_{\infty}$ algebra structures on the sheaf of sections of $E$ with an anchor map $\rho: E_{0} \to TM$ such that
\begin{enumerate}
	\item For $n=2$ and one of the entry having order 1, we have the Leibniz rule
	\begin{equation*}
	\{x, fy \}_2=f\{x,y \}_2+\rho(x)[f]y
	\end{equation*}
	where $x\in \Gamma(E_0)$, $y\in \Gamma(E)$, $f\in \CO_M$. For $n\ge 3$, all brackets $\{\cdots \}_n$ is $\CO_M$-linear. 
	\item $E$ is a dg $\CO_M$ module. In addition, $\rho\circ d^{(1)}=0$. 
\end{enumerate} 
\end{defn}

\subsection{dg manifolds}

\begin{defn}
    A {\it graded manifold}\index{graded manifold} is defined to be a locally ringed space $\mathcal{M}=(M,\CO_{\mathcal{M}})$ where $M$ is a smooth manifold, and the structure sheaf $\CO_{\mathcal{M}}$ of $\mathcal{M}$ is locally isomorphic to $\CO(U)\otimes \sym (V^*)$ for an open set $U\subset M$ and $V$ a vector space. 
\end{defn}

Here $\CO_M$ denotes the  sheaf of $\cinf$-functions on $M$. We have the following identification for positively graded manifolds,
\begin{thm}[\cite{Bat79}]
    Let $\mathcal{M} = (M,\CO_{\mathcal{M}})$ A positively graded manifold  can be realized by a graded vector bundle $E\to M$ such that $\CO_{\mathcal{M}} \simeq \Gamma (\sym E^*)$.
\end{thm}

\begin{defn}
    A {\it dg manifold}\index{dg manifold} is a $\Z$-graded manifolds $E=\bigoplus_{i\in \Z} E_i$ with a degree 1 odd homological vector field $Q$, i.e. $Q^2=0$.
\end{defn}
dg manifolds are introduced in \cite{AKSZ97}, which is called {\it Q manifolds}\index{Q mainfold}.

If we can reduced the grading from $\Z$ to $\N$, then we call $E$ is a {\it positively graded dg manifold} or $NQ$-manifold. The structure sheaf $\CO_E$ of $E$, i.e. functions on $E$, is isomorphic to $\Gamma(\sym E^*)$, where $\sym E^*$ is the graded symmetric algebra of $E^*$, i.e. if we have $e_{i_1}, e_{i_2}\in E^*$ then
\begin{equation*}
e_{i_1} \odot e_{i_2}=(-1)^{|e_{i_2}| |e_{i_2}| } e_{i_2}\odot e_{i_1} \in \sym^2 E^*
\end{equation*}
Given a function $f\in \Gamma(\sym E^*)$, we say $f$ is of \its{arity} $k$ and \its{degree} $n$ is a section of $\sum_{(\sum_{m=1}^k i_m)=n} E^*_{-i_1}\odot \cdots E^*_{-i_k}$. Then we define vector fields to be derivations on $\Gamma(\sym E^*)$. We say a vector field $X$ is of arity $n$ if it maps a function $f$ of arity $k$ to a function $X[f]$ of arity $n+k$.

Given an $L_{\infty}$-algebroid, we could construct an $NQ$-manifold by a 'dualizing' process. Note that any functions on $E$ have arity greater than or equal to 0, and it is easy to verified that any vector fields on $E$, i.e. graded derivations of $\CO_E$, have arities $\ge -1$. Given a vector field $Q$, we can decomposed it into different arities uniquely $Q=\sum_{i\ge -1} Q^{(i)}$. To see this more clearly, let's start with the case of Lie algebroids:
\begin{example}[Lie algebroid]
	First, let us consider the case of ordinary Lie algebroid. Let $E$ be a Lie algebroid over $M$ with anchor map $\rho: E \to TM$. On $\Gamma(E)$, we have the skew-symmetric bracket $[,]$. By shifting 1 degree, we can consider a symmetric bracket on $\Gamma(E[1])$ by
	\begin{equation*}
	\{x,y\}=[\tilde{x},\tilde{y}]
	\end{equation*}
	where $\tilde{x}, \tilde{y}$ are corresponding sections in $\Gamma(E)$ if $x,y \in \Gamma(E[1])$. 
	Now the functions on $E[1]$ are identified with $\Gamma(\sym E[1]^* )$. In order to construct $Q$, it suffices to define it on $C^{\infty}(M)$ and $\Gamma (A[1]^*)$.
	First, $Q[f] \in E[1]^*$, we define
	\begin{equation*}
	<Q[f],\xi>=\rho(\xi) [f]
	\end{equation*}
	for $\xi \in \Gamma(E)$. Next, $Q[f] \in \sym^2 (E[1]^*)=\bigwedge^2 (E[1]^*)$. Define
	\begin{equation*}
	<Q[\alpha], \eta\wedge \xi>= \rho(\eta)<\alpha, \xi>- \rho(\xi)<\alpha, \eta>-<\alpha, \{ \eta, \xi\}>
	\end{equation*}
	where $\alpha, \eta, \xi \in E[1]^*$. Next, we extend $Q$ to all $\sym E[1]^*$ by derivations.  For example, let $\beta \in \sym^2 E[1]^*$, then
	\begin{align*}
	<Q[\beta], x\wedge y\wedge z> =&\rho(x) <\beta, y\wedge z>-\rho(y)<\beta, x\wedge z>+\rho(z)<\beta, x\wedge y>\\
	&- <\beta,\{y,z \}\wedge x>+< \beta, \{x,z\}\wedge y>- <\beta, \{x,y\}\wedge z>
	\end{align*}

	Let $\eta \in \sym^m E[1]^*, \xi \in \sym^n E[1]^*$, then we define 
	\begin{equation*}
	Q[\eta \wedge \xi]=Q[\eta]\wedge \xi + (-1)^{|\eta|} \eta \wedge Q[\xi]
	\end{equation*} 
	Let us calculate $Q^2$. On functions, we have, 
	
	\begin{align*}
	<Q^2[f], \eta\ \wedge \xi>=&\rho(\eta)<Q[f],\xi>-\rho(\xi)<Q[f],\eta>-<Q[f], \{\eta, \xi \}>\\=& \rho(\eta)\rho(\xi)[f]-\rho(\xi)\rho(\eta)[f]-\rho(\{\eta, \xi \})[f] \\
	=& \bigg(\rho(\eta)\rho(\xi)-\rho(\xi)\rho(\eta)-\rho(\{\eta, \xi \})\bigg)[f] 
	\end{align*}
	which vanishes due to the property of the anchor map $\rho$.
	On $\sym^2 E[1]^*$, we have 
	\begin{align*}
<Q^2[\alpha], x\wedge y\wedge z> =&\rho(x) <Q[\alpha], y\wedge z>-\rho(y)<Q[\alpha], x\wedge z>+\rho(z)<Q[\alpha], x\wedge y>\\
&- <Q[\alpha],\{y,z \}\wedge x>+< Q[\alpha], \{x,z\}\wedge y>- <Q[\alpha], \{x,y\}\wedge z>\\
=& \rho(x)\bigg(\rho(y)<\alpha, z>-\rho(z)<\alpha,y>-<\alpha, \{y,z \}> \bigg)\\
&-\rho(y)\bigg(\rho(x)<\alpha, z>-\rho(z)<\alpha,x>-<\alpha, \{x,z \}> \bigg) \\
&+\rho(z)\bigg(\rho(x)<\alpha, y>-\rho(y)<\alpha,x>-<\alpha, \{x,y \}> \bigg)\\
&-\bigg(\rho(\{y,z\})<\alpha,x>-\rho(x)<\alpha, \{y,z \}>- <\alpha, \{\{y,z \},x\}> \bigg)\\
&+\bigg(\rho(\{x,z\})<\alpha,y>-\rho(y)<\alpha, \{x,z \}>- <\alpha, \{\{x,z \},y\}> \bigg)\\
&-\bigg(\rho(\{x,y\})<\alpha,z>-\rho(z)<\alpha, \{x,y \}>- <\alpha, \{\{x,y \},z\}> \bigg)\\
	\end{align*}
	Using the property of anchor map and cancellations, we get 
\begin{align*}
<Q^2[\alpha], x\wedge y\wedge z> =& <\alpha, \{\{y,z \},x\}>- <\alpha, \{\{x,z \},y\}> +<\alpha, \{\{x,y \},z\}>\\
=&<\alpha, \{\{y,z \},x\}>+ \{\{z,x \},y\}+\{\{x,y \},z\}>
\end{align*}
Hence the Jacobi identity is exactly equivalent to $Q^2=0$ on $\sym^2 E[1]^*$. Since $Q$ on higher arity terms are defined from its action on lower arity terms, we conclude that $Q^2=0$. Note that we have constructed a dga $(A^{\bullet},d)$, where  $A^{\bullet}=\Gamma(M, \sym^{\bullet} E[1]^*)$ and  $d=Q$.
\end{example}

\begin{example}[Poisson manifolds]
Recall a Poisson manifold is a smooth manifold $M$ equipped with a Poisson bracket $\{-,-\}$ satisfies Leibniz rule $\{fg,h\}=f\{g,h\}+g\{f,h\}$ and puts a Lie algebra structure on $C^{\infty}(M)$. Note that $\{f,- \}: C^{\infty}(M)\to 	C^{\infty}(M)$ is a derivation, then we can find a $\pi\in \wedge^2TM$ such that $\{f,g\}=\pi(df,dg)$.

\end{example}

Next, we shall look at the equivalence between $L_{\infty}$ algebroids  and NQ-manifolds over $\cinf$-manifolds, which is given by Voronov\cite{Vor10}.
\begin{theorem}[\cite{Vor10}]
    Let $M$ be a $\cinf$ manifolds. There is an one-to-one correspondence between $L_{\infty}$-algebroids and NQ-manifolds over $M$.
\end{theorem}
\begin{proof}
(1) Constructing a $NQ$-manifold  from an $L_{\infty}$ algebroid.

Suppose now we are given an $L_{\infty}$-algebroid structure on a dg vector bundle $\{E_{-i}, d\}_{i\ge 0}$. First notice that for any vector $X$ we can decompose $X$ into components of different arities $X=\sum_{i= -1}^{\infty}X^{(i)}$, where each $X^{(i)}$ is of the  homogeneous arity $i$. Since $Q$ is of degree 1, the $-1$ arity part which is the contraction with $\Gamma(E_{-1})$ vanishes. Hence $Q=\sum_{i= 0}^{\infty}Q^{(i)}$. 

First, the arity 0 part is given by the dual of differential, i.e. $<Q^{(0)}[\alpha],x>=(-1)^{|\alpha|}< \alpha, d^{(i)}(x)>$, where $\alpha \in \Gamma(E^*_{-i+1}), x\in \Gamma (E_{-i})$.

By analogue of formula for ordinary Lie algebroid, we define

\begin{equation*}
	<Q^{(1)}[f],\xi>=\rho(\xi) [f]
\end{equation*}
\begin{equation*}
	<Q^{(1)}[\alpha], \eta\odot \xi>= \rho(\eta)<\alpha, \xi>- \rho(\xi)<\alpha, \eta>-<\alpha, \{ \eta, \xi\}_2>
\end{equation*}
and extend the action to higher order terms by derivation. (see previous example) 

For arities $i\ge 2$, since all $\{\cdots \}_i$'s are $\mathcal{O}_M$-linear, we define

$Q^{(i)}=\{\cdots \}_i^*= E^* \to \sym^{i+1} (E^*)$ for $i\ge 2$. It follows directly that 
$Q^{(i)}$'s are $\mathcal{O}_M$ linear for $i\ge 2$.

Next, we want to verify that $Q$ is homological. Clearly $Q$ is of degree 1 by our construction. Expanding $Q^2$ gives
\begin{equation*}
Q^2=Q^{(0)}\circ Q^{(0)} + \big(Q^{(0)}\circ Q^{(1)}+Q^{(1)}\circ Q^{(0)}\big)+ Q^{{(1)}}\circ Q^{{(1)}}+\cdots=\sum_{k=0}^{\infty}\sum_{i+j=k} Q^{(i)}\circ Q^{(j)}
\end{equation*}
Let us look at first few terms. First, we have $Q^{(0)}\circ Q^{(0)}=0$ since $d^2=0$. Next, let us consider  $\big(Q^{(0)}\circ Q^{(1)}+Q^{(1)}\circ Q^{(0)}\big)$.j
\begin{align*}
<Q^{(0)}\circ Q^{(1)}(\alpha), x\odot y>=& <Q^{(1)}(\alpha), (-1)^{|\alpha|}d(x\odot y)>\\
=&<Q^{(1)}(\alpha), (-1)^{|\alpha|}(dx\odot y+(-1)^{|x|}x\odot dy)>\\
=&(-1)^{|\alpha|}\Big(\rho(dx)<\alpha, y>- \rho(y)<\alpha, dx>-<\alpha, \{ dx, y\}>\Big)\\
&+(-1)^{|\alpha|+|x|}\Big(\rho(x)<\alpha, dy>- \rho(dy)<\alpha, x>-<\alpha, \{ x, dy\}>\Big)
\end{align*} $+Q^{(1)}\circ Q^{(0)}\big)$ exact means
On the other hand,
\begin{align*}
<Q^{(1)}\circ Q^{(0)}(\alpha), x\odot y>=&\rho(x)<Q^{(0)}(\alpha), y>- \rho(y)<Q^{(0)}(\alpha), x>-<Q^{(0)}(\alpha), \{ x, y\}>\\
=&(-1)^{|\alpha|}\Big(\rho(x)<\alpha, dy>- \rho(y)<\alpha, dx>-<\alpha, d\{ x, y\}>\Big)
\end{align*}
Note that $|\alpha|=|x|+|y|-1$. Since the anchor map is nontrivial only on $\Gamma(E_{-1})$. Combined with the fact that $\rho\circ d=0$, we have, 
\begin{align*}
<\big(Q^{(0)}\circ Q^{(1)}+Q^{(1)}\circ Q^{(0)}\big)(\alpha), x\odot y>=&(-1)^{|x|+|y|}<\alpha, \{x,dy\}>+(-1)^{|y|}<\alpha, \{x,dy\}>\\
&+(-1)^{|x|+|y|}<\alpha, \{x,y\}>\\
=&(-1)^{|x|+|y|}<\alpha, d\{x,y\}+(-1)^{|x|}\{dx,y \}+ \{x,dy\}>
\end{align*}
Hence $(Q^{(0)}\circ Q^{(1)}+Q^{(1)}\circ Q^{(0)}=0$ follows from the Leibniz rule $d\{x,y\}+(-1)^{|x|}\{dx,y \}+ \{x,dy\}=0$. It follows that all higher arities term of $Q^2$ are 0 due to general Jacobi identities.

(2) Constructing an $L_{\infty}$ algebroid from a $NQ$-manifold.

Let $E=(E_{-i})_{i\ge 1}$ be an $NQ$-manifold over $M$. We want to construct an $L_{\infty}$ algebroid structure on $E$.  First, notice that given any section $e\in \Gamma (E)$, then we can identify it as a constant vector field $\partial_{e}$ on $E$ by letting $\partial_{e}(\epsilon)=<\epsilon,e>$ for $\epsilon\in \Gamma(E^*)$. Note that here we mean $\partial_{\alpha}$ is a derivation on $\Gamma \big(\sym E^* \big)$. We denote this map by $i: \Gamma(E) \to \mathfrak{X}_{const}(E)$, where $\mathfrak{X}_{const}(E)$ denotes the vector fields on $E$ which is constant on the fiber. Let $(x_i)$, $(\phi_j^k)_j$ be local coordinates of $M$ and $E_{-k}^*$'s. Then locally we can write any vector field $X$ as 
\begin{equation*}
X=\sum_{i=1}^n v^i(x) \frac{\del }{\del x_i}+ \sum_{k=1}^{\infty}\sum_{j=1}^{\dim E_{-k}}f_{j}^k(x,\phi) \frac{\del}{\del \phi_j^k}
\end{equation*}
Let $\pi$ be the operator which projects any vector field $X$ to $X'$ which is constant on fiber and equals to $X$ on the zero locus of the fibers of $E^*$. Hence locally, $\pi$ looks like
\begin{equation*}
\bigg(\sum_{i=1}^n v^i(x) \frac{\del }{\del x_i}+ \sum_{k=1}^{\infty}\sum_{j=1}^{\dim E_{-k}}f_{j}^k(x,\phi) \frac{\del}{\del \phi_j^k}\bigg) \mapsto \bigg(\sum_{i=1}^n v^i(x) \frac{\del }{\del x_i}+ \sum_{k=1}^{\infty}\sum_{j=1}^{\dim E_{-k}}\tilde{f}_{j}^k(x) \frac{\del}{\del \phi_j^k} \bigg)
\end{equation*}
where $\tilde{f}_{j}^k(x)= {f}_{j}^k(x,0)$.
By previous identification, we can regard the image of $\pi$ as sections of $\Gamma (E)$. In fact, $i^{-1}\circ \pi: \mathfrak{X}(E) \to \Gamma(E)$ gives the desired map.

Hence, we define the anchor map $\rho:\Gamma(E_{-1})\to \Gamma(TM)$ by $<Q[f],x>=\rho(x)[f]$. Denote the one bracket $\{-\}_1$ by $d$. Define $d\alpha=i^{-1}\pi([Q,\del_{\alpha}])$ for $\alpha \in \Gamma(E)$. For higher brackets, we use Voronov's higher derived bracket formula and define 
\begin{equation*}
\{ \alpha_1, \cdots ,\alpha_n\}_n=i^{-1}\circ \pi \Big([\cdots,[[Q,\del_{\alpha_1}],\del_{\alpha_2}], \cdots ]\Big)
\end{equation*}
By the property of derived bracket, we have $J^n_Q(a_1,\cdots, a_n)=\{a_1,\cdots ,a_n \}_{n, Q^2}$ where $\{-\}_{\cdots, Q^2}$ is the $n$-th derived bracket induced by $Q^2$. Since $Q^2=0$, all Jacobiator vanish and hence we get a $L_{\infty}$ algebroid structure.

\end{proof}

\begin{rem}
	Note that there exists a map $\rho^*\circ d_{dR}: C^{\infty} (M) \to \Gamma(T^*M)\to \Gamma (E^*_{-1})$. Hence, we have the following complex
	\begin{equation*}
	\cdots \stackrel{Q^{(0)}}{\longleftarrow} \Gamma(E^*_{(-2)})\stackrel{Q^{(0)}}{\longleftarrow}\Gamma(E^*_{(-1)})\stackrel{\rho^*\circ d_{dR}}{\longleftarrow}C^{\infty} (M)
	\end{equation*} 
\end{rem}

\subsection{Lie algebroid representations}
\begin{defn}
	Let $A\stackrel{\rho}{\to} TM$ be a Lie algebroid over $M$. A representation of $A$ is a pair $(E,\nabla)$ such that $E$ is a vector bundle over $M$ and $\nabla:\Gamma(A)\times \Gamma(E)\to \Gamma(E)$ is a flat $A-$connection on $E$.
	Let $\Omega(A,E)$ be the space of $E$-valued differential forms over $\Omega(A)$, the representation $(E,\nabla)$ is equivalent to a square zero differential $d_{\nabla}$.  
\end{defn}

Note that the differential $d_{\nabla}$ is given by the usual Koszul formula
\begin{align*}
d_{\nabla}(\omega)(\alpha_1,\cdots,\alpha_{n+1})=&\sum_{i=1}^{n+1} (-1)^{i+1} \nabla_{\alpha_i}\omega (\alpha_1,\cdots,\hat{\alpha}_i, \cdots\alpha_{n+1})\\
&+\sum_{1\le i<j\le  n+1} (-1)^{i+j} \omega ([\alpha_i, \alpha_j], \alpha_1, \cdots, \hat{\alpha}_i, \cdots, \hat{\alpha}_j, \cdots, \alpha_{n+1}  )
\end{align*}
\subsection{Lie algebroid cohomology}
\begin{defn}
	The Lie algebroid cohomology groups $H^{\bullet}(A,E)$ with values in representation $(E,\nabla)$ of the cohomology groups associated to the complex $\big(\Omega(A,E), d_{\nabla}\big)$.
\end{defn}

\begin{prop}
$H^1(A,E)=E^A=\{x\in E\ |\  \nabla_{a}x=0\  \forall a\in A \} $.
\end{prop}

\begin{prop}
	$H^2(A,M)=\der(A,E)/\ider(A,E)$.
\end{prop}

\begin{prop}
	Given a Lie algebroid $A\stackrel{\rho}{\to} TM$ with a representation $(E, \nabla)$, with an $(n+2)$-cocycle $\omega_{n+2}\in \Omega^{n+2}(A,E)$ where $n\ge 1$, then we can associate them a $L_{\infty}$-algebroid with only nontrivial terms concentrated in degree $0$ and $-n$ with zero differential. Conversely, for any $L_{\infty}$-algebroid with previous properties, we can construct a Lie algebroid with representations, i.e. a quadrupole $(A\stackrel{\rho}{\to} TM,E, \nabla, \omega_{n+2}).$ 
\end{prop}

\begin{proof}($\Rightarrow$) Suppose we are given an $L_{\infty}$-algebroid $\big((A_{-i})_{i\ge 1},\rho \big)$ with only nontrivial terms $A_{-n}$ and $A_{-1}$. Define $A=A_{-1}$ with anchor $\rho: A\to TM$. Jacobi identity holds since $d$ is trivial, hence $A\to TM$ forms a Lie algebroid. Next, we define $E=A_{-n}$ as a vector bundle over $M$. We can construct a representation $\nabla:\Gamma(A)\otimes \Gamma(E)\to \Gamma(E)$ through the Lie bracket. In fact, define $\nabla_{a}s=[a,s]$, where $[-,-]$ is the 2-bracket in the $L_{\infty}$ structure. From the Leibniz rule of the anchor map, we get 
	\begin{equation*}
	\nabla_{a}(fs) =[a,fs]=f[a,s]+\rho(a)(f)s=f\nabla_a(s)+L_{\rho(a)}(f)(s)
	\end{equation*}
	for $s\in \Gamma(E), a\in \Gamma(a), f\in C^{\infty}(M)$, 
	and similarly
	\begin{equation*}
	\nabla_{fa}=[fa,s]=f[a,s]=f\nabla_a (s)
	\end{equation*} Hence $(E,\nabla)$ gives a representation of $A$. Next, let us look at the $(n+2)$-bracket $l_{n+2}:\Gamma(A)^{\otimes(n+2)}\to \Gamma(E)$. We want to construct a $(n+2)$ cycle from $l_{n+2}$
	The homotopy Jacobi identity reads
	\begin{equation*}
	\sum_{{i,j \in \mathbb{N}} \atop {i+j = n+3}} 
	\sum_{\sigma \in UnShuff(i,j)}
	\chi(\sigma,v_1, \cdots, v_{n})
	(-1)^{i(j-1)}
	l_{j} \left(
	l_i \left( v_{\sigma(1)}, \cdots, v_{\sigma(i)} \right),
	v_{\sigma(i+1)} , \cdots , v_{\sigma(n+2)}
	\right)
	= 0
	\end{equation*}
	\end{proof}
where the only nontrivial $l_i$'s are $i=2, n+2$. Note that all terms $v_i$'s are of degree 0. Hence, we can break the summation into $(n+2,1)$-unshuffle $\sigma$ and $(2,n+1)$-unshuffle $\tau$.

\begin{align*}
0=& \sum_{\tau} \chi (\tau)l_{n+2}([v_{\tau(1)}, v_{\tau(2)}],v_{\tau(3)}, \cdots, v_{\tau(n+3)} )\\
&+\sum_{\sigma} \chi(\sigma) [l_{n+2}(v_{\sigma(1)},\cdots,v_{\sigma(n+2) } ), v_{\sigma(n+3)} ]\\
=& \sum_{i<j} (-1)^{i+j+1}l_{n+2}([v_{i}, v_{j}],v_{1}, \cdots, \hat{v}_i, \cdots, \hat{v}_j, \cdots, v_{n+3} )\\
&+\sum_{i} (-1)^{n+3-i}(-1)^{n+2} [l_{n+2}(v_{1},\cdots,v_{i-1}, v_{i+1},\cdots, v_{n+3}), v_{i} ]\\
=& -\sum_{i<j} (-1)^{i+j}l_{n+2}([v_{i}, v_{j}],v_{1}, \cdots, \hat{v}_i, \cdots, \hat{v}_j, \cdots, v_{n+3} )\\
&-\sum_{i} (-1)^{i+1} [l_{n+2}(v_{1},\cdots,v_{i-1}, v_{i+1},\cdots, v_{n+3}), v_{i} ]\\
&=-d_{\nabla}l_{n+2}(v_1,\cdots,v_{n+3})
\end{align*}

Here we used the fact that there are $(n+3)$ $(n+2,1)$ unshuffles each of which has sign $(-1)^{n+3-i}$. Similarly, there are $n+3$ $(2,n+1)$ unshuffles each of which has sign $(-1)^{i+j+1}$. Hence, we have shown that $l_{n+2}$ is an $(n+2)$ cocycle.

$(\Leftarrow)$

Suppose now we are given $(A\stackrel{\rho}{\to} TM,E, \nabla, \omega_{n+2})$. We construct a dg $O_M$-module $F=\bigoplus_{i\in \Z}F_i$ with only two nontrivial terms $F_0=A$ and $F_n=E$ with zero differential. For brackets, we extend the 2-bracket comes from the Lie algebroid $A$ and the cocycle $\omega_{n+2}$. In fact, we can extend $[-,-]: \Gamma(F_i)\otimes \Gamma(F_j)\to \Gamma(F_{i+j})$ by
\begin{equation*}
[a,x]=\nabla_a(x)=-[x,a]
\end{equation*}
and 
\begin{equation*}
[x,y]=0
\end{equation*}

for $a\in \Gamma(a)$, $x,y\in \Gamma(E)$. Define $l_{n+2}=\omega_{n+2}$ and $l_i=0$ for $i\not=n+2, 2$. Thus, we get an $L_{\infty}$ structure on $F$.

\begin{prop}
	Homotopy equivalent $L_{\infty}$-algebroids of the forms in the previous proposition give cohomologous cocycles. 
\end{prop}

\subsection{Deformations and obstructions}
We want to study the deformation of a Lie algebroid $A=(A,\rho,[-,-])$. 
\begin{defn}
	A multiderivation of degree $n$ on a vector bundle $E\to M$ is defined to be a skew-symmetric multilinear map  $D:\Gamma(E)^{\otimes (n+1)}\to \Gamma(E)$ which is a derivation in each entry. Hence, for any $D\in \der^n(A)$, we have an associated map $\sigma_D: \Gamma^{\otimes n}\to \Gamma(TM)$ which is called the symbol of $D$ and satisfies
	\begin{equation*}
	D(s_1,\cdots,fs_n)=fD(s_1,\cdots, s_n)+\sigma_D(s_1,\cdots,s_{n-1})(f)s_n
	\end{equation*}
\end{defn}

Consider the space of multiderivations $\der^n(A)$ of degree $n$ on $A$, we can form a cochain complex  $\big(\der^{\bullet}(A), \delta\big)$, where the differential $\delta$ is given by the usual Koszul formula. Note that $\der^{\bullet-1}(A) \simeq C^{\bullet}_{def}(A)$, where $C^{\bullet}_{def}(A)$ is the deformation complex associated to $A$.

On $\der^{\bullet}(A)$, we can define the Gerstenhaber bracket $[D_1,D_2]=(-1)^{pq}D_1\circ D_2-D_2\circ D_1$ where
\begin{equation*}
D_2\circ D_1(s_0,\cdots,s_{p+q})=\sum_{\tau}^{}(-1)^{\tau}D_2\big(D_1(s_{\tau(0)},\cdots, s_{\tau(p)}), s_{\tau(p+1)},\cdots, s_{\tau(p+q)} \big)
\end{equation*}
where the sum is over all $(p+1,q)$ shuffles for $D_1\in \der^p(E), D_2\in \der^q(E), s_i\in \Gamma(E)$.
\begin{prop}
	The Gerstenhaber bracket makes the cochain complex $\big(\der^{\bullet}(A), \delta\big)$ a differential graded Lie algebra.
\end{prop}
Note that the Lie bracket $m\in \der^1(A)$, hence we can write the differential $\delta$ as $\delta=[m,-]$ where the bracket is the Gerstenhaber bracket. Since $\der^{\bullet-1}(A) \simeq C^{\bullet}_{def}(A)$, we have $H^{\bullet}\big(\der^{\bullet}(A)\big) \simeq H^{\bullet+1}_{def}(A)$ as a differential graded Lie algebra with zero differential.

Given a Lie algebroid $A=(A,\rho,m)$, a deformation of $A$ is a one parameter family of Lie algebroid over an interval $I$, denoted $A_t=(A,\rho_t,m_t)$ varying smoothly with respect to $t$ such that $A_0=(A,\rho,m)$. By Crainic and Moerdijk \cite{CM04}, any deformation gives a cocycle $c_0\in C^2_{def}(A)$, whose cohomology class only depends on the equivalent class of deformations. 

Recall, the Jacobi identity reads $[m,m]=0$. First, let us consider $m'=m+\phi$, where $\psi:\Gamma(E)^{\otimes 2}\to \Gamma E$ is a skew-symmetric bilinear map. Since we require $\rho'=\rho+\psi$ satisfies the Leibniz rule, we have 
\begin{align*}
m'(\alpha, f\beta)=&fm'(\alpha, \beta)+\rho'(\alpha)(f)\beta\\
=& f\big(m+\phi\big)(\alpha,\beta)+\big((\rho+\psi)(\alpha)f \big)\beta
\end{align*} 
by deleting the Leibniz rule for $m$ and $\rho$, we get
\begin{equation*}
\phi(\alpha, f\beta)=f\phi(\alpha,\beta)+\psi(\alpha)(f)\beta
\end{equation*}
which says that $\phi$ is a derivation with symbol $\psi$, i.e. $\phi\in \der^1(A), \sigma_{\phi}=\psi$. Note that $\phi$ determines $\psi$ uniquely from the Leibniz rule.

In order for $m'$ to satisfy the Jacobi identity, we have,
\begin{align*}
[m',m']=&[m+\phi,m+\phi]\\
=& [m,m]+[m,\phi]+[\phi,m]+[\phi,\phi]\\
=&-2\big(\delta \phi-\frac{1}{2}[\phi,\phi] \big)
\end{align*} 
Hence we simply need $\delta \phi-\frac{1}{2}[\phi,\phi] =0$. A linearization of this equation is $\delta \phi=0$, which says that $\phi$ is a cocycle. In this case, we call $\phi$ an infinitesimal deformation. 

Next, consider a formal one parameter deformation $m_t=m+\sum_{i=1}^{\infty}\phi_i t^i$. We must have $[m_t,m_t]=0$. Expand the bracket, we have 
\begin{equation*}
[m_t,m_t]=t[m,\phi_1]+t^2\big([\phi_1,\phi_1]+[m,\phi_2]+[\phi_2,m]\big)+O(t^3)
\end{equation*}
Hence we have $[m,\phi_1]=\delta \phi_1=0$, which say that $\phi_1$ is a cocycle. The second term gives $\delta \phi_2-\frac{1}{2}[\phi_1,\phi_1]=0$. Note that $\delta\big( \frac{1}{2}[\phi_1,\phi_1]\big)=\frac{1}{2}[\delta \phi_1,\phi_1]\big)-\frac{1}{2}[\phi_1,\delta\phi_1]\big)=0$ since $\delta$ is a graded derivation on $\der^{\bullet}(A)$. Hence, $\frac{1}{2}[\phi_1,\phi_1]\big)\in \der^2(A)$ is a cocycle. Therefore, the equation gives that $\frac{1}{2}[\phi_1,\phi_1]\big)$ has to be a coboundary, i.e. $[\frac{1}{2}[\phi_1,\phi_1]\big)]=0\in H^2\big(\der^{\bullet}(A) \big)$. Hence, the space $H^2\big(\der^{\bullet}(A)\simeq H^3_{def}(A)$ is the space of {\it obstructions} to form a one parameter family of deformations with first order term $\phi_1$.

Now suppose that we have shown that $m_t=m+\sum_{i=1}^{\infty}\phi_i t^i$ satisfies Jacobi identity up to order $n$, i.e. all terms in the expansion.


\part{Homotopy theory of derived Lie $\infty$-groupoids}

We are going to study the homotopy theory of Lie $\infty$-groupoids over various derived geometric spaces. Some of these categories have homotopical structure, i.e. $\infty$-categories or model categories, which permits us to work homotopically. Others do not have good homotopy theory, and even worse than that, they usually lack of many limits, for example, pullbacks along arbitrary morphisms. Hence, we are breaking the derived Lie $\infty$-groupoids in the following two kinds of categories:
\begin{enumerate}
    \item Homotopical categories with all finite homotopy limits, which includes:
    \begin{itemize}
        \item $\dmfd$, the category of derived manifolds.
        \item $\dan_k$, the category of derived $k$-analytic spaces.
    \end{itemize}
    \item Categories without all finite limits:
    \begin{itemize}
        \item $\dban$, the category of derived Banach manifolds.
    \end{itemize}
\end{enumerate}

We will construct explicit homotopy theory on these categories, which breaks down to the above two cases. 

For the first case, we will show that they form {\it homotopy descent categories}, and derived Lie $\infty$-groupoids these categories have {\it category of fibrant objects} structure. 

For the second case, though we don't have all finite limits, we can take advantage of the Yoneda embedding $\y: \CC \to \psh(\CC)$ which naturally extends to $\y: s\CC \to s\psh(\CC)$, then compute limits in $s\psh(\CC)$ and show representabilities. Hence, we can equip these categories an {\it incomplete category of fibrant objects} structures.

\section{$\infty$-presheaves and $\infty$-stacks}
\subsection{Simplicial presheaves}

\subsubsection{Grothendieck pretopology}

Let $\CC$ be a category, we want to define presheaves and sheaves on $\CC$. Recall that, for a topological space $S$, we define sheaves on $S$ using gluing data from an open cover of the topological space. Hence, we want to define a 'topology' on a category, which is the Grothendieck topology.

\begin{defn}
	Let $\CC$ be a category with coproducts, and a terminal object $\ast$. A {\it Grothendieck pretopology}\index{Grothendieck pretopology} $\mathcal{T}$ on $\CC$ is a collection of morphisms called {\it covers} (or {\it  covering families}) satisfies:
	\begin{enumerate}
		\item each object $X\in \CC$ has a collection of covers $\{U_i\to X \}$;
		\item isomorphisms are covers;
		\item pullbacks of covers are covers;
		\item composition of covers are covers;
		\item the canonical map $X\to \ast $ is a cover.
	\end{enumerate}
	\end{defn}

	For simplicity, we will simply say pretopology for Grothendieck pretopology if there is no confusion. Grothendieck pretopology is also called {\it basis for a Grothendieck topology}. As the name suggests, each Grothendieck pretopology generates a {\it Grothendieck topology}\index{Grothendieck topology}. 
	\begin{defn}
	    A {\it Grothendieck topology} $\tau$ on a category $\CC$ consists of the following data:
	    \begin{enumerate}
	        \item for any object $x\in \CC$, there is a family $\cov(x)$ of covering sieves over $x$, i.e. subfunctors of the representable functor $\y_x = \Hom(-, x)$.
	        \item (Stability under base change) For any morphism $f:x\to y$ in $\CC$ and $u\in \cov(X)$, we have $f^*(u) = u\times_{\y_x} \y_y$.
	        \item (Local character condition) Let $x\in \CC$, $u\in \cov(x)$, and $v$ be any sieve on $X$. If for all $y\in \CC$ and $f\in u(y)$, we have $f^*(v) \in \cov(y)$, then $v\in \cov(x)$.
	    \end{enumerate}
	\end{defn}
	Given a Grothendieck pretopology $\mathcal{T}$, the Grothendieck topology $\tau$ generated from $\mathcal{T}$ is that for which a sieve ${S_i\to U}$ is covering if it contains a covering family of morphisms.
	We call a category with Grothendieck topology a {\it (Grothendieck) site}\index{site}. For simplicity, we will also call a category with pretopology a site, by which we mean the site generated by the pretopology.

Now consider a category with pretopology, we can define the category of presheaves $\psh(\CC)$ on $\CC$ consists of contravariant functors $\CC \to \set$. 
\begin{defn}
       A presheaf $F \in \psh(\CC)$ is a sheaf if $F(X)$ is the limit of the diagram
    $$ F(U)\rightrightarrows F(U\times_X U)$$.
\end{defn}
We denote the category of sheaves on $\CC$ by $\sh(\CC)$. The inclusion functor $\iota: \sh(\CC)\to \psh(\CC)$ has an exact left adjoint functor, $s:\psh(\CC) \to \sh(\CC)$ which is called the {\it associated sheaf functor} (or {\it sheafification functor}).
Now we define $s\psh(\CC)$ to be the category of simplicial objects in $\psh(\CC)$. Note that we can also define $s\psh(\CC)$ as contravariant functor from $\CC$ to $s\set$. We can endow $s\psh(\CC)$ a model structure by the following data:
\begin{enumerate}
    \item A morphism $f:F \to G$ in $s\psh(\CC)$ is called a {\it global fibration}{\index{global fibration}} if, for any $x\in \CC$, the induced morphism $F(x)\to G(x)$ is a Kan fibration of simplicial sets.
    \item A morphism $f:F \to G$ in $s\psh(\CC)$ is called a {\it global equivalence}{\index{global equivalence}} if, for any $x\in \CC$, the induced morphism $F(x)\to G(x)$ is a weak equivalence of simplicial sets.
    \item The cofibrations are defined through the standard lifting property.
\end{enumerate}
We call this model structure the {\it global model structure}{\index{global model structure!for simplicial presheaves}}(\cite{Jar87}) for simplicial presheaves.

Given a simplicial presheaf $F:\CC^{\op} \to s\set$, we define a presheaf $\pi^{\psh}_0(F):C^{\op} \to \set$ by sending any $x\in \CC$ to $\pi_0(F(x))$. Similarly, for any $x\in \CC$ and any 0-simplex $s\in F(X)_0$ we define presheaves of groups on $\CC/x$
$$
\pi^{\psh}_i(F, s): (\CC/x)^{\op} \to \gp
$$
by sending $f:y\to x$ to $\pi_i(F(y), f^*(s))$.

\begin{defn}
    Given a simplicial presheaf $F:\CC^{\op} \to s\set$, we define the {\it homotopy sheaves}\index{homotopy sheaves} of $F$ to be the sheafification of $\pi^{\psh}_0(F)$ and $\pi^{\psh}_i(F, s)$ for $i \ge 1$, which we denote by $\pi^{}_0(F)$ and $\pi^{}_i(F, s)$ respectively.
\end{defn}
Using homotopy sheaves, we can refine the global model structure as follows:
\begin{enumerate}
    \item A morphism $f:F \to G$ in $s\psh(\CC)$ is called a {\it local equivalence}{\index{local equivalence}} if it satisfies:
    \begin{itemize}
        \item The induced morphism $\pi_0(F)\to \pi_0(G)$ is an isomorphism of sheaves.
        \item For any $x\in \CC$, any $s\in F(x)_0$, and $i\ge 1$, the induced morphism $\pi_i(F,s)\to \pi_i(F', f(s))$ is an isomorphism of sheaves on $\CC/x$.
    \end{itemize}
    \item The {\it local cofibration}\index{local cofibration} is defined as the same as the global cofibration.
    \item The fibrations are defined through the standard lifting property.
\end{enumerate}
This model structure is called the {\it local model structure}\index{local model structure!} for simplicial presheaves. In \cite{DHI04}, we have an easy characterization of fibrant object in the local model structure. 
\begin{defn}
    Let $x\in \CC$, we define a {\it hypercovering}\index{hypercovering} of $x$ to be a simplicial presheaf $H$ with a morphism $H\to x$ such that,
    \begin{enumerate}
        \item For each $n$, $H_n$ is a disjoint union of representable presheaves.
        \item For each $n$, the morphism of presheaves
        \begin{equation*}
            H_n\simeq \Hom(\Delta[n], H)\to \Hom(\del \Delta[n], H) \times_{\Hom(\del \Delta[n], x)} \Hom( \Delta[n], x) 
        \end{equation*}
    \end{enumerate}
    
\end{defn} 
Let $f\in s\psh(\CC)$ and $H\to x$ a hypercovering of $x\in \CC$, we can construct an augmented cosimplicial diagram
$$
F(x) \to ([n] \mapsto F(H_n))
$$
\begin{thm}[\cite{DHI04}]
    An object $F\in s\psh(\CC)$ is fibrant in the local model structure if and only if it satisfies:
    \begin{enumerate}
        \item For any $x\in \CC$, $F(x)$ is fibrant.
        \item For any $x\in \CC$ and any hypercovering $H\to x$, the natural morphism
        \begin{equation*}
            F(x)\to \hocolim_{[n]\in \Delta}F(H_n)
        \end{equation*}
        is an equivalence of simplicial sets.
    \end{enumerate}
\end{thm}
The first condition is rather anodyne, whereas the second one is similar to the condition to be a sheaf. In fact, if $F$ is an ordinary presheaf (considered as a constant simplicial presheaf), then the second condition simplifies exactly to the sheaf condition. We will call an object $F\in s\psh(\CC)$ satisfying the second condition above a {\it stack}\index{hypersheaves} over $\CC$, which are also called {\it $\infty$-stack}\index{stack!$\infty$-} or {\it hypercomplete $\infty$-sheaves}\index{hypercomplete $\infty$-sheaves} in \cite{Lur09b}, and {\it stacks} in \cite{TV02}. We call the homotopy category $\ho(s\psh(\CC))$ the {\it homotopy category of hypersheaves} on the site $(\CC, \tau)$.

\subsubsection{$(\infty,1)$-Grothendieck topology}

In this section, we will generalize Grothendieck topology to a more general setting, which consider the underlying category has already had some homotopical structure. Roughly speaking, we will define a 'simplicial' Grothendieck topology on a 'simplicial' category.

We define an {\it ($\infty$,1)-Grothendieck topology}\index{Grothendieck topology!($\infty$,1)-} $\mathcal{T}$ on an $\infty$-category $\mathsf{C}$ consists of data such that for any object $c$ in $\mathsf{C}$, there is a collection of sieves, called {\it covering sieves}, such that
\begin{enumerate}
	\item For each $c\in \mathsf{C}$, the overcategory $\mathsf{C}/c$ is a covering sieve, i.e. the monomorphism $\id:\y(c) \to \y(c)$ is a cover.
	\item Pullback of a covering sieve is a sieve.
	\item For a covering sieve $s$ on $c\in \mathsf{C}$ and $t$ any sieve on $c$, if $f^*t$ is a covering sieve for all $f\in s$, then $t$ is a covering sieve.
\end{enumerate}

Equivalently, we have the following characterization
\begin{thm}[\cite{Lur09a}]
     The data of an $(\infty, 1)$-Grothendieck topology is given by the data of an ordinary Grothendieck topology on $\ho(\CC)$. Consider a property $\mathbf{P}$ of morphisms in $\ho(\CC)$, we say a morphism $f$ in $\CC$ satisfies $\mathbf{P}$ if its image in $\ho(\CC)$ satisfies $\mathbf{P}$.
\end{thm}

\begin{defn}
	An $\infty$-category equipped with an ($\infty$,1)-Grothendieck topology $\mathcal{T}$ is called an {\it $(\infty,1)$-site\index{site!$(\infty,1)$-}}.
\end{defn}

\begin{defn}[\cite{Lur09a}] 
	A simplicial object in an $\infty$-category $\mathsf{C}$ is defined to be an $(\infty,1)$-functor $X:\Delta^{op}\to \mathsf{C}$. We denoted the corresponding $\infty$-category to be $\mathsf{C}^{\Delta^{op}}=\funi(\Delta^{op}, \mathsf{C})$.
\end{defn}
\subsubsection{$\infty$-Yoneda embedding}
Let $\mathsf{C}$ be an $\infty$-category, then an {\it $(\infty,1)$-presheaf} on $\mathsf{C}$ is an $(\infty, 1)$-functor $F: \mathsf{C}^{op} \to \igpd$. Denote the $\infty$-category of $(\infty,1)$ sheaves on $\mathsf{C}$ by
$\pshi(\mathsf{C})$ given by
\begin{equation*}
\pshi(\mathsf{C})= \funi(\mathsf{C}^{op}, \igpd )
\end{equation*} 

\begin{defn}(($\infty,1)$-Yoneda embedding)
	Let $\mathsf{C}$ be an $\infty$ category. We define the {\it ($\infty,1)$-Yoneda embedding}\index{Yoneda embedding!$\infty,1)$-} of $\mathsf{C}$ to be the $(\infty,1)$-functor $$\y: \mathsf{C} \to \pshi(\mathsf{C})$$
	by $\y(X)=\map_{\mathsf{C}}(-,X): \mathsf{C}^{op} \to \igpd$. This map is fully faithful. 
\end{defn}

\begin{defn}
	A {\it sieve}\index{sieve} in an $\infty$-category $\mathsf{C}$ is a full sub-$\infty$-category $\mathsf{D}$ such that $\mathsf{D}$ is closed under precomposing morphisms in $\mathsf{C}$. A {\it sieve} on an object $c\in \mathsf{C}$ is a sieve in $\mathsf{C}/c$. This is equivalent to say a sieve on $c$ is an equivalent class of monomorphisms $\{U\to \y(c) \}$ in $\pshi(\mathsf{C})$.
\end{defn}

Let $s$ be a sieve, and $f\in \Hom_{\mathsf{C}}(d,c)$, we define the {\it pullback sieve $f^*s$ on $d$} to be all morphisms to $d$ such that for any $g\in f^*s$, $f\circ g$ is equivalent to a morphism in $s$.
Denote the $\infty$-category of $(\infty,1)$-presheaves on $\dm$ by $\pshi(\dm)$. Note that $\pshi(\dm)$ is naturally simplicially enriched.

\begin{prop}[\cite{Lur09a}]
	Let $\CC$ be an $\infty$-category. Consider $X\in \CC$, denote $\y_X \in \pshi(\CC)$ defined by
	\begin{equation*}
	\y_X(U) := \map_{\CC}(U,X) \in \igpd
	\end{equation*}
	
	Then for any $(\infty,1)$-presheaf $F$ on $\CC$, there is a canonical isomorphism of $\infty$-groupoids
	\begin{equation*}
	F(X) \simeq \map_{\pshi(\CC)} (\y_X, F)
	\end{equation*}
\end{prop}  

Then we have an embedding $\y: \CC \to \pshi(\CC)$
Denote the $\infty$-category of simplicial $(\infty,1)$-presheaves on $\dm$ by $s\pshi(\dm)$.

\section{Derived Lie $\infty$-groupoid}
We want to consider the $\infty$-groupoid objects in the $\infty$-category of derived manifolds $\dm$. First, we want to study simplicial derived manifolds. For simplicity, we denote $\dm$ for either one of $\dmfd$, $\dan$, and $\dban$, and call them derived manifolds unless otherwise specified.

Let $(\dm, \mathcal{T}_{ss})$ be the category of derived manifold equipped with smooth surjection pretopology. Let $X_{\bullet}: \Delta^{op}\to \dm$ be a simplicial object in derived manifolds.

\begin{defn}
	A simplicial map $p:X_{\bullet}\to Y_{\bullet}$ is a {\it Kan fibration} if for each horn inclusion $\Lambda^i[n] \subset \Delta[n]$, the matching map
	\begin{equation*}
	X(\Delta[n])\longrightarrow X(\Lambda^i[n])\times_{Y(\Lambda^i[n])} Y(\Delta[n])
	\end{equation*}
	is a cover for all $n\ge 1, 1\le 0\le i\le n$.
\end{defn}
Recall that
\begin{defn}
	A simplicial set $X_{\bullet}$ is an {\it $\infty$-groupoid}\index{$\infty$-groupoid! in simplicial sets} if the canonical map $X_{\bullet}\to *$ is a Kan fibration.
\end{defn}
\begin{defn}
	Let $\mathsf{C}, \mathsf{D}$ be two $\infty$-categories (quasi-categories), and {\it $(\infty,1)$-functor}\index{functor!$(\infty,1)$-} $F:\mathsf{C}\to \mathsf{D}$ is a $s\set$ morphism of the underlying simplicial sets.
\end{defn}
In general, a simplicial object in an $\infty$-category $\mathsf{C}$ is defined to be an $(\infty,1)$-functor $\Delta^{op}\to \mathsf{C}$ . We define the $\infty$-category of {\it simplicial derived manifolds}\index{derived manifold!simplicial} to be the $\infty$-category of $(\infty,1)$ functors
\begin{equation*}
s\dm=\dm^{\Delta^{op}}=\funi(\Delta^{op}\to \dm)
\end{equation*}

Let $\rhom_{\Delta}(-,X): s\set^{\op} \to \dm$ be the homotopy right Kan extension of $X_{\bullet}\in s\dm$. Note that since $\dm$ has all finite homotopy limits, $\rhom_{\Delta}(K, X_{\bullet}) \in \dm$ for finite $K\in s\set$. 

\begin{defn}
    Let $K$ be a simplicial set, we define the {\it homotopy $K$-matching object}\index{homotopy $K$-matching object} in $\dm$  to be $M^h_K X_{\bullet} = \rhom_{\Delta}(K,X_{\bullet}).$
\end{defn}

Let $T$ be a finite simplicial set and $S \hookrightarrow T$ an inclusion of simplicial subset. Let $f: X_{\bullet} \to Y_{\bullet}$ be a morphism between simplicial derived manifolds, then we denote the $\Hom(S \hookrightarrow T, f)$ to be the homotopy fiber product
$M^h_S X_{\bullet} \times^h_{M^h_S Y_{\bullet}} M^h_T Y_{\bullet} $.

There exists a fully faithful embedding $\y: s\dm\to s\sh(\dm)$.
\begin{defn}
	Let $X_{\bullet}$ be a simplicial derived manifold, then define $\y X_{\bullet}$ to be the representable simplicial sheaf such that $\y X(U)_n=\map_{\dm} (U,X_n)$.
\end{defn}

Let $X_{\bullet}$, $Y_{\bullet}$ be two simplicial derived manifolds, a map $f: X_{\bullet}\to Y_{\bullet}$ is a {\it Kan fibration}\index{Kan fibration} if the matching map
\begin{equation*}
X_k\to  M^h_{\Lambda^i[k]} X\times_{M^h_{\Lambda^i[k]}Y} Y_k
\end{equation*}
is a cover for all $0\le i \le k$ and $k\ge 1$. If the above matching map are all isomorphisms, we call $f$ a {\it unique Kan fibration}. If the above matching map are covers for all $1\le k \le n$ and are isomorphisms for all $k > n$, then $f$ is an $n$-{\it Kan fibration}. A Kan fibration $f$ is a {\it smooth Kan fibration} if the map restriction to the 0-simplices $f_0: X_0 \to Y_0$ is a cover.

\begin{defn}
    We call $X_{\bullet}$ a {\it derived Lie $\infty$-groupoid}\index{derived Lie $\infty$-groupoid} if the canonical map $ X_{\bullet}\to \ast$ is a Kan fibration.
\end{defn}
Denote the $\infty$-category of derived Lie $\infty$-groupoid by $\ilgpd$.

Let $X_{\bullet}\in \ilgpd$. If the Kan fibration $X_k \to X(\Lambda^i[k])$ is an equivalence for $k> n, 0\le i \le k$, then we say $X_{\bullet}$ is a {\it derived Lie $n$-groupoid}.

\begin{rem}
    This definition of derived Lie $\infty$-groupoids roughly corresponds to homotopy hypergroupoids in pseudo-model categories in \cite{Pri13}. \cite{BG17} defined geometric $\infty$-category in a descent category. \cite{RZ20} defined Lie $\infty$-groupoids in Banach manifolds in a similar fashion. 
\end{rem}

Similarly, we can define a geometric $\infty$-category in $\dm$.
\begin{defn}
    Let $X_{\bullet}$ be a simplicial derived manifold. $X$ is a {\it derived Lie $\infty$-category}\index{{\it derived Lie $\infty$-category}} if for each $0<i<k$ and $k\ge 1$.
    $$ X_k\to Y_k \times_{Y(\Lambda^i[k])} X(\Lambda^i[k]) $$ is a cover. If the above matching map are all isomorphisms $k > n$, then $X_{\bt}$ is a {\it derived Lie $n$-category}.
\end{defn}

A map between derived Lie $\infty$-groupoids $f: X_{\bullet} \to Y_{\bullet}$ is called a {\it hypercover} if 
$$X_k \to Y_k \times_{Y(\del \Delta^i[k])} X(\del \Delta^i[k]) $$
are covers for all $k$. This definition is roughly an acyclic Kan fibration.

Next, we define a homotopy version of descent category in \cite{BG17}.

\begin{rem}
    Note that in our definition, derived Lie $\infty$-category is irrelevant to the derived $\infty$-category of dg-modules which is an enhancement of the classical derived category. From now on, we shall only call Lie $\infty$-category in some specific derived spaces.
\end{rem}

\cite{BG17} defines the {\it descent category}\index{descent category}, which is a weaker notion than pretopology but with finite completeness assumption. For homotopical categories, we have a natural extension of this notion:

\begin{defn}
	Let $\CC$ be a homotopical  category with a subcategory called covers.  We call $\CC$ is a {\it homotopy descent category}\index{descent category!homotopy}, if the following axioms are satisfied:
	\begin{enumerate}
	    \item $\CC$ has finite homotopy limits;
	    \item pullback of a cover is a cover;
	    \item if $f$ is a cover and $gf$ is a cover, then $g$ is a cover.
	\end{enumerate}
\end{defn}

 By lemma \ref{lem: maps}, $\dmfd$ with smooth surjections or {\'e}tale maps  is a homotopy descent category. By lemma, $\dan_k$ with  subjective submersion or {\'e}tale maps is a homotopy descent category. $\dban$ does not satisfy this axiom, hence we want to develop other tools to fix it. 


\subsection{Points in Grothendieck topology} 
We will define a new tool in a Grothendieck topology which allows us to test properties infinitesimally. In a topological space, we consider a sequence of open neighborhoods of a given point $x$, and then watching the behavior of some geometric objects over these neighborhoods. We want to construct a similar notion for a category with a Grothendieck topology.

\begin{defn}[\cite{RZ20}]
	Let $(\CC, \mathcal{T})$ be a category of equipped with a Grothendieck pretopology.. 
	\begin{itemize}
		\item A {\it point}\index{point} is a functor $\p: \sh (\CC) \to \set $ which preserves finite limits and small colimits.
		\item $(\CC, \mathcal{T})$ is said to have {\it enough points}\index{enough points} if there exists a collection of points $\{\p_i\}_{i\in I}$ such that a sheaf morphism $\phi: F\to G$ is an isomorphism if and only if $\p_*(\phi): \p(F)\to \p(G)$ is an isomorphism of sets for all $\p\in \{\p_i\}_{i\in I}$. In this case, we say $\{\p_i\}_{i\in I}$ is {\it jointly conservative}\index{jointly conservative} with respect to $(\CC, \mathcal{T})$.
	\end{itemize}
\end{defn}
Points were originally introduced in topos theory as an adjuntion
$$
x : \set \stackrel{\overset{x^*}{\leftarrow}}{\underset{x_*}{\to}} \CC
$$
between the base topos $\set$ to $\CC$. Given an object $c\in \CC$, we call $x^*(c)$ the {\it stalk}\index{stalk}. We won't need this level of generality, so we stick with the notion in \cite{RZ20} which is sufficient for our construction.

Similarly, we define points in  a homotopical category with an $\io$-Grothendieck pretopology.

\begin{defn}[\cite{RZ20}]
	Let $(\CC, \mathcal{T})$ be a homotopical category with an $\io$-Grothendieck pretopology.
	\begin{itemize}
		\item A {\it point} is a functor $\p: \sh (\CC) \to \set $ which preserves finite homotopy limits and small homotopy colimits.
		\item $(\CC, \mathcal{T})$ is said to have {\it enough points} if there exists a collection of points $\{\p_i\}_{i\in I}$ such that a sheaf morphism $\phi: F\to G$ is an isomorphism if and only if $\p_*(\phi): \p(F)\to \p(G)$ is an isomorphism of sets for all $\p\in \{\p_i\}_{i\in I}$. In this case, we say $\{\p_i\}_{i\in I}$ is {\it jointly conservative} with respect to $(\CC, \mathcal{T})$.
	\end{itemize}
\end{defn}
For simplicity, we will usually refer $\io$-Grothendieck pretopology simply as Grothendieck pretopology when the underlying category is a homotopical category and there is no other confusion.

    Let $X_{\bt}$ be a simplicial object in $(\CC, \mathcal{T})$. Given a point $\p: \sh(\CC)\to \set$, there is a natural extension of $\p$ to a map $\sh(s\CC)\to s\set$ by
    \begin{equation*}
        \p X_n = \p(\y X_n)
    \end{equation*}
    and all structure maps are images under the Yoneda embedding.

Consider $(\dmfd, \mathcal{T}_{ss} )$ the category of derived manifold equipped with smooth surjection pretopology and $(\dban, \mathcal{T}_{ss} )$ with surjection submersion pretopology. We want to show that both these sites have enough points.

\begin{prop}
    Let $G_{\bt}$ be a derived Lie $\infty$-groupoid in $(\dmfd, \mathcal{T}_{ss} )$ or $(\dban, \mathcal{T}_{ss} )$, and $K_{\bt}$ a finitely generated simplicial set. Then there exists a unique natural isomorphism
    $$
        \p\Hom(K, X)\simeq \Hom_{s\set}(K, \p G_{\bt})
    $$
    for each $\p$.
\end{prop}
\begin{proof}
    First, we have $\p \Hom(K_{\bt}, X_{\bt}) \simeq \p\Hom(K_{\bt}, \y X_{\bt})$. Since $K_{\bt}$ is finitely generated, $\Hom(K_{\bt}, \y X_{\bt})$ is a finite limit. Using the fact that $\p$ preserves finite limits, we get directly $\p\Hom(K_{\bt}, \y X_{\bt}) \simeq \Hom_{s\set}(K_{\bt}, \p X_{\bt})$.
\end{proof}

\begin{defn}
    Consider $(\dm, \mathcal{T})$ be either $(\dmfd, \mathcal{T}_{ss} )$ or $(\dban, \mathcal{T}_{ss} )$. Let $F,G\in \sh(\dm)$. We say a sheaf morphism $\phi: F \to G$ is a {\it local surjection} \index{local surjection} if, given any object $X\in \dm$ and $y\in G(X)$, there exists a cover $f: U\to X$ such that $f^*y$ lies in the image of $\phi_U: F(U)\to G(U)$. 
\end{defn}

\begin{defn}
Consider $(\dm, \mathcal{T})$ be either $(\dmfd, \mathcal{T}_{ss} )$ or $(\dban, \mathcal{T}_{ss} )$. Let $\{\p_i\}_{i\in I}$ be a collection of jointly conservative points of $(\dm, \mathcal{T})$. We say $\phi$ is a {\it stalkwise surjection}\index{stalkwise surjection} with respect to $\{\p_i\}_{i\in I}$ if, for all $\p_i$,
$\p_i(\phi): \p_i(F) \to \p_i(G)$ is surjective.
\end{defn}

\begin{prop}
Let $(\dm, \mathcal{T})$ be either $(\dmfd, \mathcal{T}_{ss} )$ or $(\dban, \mathcal{T}_{ss} )$.
We have 
\begin{enumerate}
    \item Local surjections of sheaves on $(\dm, \mathcal{T})$ are epimorphisms.
    \item Epimorphism of sheaves on $(\dm, \mathcal{T})$ are stalkwise surjection with respect to any collection of jointly conservative points.
    \item Let $f$ be a cover, then $\y(f)$ is a stalkwise surjection with respect to any collection of jointly conservative points.
\end{enumerate}

\end{prop}

\begin{proof}
Let $\phi: F \to G$ be a morphism in $\sh(\dm)$.

(1) Suppose $\phi$ is a local surjection. Consider two morphism $\alpha, \beta: G\to H$ in $\sh(\dm)$, such that $\alpha\circ \phi = \beta\circ \phi$, we want to show that   
\end{proof}

\begin{defn}
Let $\{\p_i\}_{i\in I}$ be a collection of jointly conservative points of $(\dm, \mathcal{T})$. 
    A morphism $\psi: X_{\bt} \to Y_{\bt}$  of derived Lie $\infty$-groupoid in $(\dm, \mathcal{T})$ is a {\it stalkwise weak equivalence} \index{stalkwise weak equivalence} if for any $\p_i \in \{\p_i\}_{i\in I}$, the induced map $\p\psi: \p X_{\bt} \to \p Y_{\bt}$ is a weak equivalence of simplicial sets.
\end{defn}
If a morphism of derived Lie $\infty$-groupoids is both a stalkwise Kan fibration and stalkwise weak equivalence, then we call it a {\it stalkwise acyclic fibration}\index{acyclic fibration!stalkwise}.

\begin{prop}
Let $\{\p_i\}_{i\in I}$ be a collection of jointly conservative points of $(\dm, \mathcal{T})$. A morphism $\psi: X_{\bt} \to Y_{\bt}$  of derived Lie $\infty$-groupoids in $(\dm, \mathcal{T})$ is a stalkwise acyclic fibration if and only if
    $$X_k \to Y_k \times_{Y(\del \Delta^i[k])} X(\del \Delta^i[k])$$
    are stalkwise surjections for all $k\ge 0$.
\end{prop}


\begin{cor}
    A hypercover of derived Lie $\infty$-groupoids is both a Kan fibration and a stalkwise weak equivalence.
\end{cor}
\begin{prop}
	$(\dmfd, \mathcal{T}_{ss} )$ has enough points. In fact, let $M\in \dm$ and $x\in M$, define
	\begin{equation}
	 \p_x=\colim_{U^{\aff, open}\subset M} F(U)
	\end{equation} 
	where each $U^{\aff, open}$ is an affine open derived manifold and $U\to M$ is a cover, then $\{\p_x \}$ is a jointly conservative collection of points.
\end{prop}
\begin{proof}
	First, note that $\p_x$ is a filtered colimit for any $x\in M$, hence it preserves finite limits and small colimits. Let $\phi: F \to G$ be a sheaf morphism, and suppose $(\p_x)_*(\phi): \phi_x F\to \phi_x G$ is an isomorphism in $\set$.  
	
	First we show $\phi$ is injective. Let $M\in \dm$ and $f,g\in F(M)$ with $\phi_M(f)=\phi_M(g)\in G(M)$. Note that since $(\p_x)_*(\phi)$ is an injection, we have if $f\in F(U_1), g\in F(U_2)$ and $\p_x(\phi)(\overline{f})=\p_x(\phi)(\overline{g})$ where $x\in U_1\cap U_2$, then there exists $U_{12}\subset U_1\cap U_2$ containing $x$ such that $i^*_1 f=i^*_2g$ where $i_n: U_{12}\to U_n, n=1,2$ are inclusions. Let $x\in M$, then there exists an affine open derived manifold $U_x$, i.e. $U_x=\spec A_x$ for some $A_x\in \cialg$ such that 
	\begin{equation*}
	\coprod_{x\in M} U_x \stackrel{(i_x)}{\longrightarrow} M
	\end{equation*} 
	is a cover, where each $i_x: U_x \to M$ is an inclusion. By pulling back along each $i_x$,
	\begin{equation*}
	\phi_{U_x}(i_x^*f)=\phi_{U_x}(i_x^*g)
	\end{equation*}
    implies that
	\begin{equation*}
	(\p_{x})_*(\phi_{U_x})(\overline{i_x^*f})=(\p_{x})_*(\phi_{U_x})(\overline{i_x^*g})
	\end{equation*}
	From previous observation, for each $U_x$, there exists $i':U_x'\subset U_x$ such that $(i')^*f=(i')^*g$. Now since $F$ is a sheaf and $\coprod_{x\in M} U_x \stackrel{(i_x)}{\longrightarrow} M$ is a cover, we have $f=g$.
	
	Next, we show  $\phi$ is surjective. Let $g\in G(M)$ and the pullback along $i_x$ is $i^*_x(g)\in G(U_x)$. Since $(\p_x)_*(\phi)$ is surjective, there exists a $j_x:U_x'\subset U_x$ and $f_x\in F(U_x')$ such that $\phi(f_x)=j_x^*i_x^* (g) \hookrightarrow M$. Now consider the fiber product
	\begin{equation*}
	\coprod_{x\in M} U_x' \times_M \coprod_{x\in M} U_x'\quad  \overset{p_1}{\underset{p_2}{{\rightrightarrows}}}\quad M
	\end{equation*}
	Observe that
	\begin{equation*}
	\phi(p_1^*(f_x))=p_1^*(\phi(f_x))=p_1^*(j^*_xi_x^*(g) )
	\end{equation*}
	and $p_1^*(j^*_xi_x^*(g) )=p_2^*(j^*_xi_x^*(g) )$ since $G$ is a sheaf. 
	Therefore, 
	\begin{equation*}
	p_2^*(j^*_xi_x^*(g) )=p_2^*(\phi(f_x)=\phi(p_2^*(f_x))
	\end{equation*}
	By the injectivity of $\phi$ from the first part, we have $p_1^*(f_x)=p_2^*(f_x)$. Since $F$ is sheaf, there exists a global section $f\in F(M)$ such that $f\big|_{U_x'}=f_x$. Hence $j_x^*\phi(f)=\phi(f_x)=j^*_xi_x^*(g)$, which implies that $\phi(f)=g$.
	
	Finally, suppose $\phi$ is a sheaf isomorphism, then it is obvious that each $\p_x$ gives isomorphism of sets.
	\end{proof}

\begin{cor}
	$(\dban, \mathcal{T}_{ss} )$ has enough points. In fact, let $M\in \dm$ and $x\in M$, define
	\begin{equation}
	 \p_x=\colim_{U^{ open}\subset M} F(U)
	\end{equation} 
	where each $U^{, open}$ is an open derived Banach manifold and $U\to M$ is a cover, then $\{\p_x \}$ is a jointly conservative collection of points.
\end{cor}
\begin{rem}
	 By our construction, each $\p_x$ is local. hence only depends on affine $\{U_x=\spec A_x\}$ which contains $x$.  
\end{rem}
\begin{proof}
Similar to the case of derived manifold, we just need to replace the affine opens to be open balls.
\end{proof}


\subsection{Locally stalkwise pretopology}

\begin{defn}
Let $(\mathsf{C},\mathcal{T}, \{\p_i\}_{i\in I})$ be a site with enough points. A morphism $F\stackrel{g}{\to} \y Y$ in $\sh(\mathsf{C})$ is a {\it local stalkwise cover}\index{local stalkwise cover} iff there exists an object $X\in \mathsf{C}$ and a stalkwise surjection $\y X \stackrel{f}{\to} F$ such that $g\circ f $ is a cover.
\end{defn}

\begin{defn}[\cite{RZ20}]
The pretopology on a site $(\mathsf{C},\mathcal{T})$ is a {\it locally stalkwise pretopology}\index{locally stalkwise pretoplogy} if it satisfies
\begin{enumerate}
	\item Let $g, f$ be morphisms in $\CC$. If $g\circ f$ is a cover and $\y(f)$ is a stalkwise surjection in $\sh(\mathsf{C})$ with respect to a joint conservative collection of points $\{\p_i\}_{i\in I}$, then $g$ is a cover.
	\item Let $X\stackrel{q}{\to} Y$ and $Z\stackrel{p}{\to} Y$ be two morphisms in $\mathsf{C}$. Suppose $\y (q)$ is a stalkwise surjection with respect to $\{\p_i\}_{i\in I}$ and the base change $X\times_Y Z \stackrel{\tilde{p}}{\to} Z$ is a local stalkwise cover, then $p$ is a cover.  
\end{enumerate}
\end{defn}
The most important property of a locally stalkwise pretopology is that it allows us to characterize Hypercovers by Kan fibrations and stalkwise weak equivalences.
\begin{prop}[\cite{RZ20}]
    Consider a category with pretopology $(\CC, \mathcal{T})$ equipped with a locally stalkwise pretopology with respect to a jointly conservative collection of points $\{\p\}_i$. Let $f: X_{\bt} \to Y_{\bt}$ be a morphism of Lie $\infty$-groupoids in $(\CC, \mathcal{T})$, then the followings are equivalent:
    \begin{enumerate}
        \item $f$ is a Kan fibration and a stalkwise weak equivalence with respect to $\{\p\}_i$.
        \item $f$ is a Kan fibration and a stalkwise weak equivalence with respect to any jointly conservative collection of points of $(\CC, \mathcal{T})$.
        \item $f$ is a hypercover.
    \end{enumerate}
\end{prop}
\begin{proof}
    See \cite[Proposition 6.7]{RZ20}.
\end{proof}

We will show the pretopologies of both  $(\dmfd,\mathcal{T}_{ss})$, $(\dan, \mathcal{T}_{ss})$, and $(\dban,\mathcal{T}_{ss})$ are locally stalkwise.
\begin{lem}
	Let $f:X\to Y$, $g:Y\to Z$ be morphisms in $\dmfd$. Suppose $g\circ f$ is a smooth surjection and $f$ is surjective, then $g$ is also a smooth surjection.
\end{lem}
\begin{proof}
Let $y\in Y$. Denote $z=f(y)$ and $x\in f^{-1}(y)$ any preimage of $y$. Since $g\circ f$ is smooth, there exists a local section $\sigma_{UV}:U\to V$ where $U$ and $V$ are affine open neighborhood of $z$ and $x$, i.e. $\sigma_{UV}(z)=x$ and $(g\circ f)\circ \sigma_{UV} =\id_U$. Let $W=g^{-1} U$, then $\sigma =f\circ \sigma_{UV}: U \to W$ is a local section such that $\sigma(z)=y$. $g$ is clearly surjective. 
	\end{proof}
	\begin{rem}
	    The above proof also works for \'{e}tale topology (i.e. \'{e}tale maps as covers) since we have local lifting property for  \'{e}tale maps as well.
	\end{rem}
\begin{cor}
    Let $f:X\to Y$, $g:Y\to Z$ be morphisms in $\dban$ or $\dan$. Suppose $g\circ f$ is a surjective submersion and $f$ is surjective, then $g$ is also a surjective submersion.
\end{cor}
\begin{proof}
    Similar to previous proof.
\end{proof}

\begin{lem}
	Let $f:X\to Y$, $g:Y\to Z$ be morphisms in $\dmfd$. Suppose $g\circ f$ is a smooth surjection and $f$ is a stalkwise surjection, then $g$ is also a smooth surjection.
\end{lem}
\begin{proof}Since $f$ is stalkwise surjective, it has to be surjective.
\end{proof}
\begin{rem}
    Again, the result still holds if we replace surjective submersion by  \'{e}tale maps.
\end{rem}
\begin{cor}
Let $f:X\to Y$, $g:Y\to Z$ be morphisms in $\dban$ or $\dan$. Suppose $g\circ f$ is a surjective submersion and $f$ is a stalkwise surjection, then $g$ is also a surjective submersion.
\end{cor}
\begin{prop}
	$(\dmfd,\mathcal{T}_{ss})$ and $(\dan,\mathcal{T}_{ss})$ are categories with locally stalkwise pretopology.
\end{prop}

\begin{proof} We will prove the case for  $(\dmfd,\mathcal{T}_{ss})$, and the case for $(\dan,\mathcal{T}_{ss})$ is similar. 

It suffices to verify the second axiom. Pick $z\in Z$ and denote $p(z)=y$. Since $\y(q)$ is a stalkwise surjection, for any $y\in Y$, we can find an affine open neighborhood $O_y$ of $y$ such that there exists a $x\in q^{-1}(y) \subset X$ and an affine open neighborhood $O_y$ of $y$ such that $s_{yx}:O_y \to O_x$ is a local section of $q|_{O_x}$ such that $s_{yx}(y)=x$. Let $O_z=p^{-1}(O_y)$, then we have a pullback diagram
\begin{center}
	\begin{tikzcd}
	O_x\times_{O_y} O_z \arrow[r, "p'"]\arrow[dr, phantom, "\ulcorner", very near start] \arrow[d,"q'"] & O_x \arrow[d, "q"] \\
	O_z \arrow[r,"p"]           & O_y         
	\end{tikzcd}
\end{center}
where $O_x\times_{O_y} O_z= p'^{-1}(O_x)=q'^{-1}O_z$ by construction. Note that both $q|_{O_x}$ and $q'|_{O_x\times_{O_y} O_z}$ are smooth surjections, In particular, we can shrink both of them to make $q$ a projection when restricts to $O_x$. By construction, we can find a $w\in q'^{-1}(z)\subset O_x\times_{O_y} O_z$ such that $p'(w)=x$ and a local section $s_{zw}$ of $q'|_{O_x\times_{O_y} O_z}$ with $s_{zw}(z)=w$. Since $p'$ is a locally stalkwise cover, there exists some $U\in \dmfd$ with a map $f: U \to O_x\times_{O_y} O_z $ such that $ p'\circ f$ is a smooth surjection. Hence, we can find a section $s_{xu}:O_x \to O_u$ where $O_u = ( p'\circ f)^{-1} O_x$. Now $s_{xw}=f\circ s_{xu}: O_x \to O_w$ is the desired section of $p'$. To construct the section of $p$, we just need to take $s_{yz}=q'\circ f\circ s_{xu}\circ s_{yx}$.
	\end{proof}
	
	\begin{rem}
 	    The key of the proof is the 'inverse function theorem' for smooth surjections of derived manifolds. Hence, it also shows that the result hold for \'{e}tale topology and other topology which satisfies the 'inverse' function theorem. For more about inverse function theorem in derived manifolds, see \cite[Proposition 6.2.1]{Nui18}.
	\end{rem}
	
	\begin{prop}
	$(\dban,\mathcal{T}_{ss})$ is a category with locally stalkwise pretopology.
\end{prop}
\begin{proof}
    The case for the category of (ordinary) Banach manifolds is shown in \cite[Proposition 6.12]{RZ20}. The proof for derived case follows mostly from the ordinary case. The only thing different from the above proof is that we need to take care of the representability issues in Banach manifolds. In this paper we do not construct homotopy structures on $\dban$ which will be developed in future work. When we compute fiber product in $\dban$, we compute pushout in the dga's, hence we only need to care about the representability of the degree 0 terms, which follows from the \cite[Lemma 6.9]{RZ20} and \cite[Lemma 6.11]{RZ20}. Again, the core of the proof is the inverse function theorem for submersion.
\end{proof}


\subsection{Collapsible extensions }
In this section, we will study a special class of simplicial maps, which will be used heavily later when we prove some representability results on simplicial sheaves.

\begin{defn}
    Let $T_{\bt}$ be a finitely generated simplicial set and $S_{\bt}$ a simplicial subset.  The inclusion map $\iota: S_{\bt} \to T_{\bt}$ is called a {\it collapsible extension}\index{collapsible extension} if and only if it can decomposed as a sequence of inclusion  maps
    \begin{equation*}
        S_{\bullet} = S^0_{\bt} \hookrightarrow
        S^1_{\bt} \hookrightarrow \cdots \hookrightarrow S^l_{\bt} = T_{\bt}
    \end{equation*}
    i.e. for each $i$, $S^i = S^{i-1}\sqcup_{\Lambda^j[m]} \Delta[m]$ for some horn $\Lambda^j[m]$ and $m>0$. If $T_{\bt}$ is a collapsible extension of a point, we say it is {\it collapsible}.
\end{defn}
So roughly speaking, collapsible extension is a sequence of filling some horns. We can also define similar maps which fill in boundaries.

\begin{defn}
    Let $T_{\bt}$ be a finitely generated simplicial set and $S_{\bt}$ a simplicial subset.  The inclusion map $\iota: S_{\bt} \to T_{\bt}$ is called a {\it boundary extension}\index{boundary extension} if and only if it can decomposed as a sequence of inclusion  maps
    \begin{equation*}
        S_{\bullet} = S^0_{\bt} \hookrightarrow
        S^1_{\bt} \hookrightarrow \cdots \hookrightarrow S^l_{\bt} = T_{\bt}
    \end{equation*}
    i.e. for each $i$, $S^i = S^{i-1}\sqcup_{\del \Delta[m]} \Delta[m]$ for some horn $\Lambda^j[m]$ and $m>0$. 
\end{defn}

An obvious result is
\begin{lem}[]
    The inclusion of any face $\Delta[k] \to \Delta[n]$ is a collapsible extension for $0 \le k \le n$.
\end{lem}
\begin{proof}
    See \cite[Lemma2.44]{Li15}.
\end{proof}
Next, we will see how collapsible extension relate to representablity of Lie $\infty$-groupoids.

\begin{lem}\label{c1}
    Let $X_{\bt}$ be a Lie $\infty$-groupoid. Suppose $S_{\bt}\hookrightarrow T_{\bt}$ is a collapsible extension. If $\Hom(S_{\bt}, X_{\bt})$ is representable, then $\Hom(T_{\bt}, X_{\bt})$ is also representable, and the induced map 
    $$
    \Hom(T_{\bt}, X_{\bt}) \to \Hom(S_{\bt}, X_{\bt})
    $$
    is a cover.
\end{lem}
\begin{proof} This is \cite[Lemma 3.7]{RZ20}. 
    Let $S_{\bullet} = S^0_{\bullet} \hookrightarrow
        S^1_{\bullet} \hookrightarrow \cdots \hookrightarrow S^l_{\bullet} = T_{\bt}$ be the collapsible extension. Since covers are closed under composition, we can just restrict to the case of one inclusion. Let $T_{\bt} = S_{\bt} \sqcup_{\Lambda^j[m]} \Delta[m]$. Applying $\Hom(-, X_{\bt})$
        \begin{equation*}
            \Hom(T_{\bt}, X_{\bt}) = \Hom(S_{\bt}, X_{\bt}) \times_{\Hom(\Lambda^j[m], X_{\bt})}  \Hom(\Delta[n], X_{\bt})
        \end{equation*}
        Since $X$ is a Lie $\infty$-groupoid, $\Hom(\Delta[n], X_{\bt})\to \Hom(\Lambda^j[m], X_{\bt})$ is a cover between representable sheaves. Therefore, by axioms of pretopology, we get $\Hom(T_{\bt}, X_{\bt})$ is representable and $
    \Hom(T_{\bt}, X_{\bt}) \to \Hom(S_{\bt}, X_{\bt})$ is a cover.
\end{proof}
\begin{rem}
    For $X_{\bt}$ being a Lie $n$-groupoid, and $T_{\bt} = S_{\bt} \sqcup_{\Lambda^j[m]} \Delta[m]$ with $m > n, 0\le j \le m$, and suppose $\Hom(S_{\bt}, X_{\bt})$ is representable, then $
    \Hom(T_{\bt}, X_{\bt}) \to \Hom(S_{\bt}, X_{\bt})$  is actually an isomorphism.
\end{rem}
Next, we consider the representability of sheaves.
\begin{lem}
\label{c2}
    Let $S\subset \Delta[n]$ be a collapsible simplicial subset, $X_{\bt}$ a simplicial manifold, and $Y_{\bt}$ a Lie $\infty$-groupoid. If $f: X_{\bt}\to Y_{\bt}$ is a morphism satisfies $\kan(m, j)$ for $m < k$ and $0 \le j \le m$, then the sheaf on $\dm$ $\Hom(S_{\bt} \hookrightarrow \Delta[k], X_{\bt} \stackrel{f}{\to} Y_{\bt})$ is representable.
\end{lem}
\begin{proof}
    This is \cite[Lemma 3.9]{RZ20}. Consider $S^0 = * \hookrightarrow \Delta[k]$. Note that $\Hom(* \hookrightarrow \Delta[k], X_{\bt} \stackrel{f}{\to} Y_{\bt})$ is represented by $X_0 \times_{Y_0} Y_k$. By previous two lemmas, we see $Y_k \to Y_0$ is a cover.
\end{proof}
\begin{cor}
    For same assumption as above, $\Hom(S_{\bt} \hookrightarrow \Delta[k], X_{\bt} \stackrel{f}{\to} Y_{\bt})$ is representable.
\end{cor}
\begin{proof}
    Applying previous lemma to the horn $\Lambda^j[m]$ which is collapsible for all $j$'s.
\end{proof}

\begin{cor}
    Suppose $X_{\bt}\to *$ satisfies $\kan(m, j)$ for $1 \le m <k$, then $\Hom(\Lambda^j[m], X_{\bt})$ is representable.
\end{cor}

\begin{lem}
    Let $f: \Lambda^i[1] \to \Delta[1]$ for $i= 0, 1$ be the standard inclusion. If $\iota: S_{\bt}\to T_{\bt}$ is a boundary extension, then the induced map
    $$
    (S_{\bt} \otimes \Delta[1]) \sqcup_{S_{\bt} \otimes \Lambda^i[1] }(T_{\bt} \otimes \Lambda^i[1]) \to (T_{\bt} \otimes \Delta[1])
    $$
    is a collapsible extension.
\end{lem}
\begin{proof}
    See \cite[Lemma 3.3.3]{Hov07} for the case of $\iota$ being the standard inclusion $\del \Delta[n] \to \Delta[n]$. Suppose $F:s\set \times s\set \to s\set$  is a co-continues functor and 
    
\begin{center}
	\begin{tikzcd}
	X_{\bullet} \arrow[r, ""] \arrow[d, ""] & Y_{\bullet}
	\arrow[d, ""] \\
	S_{\bullet} \arrow[r, ""] \arrow[r, ""]      &  T_{\bullet}   \arrow[ul, phantom, "\ulcorner", very near start]              
	\end{tikzcd}
\end{center}
is a pushout square of simplicial set, then applying the same technique in \cite[Lemma 2.42]{Li15}, we have
\begin{center}
	\begin{tikzcd}
	 F(X_{\bt}, \Delta[1]) \sqcup_{F(X_{\bt}, \Lambda^i[1]) }F(Y_{\bt}, \Lambda^i[1]) \arrow[r, ""] \arrow[d, ""] & F(Y_{\bt}, \Delta[1])
	\arrow[d, ""] \\
	F(S_{\bt}, \Delta[1]) \sqcup_{F(S_{\bt}, \Lambda^i[1]) }F(T_{\bt}, \Lambda^i[1]) \arrow[r, ""] \arrow[r, ""]      &  F(T_{\bt}, \Delta[1]) \arrow[ul, phantom, "\ulcorner", very near start]              
	\end{tikzcd}
\end{center}
is also a pushout square. Now take $F$ to be the product and proceed by induction.
\end{proof}

\begin{rem} Since collapsible extensions are boundary extensions, 
    replacing the boundary extension assumption in the previous lemma by collapsible extension, the result still holds.
\end{rem}
\begin{lem}
    The inclusion $\Lambda^j[n]\times \Delta[1] \to \Delta[n]\times \Delta[1]$ is a collapsible extension.
\end{lem}
\begin{proof}
    Regard $\Lambda^j[n]\times \Delta[1] \to \Delta[n]\times \Delta[1]$ as a composition
    \begin{equation*}
        \Lambda^j[n]\times \Delta[1] \to (\Lambda^j[n] \times \Delta[1]) \sqcup_{\Lambda^j[n] \times \Lambda^i[1] }(\Delta[n] \times \Lambda^i[1])\to \Delta[n]\times \Delta[1]
    \end{equation*}
    It is clear that the first map is collapsible. The second map is also collapsible by the previous lemma.
\end{proof}

\begin{lem}\label{c5}
    The inclusion 
    $$(\Lambda^j[n] \times \Delta[1]) \sqcup_{\Lambda^j[n] \times \del\Delta[1] }(\Delta[n] \times \del\Delta[1])\to \Delta[n]\times \Delta[1]$$
    is a collapsible extension.
\end{lem}
\begin{proof}
    See \cite[Appendix A]{RZ20}.
\end{proof}
\section{Homotopy theory of derived Lie $\infty$-groupoids}

\subsection{Category of fibrant objects}
Category of fibrant objects (CFO), also known as Brown category, is a weaker notion of a Quillen model category which still allow us to perform many operations in homotopy theory.

\begin{defn}[\cite{Bro73},\cite{BG17}, \cite{RZ20}]

Let $\CC$ be a small category, we say that $\CC$ is a {\it category of fibrant objects}\index{category of fibrant object} (CFO) such that there exists two distinguished subcategories $\mathcal{W}$ and $\mathcal{F}$ called {\it weak equivalences}\index{weak equivalence!in a category of fibrant objects} and {\it fibrations}\index{fibration!in a category of fibrant objects} respectively, and it satisfies the following conditions

\begin{enumerate}
	\item $\CC$ has all finite products, and in particular a terminal object $\ast$.
	\item Pullback of a fibration along arbitrary morphisms exist, and it is also a fibration.
	
	\item The morphisms which sit in both $\mathcal{W}$ and $\mathcal{F}$ are call {\it acyclic fibrations}\index{acyclic fibrations!in a category of fibrant objects}. The Pullbacks of acyclic fibrations are acyclic fibrations.
	\item  Weak equivalences satisfy 2-out-of-3, and contain all isomorphisms.
	\item  Composition of fibrations are fibrations, and all isomorphisms are fibrations.
	\item Given any object $B$, there exists a {\it path object}\index{path object} $B^{\Delta[1]}$ that fits into the diagram
	\begin{equation*}
	B\stackrel{\sigma}{\longrightarrow} B^{\Delta[1]} \stackrel{(d_0, d_1)}{\longrightarrow} B \times B
	\end{equation*}
	where $\sigma$ is a weak equivalence and $(d_0, d_1)$ is a fibration, and the composition $B\to B\times B$ is the diagonal map. 
	\item For any objects $B$, the canonical map $B \to \ast$ is a fibration, i.e. all objects are {\it fibrant}.
\end{enumerate}

\end{defn}
\begin{example}
    Fibrant model categories are trivial examples of categories of fibrant objects, for example:
    \begin{itemize}
        \item $\Top$ with the  Quillen model structure.
        \item $\ch_k^{\ge 0}$ with projective model structures.
        \item $\Mod_A^{\ge 0}$ with projective model structures.
    \end{itemize}
\end{example}
\begin{example}
    The next simple examples are restriction of model categories to their fibrant objects, for example
    \begin{itemize}
        \item The subcategory of $s\set$ consisting of Kan complexes, which we call the category of $\infty$-groupoids $\igpd$.
    \end{itemize}
\end{example}

\begin{example}[Simplicial sheaves]
    Let $(\CC, \mathcal{T})$ be a site with enough points. For example, take $\CC = \open(X)$ the category of open subsets for a topological space $X$. Then the category of simplicial sheaves on $\CC$  whose stalks are Kan complexes form a category of fibrant objects. Hence, this gives a model for the homotopy category of $\infty$-stacks over $\CC$. This is a motivating example in \cite{Bro73} to introduce categories of fibrant objects. 
\end{example}

\begin{example}[$C^*$-algebras]
    Let $\cstar$ be the category of $C^*$-algebras. \cite{Sch84} construct a category of fibrant objects structure on  $\cstar$ as follows. 
    
    Denote $\pi_0 \cstar$ the ordinary homotopy category of $\cstar$, i.e. the same objects as $\cstar$ with homotopy classes of maps in $\cstar$. We say a map $f:A\to B$ is a {\it homotopy equivalence}\index{homotopy equivalence!$C^*$-algebras} if $\pi_0(f)$ is invertible in $\pi_0 \cstar$. A map $f:A\to B$ is called a {\it a Schochet fibration}\index{Schochet fibration} if its induced map
    $$
    f_*:\Hom_{\cstar}(C, A) \to \Hom_{\cstar}(C, B) 
    $$
    has the path lifting property for all $C\in \cstar$. 
    
     $\cstar$ with homotopy equivalences as weak equivalences and Schochet fibrations as fibrations is a category of fibrant objects. 
     
     \cite{UUY} construct another category of fibrant objects structure the category of separable $C^{*}$-algebras $\cstar^{\sep}$ forms a category of fibrant objects with weak equivalences the $KK$-equivalences and fibrations the Schochet fibrations, whose homotopy category $\ho(\cstar^{\sep})$ is equivalent to the $KK$-category of Kasparov \cite{Kas07}. This implies that Kasparov’s $KK$-category is a stable triangulated category.
\end{example}

\begin{example}[Behrend-Liao-Xu derived manifolds]
    \cite{BLX21} develops a theory of derived manifolds using bundles of curved $\li[1]$-algebras. They construct a category of fibrant objects on their category of derived manifolds as follows:
    \begin{itemize}
        \item A morphism is a weak equivalence \index{weak equivalence!BLX derived manifolds} if:
        \begin{enumerate}
            \item It induces a bijection on classical loci.
            \item Its linear part induces a quasi-isomorphism on tangent complexes at all classical points.
        \end{enumerate}
        \item A morphism is a fibration \index{fibration!BLX derived manifolds}if:
            \begin{enumerate}
                \item  The underlying morphism of manifolds is a submersion,
                \item The linear part of the morphism of $\li[1]$-algebras is levelwise surjective.
            \end{enumerate}
    \end{itemize}
\end{example}
Sometimes we want to deal with categories which do not contain all finite limits, but we still want to do homotopy theory on it. It turns out that we can loosen the limits criteria sometimes, and consider {\it incomplete category of fibrant objects} (iCFO), where we do not assume all pullbacks of fibrations exists, and only for those pullbacks exist, the pullbacks are still fibrations. In summary,
\begin{defn}[\cite{RZ20}]
    We say a category $\CC$ is an {\it incomplete category of fibrant objects} (iCFO), if it satisfies the conditions (3)-(7) of categories of fibrant objects, and we replace (2) by
\begin{itemize}
    \item If the pullback of a fibration exists, then it is a fibration.
\end{itemize}
\end{defn}

\begin{example}
    As a prototypical example, \cite{RZ20} shows that Lie $\infty$-groupoids in $\ban$ with surjective submersion pretopology is an incomplete category of fibrant objects.
\end{example}
\subsubsection{Homotopical algebra for categories of fibrant objects}
CFO allows us to perform explicit homotopical operations, for example, compute $\io$-limits explicitly.

Recall that in a homotopical category, the homotopy pullback of two maps $\f: A\to C$, $g: B\to C$ are defined as a universal object $X$ such that the diagram
\begin{center}
    \begin{tikzcd}
X \arrow[r] \arrow[d] & {A} \arrow[d, "f"] \\
{B} \arrow[r, "g"]           & {C}          
\end{tikzcd}
\end{center}
commutes up to homotopy. Thanks to the existence of path object, we can compute homotopy pullbacks explicitly and easily in categories of fibrant objects.

\begin{thm}[\cite{Bro73}]\label{hfiber}
Let $\CC$ be a category of fibrant object, then the {\it homotopy pullback}\index{homotopy pullback!category of fibrant objects}(or {\it homotopy fiber product} $A\times^h_C\times B$ of two maps $\f: A\to C$, $g: B\to C$ is presented by $A\times_C C^I \times_C B$, i.e. the ordinary limit of
    \begin{center}
        \begin{tikzcd}
{A\times_C C^I \times_C B} \arrow[rr] \arrow[dd] &                        & {A} \arrow[d, "f"] \\
                         & {C^I} \arrow[d, "d_0"] \arrow[r, "d_1"] & {C}           \\
{B} \arrow[r, "g"]             & {C}                     &             
\end{tikzcd}
    \end{center}
    Moreover, the projection map $\pi: A\times_C C^I \times_C B \to A$ a fibration. If in addition $v:B\to C$ is a weak equivalence, then $\pi$ is an acyclic fibration.
\end{thm}

Another useful property of CFO is that we have a nice simplification of homotopy mapping space. We follow the construction in \cite[Section 3.6.3]{NSS12}.
\begin{defn}
    Let $\CC$ be a category of fibrant objects. Let $X, Y \in \CC$ be two objects. Define a category $\cocycle(X,Y)$ by
   \begin{enumerate}
       \item Objects are spans, i.e. diagrams of the following form
       $$X \stackrel{\simeq}{\longleftarrow}
  A
  \stackrel{}{\longrightarrow}
  Y$$ where the left morphism is an acyclic fibration.
  \item Morphisms are given by commutative diagrams of the following form
  \begin{center}
      \begin{tikzcd}
  & A_1 \arrow[ld, "\simeq", two heads] \arrow[dd, "f"] \arrow[rd, "g_1"] &   \\
X &                                                                  & Y \\
  & A_2 \arrow[lu, "\simeq", two heads] \arrow[ru, "g_2"]            &  
\end{tikzcd}
  \end{center}
  Note that the map $f:A_1\to A_2$ is necessarily a weak equivalence by 2-out-of-3. 
   \end{enumerate}
\end{defn}

\begin{thm}[\cite{NSS12}]
    Let $\CC$ be a category of fibrant objects. Let $X, Y \in \CC$ be two objects. Given any objects $X, Y\in \CC$, the canonical inclusions
    \begin{equation*}
        \mathcal{N} \cocycle(X, Y)\to \mathcal{N} \wcocycle(X, Y) \to L^H\CC(X,Y)
    \end{equation*}
    weak equivalences, where $L^H \CC$ is the hammock localization of $\CC$.
\end{thm}
Hence, we can compute the homotopy mapping spaces (or derived Hom space in \cite{NSS12}) of a category of fibrant objects simply using its category of spans $\cocycle(\CC)$ or $\wcocycle(\CC)$. As an easy consequence of this theorem, the homotopy fiber product we get in Theorem \ref{hfiber} presents the correct $\io$-limit.

\subsubsection{Fibrations in derived Lie $\infty$-groupoids}
In this section, we will prove some basic properties of Kan fibrations and hypercovers of derived Lie $\infty$-groupoids in various categories. For simplicity, we use $\dm$ to denote either $\dmfd$, $\dan$, or $\dban$. The proofs will work for any of these categories unless specified explicitly.
\begin{prop}
Let $f: X_{\bullet}\to Y_{\bullet}$, $g: Y_{\bullet}\to Z_{\bullet}$ be Kan fibrations between derived Lie $\infty$-groupoids in $\dm$, then $g\circ f $ is also a Kan fibration.\label{fib}
\end{prop}
\begin{proof}
    We want to show the induced map $X_{k}\to M^h_{\Lambda^i[k]} X\times_{M^h_{\Lambda^i[k]}Z} Z_k $ is a cover. We have the following commutative diagram.
    
    \begin{center}
        \begin{tikzcd}
{M^h_{\Lambda^i[k]} X\times_{M^h_{\Lambda^i[k]}Y} Y_k} \arrow[r, "g_*"] \arrow[d, "pr_2"] & {M^h_{\Lambda^i[k]} X\times_{M^h_{\Lambda^i[k]}Z} Z_k} \arrow[d, "f_*"] \arrow[r] \arrow[r, "pr_1"] & {M^h_{\Lambda^i[k]} X} \arrow[d, "f_*"] \\
{Y_k} \arrow[r, "\xi"]           & {M^h_{\Lambda^i[k]} Y \times_{M^h_{\Lambda^i[k]}Z} Z_k} \arrow[d, "\psi"] \arrow[r, "pr_1"]           & {M^h_{\Lambda^i[k]} Y} \arrow[d, "g_*"] \\
                       & {Z_k} \arrow[r, "\iota_*"]                     & {M^h_{\Lambda^i[k]} Z}          
\end{tikzcd}
    \end{center}
    
    The bottom square and the composition of bottom and middle squares are pullbacks, hence the middle one is as well. The composition of left and middle squares are pullbacks, hence the left square is also a pullback. Therefore, $g^*$ is a cover since $\xi$ is. Hence, the composition $X_k \to {M^h_{\Lambda^i[k]} X\times_{M^h_{\Lambda^i[k]}Y} Y_k} \to {M^h_{\Lambda^i[k]} X\times_{M^h_{\Lambda^i[k]}Z} Z_k}$ is a cover.
\end{proof}

\begin{cor}
    The composition of $n$-Kan fibrations are $n$-Kan fibrations. In particular, composition of unique Kan fibrations are unique Kan fibrations.
\end{cor}

\begin{lem}
    Hypercovers between derived Lie $\infty$-groupoids in $\dm$ are Kan fibrations.
\end{lem}
\begin{proof}
    Apply the canonical inclusion $\Lambda^i[k] \to \del \Delta[k]$.
\end{proof}

\begin{lem}
    Let $f: X_{\bullet}\to Y_{\bullet}$, $g: Y_{\bullet}\to Z_{\bullet}$ be hypercovers between derived Lie $\infty$-groupoids in $\dm$, then $g\circ f $ is also a hypercover.
\end{lem}
\begin{proof}
    Similar to the case of Kan fibrations by replacing homotopy matching space of $\Lambda^{i}[k]$ by homotopy matching space of $\del \Delta[n]$.
\end{proof}

\begin{prop}\label{pb1}
	Let $f:X_{\bullet}\to Y_{\bullet}$ be a Kan fibration between derived Lie $\infty$-groupoids in $\dmfd$ or $\dan$. Then the  pullback of $f$ along any morphisms
	$g: Z_{\bullet}\to Y_{\bullet}$ exists and $h: X_{\bullet} \times_{Y_{\bullet}} Z_{\bullet} \to Z_{\bullet}$ is a Kan fibration.
\end{prop}

\begin{proof}
The pullback has $n$-simplices $X_n \times^h_{Y_n} Z_n \in \dm$, hence is a simplicial derived manifold. We have the following commutative diagram
    \begin{center}
        \begin{tikzcd}
                                     &  & {X \times^h_{Y} Z} \arrow[rd] \arrow[d] &              \\
{\Lambda^i[n]} \arrow[rru] \arrow[rrr] \arrow[d] &  & {Z} \arrow[rd]           & {X} \arrow[d] \\
{\Delta[n]} \arrow[rrr] \arrow[rru]           &  &                         & {Y}          
\end{tikzcd}
    \end{center}
Hence we have a pullback diagram
\begin{center}
    \begin{tikzcd}
{X_n \times^h_{Y_n} Z_n} \arrow[r] \arrow[d, "\psi"] & {X_n} \arrow[d, "\xi"] \\
{ \big(M^h_{\Lambda^i[n]} {(X \times^h_{Y} Z)}\big) \times^h_{M^h_{\Lambda^i[n]}Z} Z_n} \arrow[r]           & {M^h_{\Lambda^i[n]}X} \times^h_{M^h_{\Lambda^i[n]}Y} Y_n          
\end{tikzcd}
\end{center}
Hence $\psi$ is a cover since $\xi$ is. Hence $h$ is a fibration. Note that 
$$
 \big(M^h_{\Lambda^i[n]} (X \times^h_{Y} Z)\big) \times^h_{M^h_{\Lambda^i[n]}Z} Z_n \to \big(M^h_{\Lambda^i[n]} {(X \times^h_{Y} Z)}\big)
$$ is a cover since it  is the pullback of $Z_n \to M^h_{\Lambda^i[n]}Z $  which is a cover.
    
\end{proof}
The above proposition generalizes to the case where $X_{\bullet}, Y_{\bullet}, Z_{\bullet}$ are Lie $n$-groupoids. It is easy to show that in this case the fiber product is also a Lie $n$-groupoid.

\begin{rem}
    Clearly this won't work for derived Banach manifolds due to lacking of limits. We will prove later that once pullback of a fibration exists, then it is a fibration, which is a key component in the iCFO structure on  derived Lie $\infty$-groupoids in $\dban$. 
\end{rem}
\begin{prop}\label{pb2}
	Let $f:X_{\bullet}\to Y_{\bullet}$ be a hypercover between derived Lie $\infty$-groupoids in $\dmfd$ or $\dan$. Then the pullback of $f$ along any morphism 
	$g: Z_{\bullet}\to Y_{\bullet}$ exists and $h: X_{\bullet} \times_{Y_{\bullet}} Z_{\bullet} \to Z_{\bullet}$ is a hypercover.
\end{prop}
\begin{proof}
    By similar argument as above, we have a pullback diagram
\begin{center}
    \begin{tikzcd}
{X_n \times^h_{Y_n} Z_n} \arrow[r] \arrow[d, "\psi"] & {X_n} \arrow[d, "\xi"] \\
 \big(M^h_{\Delta[n]} (X \times^h_{Y} Z)\big) \times^h_{M^h_{\Delta[n]}Z} Z_n \arrow[r]           & M^h_{\Delta[n]}X \times^h_{M^h_{\Delta[n]}Y} Y_n          
\end{tikzcd}
\end{center}
Hence $\psi$ is a cover and $h$ is then a hypercover.

\end{proof}

Next, we will show simplicial derived manifolds and simplicial derived $k$-analytic spaces also form a homotopy descent category.

\begin{prop}
Both $s\dmfd$ and $\dan$ are homotopy descent categories with hypercovers as covers.
\end{prop}
We have shown pullbacks of hypercovers are hypercovers, we just need to verify the last criteria.

\begin{lem}
Let $f: X_{\bullet}\to Y_{\bullet}$, $g: Y_{\bullet}\to Z_{\bullet}$ be morphisms of  $s\dmfd$ or $s\dan$. Suppose $f$ and $g\circ f $ hypercovers, then $g$  is also hypercover.
\end{lem}

\begin{proof}
We have the following commutative diagram
    \begin{center}
        \begin{tikzcd}
{X_n} \arrow[r, "\alpha"] \arrow[rd, "\beta"] & {M^h_{\Lambda^i[n] X }}\times^h_{M^h_{\Lambda^i[n]  }Y} Y_n \arrow[r] \arrow[d] & {Y_n} \arrow[d, "\gamma"] \\
                        & {M^h_{\Lambda^i[n] X }}\times^h_{M^h_{\Lambda^i[n]  }Z} Z_n \arrow[r]           & {M^h_{\Lambda^i[n] Y }}\times^h_{M^h_{\Lambda^i[n]  }Z} Z         
\end{tikzcd}
    \end{center}
    
\end{proof}

\subsection{Path object}
\begin{prop}\label{path}
	Let $X_{\bullet}$ be a derived Lie $\infty$-groupoid in $\dm$, then there exists a path object $X^{\Delta[1]}_{\bullet}$, that is, we have a factorization
	\begin{equation*}
	X_{\bullet}\stackrel{s^{*}_0}{\longrightarrow} X^{\Delta[1]}_{\bullet} \stackrel{(d^{*}_0, d^*_1)}{\longrightarrow} X_{\bullet} \times X_{\bullet}
	\end{equation*}
	which is a factorization of the diagonal map $X_{\bullet}\to X_{\bullet}\times X_{\bullet}$ into a stalkwise weak equivalence $s^{*}_0$ followed by a Kan fibration $(d^{*}_0, d^*_1)$.
\end{prop}

First, we want to look at the sheaf level. 
\begin{lem}
	Let $X_{\bullet}$ be a derived Lie $\infty$-groupoid in $\dm$, then we have a factorization
	\begin{equation*}
	\BR \y X_{\bullet}\stackrel{s^{*}_0}{\longrightarrow} (\BR\y X_{\bullet})^{\Delta[1]} \stackrel{(d^{*}_0, d^*_1)}{\longrightarrow} \BR\y X_{\bullet} \times \BR\y X_{\bullet}
	\end{equation*}
	where $s^*_0$ is a stalkwise weak equivalence and $(d^{*}_0, d^*_1)$ is a stalkwise Kan fibration.
\end{lem}
\begin{proof}
	Note that for any $\p\in \po$, $\p X$ is a Kan complex in $s\set$. Applying $\p$ to the previous diagram we have
	\begin{equation*}
	\p X_{\bullet}\stackrel{\p (s^{*}_0)}{\longrightarrow} (\p X_{\bullet})^{\Delta[1]} \stackrel{\p(d^{*}_0, d^*_1)}{\longrightarrow} \p X_{\bullet} \times \p X_{\bullet}
	\end{equation*}
	where $$\big(\p  X_{\bullet}\big)_n^{\Delta[1]}=\big(\p (\y X_n)\big)^{\Delta[1]}= \Hom_{s\set} (\Delta[n] \times \Delta[1], \p  X_{\bullet})$$
	which is the path object in $s\set$ for $\p X_{\bullet}$.
	\end{proof}

\begin{lem}
	 $(\y X_{\bullet})^{\Delta[1]}\in s\sh(\dm)$ is a representable simplicial presheaf which is represented by a derived Lie $\infty$-groupoid $X^{\Delta[1]}_{\bullet}$.
\end{lem}
\begin{proof}First note that by the simplicial structure of simplicial sheaves
	\begin{align*}
	(\mathbf{R}\y X_{\bullet})^{\Delta[1]} (U)_n&\simeq\Hom_{s\set}(\Delta[n]\times \Delta[1], \BR \y X_{\bullet}(U) )\\
	&\simeq \Hom_{s\dm}\big((\Delta[n]\times \Delta[1])\otimes U, X_{\bullet} \big)\\
	&\simeq \Hom_{s\dm}\big((U, \rhom_{s\set} (\Delta[n]\times \Delta[1],X_{\bullet} )\big)
	\end{align*}
	In order for $(\BR\y X_{\bullet})^{\Delta[1]}$ to be representable, we to show that $\rhom_{s\set} (\Delta[n]\times \Delta[1],X_{\bullet} )$ is a derived manifold.
	
	Note that $\Delta[n]\times \Delta[1]$ has a canonical decomposition into $(n+1)$ $n+1$-simplices $\Delta[n]$, hence we have 
	$$  \BR\Hom (\Delta[n]\times \Delta[1],X_{\bullet} ) \simeq X_{n+1}   \times_{d_1,d_1}^h X_{n+1}\times_{d_2,d_2}^h\cdots \times_{d_n,d_n}^h X_{n+1}$$ which is a derived manifold. Hence, $(\BR \y X_{\bullet})^{\Delta[1]}$ is represented by a simplicial derived manifold,  and we denote it by $X_{\bullet}^{\Delta[1]}$.
	
	Next, we want to prove that $X_{\bullet}^{\Delta[1]}$ is a derived Lie $\infty$-groupoid, i.e. $X_{\bullet}^{\Delta[1]}(\Delta[n]) \to  X_{\bullet}^{\Delta[1]}(\Lambda^i[n])$ is a cover.
	
	First we have 
	$$M^h_{\Lambda^i[n]} X_{\bullet}^{\Delta[1]}=\rhom_{s\set}(\Lambda^i[n], X_{\bullet}^{\Delta[1]}) $$
	which is characterized by the sheaf $$U\mapsto \Hom_{s\dm}(\Lambda^i[n] \otimes U, X_{\bullet}^{\Delta[1]}) $$
	Since
	\begin{align*}
	\Hom_{s\dm}(\Lambda^i[n] \otimes U, X_{\bullet}^{\Delta[1]}) \simeq \Hom_{s\set} (\Lambda^i[n], \BR \y\big(X_{\bullet}^{\Delta[1]})(U)  \big)
	\end{align*}
	note that by construction  $\BR \y\big(X_{\bullet}^{\Delta[1]})\simeq (\BR \y X_{\bullet})^{\Delta[1]}$
	we get 
	\begin{align*}\Hom_{s\set} (\Lambda^i[n], \BR\y\big(X_{\bullet}^{\Delta[1]})(U)  \big) &\simeq \Hom_{s\set} \big(\Lambda^i[n]\times \Delta[1], \BR\y X_{\bullet}(U) \big) \\
	&\simeq \Hom_{s\dm} \big((\Lambda^i[n]\times \Delta[1])\otimes U, X_{\bullet} \big) 
	\end{align*}
	Hence, we have $M^h_{\Lambda^i[n]} X_{\bullet}^{\Delta[1]}\simeq \rhom_{s\set} \big((\Lambda^i[n]\times \Delta[1]),  X_{\bullet}\big)$.
	
	We will show that $\Hom (\Delta[n]\times \Delta[1],X_{\bullet} )\to \Hom \big((\Lambda^i[n]\times \Delta[1]),  X_{\bullet}\big)$ is a cover by induction. Since  $\Hom(\Lambda^i[1]\times \Delta[1]X_{\bullet})\simeq \hom(\Delta[1], X_{\bullet})\simeq X_1$ which is representable. Combining this with the fact that $\Lambda^i[1]\times \Delta[1]\hookrightarrow \Delta[1]\times \Delta[1]$ is a collapsible extension and $X_{\bullet}$ is an $\infty$-groupoid object, the base case holds. Now consider $n>1$ and $\kan(k,i)$ holds for all $k<n, 0\le i\le k$. Since $\Lambda^i[n] \hookrightarrow \Delta[n]$ is a collapsible extension, we have $\Hom(\Lambda^i[n], X_{\bullet}^{\Delta[1]})$ is representable. Therefore, $\Hom (\Delta[n]\times \Delta[1],X_{\bullet} )\to \Hom \big((\Lambda^i[n]\times \Delta[1]),  X_{\bullet}\big)$ is a cover. 
	\end{proof}

\begin{rem}
    If $X_{\bt}$ happens to be a Lie $k$-groupoid for some $k<\infty$, we can actually show that $X_{\bt}^{\Delta^1}$ is also a Lie $k$-groupoid.
\end{rem}
Now by our previous construction, we have a factorization $X_{\bullet}\stackrel{s^{*}_0}{\longrightarrow} X^{\Delta[1]}_{\bullet} \stackrel{(d^{*}_0, d^*_1)}{\longrightarrow} X_{\bullet} \times X_{\bullet}$ with 
$d^*_i \circ s_0^*=\id_{X_{\bullet}}$, where $s^*_0:X_{\bullet} \simeq X_{\bullet}^{\Delta[0]}\to X_{\bullet}^{\Delta[1]}$ is induced by $s_0: \Delta[1]\to \Delta[0] \subset \partial \Delta[1]$, and $(d^{*}_0, d^*_1):  X_{\bullet}^{\Delta[1]} \to  X_{\bullet} \times X_{\bullet}\simeq X_{\bullet}^{\Delta[0]} \times X_{\bullet}^{\Delta[0]}$ is induced by $(d_0, d_1): \partial \Delta[1] \to \Delta[1]$. Note that $s^*_0$ is a stalkwise weak equivalence by definition. Hence, we simply need to show  $f=(d^{*}_0, d^*_1)$ is a Kan fibration, i.e. for all $n >= 1$ and $0 \le j\le n$, the morphisms of sheaves
\begin{equation*}
    X_n^{\Delta[1]}\stackrel{(\iota^*_{n,j}, f_*)}{\longrightarrow} \BR\Hom(\iota_{n,j}, f)
\end{equation*}
can be represented by a cover. By the isomorphism
\begin{equation*}
    \BR\Hom(L, X^K) \simeq \BR\Hom(L\times K, X) 
\end{equation*}
for any finitely generated simplicial set $K$ and, we have 
\begin{align*}
    \BR\Hom(\iota_{n,j}, f) =& \BR\Hom(\Lambda^j[n] \stackrel{\iota_{n,j}}{\rightarrow} \Delta[n],X^{\Delta[1]}_{\bullet} \stackrel{f}{\longrightarrow} X_{\bullet} \times X_{\bullet} )\\
    =& M^h_{\Lambda^j[n]}X^{\Delta[1]}_{\bullet} \times_{M^h_{\Lambda^j[n]} (X_{\bullet} \times X_{\bullet})} M^h_{\Delta[n]} (X_{\bullet} \times X_{\bullet})\\
    =& \BR \Hom(\Lambda^j[n]\times \Delta[1], X_{\bt}) \times_{\BR \Hom(\Lambda^j[n]\times \del \Delta[1], X_{\bt})}  \BR \Hom( \Delta[n]\times \del \Delta[1], X_{\bt})\\
        =& \BR\Hom\big( (\Lambda^j[n] \times \Delta[1]) \sqcup_{\Lambda^j[n] \times \del\Delta[1] }(\Delta[n] \times \del\Delta[1]), X_{\bt}\big)
\end{align*}

Hence, the $\kan(n,j)$ condition simplifies to 
\begin{equation*}
   \BR \Hom(\Delta[n]\times \Delta[1], X_{bt}) \stackrel{(\iota^*_{n,j}, f_*)}{\longrightarrow} \BR\Hom\big( (\Lambda^j[n] \times \Delta[1]) \sqcup_{\Lambda^j[n] \times \del\Delta[1] }(\Delta[n] \times \del\Delta[1]), X_{\bt}\big)
\end{equation*}
being a cover. 
 Since $(\Lambda^j[n] \times \Delta[1]) \sqcup_{\Lambda^j[n] \times \del\Delta[1] }(\Delta[n] \times \del\Delta[1])\to \Delta[n]\times \Delta[1]$ is a collapsible extension by Lemma \ref{c5}, it suffices to show that $$ \BR\Hom(\iota_{n,j}, f)=\BR\Hom\big( (\Lambda^j[n] \times \Delta[1]) \sqcup_{\Lambda^j[n] \times \del\Delta[1] }(\Delta[n] \times \del\Delta[1]), X_{\bt}\big)$$ is represented by a cover  for all $n,j$ and then applying Lemma \ref{c1} we are done.
 
 First consider $n = 1$. We have a pullback square
 \begin{center}
	\begin{tikzcd}
{ \BR\Hom(\iota_{n,j}, f)} \arrow[r] \arrow[d]    \arrow[dr, phantom, "\ulcorner", very near start]     & X_1\times X_1 \arrow[d, "{(d_j, d_j)}"] \\
X_1 \arrow[r, "{(d_0, d_1)}"] & X_0\times X_0                          
\end{tikzcd}
\end{center}
$(d_j, d_j)$ is clearly a cover, hence the pullback exists and ${ \BR\Hom(\iota_{n,j}, f)}$ is representable. For higher $j, n$ we can apply Lemma \ref{c5} and proceed by induction.

\begin{rem}
    The canonical identification of $\Delta[n]\times \Delta[1] = \cup_{1\le k\le n} x_k$ as $n+1$ $(n+1)$-simplices is constructed by defining $x_k$ to be the $(n+1)$-simplex generated by the following points
    $$
    \big\{(0,0), (1,0),\cdots, (k,0), (k,1), (k+1, 1),\cdots,(n,1)\big\}
    $$
    Let's look at first few examples. For $n=0$, we have $\Delta[0]\times \Delta[1]\simeq{\Delta[1]}$ which is a single 1-simplex. For $n=1$, the decomposition of $\Delta[1]\times \Delta[1]$ is 
    \begin{center}
        \begin{tikzcd}
{(0,1)} \arrow[r]                      & {(1,1)}           \\
{(0,0)} \arrow[r] \arrow[u] \arrow[ru, red] & {(1,0)} \arrow[u]
\end{tikzcd}
    \end{center}
    where $x_0$ is the 2-simplex generated by $\{(0,0), (0,1), (1,1)\}$ which corresponds to the upper 2-simplex, and $x_1$ is then the 2-simplex generated by $\{(0,0), (1,1), (1,1)\}$ which corresponds to the lower 2-simplex. The fiber product quotients out the edge $(0,0)\to (1,1)$.
    
    Similarly, for $n=2$. we decompose $\Delta[2]\times \Delta[1]$ as three 3-simplices
    
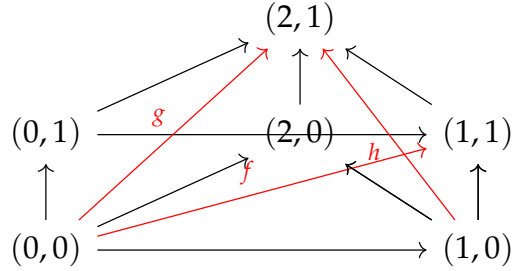
\begin{center}
      \begin{tikzcd}
                                                                                &  & {(2,1)}           &                                                                     \\
{(0,1)} \arrow[rru] \arrow[rrr]                                                 &  & {(2,0)} \arrow[u] & {(1,1)} \arrow[lu]                                                  \\
{(0,0)} \arrow[rrr] \arrow[u] \arrow[rru] \arrow[rrru, red, "f"] \arrow[rruu,red,  "g"] &  &                   & {(1,0)} \arrow[lu] \arrow[u] \arrow[u] \arrow[lu] \arrow[luu, red, "h"]
\end{tikzcd}
    \end{center}
    where $x_0$ is the 2-simplex generated by $\{(0,0), (0,1), (1, 1), (2,1)\}$ which corresponds to the 2-simplex bounded by the diagonal edges $f$ and $g$, $x_1$ is the 2-simplex generated by $\{(0,0), (1,0), (1,1), (2,1)\}$ which corresponds to the 2-simplex bounded by the diagonal edges $h$ and $f$, and 
    $x_2$ is the 2-simplex generated by $\{(0,0), (0,1), (0,2), (2,1)\}$ which corresponds to the 2-simplex bounded by the diagonal edges $h$ and $g$.
\end{rem}

\begin{rem}
    Note that for constant simplicial objects without underlying homotopy theory, the path objects are trivial. For example, take $M$ be a Banach manifold and consider $M_{\bt}$ with $M_i=M$ and all structure maps are identities. Let's look at $M_{\bt}^{\Delta[1]}$. The 0-simplex is simply $\Hom(\Delta[0]\times \Delta[1]\to M_{\bt})\simeq \Hom(\Delta[1]\to M_{\bt})\simeq M$. The 1-simplex is 
    \begin{align}
        \Hom(\Delta[0]\times \Delta[1]\to M_{\bt})  &\simeq M_2\times_{d_1,d_1} M_2\\
        &\simeq M_2\times_{M_1} M_2 \simeq M
    \end{align}
    By similar computation we see that $M_{\bt}^{\Delta[1]}$ is just the $M$ itself. This justifies that we do not suppose any homotopy theory on Banach manifolds.
\end{rem}
 
 \begin{thm}
    Let $(\dm, \mathcal{T})$ be a category with pretopology, then the category of derived Lie $\infty$-groupoids in $(\dm, \mathcal{T})$, $\ilgpd_{\dm}$,  carries a category of fibrant object structure, where fibrations are Kan fibrations, and weak equivalences are stalkwise weak equivalences.
 \end{thm}
 By our construction, the result can be adapted to any homotopy descent categories.
 \begin{cor}
    Let $\CC$ be a homotopy descent category, then the category of derived Lie $\infty$-groupoids in $\CC$ is a category of fibrant objects.
 \end{cor}
 First, let's verify axiom (4), (5).
 \begin{lem}For a homotopy descent category with pretopology $(\dm, \mathcal{T})$, $\ilgpd_{\dm}$ satisfies:
     \begin{itemize}
         \item All isomorphisms in $\ilgpd_{\dm}$ are both weak equivalences and Kan fibrations.
         \item Weak equivalences satisfy 2-out-of-3.
         \item Composition of fibrations are fibrations.
     \end{itemize}
 \end{lem}
 \begin{proof}
    (1) is trivial. For (2), note that weak equivalence of simplicial set satisfies 2-out-of-3, and by the construction of stalkwise weak equivalences, it is obvious that it also satisfies 2-out-of-3. (3) follows from Proposition \ref{fib}.
 \end{proof}
 (2) and (3) follows from Proposition \ref{pb1} and Proposition \ref{pb2} respectively. (1) and (7) are trivial by our construction. Therefore, we finish the proof.


\subsection{iCFO for incomplete derived spaces}
For incomplete or non-small categories, the previous technique won't work. The main issue is that due to lack of limits, pullback might not exist in general. Hence, we can only construct an iCFO structure. We will take advantage of the technique developed in \cite{RZ20} which looks at category with locally stalkwise pretopology.

\begin{thm}
    Given an incomplete category with locally stalkwise pretopology $(\dm, \mathcal{T})$, then the category of derived Lie $\infty$-groupoids in $(\dm, \mathcal{T})$ carries a category of fibrant object structure, where fibrations are Kan fibrations, and weak equivalences are stalkwise weak equivalences.
 \end{thm}
 Note that we only need to verify (2) and (3), and all (4)-(7) follows from the homotopy descent category case.
 \begin{prop}
    Let $f: X_{\bt} \to Y_{\bt}$ and $g: Z_{\bt} \to Y_{\bt}$ be two morphisms in $\ilgpd_{\dm}$, where $f$ is a Kan fibration and $g$ is arbitrary. Suppose the pullback $Z_0\times_{Y_0} X_0$ exists in $\dm$, then we have
    \begin{enumerate}
        \item All $Z_k\times_{Y_k} X_k$ exists for $k >= 1$, and the induced map $Z_{\bt} \times_{Y_{\bt}} X_{\bt}\stackrel{p_f}{\rightarrow} Z_{\bt}$ is a Kan fibration in $s\dm$.
        \item $Z_{\bt} \times_{Y_{\bt}} X_{\bt}$ is a derived Lie $\infty$-groupoid in $\dm$.
    \end{enumerate}\end{prop}
    
    \begin{proof}
        First, let's look at (1). We want to show that the morphism of sheaves
        \begin{equation*}
            Z_n \times_{Y_n} X_n \stackrel{(\iota^* , {p_{f}}_*)}{\longrightarrow} \Hom(\Lambda^j[n]\stackrel{\iota}{\to} \Delta[n], Z_{\bt} \times_{Y_{\bt}} X_{\bt}\stackrel{p_f}{\rightarrow} Z_{\bt})
        \end{equation*}
        is represented by a cover. Since $M_{K_{\bt}}(-): s\dm \to s\sh(\dm)$ preserves limits, we have
        \begin{align*}
            \Hom(\iota, p_f) =& M_{\Lambda^j[n]}(Z_{\bt} \times_{Y_{\bt}} X_{\bt}) \times_{M_{\Lambda^j[n]}Z_{\bt}} Z_n\\
            \simeq& \big(M_{\Lambda^j[n]}Z_{\bt} \times_{M_{\Lambda^j[n]}Y_{\bt}} M_{\Lambda^j[n]}X_{\bt}\big) \times_{M_{\Lambda^j[n]}Z_{\bt}} Z_n
        \end{align*}
        On the other hand, we also have a composition of pullbacks
        \begin{center}
        \begin{equation}\label{dia1}
            \begin{tikzcd}
{ \Hom(\iota, p_f)} \arrow[r, "\pr_1"] \arrow[d, "\pr_2"] \arrow[dr, phantom, "\ulcorner", very near start]& Z_n \arrow[d, "\iota^*"]  \\
{M_{\Lambda^j[n]}Z_{\bt} \times_{M_{\Lambda^j[n]}Y_{\bt}} M_{\Lambda^j[n]}X_{\bt}} \arrow[r, "{p_f}_*"] \arrow[d, "\pr_3"]      \arrow[dr, phantom, "\ulcorner", very near start]& {M_{\Lambda^j[n]}Z_{\bt}} \arrow[d, "g_*"] \\
{M_{\Lambda^j[n]}X_{\bt}} \arrow[r, "f_*"]                & {M_{\Lambda^j[n]}Y_{\bt}}          
\end{tikzcd}
\end{equation}
        \end{center}
        Therefore we have 
        \begin{equation}\label{iso1}
            \Hom(\iota, p_f) \simeq {M_{\Lambda^j[n]}X_{\bt}}\times_{M_{\Lambda^j[n]}Y_{\bt}} Z_n
        \end{equation}
    
        We also have another compositions of pullbacks
        
        \begin{center}
        \begin{equation}
            \begin{tikzcd}
Z_n\times_{y_n} X_n \arrow[r] \arrow[d, "{(\iota^*, {p_{f}}_*)}"]\arrow[dr, phantom, "\ulcorner", very near start] & {X_n} \arrow[d, "{(\iota^* , f_*)}"]           &              \\
{\Hom(\iota, p_f)} \arrow[r, "{p_{f}}_*"] \arrow[d]\arrow[dr, phantom, "\ulcorner", very near start] & {Y_n \times_{M_{\Lambda^j[n]}Y_{\bt}} M_{\Lambda^j[n]}X_{\bt}}       \arrow[d] \arrow[r]\arrow[dr, phantom, "\ulcorner", very near start] & {M_{\Lambda^j[n]}X_{\bt}} \arrow[d, "f_*"] \\
{Z_n} \arrow[r, "g_*"]           & {Y_n} \arrow[r, "\iota^*"]           & {M_{\Lambda^j[n]}Y_{\bt}}          
\end{tikzcd}
 \end{equation}
        \end{center}
       
    Now by assumption we see $(\iota^* , f_*)$ is a cover, so we just need to show that $\Hom(\iota, p_f) \simeq {M_{\Lambda^j[n]}X_{\bt}}\times_{M_{\Lambda^j[n]}Y_{\bt}} Z_n$ is representable, then $(\iota^*, {p_{f}}_*)$ is a cover and $Z_n\times_{y_n} X_n$ is representable. 
    
    We shall proceed by induction. Consider $n = 1$. By the top square in diagram \ref{dia1}, and replacing $\Hom(\iota, p_f)$ by the isomorphism \ref{iso1} from the whole square, we have a pullback square
    \begin{center}
            \begin{equation}
                \begin{tikzcd}
    X_0\times_{Y_0} Z_1 \arrow[r, "\pr_1"] \arrow[d, "\pr_2"] \arrow[dr, phantom, "\ulcorner", very near start]& Z_1 \arrow[d, "\iota^*"]  \\
    Z_{0} \times_{Y_0} X_{0} \arrow[r, "{p_f}_*"] ]& Z_{0}
    \end{tikzcd}
    \end{equation}
            \end{center}
    Note that $\iota^* = d_j: Z_1\to Z_0$ which is a cover. Since $ Z_{0} \times_{Y_0} X_{0} $ is representable, so is $X_0\times_{Y_0} Z_1 $ and $\pr_2$ is a cover. Hence, $X_1\times_{Y_1} Z_1$ is representable and $p_f$ satisfies Kan condition for $n=1.$
    
    Now suppose $p_f$ satisfies Kan condition $\kan(m, j)$ for $1\le m < n$ and $1 \le j \le m$, we want to show the Kan condition holds for $m=n$ and $1 \le j \le m$ as well. Applying Lemma \ref{c2} with $S_{\bt} = \Lambda^j[n]$ and $T_{\bt} = \Delta[n]$, we get that $\Hom(\iota_{n,j}, p_f)$ is representable and $p_f$ satisfies $\kan(n,j)$ for all $0\le j \le n$. Therefore, $p_f$ is a Kan fibration. 
    
    Finally, let's look at (3). By assumption $z:Z_{\bt}\to *$ is a Kan fibration, then by Proposition \ref{fib}, $z\circ p_f: Z_{\bt}\times_{Y_{\bt}} X_{\bt}\to *$ is also a Kan fibration. Hence, $Z_{\bt}\times_{Y_{\bt}} X_{\bt}$ is a Lie $\infty$-groupoid.
    \end{proof} 


 \part{Homotopical algebra of derived  Lie $\infty$-groupoids and algebroids}
 In this chapter,
 we study homotopical algebras for derived Lie $\infty$-groupoids and algebroids and study their homotopy-coherent representations, which we call $\infty$-representations. We relate $\infty$-representations  of $\li$-algebroids to (quasi-) cohesive modules developed by Block, and $\infty$-representations  of Lie $\infty$-groupoids to $\infty$-local system introduced by Block-Smith. Then we apply these tools in studying singular foliations and their characteristic classes. We then construct Atiyah classes for $\li$-algebroids pairs.
 
\subsection{Semi-model categories}
\begin{defn}[\cite{Nui19}, \cite{WY18}]
    Let $\CC$ be a bicomplete category. We say $\CC$ is a (left) {\it semi-model category}\index{semi-model category} if it is equipped with wide subcategories of weak equivalences\index{weak equaivalence!in a semi-model category} $\mathcal{W}$, cofibrations\index{cofibration!in a semi-model category} $\mathcal{C}$, and fibrations\index{fibration!in a semi-model category} $\mathcal{F}$, which satisfy the following data:
    \begin{enumerate}
        \item The weak equivalences satisfy 2-out-of-3.
        \item The weak equivalences, fibrations, and cofibrations are stable under retracts.
        \item The cofibrations have the left lifting property with respect to the trivial fibrations. The trivial cofibrations with cofibrant domain (i.e. with a domain $X$ for which the map $\empty \to X$ is a cofibration) have the left lifting property with respect to the fibrations.
        \item Every map can be factored functorially into a cofibration followed by a trivial fibration. Every map with cofibrant domain can be factored functorially into a trivial cofibration followed by a fibration.
        \item The fibrations and trivial fibrations are stable under transfinite composition, product, and base change.
    \end{enumerate}
\end{defn}
\begin{rem}
    It is also possible to define a semi-model category through specific adjunction to model categories. For more details, see \cite{WY18}.
\end{rem}
\begin{rem}
    In a semi-model category, only the cofibrations and trivial fibrations are determined by each other via the lifting property, which implies that a semi-model structure is only determined by its weak equivalences and fibrations.
\end{rem}
\section{Derived  $L_{\infty}$-algebroid}

Let $A$ be a (unital) commutative dga over characteristic 0. The {\it tangent module} $T_A$ associated to $A$ is defined by the space of graded $k$-derivations $\der_k(A, A)$. Note that $T_A$ is both a dg-$A$-module and dg-Lie algebra over $k$.

\begin{defn}
	Let $\g$ be a dg-$A$-module.
	We say $\g$ is a {\it dg-Lie algebroid} \index{dg-Lie algebroid}  over $A$ if it also has a dg-Lie algebra structure over $k$ and with anchor map $\rho: \g \to T_A$ satisfies
	\begin{enumerate}
		\item $\rho$ is a map of dg-$A$-modules;
		\item $\rho$ is a map of dg-Lie algebras;
		\item The following Leibniz rule holds
		$$ [x, a\cdot y ] = (-1)^{|x|\cdot |a|}[x,y] + \rho(x)(a) y $$
		for $a \in A$, $x, y\in \g$.
	\end{enumerate}
\end{defn}
Morphisms between two  dg-$A$-algebroids over $A$ are $A$-linear morphisms of dg-$A$-module over $T_A$ which preserves the Lie brackets.

\begin{defn}
    Let $\g$ be a dg-$A$-module.
    We say $\g$ is a {\it $L_{\infty}$-algebroid} \index{$L_{\infty}$-algebroid}  over $A$ if it also has a $L_{\infty}$-algebra structure over $k$ and with an anchor map $\rho: \g \to T_A$ which satisfies
    \begin{enumerate}
    	\item $\rho$ is a map of dg-$A$-modules;
    	\item $\rho$ is a map of $L_{\infty}$-algebras;
    	\item The following Leibniz rule holds
    	\begin{align}
    	[x, a\cdot y ] &= (-1)^{|x|\cdot |a|}a[x,y] + \rho(x)(a) y\\
    	[x_1, \cdots, x_{n-1},  a\cdot x_n ] &= (-1)^{n|a|}(-1)^{a(|x_1|+\cdots+|x_n|)}a[x_1, \cdots, x_n]\quad n\ge 3 
    	\end{align}
    	for $a\in A$, $x, y, x_1, \cdots, x_n \in \g$.
    \end{enumerate}
\end{defn}
\begin{example}
	Let $A$ be an ordinary $k$-algebra and $\g$ is $\ge 0$-graded, then an $L_{\infty}$-structure on $\g$ is equivalent to a differential on the Chevalley-Eilenberg algebra $\sym_A(\g[1])^{\vee}$.
	
	In particular, if $A$ is a (dg-) $\cinf$-ring, we call $\g$ a (derived)NQ-(super)manifold.
\end{example}

\begin{example}[Action $L_{\infty}$-algebroids]
    Let $\g$ be an $L_{\infty}$-algebra, and $\rho: \g \to T_A$ be a map of $L_{\infty}$-algebra over $k$. We can equip $A\otimes \g$ with a structure of $L_{\infty}$-algebroid by extending $\rho$ to an $A$-linear map, and brackets are given by
    \begin{align*}
        [a\otimes x, b\otimes y] =& \pm ab\otimes [x, y], + a\cdot \rho(x)(b) \otimes y - (\pm)b\cdot \rho(y)(a) \otimes x\\
        [a_1 \otimes x_1, \cdots, a_n \otimes x_n] =& a_1\cdots a_n[x_1, \cdots, x_n]
    \end{align*}
    Here $\pm$ is the Koszul sign for the grading.
\end{example}

\begin{example}[Singular foliations]
    Consider $A = \cinf(M)$. Let $\CF$ be a singular foliation which admits a resolution by a complex of vector bundles $E_{\bullet}$, then we are able to construct an $L_{\infty}$-algebroid structure on $E_{\bullet}$ \cite{LLS20}, which is called the {\it universal $L_{\infty}$-algebroid} of the singular foliation $\CF$.
\end{example}

Apparently there are more choices of defining a sub-$L_{\infty}$-algebroids due to the homotopical nature of $L_{\infty}$-algebroids. We will use the following definition throughout this paper.
\begin{defn}
Let $\g$ be an $L_{\infty}$-algebroid over $A$, then we define a {\it sub-$L_{\infty}$-algebroid}\index{$L_{\infty}$-algebroid!sub} (or simply {\it subalgebroid}) of $\g$ to be a sub-$A$-module of the kernel of the anchor map which is also closed under the brackets and the differential. Later we will see that, $\g/\h$ inherits an $L_{\infty}$-algebroid structure, which plays the role of 'normal bundle' of $\h$. 
\end{defn}
We have two differential type of morphisms of $L_{\infty}$-algebroid.
\begin{defn}
	A {\it (strict) morphism}\index{$L_{\infty}$-algebroids!morphism} between $L_{\infty}$-algebroids is a $A$-linear morphisms of dg-$A$-module over $T_A$ which preserves the $L_{\infty}$ structure.
\end{defn}

In differential geometry, we often work with a weaker type of morphisms. 
\begin{defn}
An $L_{\infty}$-morphism $\g \leadsto \h$  between $L_{\infty}$-algebras is a {\it twisting cochain}\index{twisting cochain}
$$
\overline{C}_{\bullet}(\g) \to \h[1]
$$
or, equivalently, a map of commutative dg-coalgebra $\overline{C}_{\bullet}(\g) \to \overline{C}_{\bullet}(\h)$.
\end{defn}

\begin{defn}
An {\it $L_{\infty}$-morphism}\index{\it $L_{\infty}$-morphism} $\g \leadsto \h$ between $L_{\infty}$-algebroids is an $L_{\infty}$-morphism of  $L_{\infty}$-algebras $\tau: \g \to \h$ such that
\begin{enumerate}
    \item the composition $\rho_{\h}: \overline{C}_{\bullet}(\g) \to \h[1] \to T_A[1]$ first takes the quotient by $\sym_k^{\ge 2} \g[1] \subset \overline{C}_{\bullet}(\g)$ and then applies the anchor map $\rho$ of $\g$ to the remaining of $\g[1]$.
    \item the map of graded vector spaces $\tau:\sym_k^{\ge 1}\g[1] \to \h[1]$ descends to a graded $A$-linear map.
\end{enumerate}
\end{defn}
\begin{rem}
    One motivation for $L_{\infty}$-morphisms comes from differential geometry and mathematical physics. Consider $A$ to be a ordinary algebra, and the dg-$A$-module underlying the $L_{\infty}$-algebroids $\g$ is a nonnegatively graded complexes of finitely generated projective $A$-module, then an $L_{\infty}$-morphism $\g \leadsto \h$ is equivalent to a map of cdga's $\sym_A (\g[1])^{\vee} \to \sym_A (\h[1])^{\vee}$. In fact, this is the same as the morphisms of NQ-manifolds or dg-manifolds. For example $A = \cinf{(M)}$ and $\g$ is a complex of finite dimensional vector bundles.
\end{rem}

\begin{thm}[\cite{Nui19}]
Let $A$ be a fixed cdga. There is a 
right proper, tractable semi-model structures on the category of derived $L_{\infty}$-algebroid $\lalgd_A$, in which a map is a weak equivalence (resp. a fibration) if and only if it is a quasi-isomorphism (a degreewise surjection). 
\end{thm}

Let $A$ be a dg-$\cinf$-ring (or simply a cdga over $\R$). Denote $A^c$ the cofibrant replacement of $A$. The $\infty$-category associated to the semi-model category of $\li$-algebroids over $A$ is 
\begin{equation*}
    \lalgdi_A = \lalgd_{A^c}[W^{-1}].
\end{equation*}
Note that here we need to pass to the cofibrant replacement of $A$ since the tangent module $T_A$ is only homotopy invariant when $A$ is cofibrant.

\begin{prop}
There is combinatorial stable model category structures on the category of derived $L_{\infty}$-algebroid $\lalgd_A$, which presents the $\infty$-category $\lalgdi$.
\end{prop}

\begin{defn}
   Let $\g$ be a derived $L_{\infty}$-algebroid. Define $\Omega^1_{\g}$ to be the sheaf of smooth 1-form of $\CO_{\g}$. Note that $\Omega^1_{\g}$ is actually a sheaf of chain cochain complexes. We write $\Omega^p_{\g} = \bigwedge^p_{\CO_{\g}}\Omega^1_{\g}$. Denote the chain cochain complexes of global sections $\Gamma(\g, \Omega^p_{\g})$ by $\Omega^p_{\cinf(g)}$.
   
    Define  $T_{\g}$ to be the sheaf $\Hom_{\CO_{\g}}(\Omega^1_{\g}, \CO_{\g})$, and let $T_{\cinf(\g)}$ denote the chain cochain complexes of the global section $\Gamma(\g, T_{\g})$.
\end{defn}

\begin{rem}
    Here $\Omega^1_{\g}$ is not the module of K\"{a}hler differential of the cdga $\cinf(\g)$ since we need the derivation $d: \cinf(\g) \to \Omega^1_{\g}$ to be a $\cinf$-derivation when restricted to $A^0$.
\end{rem}

\begin{defn}
Let $\g$ be a derived $L_{\infty}$-algebroid over $A$. We define the de Rham complex $\dr(\g)$ to be the product of total cochain complex of the triple complex
$$
    \cinf(\g) \stackrel{d^{\dr}}{\to} \Omega^1_{\cinf(\g)}  \stackrel{d^{\dr}}{\to} \Omega^2_{\cinf(\g)}  \stackrel{d^{\dr}}{\to} \cdots
$$
therefore, $\dr(\g)^m = \prod_{i + j - k = m} (\Omega^i_{\cinf(\g)})^j_k$ with total differential $D = \delta + d^{\CE} + d^{\dr}$, where $d^{\CE}$ is the Chevalley-Eilenberg differential on $\g$ and $d^{\dr}$ is the de Rham differential defined above.

\end{defn}
\subsection{Perfect complexes}

\begin{defn}
    Let $M$ be a derived manifold, a quasi-coherent sheaf $\E$ is a {\it perfect complex}\index{perfect complex} if it is locally finitely presentable, i.e for every $x\in M$, there exists an open neighborhood $U$ such that $\E|_U$ can be obtained from the structure sheaf $\CO|_U$ be finite limits and colimits. We denote $\perf(M) \subset \QCoh(M)$ the category of perfect complexes.
\end{defn}
Let $M$ be a smooth manifold, then the most common perfect complexes are finite chain complexes of vector bundles. On the other hand, any perfect complex on $M$ is locally quasi-isomorphic to finite chain complexes of vector bundles. For compact manifolds, global resolutions exist.

\begin{defn}
    Let $\E \in \perf(M)$. We say $\E$ has {\it Tor-amplitude} contained in $[a,b]$ if the associated sheaf $\E\otimes_{\CO_M} \pi_0(\CO_M)$'s homotopy sheaves vanish outside degree $[a,b]$. We denote the subcategory of perfect complexes with Tor-amplitude contained in $[a, b]$ by $\perf^{[a,b]}(M)$.
\end{defn}
\section{$\infty$-Lie differentiation}
\begin{defn}
    Let $A$ be a commutative cosimplicial dga, we construct its {\it normalization}\index{normalization} as follows.
    
    First, we define a cochain complex $N^{\bt} A$ where
    \begin{equation*}
        N^mA = \{ a\in A^m: \sigma^j \in A^{m-1}, 0\le j\le m\}
    \end{equation*}
    with differential $d = \sum_i (-1)^i\del_i$. Next we define an associative product $\cup$ similar to the usual Alexander-Whitney product on $N^{\bt} A$ by
    \begin{equation*}
        a \cup b = (\del^{[m+1, m+n]}a)\cdot (\del^{[1, m]} b)
    \end{equation*}
    for $a\in N^m A, b \in N^n B$.
    
    Now we define a commutative cochain algebra $D^{\bt} A$ as the quotient of $N^{\bt} A$ by
    \begin{equation*}
        (\del^I a) \cdot (\del^J b) \simeq\begin{cases}
         (-1)^{(J, I)}(a\cup b) & a\in A^{|J|}, b\in A^{|I|}\\
         0 &\text{otherwise}
        \end{cases}
    \end{equation*}
    For each disjoint sets (possibly empty) $I, J$. Here $(-1)^{(J, I)}$ denote the sign of permutation of integers $I\sqcup J$ which sends first $|I|$ elements to $I$ (in order) and the left to $T$ (in order).
    \end{defn}
    We can easily see that $D^{\bt}$ is a functor from the category of cosimplicial cdga's to the category of stacky cdga's. In fact, we have:
    \begin{prop}[\cite{Pri17, Lemma 3.5}]
        $D^{\bt}$ is a left Quillen functor from the Reedy model structure on cosimplicial cdga's to the model structure on stacky cdga's in Lemma \ref{stacky}.
    \end{prop}

    Next, we can define a generalization of constructing Lie algebroids from Lie groupoids.
    \begin{defn}
        Let $G_{\bt}$ be a derived Lie $\infty$-groupoid, define the normalization $NG$ to be a derived $L_{\infty}$-algebroid with structure sheaf $D^{\bt}\big((\sigma^0)^{-\bt} \CO_G\big)$, where $(\sigma^0)^{-\bt} \CO_G$ is the cosimplicial sheaf with 
        \begin{equation*}
            \big((\sigma^0)^{-\bt} \CO_G\big)^m = \big((\sigma^0)^m\big)^{-1}\CO_{G_m}
        \end{equation*}
    \end{defn}
    
    \begin{lem}[\cite{Pri20a}]
        Let $G_{\bt}$ a derived Lie $n$-groupoid, then its normalization is an $L_{\infty}$-algebroid with degree $\le n$. 
    \end{lem}
    \begin{proof}
    Similar to the standard Dold-Kan construction. See \cite[Lemma 4.17]{Pri20a}.
    \end{proof}

\subsection{Tangent complex}

Let $A\in \cialg^{\dg}$ be a dg $\ci$ ring and $E\in\Mod_A^{\dg}$ be a dg module over $A$. A {\it multiderivation} of degree $n$ is a graded symmetric multilinear map
$D: E^{\otimes(n+1)} \to E$
which is a derivation in each variable, i.e. there is a {\it symbol} map $\sigma_{D}:E^{\otimes(n+1)} \to T_A$ such that the following graded Leibniz rule holds
\begin{equation}
D(s_0,s_1,\cdots fs_n) = f  D(s_0,s_1,\cdots s_n) + \sigma_{D}(s_0,s_1,\cdots s_n)(f)
\end{equation} 

\begin{lem}
	Let $D \in \der^n(E)$, then we have the short exact sequence of dg modules
	\begin{equation}
	0 \to \sym^{n+1} E^{\vee} \otimes E \to \der^n E \to \sym^n E^{\vee} \otimes T_A \to 0 
	\end{equation}
	and $\der^n(E)=0$ for $n> \rk E$. 
\end{lem}

\section{Cohesive modules over stacky dga}
\subsection{Stacky dga}

Though the category of differential graded algebra suffices for most of our work, sometimes it is still necessary to consider a (cohomologically graded) cdga $A^{\bullet}$ (for example, Chevalley-Eilenberg algebra of a (derived) $L_{\infty}$-algebroid), where $A^0$ is a (homologically graded) dga. Hence, we are looking for a specific kind of double complex where the homological dga and cohomological dga structures are compatible. 
\begin{defn}
Define a chain-cochain complex $V$ over $k$ to be a bigraded $k$-vector space equipped with two differentials $d: V^i_j \to V_j^{i+1}$ and $\delta: V^i_j \to V_{j-1}^j$ such that $(d+\delta)^2 = d \delta + \delta d = 0.$
\end{defn}
There is an obvious tensor product $\otimes$ in the category defined above, which allows us to define the algebra structure on it.
\begin{defn}[\cite{Pri17}]
    A {\it stacky dga}\index{ stacky dga} $A$ is a chain-cochain complex $A^{\bullet}_{\bullet}$ equipped with a commutative product $A\otimes A \to A$ and a unit $k \to A$. We can regard all chain complexes as chain-cochain complexes by $V = V^0_{\bullet}$. Given a chain dga $R$, a stacky dga $A$ over $R$ is given by a map of stacky dga $R \to A$. If in addition  $A$ is graded commutative, then we say $A$ is a {\it stacky cdga}\index{stacky cdga}.

    As the name suggests, 'stacky' means these dga's are enhanced in the 'stacky' direction, i.e. they are not only model for derived spaces, but also  derived (infinitesimal) stacks. We denote $\dgcalg_R^{\st}$ the category of stacky cdga's over $R$, and $\dcalg_R^{\st}$ the full subcategory consists of stacky cdga's which are concentrated in non-negative cochain degrees.
\end{defn}

\begin{example}[derived $L_{\infty}$-algebroids]
     A large class of stacky cdga's is given by derived $L_{\infty}$-algebroids. Recall that a derived $L_{\infty}$-algebroid $\g$ over a derived manifold $(X, \CO_X)$ is given by an $L_{\infty}$-algebroid structure over the dga $
     A = \Gamma(\CO_X)$. The chain part is given by the derived direction
     $$
     \cdots \CO_{X,2}\stackrel{\delta}{\to} \CO_{X,1}\stackrel{\delta}{\to} \CO_{X,0} = \cinf(X)
     $$
     and the cochain part is given by the stacky direction
     $$
     \cinf(X) = (\sym \g^{\vee}[-1])^0 \stackrel{d}{\to} (\sym \g^{\vee}[-1])^1 \stackrel{d}{\to} (\sym \g^{\vee}[-1])^2 \stackrel{d}{\to} \cdots
     $$
     
\end{example}
\begin{example}[BRST complex for coisotropic reduction]\label{brst}
    BRST complexes, introduced in \cite{BRS75}\cite{Tyu75}, is a tool in mathematical physics to describe both the homotopy quotients and homotopy intersections.

    Let $(M, \omega)$ be a symplectic manifold of dimension $2n$, and $M_0$ a coisotropic submanifold of codimension $k$, i.e. $(T_p M_0)^{\perp} \subset T_p M_0$ for all $p \in M_0$. For simplicity, we assume that $M_0$ has a trivial normal bundle. Now we can write $M_0$ as the zero set of a smooth function $\phi: M \to V$, where $V$ is a vector space of codimension $k$. Pick a basis $\{e_i\}_{1\le i\le k}$ for $V$, then we can write $\phi = \sum_{i = 1}^k \phi_i e_i$, where $\phi_i \in \cinf(M)$. Since $M_0$ is a submanifold of $M$, $\phi_i$'s generate the vanishing ideal $\mathcal{I}$ of $M_0$, so 
    $$\mathcal{I} = \{\sum_i f_i \phi_i| f_i \in \cinf(M)\}$$
    Note that, since $M_0$ is coisotropic, $\mathcal{I}$ is closed under the Poisson bracket, i.e. $\{ \mathcal{I}, \mathcal{I}\} \subset \mathcal{I}$. 
    
    Now we define the BRST complex by 
    $$
    C^{p,q} = \bigwedge^p V^{\vee} \otimes \bigwedge^q V
    $$
    where we regard $V$ as a trivial vector bundle on $M$.
    Clearly, by the graded algebra structure induced from exterior product, we get a stacky dga $A^{\bullet}_{\bullet}$ with $A^q_p = C^{p, q}$, with $D = d + (-1)^i\delta$ where $d$ is the Chevalley-Eilenberg differential and $\delta$ is the Koszul differential. Note that $H^0(A) \simeq \cinf(M_0/\CF)$ as a Poisson algebra, where $\CF$ is the foliation generate by the $(T M_0)^{\perp}$. This is a prototypical example of a stacky cdga, where the underlying geometric space is the leaf space of the foliation $\CF$ on $M_0$. We can see that the derived direction is a generalization of submanifolds or subspaces, whereas the stacky direction generalizes the (homotopy) quotient or orbit space.
\end{example}

\begin{example}[de-Rham algebras of derived manifolds]
    Let $(X, \CO_X)$ be a derived manifold $(X, \CO_X)$. Denote the chain complex of the global sections of $\CO_X$ by $\cinf(X)$. Recall $\Omega_X^1$ denotes the sheaf of chain complexes of smooth 1-form on $\CO_X$, and $\Omega^P_X = \bigwedge^p_{\CO_X} \Omega^1_X$. The de Rham complex $\dr(X)$ is the product total cochain complex of the double complex,
    $$
    \cinf(X) \stackrel{d}{\to} \Omega^1_{\cinf(X)}  \stackrel{d}{\to} \Omega^2_{\cinf(X)}  \stackrel{d}{\to} \cdots
    $$
    hence $\dr(X)^m = \prod_j(\Omega^{m+j}_{\cinf(X)})_j$. $\dr(X)$ is then a stacky cdga with $D = d + (-1)^i \delta$, where $\delta$ is the differential of $\CO_X$, and the product structure comes from the exterior algebra.
\end{example}

    

\begin{defn}
    Recall that a morphism $U\to V$ between a chain-cochain complexes is a levelwise quasi-isomorphism if $U^i 
    \to V^i$ is a quasi-isomorphism for all $i \in \Z$. We call a morphism of stacky cdga's a levelwise quasi-isomorphism if the underlying chain-cochain complex is so.
\end{defn}
The following is \cite{Pri17}, Lemma 3.4.
\begin{lem}\label{stacky}
    There is a cofibrantly generated model structure on stacky cdga's over $R$ in which fibrations are surjections and weak equivalences are levelwise quasi-isomorphisms.
\end{lem}

Next, we will generalize stacky dga's to curved stacky dga's, where the integrability condition $d^2=0$ is no longer satisfied. Instead, we have a 'curvature' for each stacky dga.

\begin{defn}
    A {\it curved stacky dga}\index{curved stacky dga} is a quadruple $\mathsf{A} = (A^{\bullet}_{\bullet}, d, \delta, c)$ where $A^{\bullet}_{\bullet}\in \dcalg_k$ is a stacky dga where the cohomological degree is non-negatively graded, with a derivation $D = d + \delta$ which satisfies the usual graded Leibniz rule and
    \begin{equation*}
        D^2 (a) = [c, a]
    \end{equation*}
    for a fixed $c \in (A^{\bullet}_{\bullet})^2$ satisfying the Bianchi identity $Dc = 0$.
\end{defn} 
\cite{PP05}\cite{Blo05} discuss the case where $A$ are ordinary dga's.

For simplicity, when we write a single superscript $A^{\bt}$ for a stacky dga, we will always mean total degrees, i.e $A^n = \oplus_{p-q=n} A^p_q$.

\begin{example}[Endomorphism module of an affine derived manifold]
    
\end{example}

\begin{example}[derived $L_{\infty}$-algebroids]
    
\end{example}

\subsection{Cohesive modules over stacky dga's}
Let $\mathsf{A} = (A^{\bullet}_{\bullet}, d, \delta, c)$ be a curved stacky dga, and let $E= E^{\bullet}$ be a (right) dg-$A^0$-module.
\begin{defn}
    Let  $\mathbb{E}: E^{\bullet}\otimes_{ {A}^0} A^{\bullet}_{\bullet} \to E^{\bullet}\otimes_{{A}^0} A^{\bullet}_{\bullet}$ be a $k$-linear map of total degree one which satisfies the graded Leibniz rule
    \begin{equation*}
        \mathbb{E}(e\omega) = (\mathbb{E}(e\otimes 1))\omega + (-1)^{|e|}ed\omega,
    \end{equation*}
    then we call $\mathbb{E}$ a {\it $\Z$-connection}\index{connection!$\Z$-} on $E$.
\end{defn}

A $\Z$-connection is determined by its value on $E^{\bullet}$ part. We can write $\mathbb{E} = \mathbb{E}^0 + \mathbb{E}^1 + \mathbb{E}^2+ \cdots$, where $\mathbb{E}^k:  E^{\bullet}\to E^{\bullet - k + 1}\otimes_{{A}^0} {A}^{k}$. Clearly, $\mathbb{E}^1$ part corresponds to ordinary connections on each $E^n$, and $\mathbb{E}^k$ is $A^0$-linear. 

Note that the usual definition of curvature $\mathbb{E}^2$ will not work, since it won't be ${A}^0$-linear. Instead, we define the {\it relative curvature}\index{relative curvature} of $\mathbb{E}$ to be the operator 
\begin{equation*}
    R_{\mathbb{E}} = \mathbb{E}^2(e) + e\cdot c
\end{equation*}
where $c$ is the curvature of ${A}$. Note that $R_{\mathbb{E}}$ is then $A^0$-linear.

\begin{defn}
    Let $E = E^{\bullet}$ be a dg ${A}^0$-module (bounded in both directions) together with a flat $\Z$-connection $\mathbb{E}$, i.e. $R_{\mathbb{E}} = 0$, then we call $E$ a {\it quasi-cohesive module}\index{cohesive module!quasi-}. If $E$ is also finitely generated and projective over $A^0$ and bounded in both directions, then we call $E$ a {\it cohesive module}\index{cohesive module}. Denote the category of cohesive modules over ${A}$ to be $\Mod^{\coh}_{{A}}$, and the category of quasi-cohesive modules by $\Mod^{\qcoh}_{{A}}$. Note that $\Mod^{\coh}_{{A}}$ is the same as $\mathcal{P}_A$ in \cite{Blo05}, \cite{BS14}, \cite{BD10}, and \cite{BZ21}. For more about the theory of cohesive modules and quasi-cohesive modules, see \cite{Blo05}.
\end{defn}

We define the degree $k$ morphisms between two cohesive modules $E_1 = (E_1^{\bullet}, \mathbb{E}_1)$ and $E_2 = (E_2^{\bullet}, \mathbb{E}_2)$ to be $$
\HOM^k(E_1, E_2) =  \HOM_{A^{\bullet}}^k(E^{\bullet}_1\otimes_{ A^0} A^{\bullet}, E^{\bullet}_2\otimes_{ A^0} A^{\bullet})
$$
, i.e. the set of degree $k$ $A^{\bullet}$-linear map from $E^{\bullet}_1\otimes_{ A^0} A^{\bullet}$ to  $E^{\bullet}_2\otimes_{ A^0} A^{\bullet}$. By a similar argument as above, we have 
$$
\HOM_{A^{\bullet}}^k(E^{\bullet}_1\otimes_{ A^0} A^{\bullet}, E^{\bullet}_2\otimes_{ A^0} A^{\bullet}) = \Hom_{A^{0}}^k(E^{\bullet}_1, E^{\bullet}_2\otimes_{ A^0} A^{\bullet})
$$. We define a differential on the morphisms $d_{\HOM}: \HOM^{\bullet}(E_1, E_2) \to \HOM^{\bullet + 1}(E_1, E_2) \to$ by 
$$
d_{\HOM}(e) = \mathbb{E}_2(\phi(e)) - (-1)^{|\phi|}\phi(\mathbb{E}_1(e)).
$$
It is easy to verify that $d_{\HOM}^2 = 0$, and hence $\Mod^{\coh}_{{A}}$ is a dg-category.

Given a dg-category $\CC$, we have a subcategory $Z^0(\CC)$ which has the same objects as $\CC$ and morphisms
$$
Z^0(\CC)(x, y) = Z^0(\CC(x, y))
$$
i.e. degree 0 closed morphisms in $\CC(x, y)$. On the other hand, we can form the {\it homotopy category}\index{homotopy category} $\ho(\CC)$ which has the same objects as $\CC$ and morphisms, 
$$
\ho(\CC)(x,y) = H^0(\CC(x,y))
$$
which is the 0th cohomology of the morphism complex.

Next, we will briefly discuss the triangulated structure of cohesive modules and explore homotopy equivalences between cohesive modules. 

First, we define a shift functor. For $(E, \mathbb{E}) \in \Mod^{\coh}_{{A}}$, we set $E[1] = (E[1] = (E^{\bt+1}, -\EE)$. Next, for $(E_1, \EE_1), (E_2, \EE_2) \in \Mod^{\coh}_{{A}}$ and $\phi \in Z^0\Mod^{\coh}_{{A}}(E_1, E_2)$, we define the cone of $\phi$, $C_{\phi} = (C_{\phi}^{\bt}, \EE_{\phi})$ by
\begin{equation*}
    C_{\phi}^{\bt}=\begin{pmatrix}
    E_2^{\bt}\\
    \oplus\\
    E_1[1]^{\bt}
    \end{pmatrix}
\end{equation*}
and 
\begin{equation*}
    C_{\phi}^{\bt}=\begin{pmatrix}
    \EE_2& \phi\\
    0 & -\EE_1^{\bt}
    \end{pmatrix}
\end{equation*}
Now we have a triangle of degree 0 closed morphisms 
\begin{equation}\label{triangle}
    \E\stackrel{\phi}{\to} F\to C_{\phi}\to E[1]
\end{equation} Under this construction, $\Mod^{\coh}_{{A}}$ is pre-triangulated, and $\ho(\Mod^{\coh}_{{A}})$ is triangulated with the collection of distinguished triangles being isomorphic to form \ref{triangle}.

    A degree 0 closed morphism $\phi \in \Mod^{\coh}_{{A}}(E_1, E_2)$ is a {\it homotopy equivalence} if it induces an isomorphism in $\ho(\Mod^{\coh}_{{A}})$. We will give a simple criterion to determine whether a map is a homotopy equivalence. Consider the following decreasing filtration
\begin{equation*}
    F^k\Mod^{\coh}_{{A}}(E_1, E_2) = \{ \phi \in \Mod^{\coh}_{{A}}(E_1, E_2)| \phi^i = 0  \text{ for $i < k$}\}
\end{equation*}
\begin{lem}
    There exists a spectral sequence
    $$
    E^{pq}_0 \Longrightarrow H^{p+q}(\Mod^{\coh}_{{A}}(E_1, E_2))
    $$
    where
    $$
    E^{pq}_0 = \gr\big(\Mod^{\coh}_{{A}}(E_1, E_2)\big) = \{\phi^p \in (\Mod^{\coh}_{{A}})^{p+q}(E_1, E_2):E_1^{\bt} \to E_2^{\bt+q}\otimes_{A_0} A^p\}
    $$
    with differential $d_0(\phi^p) = \EE_2\circ \phi^p - (-1)^{p+q}\phi^p \circ \EE_1$.
\end{lem}

\begin{prop}
    A closed morphism $\phi \in (\Mod^{\coh}_{{A}})^0(E_1, E_2)$ is a homotopy equivalence if and only if $\phi^0: (E_1^{bt}, \EE_1) \to (E_2^{bt}, \EE_2)$ is a quasi-isomorphism of complexes of $A^0$-modules.
\end{prop}

\begin{proof}
    Follows from \cite[Proposition 2.9]{Blo05}.
\end{proof}

\begin{defn} The $\infty$-category $\Modi_{{A}}^{\coh}$ of cohesive modules over ${A}$ is the $\infty$-category associated to the dg-category $\Mod_{{A}}^{\coh}$ under the dg-nerve
  $$\cohmodi_{\mathsf{A}} = \nerv_{\dg}(\cohmod_{\mathsf{A}}).$$
\end{defn}
\section{$\infty$-representations}
\subsubsection{$\infty$-representations of $L_{\infty}$-algebroids}

\begin{defn}
	The {\it Chevalley-Eilenberg algebra}\index{Chevalley-Eilenberg algebra} of $\g$ with coefficients in a dg $A$-module is the dga
	$$
	\CE(\g,E)= \Hom_A(\sym_A \g[-1],E)
	$$
	with differential given by
	$$
	(\del \alpha)(X_1,X_2,\cdots,X_n) 
	$$
\end{defn}

If $E=\g$, then we denote the $\CE(\g,\g)$ simply by $\CE(\g)$.

\begin{defn}
	Let $\g\in \lalgd_A$, an {\it $\infty$-representation}\index{$\infty$-representation!of a $\li$-algebroid} of $\g$ is a dg $\mathsf{A}$-module $E$, together with a $\Z$-connection
	\begin{equation*}
	\nabla: \CE(\g)\otimes_{A} E \to \CE(\g)\otimes_{A}E
	\end{equation*}
	of total degree one which is flat and satisfies graded Leibniz rule
	$$
	\nabla(\omega\eta)= d_A(\omega)\eta + (-1)^{|\omega|}\omega \nabla(\eta)
	$$
	for all $\omega \in \CE(\g)$, $\eta \in \CE(\g)\otimes_{A}E$.
\end{defn}

In literature, $\infty$-representations are also called {\it representations up to homotopy} or {\it sh-representations}. Denote the category of $\infty$-representations of an $L_{\infty}$-algebroid $\g \in \lalgd_A$ by $\repi_{\g,A}$

\begin{prop}
    Let $\g$ be an $L_{\infty}$-algebroid over $\mathsf{A}$, and $E$ a dg-$\mathsf{A}$-module, then an $\infty$-representation of $\g$
 on $E$ is equivalent to any of the following:
 \begin{enumerate}
    \item A quasi-cohesive module structure on $E$ over $\CE(\g)$.
    \item A quasi-cohesive module structure on $E^{\vee}$ over $\CE(\g)$.
    \item A square-zero degree 1 derivation $Q \in \der(\CA(\g, E))$ extending the differential $d_{\CA(\g)}$ on $\CA(\g)$.
    
     \item An Abelian extension $\g \oplus E$ of $\g$ along $E$:
     \begin{itemize}
         \item $L$ is an $L_{\infty}$-subalgebroid;
         \item E is an ideal, i.e. $\Tilde{l}_k(E, \cdots) \subset E$, where $\Tilde{l}_k$'s are the extension of $l_k$'s of $\g$;
         \item E is Abelian, i.e. $\Tilde{l}_k$ vanishes when evaluating at more than two elements of $E$.
     \end{itemize}
     
    \item The structure of a retract diagram of $L_{\infty}$-algebroids on $\g \to \g \oplus E \to \g$, which is square zero, i.e. all brackets vanish when evaluated on at least two elements of $E$.
    
     \item A collection of operations $[x_1, \cdots, x_n, -]: E \to E$ of degree $|x_1| + \cdots |x_n| + n - 2$ for $x_1, \cdots, x_n, e$ such that 
     \begin{align*}
         [x_{\sigma(1)}, \cdots, x_{\sigma(n)}, e] &= (-1)^{\sigma}[x_1, \cdots, x_n, e] \quad \sigma \in \Sigma_n\\
         [a\cdot x_1, \cdots, x_n, e] &= (-1)^{(n-1)a} a \cdot [x_1, \cdots, x_n, e]\\
         [ x_1, \cdots, x_n, a\cdot e] &= (-1)^{(n-1)a} a \cdot [x_1, \cdots, x_n, e] \quad n \ge 2\\
         [x_1, a\cdot e] &= a\cdot [x_1, e]+ x_1(a)\cdot e. 
     \end{align*}
     Here we ignore all Koszul signs due to permutations of variables. Moreover, these brackets determines a module structure, i.e.
     \begin{equation*}
         J^{n+1}(x_1, \cdots, x_n, a\cdot e) = 0
     \end{equation*}
     for all $n \ge 0$.

     \item An $L_{\infty}$-morphism $\g \to \At(E)$.
     
 \end{enumerate}\end{prop}
 
 \begin{proof}
 (1) is apparently equivalent to the definition of the $\infty$-representation.
 
 (1) $\Leftrightarrow$ (2): This is apparent from the construction of $(E^{\bullet \vee}, \mathbb{E}^{\vee}$, where 
 $$
 (\mathbb{E}^{\vee}\phi)(e) = d(\phi(e)) - (-1)^{|\phi|} \phi(\mathbb{E}(e)).
 $$
 
  (2) $\Leftrightarrow$ (3): Obvious.
  
  (3) $\Leftrightarrow$ (4): Note that $Q^2 = 0$ implies that $\g \oplus E$ is an $L_{\infty}$-algebroid, where the natural inclusion $\g \subset \g \oplus E$ is a subalgebroid. 
  
  (i) is obvious. (ii) and (iii) follows from the fact that $Q(E^{\vee}) \subset \CA(\g, E^{\vee})$.

  (4) $\Leftrightarrow$ (5): Denote the $m$-ary bracket $[\cdots]$ by $m_k$, then we simply set
  $$
  m_k(x_1, \cdots, x_{k-1}, e) = \Tilde{l}_k (x_1, \cdots, x_{k-1}, e).
  $$ The Jacobi identities follows from the $L_{\infty}$ structure on $\g \oplus E$.
  
    (5) $\Leftrightarrow$ (6) The $\infty$-representation of $\g$ on $E$ is equivalent to the data of a twisting cochain $\tau: \overline{C}_{*}(\g) \to \End_k(E)[1]$, where $\End_k(E)$ is the endomorphism Lie algebra of $E$, and $\tau$ is given by
    $$
    \sym^n_k \g[1] \to \End_k(E)[1]: x_1\otimes \cdots \otimes x_n \to [x_1, \cdots, x_n, -]
    $$
  
  is equivalent to the data
  
  $$
  (\rho, \tau): \overline{C}_{*}(\g) \to (T_{\mathsf{A}} \oplus \End_k(E))[1]
  $$
  which is graded $\mathsf{A}$-linear and takes values in the Atiyah Lie algebroid of $E$.
 \end{proof}
Hence, we see an $\infty$-representation of an $L_{\infty}$-algebroid $\g$ on $E$ is equivalent to a cohesive module structure on $E$ over $\CE(\g)$. Therefore, we get
\begin{lem}
There exists an equivalence of dg-categories
\begin{equation*}
    \cohmod_{\CE{(\g)}}\simeq \repi_{\mathsf{A}}^{\infty}(\g)
\end{equation*}
\end{lem}
Hence we will use cohesive modules and $\infty$-representations over $L_{\infty}$-algebroids interchangeably. In particular,  we call a cohesive module $E$ over an $L_{\infty}$-algebroid $\g$ when $E$ is a cohesive module over $\CE(\g)$. For simplicity, we will also call $E$ a $\g$-module if there is no confusion. 

\subsubsection{$\infty$-representations of simplicial sets}
For any $K_{\bullet}$ be a simplicial set, let $(C_{\bullet}(K_{\bullet}), \del, \Delta)$ be the dg coalgebra of simplicial chains on $K_{\bullet}$ over $k$ with the Alexander-Whitney coproduct $\Delta$. Consider the maps $\del'(x) = \sum_{i = 1}^{n - 1} (-1)^{i} K(d_i)$ and the reduced coproduct $\Delta'(x) = \Delta(x) - x\otimes 1 - 1 \otimes x$, we get dg coalgebra structure on $C_{\bullet}(K_{\bullet})$ and on the shifted graded module $s^{-1}C_{\bullet > 0}(K_{\bullet})$. 

We want to define a functor $\Lambda: s\set \to \dgcat_k$. For $K_{\bullet}$, define a dg-category $\Lambda(K_{\bullet})$, where the objects are $K_0$, and for any $x, y \in K_0$ we construct a chain complex $$
(\Lambda(K_{\bullet})(x, y), d_{\Lambda})$$ as follows: $\Lambda(K_{\bullet})(x, y)$ is the quotient of a free $k$-module generated by monomials 
$$(\sigma_1|\cdots|\sigma_k)$$
, where each 
$
\sigma_i \in s^{-1}C_{\bullet > 0}(K_{\bullet})
$ is a generator and satisfies $\max \sigma_i = \min \sigma_{i+ 1}$, by the equivalence relations generated by 
\begin{enumerate}
    \item $$(\sigma_1|\cdots|\sigma_k) \sim (\sigma_1|\cdots|\sigma_{i - 1}|\sigma_{i+1}|\sigma_k)$$
    if $\sigma_i$ is a degenerate 1-simplex for some $1\le i \le k$ and $k \ge 2$;
    \item $(\sigma_1|\cdots|\sigma_k) \sim 0$ if $\sigma_i \in C_{n_i}(K_{\bullet})$ is a degenerate simplex for some $1\le i \le k$, $n_i \ge 2$, and $k \ge 1$. Denote the equivalence class of $(\sigma_1|\cdots|\sigma_k)$ by $[\sigma_1|\cdots|\sigma_k]$.
\end{enumerate} Compositions are given by concatenations of monomials. The differential $d_{\Lambda}$ is given by extending $-\del' + \Delta'$ as a derivation on monomials. $d_{\Lambda}$ is then well-defined on equivalence classes and satisfies $d_{\Lambda}\circ d_{\Lambda} = 0$.

\begin{defn}
    We define the {\it dg nerve} functor $\n_{\dg}: s\set \to \dgcat_k$ by setting
    $$
        \n_{\dg}(\CC) := \Hom_{\dgcat_k}(\Lambda(\Delta^n), \CC)
    $$
\end{defn}
This definition agrees with Lurie's dg-nerve functor, hence $\Lambda$ is the left adjoint to  Lurie's dg-nerve, and $\nerv_{\dg}(\CC)$ is an $\infty$-category for any dg-category $\CC$. 

Let $\Ch_k$ denote the dg-category of chain complexes over a field $k$ of characteristic 0. Let $\CC$ be a dg-category and $\n_{\dg}\CC \in s\set$ its dg nerve.  
\begin{defn}
	An {\it $\infty$-representation}\index{$\infty$-representation!of a simplicial set} of $K_{\bullet}$ valued in $\CC$ is an $\infty$-functor $\mathtt{F}: K_{\bullet} \to \n_{\dg}\CC$, i.e. a morphism between the underlying simplicial sets. Denote $\repi_{\CC }(K_{\bullet}) = \funi(K_{\bullet}, \n_{\dg}\CC)$ the $\infty$-category of $\infty$-representations of $K_{\bullet}$ valued in $\CC$. 
\end{defn} 

The $n$-simplices of $\repi_{\CC }(K_{\bullet})$ are $\funi(\Delta^n \times K_{\bullet}, \n_{\dg}\CC)\simeq \dgcat_k\big(\Lambda(\Delta^n \times K_{\bullet}), \CC \big)$.

Let's look at the structure of an $\infty$-representation of a simplicial set. 

\begin{defn}
    Let $G_{\bullet}$ be a (derived) Lie $\infty$-groupoid, then an $\infty$-representation of $G_{\bullet}$ on a dg-category is defined as an $\infty$-representation of simplicial sets and all structure maps are required to be $\cinf$. We denote the category of $\infty$-representation of a (derived) Lie $\infty$-groupoid by $\repi_{\CC }(G_{\bullet})$. 
\end{defn}

\begin{lem}
    
\end{lem}

\begin{defn}
    We define an {\it $\infty$-local system}{\it $\infty$-local system} on a (derived) Lie ${\infty}$-groupoid $G_{\bullet}$ to be an $\infty$-representation of $G_{\bullet}$ valued in $\CC$.
\end{defn}

Note that the data of an $\infty$-local system is roughly a simplicial map from the simplicial set $G_{\bullet}$ to the dg-nerve of $\Ch_k$. By a Dold-Kan type correspondence, we can characterize the data of an $\infty$-local system as a dg-map between dg-categories. 

Let $K_{\bullet}$ be a Lie $\infty$-groupoid and $\CC$ a dg-category over $k$. Fix a map $F: K_0 \to \ob \CC$, i.e. a map on 0-simplices. Define $$\CC^{i,j}_F = \{f: K_i \to \CC^j | f(\sigma) \in \CC^j(F(\sigma_{(i)}), F(\sigma_{(0)})) \}$$
and 
$$
\CC^k_F(K_{\bullet}) = \bigoplus_{i+j = k, k\ge 0} \CC^{i, j}_F
$$
Now $C_F(K_{\bullet}) = \bigoplus_k \CC^k_F(K_{\bullet})$ forms a dga with differential $\hat{\delta}$ and product $\cup$ defined by

\begin{align*}
    (\hat{\delta} f^i)(\sigma_{i + 1}) =& \sum^i_{l = 1}(-1)^{l + |f^i|} f^i(\del_l(\sigma))\\
    (f\cup g)(\sigma_k) =& \sum_{t = 0}^k (-1)^{t|g^{k - t}|} f^t(\sigma_{(0\cdots t)}) g^{k - t}(\sigma_{(t\cdots k)})
\end{align*}

\begin{defn}
We define an $\infty$-local system to be a pair $(F, f)$ with $F: K_0 \to \ob \CC$ and $f\in C_F(K_{\bullet})$ which satisfies Maurer-Cartan equation, i.e. $f\in C^1_F(K)$ and 
\begin{align}\label{mc}
    &d_{F(\sigma_{(i)})} = f^0(\sigma_{(i)})\\
    &\hat{\delta} f + f\cup f = 0.
\end{align}
\end{defn}
\begin{rem}
    By a little abuse of notation, we will refer an $\infty$-local system $(F, f)$ simply by $F$. To avoid confusion about $F(x)$ and $f(x)$ for any zero simplices $x$, we will use $F_x$ to denote the former, and $F(x)$ to denote the latter.
\end{rem}
\begin{example}
    Let's take $G_{\bt}$ to be the smooth fundamental groupoid $\Pi^{\infty}(M)$ of a manifold $M$, and $\CC = \Ch_k$. Then the data of an $\infty$-local system consists of:
    \begin{enumerate}
        \item A graded vector space $E_x = \bigoplus_i E^i_x$ for $x \in M$.
        \item A sequence of $k$-cochains $f^k\in \Hom^{1-k}(E_{\sigma_{(k)}}, E_{\sigma_{(0)}})$ for $\sigma \in \Pi^{\infty}(M)_k$, which satisfies equation \ref{mc} (note the notation is a little different here).
    \end{enumerate}
\end{example}

We can put a dg-category structure on the category of $\infty$-local systems. For two $\infty$-local systems $F, G$ over $K_{\bullet}$ valued in $\CC$, define a complex of morphisms
$$
\Loc_{\CC}^{\dg}(K_{\bullet})(F,G) = \bigoplus_{i + j = k} \{\phi: K_i \to \CC^j| \phi(\sigma) \in \CC^j(F(\sigma_{(i)}), G(\sigma_{(0)}) \} 
$$
and a differential $D$ on it

$$
D \phi = \hat{\delta}\phi + G\cup \phi  - (-1)^{|\phi|}\phi \cup F.
$$
where $\phi = \sum_{i\ge 0} \phi^i$ with total degree $|\phi| = p$, and 
$$
(\hat{\delta}(\sigma_k) = \delta\circ T = \sum_{j = 1}^{k - 1}(-1)^{j + |\phi|} \phi^{k - 1}(\del_j(\sigma_k))
$$

This yields a dg-category $\Loc^{\dg}_{\CC}$, where the composition of morphisms is given by $\cup$. Denote the corresponding $\infty$-category $\loci^{\dg}_{\CC} = \Loc^{\dg}_{\CC}(W^{-1})$.

\begin{prop}
Given a Lie $\infty$-groupoid $K_{\bullet}$ and a dg-category $\CC$. There exists an equivalence of $\infty$-categories
\begin{equation*}
    \loci^{\dg}_{\CC}(K_{\bullet}) \simeq \repi^{\infty}_{\CC}(K_{\bullet})
\end{equation*}
\end{prop}
\begin{proof}
 See \cite{Smi11} Appendix.
\end{proof}

For pre-triangulated $\CC$, we can define shift and cone on $\Loc^{\dg}_{\CC}(K_{\bt})$.
First, let's define the shift functor. Let $F\in \Loc^{\dg}_{\CC}(K_{\bt})$, we define $F[i]$ by
\begin{align*}
    &F[i]_{(x\in K_0)} = F_x[i]\\
&F[i](\sigma_k) = (-1)^{i(k - 1)} F(\sigma_k)
\end{align*}
On morphisms, we define
$$
\phi[i](\sigma_k) = (-1)^{ik} \phi
$$

Next, we define the cone. Given a morphism $\phi \in \Loc^{\dg}_{\CC}(K_{\bt})(F, G)$ of total degree $i$. Define  
\begin{align*}
    C_{\phi}: K_0&\to \ob\CC\\
    x &\mapsto F[1-i]_x \oplus F_x
\end{align*}
and $c_{\phi} \in \CC^1_{C_{\phi}}$ by
$$
c_{\phi} = \begin{pmatrix}
F[1-i] & 0\\
\phi[1-i] & G
\end{pmatrix}
$$
\begin{rem}
    Note that $(C_{\phi}, c_{\phi})$ will not be an $\infty$-local system in general unless $\phi$ is closed.
\end{rem}

\begin{defn}
    Let $\phi\in \Loc^{\dg}_{\CC}(K_{\bt})(F, G)$ be a degree 0 closed morphism, we say $\phi$ is a {\it homotopy equivalence} if it induces an isomorphism in $\ho  \Loc^{\dg}_{\CC}(K_{\bt})$.
\end{defn}
Next, we will give an easy criterion to determine whether a map is a homotopy equivalence. Consider the following decreasing filtration
\begin{equation*}
    F^k\Loc^{\dg}_{\CC}(K_{\bt})(F, G) = \{ \phi \in \Loc^{\dg}_{\CC}(K_{\bt})(F, G)| \phi^i = 0 \text{for} i < k\}
\end{equation*}
\begin{lem}
    There exists a spectral sequence
    $$
    E^{pq}_0 \Longrightarrow H^{p+q}(\Loc^{\dg}_{\CC}(K_{\bt})(F, G))
    $$
    where
    $$
    E^{pq}_0 = \gr\big(\Loc^{\dg}_{\CC}(K_{\bt})(F, G)\big) = \{\phi: K_p \to \CC^q \phi(\sigma) \in \CC^q(F(\sigma_{(p)}), G(\sigma_{(0)})\}
    $$
    with differential $d_0(\phi^p) = d_G\circ \phi^p - (-1)^{p+q}\phi^p \circ d_F$.
\end{lem}
\begin{cor}
    The $E_1$-page of the above spectral sequence is a local system valued in graded vector space in the usual sense.
\end{cor}
Now we can give the criterion we want.
\begin{prop}
    For $\CC = \Ch_k$, a closed morphism $\phi \in \Loc^{\dg}_{\CC}(K_{\bt})^0(F, G)$ is a homotopy equivalence if and only if $\phi^0: (F_x, d_F) \to (G_x, d_G)$ is a quasi-isomorphism of complexes for all $x\in K_0$.
\end{prop}

\begin{proof}
    Follows from \cite[Proposition 2.9]{Blo05}.
\end{proof}
\section{Cohomology of derived Lie $\infty$-groupoids}
\subsection{Actions of derived Lie $\infty$-groupoids}
    Let's look at the action of derived Lie $\infty$-groupoids on a general space. First, let's recall the ordinary Lie groupoid action on a manifold. Let $G_{\bt}$ be a Lie groupoid acting on a manifold $M$. The data of this groupoid action is encoded in an {\it action groupoid}\index{action groupoid} $A_{\bt}$ and a groupoid morphism $\pi: A_{\bt} \to G_{\bt}$ over a $\cinf$ map $M \to G_0$, where $A_0 = M$ and
    \begin{equation*}
        A_{1}= M\times_{G_0, t} G_1 = \{(x, g): t(g) = \epsilon(y)
    \end{equation*}
    with structure maps $s(x, g) = xg, t(x, g) = x$. In fact we have a double pullback square
    
    	\begin{center}
	\begin{tikzcd}
		
        M\times_{G_0, t} G_1 \arrow[r,"\pr_2"] \arrow[d, red, shift right=1.5ex, "s"] \arrow[d, blue, "t"]& G_1 \arrow[d, red, "s"]\arrow[d, blue, shift left=1.5ex, "t"] \\
        		M \arrow[r,"\epsilon"]
        		&G_0
	\end{tikzcd}
\end{center}
\begin{lem}
    Kan fibrations between Lie groupoids are equivalent to the data of Lie groupoid actions.
\end{lem}
\begin{proof}
    $(\Leftarrow)$ Given a Lie groupoid action $G_{\bt}$ on $M$, we get a groupoid morphism $\pi: A_{\bt} \to G_{\bt}$. It suffices to show $\pi$ is a Kan fibration. $\kan(1,0)$ is equivalent to 
    $$A_1 \to M\times_{\epsilon, G_0, t} G_1
    $$
    which is an isomorphism by construction. By applying the inverse map, $\kan(1,0)$ is also satisfied. Higher Kan conditions follows from degree 1 case.
    
    $(\Rightarrow)$ Given a Kan fibration $\pi: A_{\bt} \to G_{\bt}$, we want to show there is an action of $G_{\bt}$ on $A_0$. By the $\kan(1,0)$ condition and unique Kan conditions for $n > 1$, we see that 
    $$A_1 \to A_0\times_{\epsilon, G_0, t} G_1$$
    Hence we can define $s(x, g) = xg, t(x, g) = x$, which gives us the desired data for  action.
\end{proof}
Therefore, this inspires us to define a higher groupoid action using Kan fibrations.

\begin{defn}
    Let $G_{\bt}$ be a Lie $\infty$-groupoid, then an {\it $\infty$-action}\index{$\infty$-action} of $G_{\bt}$ is a Kan fibration $\pi: A_{\bt} \to G_{\bt}$.
    
    If $G_{\bt}$ be a Lie $n$-groupoid, an $n$-action is a $n$-Kan fibration $\pi:A_{\bt} \to G_{\bt}$ of Lie $n$-groupoids.
\end{defn}
In \cite{Li15} the 2-groupoid case is shown to be the correct definition of actions.
\subsection{Derived Lie $\infty$-groupoid Cohomology}
    \begin{defn}
        Let $G_{\bt}$ be a derived Lie $\infty$-groupoid over a dga $A$. We denote by $C^{\bt}(G_{\bt})$ the smooth cochain complex on $G_{\bt}$, where $C^k(G_{\bt})$ consists of smooth functions on $G_k$, i.e.
        \begin{equation*}
            C^k(G_{\bt}) = \CO_{G_k}
        \end{equation*}
        The differential $d = \sum_i(-1)^i d^*_i$.
    \end{defn}

    Consider a derived Lie $\infty$-groupoid $G_{\bt}$ over a cdga $A$, and $E_{\bt}\in \Mod_A^{\dg}$. We form a dg-$C^{\bt}(G_{\bt})$-module $C^{\bt}(G_{\bt}; E_{\bt})$ whose degree $k$ part is 
    \begin{equation}
        C^k(G_{\bt}; E_{\bt}) = \bigoplus_{i+j = k}\Gamma(G_i;Q_0^*E_j)
    \end{equation}
    where $Q_0$ is defined as 
    $$
    Q_i := d_{1}\circ\cdots\circ d_i: G_i \to G_0
$$ is the projection on the last vertex. 
\begin{defn}
    We define a $\Z$-connection $$
    \BE: C^{\bt}(G_{\bt}; E_{\bt}) \to C^{\bt + 1}(G_{\bt}; E_{\bt})$$
    on $C^{\bt}(G_{\bt}; E_{\bt})$ to be a $k$-linear map of total degree one  which satisfies the grade Leibniz rule
    \begin{equation*}
        \mathbb{E}(e\omega) = (\mathbb{E}(e\otimes 1))\omega + (-1)^{|e|}ed\omega,
    \end{equation*}
    for $e\in C^{\bt}(G_{\bt}; E_{\bt}), \omega \in C^{\bt}(G_{\bt})$. 
\end{defn}

\begin{defn}
    We define an {\it $\infty$-representation}\index{$\infty$-representation! of a Lie $\infty$-groupoid} of $G_{\bt}$ to be a dg-$A$-module $E_{\bt}$ together with a flat $\Z$-connection on $C^{\bt}(G_{\bt}; E_{\bt})$.
\end{defn}
    We denote the resulting
    
We equip the category of $\infty$-representation a dg structure by defining the morphism complex
    
    \begin{equation*}
        \rep(G_{\bt})(F, G)^k = \bigoplus_{i+j = k} \Gamma\big(G_i, \HOM^j(F,G)\big)
    \end{equation*}
\begin{prop}
    Let $\CC = \Mod_A^{\dg}$, then there exists an dg equivalence  between $\Loc_{\infty}^{A}(G_{\bt})$ and $\rep^A_{\infty}(G_{\bt})$.
\end{prop}

\begin{proof}
    It's easy to see that $C^{\bt}(G_{\bt}; E_{\bt})$ is a $C^{\bt}(G_{\bt})$-module generated by $$\Gamma(G_{\bt}, E_{\bt}) = \bigoplus_{i+j = k} \Gamma(G_{k}, Q_1^* E_{j})$$
    The Leibniz rule for $\BE$ implies that we have a decomposition
    \begin{equation*}
        \BE = \BE^0+\BE^1 + \BE^2 +\cdots 
    \end{equation*}
    where $\BE^i \in  \Hom\big(\Gamma(G_{\bt}, E_{\bt}), \Gamma(G_{\bt+i}, E_{\bt+ 1-i})\big)$. For $i\not=1$, since $D_i$ is $C^{\bt}(G_{\bt})$-linear, we can identify them as an element of $\Gamma\Big(G_i, \Hom^{1-i}\big(P^*_i(E_{\bt}), Q^*_0(E_{\bt})\big)\Big)$, which  is exactly the $\CC^{i,1-i}_F$ we defined in $\infty$-local systems. For $i=1$, $\BE^1$ is a derivation, which can be identified as $\BE^1 = \hat{\delta} + \omega$ for some $\omega \in \Gamma\Big(G_1, \Hom^{0}\big(P^*_1(E_{\bt}), Q^*_0(E_{\bt})\big)\Big)$ and $\hat{\delta}$ is dual to the face map. Now the conditions for $E$ to be an  $\infty$-local system are
    \begin{enumerate}
        \item  $d_{F(\sigma_{(i)})} = f^0(\sigma_{(i)})$, which means $\mathbb{E_0}\in \Gamma\Big(G_0, \Hom^{1}\big(P^*_0(E_{\bt}), Q^*_0(E_{\bt})\big)\Big)$ is exactly the differential for $E_{\bt}$.\\
        \item $\hat{\delta} f + f\cup f = 0$ means $\BE\circ \BE = 0$.
    \end{enumerate}
    which are both satisfied by the construction. The dg structures on $\Loc_{\infty}^{A}(G_{\bt})$ and $\rep^A_{\infty}(G_{\bt})$ are exactly the same.
\end{proof}

\begin{defn}
    We define the {\it differentiable cohomology}\index{differentialble cohomology} of a derived Lie $\infty$-groupoid valued in $E$ to be
    \begin{equation*}
        H^{\bt}_{\diff}(G_{\bt}; E) = H^{\bt}\big(C^{\bt}(G_{\bt}; E_{\bt}), \BE \big)
    \end{equation*}
\end{defn}


\section{Chern-Weil theory for perfect dg modules}

\subsection{Chern-Weil theory for perfect dg modules}
Let $(E^{\bullet}, \mathbb{E})$ be a dg-$A^0$ module over a dga $A^0$ with a $\Z$-connection $\mathbb{E}$. In this section, we will develop a general theory of constructing characteristic classes valued in $A$. Later we shall apply it to singular foliations.

The curvature of the $\Z$-connection $\mathbb{E}$ is defined by the usual formula 
\begin{equation*}
    R_{\mathbb{E}} = \mathbb{E}^2 = \frac{1}{2}[\mathbb{E}, \mathbb{E}] \in A(\CF, \underline{\End}(E))
\end{equation*}

In order to define the characteristic forms, we need to define a $\Z$-graded supertrace map $\str: A(M, \underline{\End}(E)) \to \mathsf{A}$. For each 
$\phi_i \in \Gamma(\End(E^i))$, define $\str(\phi_i) = (-1)^i \tr(\phi_i)$. Extend $\str$ to a ($\Z$-graded) $A$-linear map we get  $\str:  A(\CF, \underline{\End}(E)) \to A(\CF)$. Note that by construction $\str$ vanishes on $A(\CF, \underline{\End}_i(E))$ for all $i \not= 0$. 

\begin{prop}[Bianchi identity]
$\mathbb{E}R_{\mathbb{E}}^{i} = 0$ for all $i \ge 1$.
\end{prop}
\begin{proof}
It follows from $\mathbb{E}R_{\mathbb{E}}^{i} = [\mathbb{E}, \mathbb{E}^{2i}] = 0$.
\end{proof}

\begin{lem}
$d_{\mathsf{A}}\str(R_{\mathbb{E}}^{i}) = 0$.
\end{lem}

\begin{proof}
\begin{equation*}
    d_{\mathsf{A}}\str(R_{\mathbb{E}}^{i}) = \sum_i \str(R_{\mathbb{E}}^{j - 1}[\mathbb{E}, R_{\mathbb{E}}]R_{\mathbb{E}}^{i-j})
\end{equation*}
where each summand is zero by the Bianchi identity. Hence, $\str(R_{\mathbb{E}}^{i})$ is closed.
\end{proof}
Let $f(z)$ be a convergent formal power series in $z$, then $f(R_{\mathbb{E}})$ is an element of $A^{even}(M, \underline{\End}(E))$, defined by
\begin{equation*}
    f(R_{\mathbb{E}}) = \sum \frac{f^{(k)}(0)}{k!} (\mathbb{E}^2)^k
\end{equation*}
Applying supertrace map to $f(R_{\mathbb{E}})$ we get an element in $\mathsf{A}$ which is a combination of even elements. We will call this element the {\it $\mathsf{A}$-characteristic form}, or simply characteristic form if there are no confusions, of $\mathbb{E}$ corresponding to $f(z)$.
\begin{prop}
Given a perfect $A^0$-module with $\Z$-connection $(E, \mathbb{E})$ over a dga $\mathsf{A}$. Then 
\begin{enumerate}
    \item The characteristic form $\str(f(R_{\mathbb{E}}))$ is a closed element of even degree.
    \item  (Transgression formaula) If $\mathbb{E}_t$ is a smooth 1-parameter family of $\Z$-connection on $E$, then 
    \begin{equation}
        \frac{d}{dt}\str (f(R_{\mathbb{E}_t})) = d \str\Big(\frac{d \mathbb{E}_t}{dt} f'(R_{\mathbb{E}_t})\Big)
    \end{equation}
    \item The cohomology class of  $\str(f(R_{\mathbb{E}}))$ in $H^{\bullet}(\mathsf{A})$ is independent of the choice of $\mathbb{E}$. 
    
\end{enumerate}

\end{prop}

\begin{proof}

We need the following lemma
\begin{lem}
For any $\alpha \in A(M, \underline{\End}(E))$, we have $d_{\mathsf{A}}(\str \alpha) = \str([\mathbb{E}, \alpha])$.
\end{lem}
\begin{proof}
Locally, we can write $\mathbb{E} = d + \omega$, so $\str([\mathbb{E}, \alpha]) = \str([d, \alpha]) + \str([\omega, \alpha])$. The second term vanishes by the definition of $\str$, and the first terms equals $\str(d\alpha)$.
\end{proof}
By this lemma, we have 
\begin{equation*}
    d \str(f(R_{\mathbb{E}})) = \str([\mathbb{E}, f(\mathbb{E}^2)]) = 0.
\end{equation*}
Hence, $\str(f(R_{\mathbb{E}}))$ is closed. The degree of it is clearly even.
We show prove (2) and (3) together. 
Let $\mathbb{E}_1$ and $\mathbb{E}_2$ be two $\Z$-connection over $\mathsf{A}$ on $E^{\bullet}$. $\mathbb{E}_1 -\mathbb{E}_2  $ is $\mathsf{A}$-linear hence we can write $\mathbb{E}_1 -\mathbb{E}_2 = \mathbb{E}'$ for some $\mathbb{E}' \in A(M, \underline{\End}(E))_1$. Set $\mathbb{E}_t = \mathbb{E}_2 + t \mathbb{E}'$. Note that $\mathbb{E}_t$ is a $\Z$-connection for any $t \in [0, 1]$. 
Let $\mathbb{E}_t$ be a smooth 1-parameter family of $\Z$-connections on $E$. The curvature of $\mathbb{E}_t$ is 
\begin{equation*}
    R_{\mathbb{E}_t} = (\mathbb{E}_2 + t \mathbb{E}')^2 = \mathbb{E}_2^2 + t[\mathbb{E}_2, \mathbb{E}'] + \frac{1}{2}t^2[\mathbb{E}', \mathbb{E}']
\end{equation*}

The deformation of curvature is

\begin{align*}
    \frac{d}{dt}R_{\mathbb{E}_t} = [\mathbb{E}_2, \mathbb{E}'] + t[\mathbb{E}', \mathbb{E}'] = [\mathbb{E}_t, \mathbb{E}']
\end{align*}
We will prove a more general lemma first.
\begin{lem}
Let $\alpha_t \in A^{even}(M, \underline{\End}(E))$ be a smooth family of forms of even total degree, then
\begin{equation*}
    \frac{d}{dt} \str(f(\alpha_t)) = \str(\frac{d}{dt}(\alpha_t) f'(\alpha_t))
\end{equation*}
\end{lem}
\begin{proof}
It suffices to consider $f$ being monomials. Consider $f(z) = z^n$. Then
\begin{equation*}
    \frac{d}{dt} \str(\alpha_t^n) = \str\Big(\sum_{i = 0}^{n - 1}\alpha_t^i\Big(\frac{d}{dt}\alpha_t\Big)\alpha_t^{n - i - 1} \Big) = n \str
    \Big(\Big(\frac{d}{dt}\alpha_t\Big)\alpha_t^{n- 1}\Big)
\end{equation*}

\end{proof}
Apply this lemma to $\alpha_t = \mathbb{E}_t^2$,
\begin{align*}
    \frac{d}{dt} \str(f(\mathbb{E}_t^2)) =& \str\bigg(\frac{d\mathbb{E}_t^2}{dt} f'(\mathbb{E}_t^2)\bigg) \\
    =&\str\bigg(\Big[\mathbb{E}_t,  \Big(\frac{d}{dt}\mathbb{E}_t^2\Big)f'(\mathbb{E}_t^2)\Big]\bigg)\\
    =& d \str\bigg(\frac{d\mathbb{E}_t}{dt} f'(\mathbb{E}_t^2)\bigg)
\end{align*}
This proves (2). 

Now we integrate the transgression formula with respect to $t$, then we get
\begin{equation}
    \str(f(R_{\mathbb{E}_1})) - \str(f(R_{\mathbb{E}_2})) = d\int_0^1 \str\big(\mathbb{E}' f'(\mathbb{E}_t^2)\big) dt
\end{equation}
Note that $\frac{d\mathbb{E}_t}{dt} = \mathbb{E}'$. Hence, the cohomology class of $\str(f(R_{\mathbb{E}_1}))$ and $\str(f(R_{\mathbb{E}_0}))$ in $H^{\bullet}(\mathsf{A})$ are the same.

\end{proof}
\begin{defn}
    Given a dg-$A^0$ module over a dga $A^0$ with a $\Z$-connection $\mathbb{E}$  $(E, \mathbb{E})$, we define the {\it $\mathsf{A}$-Pontryagin algebra of $E$} \index{Pontryagin algebra}
    \begin{equation*}
        \pont^{\bullet}_{\mathsf{A}} \subset H^{\bullet}(\mathsf{A})
    \end{equation*}
    to be the subalgebra generated by 
    \begin{equation*}
        \sigma^i_{\mathsf{A}}(E) = [ \str(R_{\mathbb{E}}^{i})] \in H^{2i}(\mathsf{A})
    \end{equation*}
    and we call $\sigma^i_{\mathsf{A}}(E)$ the {\it $\mathsf{A}$-Pontryagin character of $E$}\index{Pontryagin character}. 
\end{defn}

\subsection{$L_{\infty}$-pairs over a dga}

In this section, we will define $L_{\infty}$-pairs over a dga. As a special case, given a regular foliation $F$, then $(T_M, \CF)$ is an $L_{\infty}$-pair over $\cinf(M)$.

\begin{defn}
Let $\g$ be a $L_{\infty}$-algebroid over a dga $\mathsf{A}$, and $\h \subset \g$ a subalgebroid. Note the brackets $\{\lambda_i\}_i$ of $\h$ is the restriction of the brackets $\{\lambda_i\}_i$. We call $(\g, \h)$ an {\it $L_{\infty}$-pair}\index{$L_{\infty}$-pair}.
\end{defn}

Note that the inclusion map $\iota: \h \to \g$ gives $\iota^{\vee}:\g^{\vee} \to \h^{\vee}$, consequently we have a surjective morphism of dga $\iota^{\vee}:\sym \g^{\vee}[-1] \to \sym \h^{\vee}[-1]$
\begin{example}
    For a regular foliation $(M, \CF)$, $(T_M, \CF)$ is an  {\it $L_{\infty}$-pair}.
\end{example}
\begin{example}
    Let $\g$ be a Lie algebroid and $\h \subset \g$ a subalgebroid, then  $(\g, \h)$ is an {\it $L_{\infty}$-pair}, which is called a {\it Lie pair}. As a special case, if $\g$ is a Lie algebra and $\h \subset \g$ a Lie subalgebra, then 
    $(\g, \h)$ is an {\it $L_{\infty}$-pair} over a point. 
\end{example}

Now, for a Lie pair $(\g, \h)$, we denote the quotient $\g / \h$ by $N$.

\begin{lem}
There is an $\infty$-representation over $\h$ which gives $N$ an $\h$-module structure.
\end{lem}

\begin{proof}
There is an exact sequence of $dg$-$\mathsf{A}$-modules
\begin{equation*}
    0\longrightarrow \h \stackrel{\iota}{\longrightarrow} \g \stackrel{p}{\longrightarrow} N \to 0
\end{equation*}
The $\h$-module structure 
\begin{equation*}
    [x_1, \cdots, x_{n-1}, y]: \sym^{n-1}\g \otimes N \to N, n \ge 1
\end{equation*}is given by
\begin{equation*}
    [x_1, \cdots, x_{n-1}, y] = p \circ l_k(x_1, \cdots, x_{n-1}, y')
\end{equation*}
for any $y' \in \g$ such that $p(y') = y$. These brackets are well-defined, since $\h$ is a subalgebroid. By construction, $[\cdots]$ is an $\infty$-representation of $\h$ on $N$.
\end{proof}
\begin{lem}
$N^{\vee}$ is also an $\h$-module.
\end{lem}
\begin{rem}
    Note that $N^{\vee} = (\g / \h)^{\vee} \simeq \h^{\perp} = \ker (\iota^{\vee}:\g^{\vee} \to \h^{\vee})$. We denote the $\h$-module structure of $N^{\vee}$ by $(\h^{\vee}, \mathbb{E}_{\h}^{\h^{\perp}})$, where $\mathbb{E}_{\h}^{\h^{\perp}} = d_{\CO(\h)} + D^{\h^{\perp}}$.
\end{rem}
Let $d^{\dr}_{\g}$ and $d^{\dr}_{\h}$ be the algebraic de Rham operator on $\CO(\g)$ and $\CO(\h)$ respectively.
Define an operator $J: \CO(\g) \to A(\h, \g^{\vee})$ by $J = (\iota^{\vee} \otimes 1 \circ d^{\dr}_{\g}$, i.e. the following diagram commutes

\[\begin{tikzcd}
	{\CO(\g)} & {A(\g, \g^{\vee})} \\
	& {A(\h, \g^{\vee})}
	\arrow["{d_{\g}^{\dr}}", from=1-1, to=1-2]
	\arrow["J"', from=1-1, to=2-2]
	\arrow["{\iota^{\vee} \otimes 1}", from=1-2, to=2-2]
\end{tikzcd}\]
Hence $(1\otimes \iota^{\vee}) \circ J = \iota^{\vee} \circ d_{\g}^{\dr} = d_{\h}^{\dr} \circ \iota^{\vee}$.
It is obvious that $J$ is a derivation on the $\CA(\g)$-bimodule $\CA(\h, \g^{\vee})$, i.e. for all $\omega, \omega' \in \CA(\g)$,
$J(\omega \odot \omega') = \iota^{\vee}(\omega) \odot J(\omega') + J(\omega) \odot \iota^{\vee}(\omega') = \iota^{\vee}(\omega) \odot J(\omega') +  (-1)^{|\omega||\omega'|} \iota^{\vee}(\omega')\odot J(\omega)$. 

Now we get a map
\begin{equation}
    J\otimes 1: \CA(\g, \underline{\End}(E)) \to \CA(\h, \g^{\vee} \otimes \underline{\End}(E))
\end{equation}
by setting

\begin{equation}
    (J\otimes 1)(\phi \circ \psi) = (\iota^{\vee} \otimes 1)(\phi) \circ (J\otimes 1)(\psi) + (J\otimes 1)(\phi) \circ (\iota^{\vee} \otimes 1)(\psi)
\end{equation}
for all $\phi, \psi \in \CA(\g, \underline{\End}(E))$.

\begin{lem}
Let $\mathbb{E}_{\h}^{\h^{\perp}} = d_{\CA(\h)} + D^{\h, \h^{\perp}}$ be the $\h$-module structure on $N^{\vee} \simeq \h^{\perp}$, then for $\omega \in \ker (\iota^{\vee})$, we have $J(\omega) \in \CA(\h, \h^{\perp})$ and 
\begin{equation}
    \mathbb{E}_{\h}^{\h^{\perp}}(J(\omega)) = J(d_{\CA(\g)}(\omega))
\end{equation}

\end{lem}
\begin{proof}
First, we will show that $$D^{\h, \h^{\perp}}(\xi) = J(d_{\CA(\g)}\xi)$$ for all $\h^{\perp}$.
For $k = 2$
For $k >= 3$. 
\begin{align*}
    <D^{\h, \h^{\perp}}(\xi), p(l)>(a_1\odot \cdots \odot a_k) =& (-1)^{|\xi| + 1}<\xi, D^{\h, \h}(p(l))(a_1\odot \cdots \odot a_k)>\\
    =& (-1)^{|\xi| + |l|_{*_k} + k}<\xi, m_{k+1}(a_1\odot \cdots \odot a_k, p(l))>\\
    =&(-1)^{|\xi| + |l|_{*_k} + k}<\xi, p\circ l_{k+1}(a_1\odot \cdots \odot a_k, p(l))>\\
    =&(-1)^{k + 1}<(\iota^{\vee} \circ d^{\dr}_{\g} \circ l_{k+1}^{\vee})(\xi), p(l))>(a_1\odot \cdots \odot a_k)\\
    =& <(J\circ d_{\CA(\g)}(\xi), p(l)>(a_1\odot \cdots \odot a_k)
\end{align*}
where $*_k$ denotes $\sum_{i = 1}^k |a_i|$.

Now let us prove the proposition. It suffices to consider the elements of the form $\omega \odot \xi \in \CA(\g)\odot \h^{\perp}$. Applying the previous equation, we get
\begin{align*}
    J\circ d_{\CA(\g)}(\omega \odot \xi ) =& J\big( d_{\CA(\g)}(\omega) \odot \xi + (-1)^{|\omega|}\omega \odot d_{\CA(\g)}\xi\big) \\
    =& \iota^{\vee}(d_{\CA(\g)}(\omega) \odot \xi + J(d_{\CA(\g)}(\omega))\odot \iota^{\vee}(\xi) + \\
    & (-1)^{|\omega|}\big(J(\omega) \odot \iota^{\vee}(d_{\CA(\g)}\xi) + \iota^{\vee}(\omega) \odot J(d_{\CA(\g)}\xi) \big)\\
    =& d_{\CA(\g)}(\iota^{\vee}(\omega))\odot \xi + (-1)^{|\omega|} \iota^{\vee}(\omega)\odot D^{\h, \h^{\perp}}(\xi) = \mathbb{E}_{\h}^{\h^{\perp}}(J(\omega\odot \xi))
\end{align*}
\end{proof}

\subsection{Characteristic classes of singular foliations}

Let $(M, \CF)$ be a perfect singular foliation, i.e. $\CF$ is a perfect module. Let $(E_{\bullet}, \mathbb{E})$ be a cohesive module resolves $\CF$, i.e \begin{equation*}
0\to  E_{-n} \stackrel{d_n}{\longrightarrow}  \cdots \stackrel{d_2}{\longrightarrow} E_{-1} \stackrel{d_1}{\longrightarrow} E_{0} \stackrel{\rho}{\longrightarrow} \CF \longrightarrow 0
\end{equation*}

By similar method for holomorphic singular foliation, we can construct an $L_{\infty}$-algebroid structure on $E_{\bullet}$. Let $\CE(E_{\bullet}) = \sym^{\bullet} E^{\vee}_{\bullet}[-1]$ be Chevalley-Eilenberg algebra of the $L_{\infty}$-algebroid $E_{\bullet}$. 
\begin{prop}
Let $N\CF = T_M/\CF$ be the normal sheaf of the singular foliation $\CF$ which is perfect, then $N\CF$ is also perfect, and we have a normal complex
\begin{equation*}
0\to  E_{-n} \stackrel{d_n}{\longrightarrow}  \cdots \stackrel{d_2}{\longrightarrow} E_{-1} \stackrel{d_1}{\longrightarrow} E_{0} \stackrel{\rho^{}}{\longrightarrow} T_M \longrightarrow N\CF \longrightarrow 0
\end{equation*}
which resolves $N\CF$ and carries a Bott $\Z$-connection $\mathbb{B}$.
\end{prop}
\begin{proof}
     Directly follows from the $L_{\infty}$-pair $(T_M, E_{\bullet})$.
\end{proof}

\begin{cor}
Let $\CF$ be a perfect singular foliation and let $E_{\bullet}$ be an $L_{\infty}$-algebroid which resolves $\CF$. Then there exists an $L_{\infty}$-algebroid structure on $E_{\bullet}\oplus(E_{\bullet}[1]\to T_M)$ which is quasi-isomorphic to the tangent module $T_M$.
\end{cor}

\subsection{Atiyah class for $L_{\infty}$-algebroids}

In this section, we will construct the {\it Atiyah class}\index{Atiyah class} for an $L_{\infty}$-pair $(\g, \h)$. Let $(E, \mathbb{E}_{\h}^E)$ be a cohesive module over $\h$. Here $\mathbb{E}_{\h}^E = d_{\CO(\h)} + D^{\h, E}$ where $D^{\h, E}$ corresponds to the $\h$-module structure on $E$ which is an $\CO(\h)$-linear map $\CO(\h) \otimes E\to \CO(\h) \otimes E$.

Recall that $\h^{\perp} = N^{\vee}$ is also an $\h$-module, hence $\h^{\perp}\otimes \underline{\End}(E)$ inherits a $\infty$-representation over $\h$, with the $\Z$-connection defined by

\begin{equation}
    \mathbb{E}^{\h^{\perp}\otimes \underline{\End}(E)}_{\h} = d_{\CO(\h)} + D^{\h^{\perp}} + [D^{\h, E}, -]
\end{equation}
Denote the cohomology of the complex $(\h^{\perp}\otimes \underline{\End}(E), \mathbb{E}_{\h}^{\h^{\perp}\otimes \underline{\End}(E)})$ by $H^{\bullet}(\h, \h^{\perp}\otimes \underline{\End}(E))$.

By the surjectivity of the map $\iota^{\vee}$, we can lift $D^{\h, E} \in (\CO(\h) \otimes \underline{\End}(E))_1$ to an element $D^{\g, E} \in (\CO(\g) \otimes \underline{\End}(E))_1$. We get a $\Z$-connection $\mathbb{E}_{\g}^E = d_{\CO(\g)} + D^{\g, E}$. 

Note that $(E, \mathbb{E}_{\g}^E)$ is not necessarily a cohesive module, i.e. the curvature $R_{\mathbb{E}_{\g}^E}$ might not vanish.

We can easily calculate $R_{\mathbb{E}_{\g}^E} = d_{\CO(\h)} \circ D^{\g, E} + (D^{\g, E})^2$.

We have the following commutative diagram
\[\begin{tikzcd}
	E & {\CO(\g)\otimes E} & {\CO(\h)\otimes \g^{\vee}\otimes E} \\
	E & {\CO(\h)\otimes E} & {\CO(\h)\otimes \h^{\vee}\otimes E}
	\arrow["{R_{\mathbb{E}_{\g}^E}}", from=1-1, to=1-2]
	\arrow["{J\otimes 1}", from=1-2, to=1-3]
	\arrow["{=}"', from=1-1, to=2-1]
	\arrow["{R_{\mathbb{E}_{\h}^E}}"', from=2-1, to=2-2]
	\arrow["{\iota^{\vee}\otimes 1}"', from=1-2, to=2-2]
	\arrow["{1\otimes \iota^{\vee}\otimes 1}", from=1-3, to=2-3]
	\arrow["{d^{\dr}_{\h} \otimes 1}"', from=2-2, to=2-3]
\end{tikzcd}\]

which implies that
\begin{equation*}
    (1\otimes \iota^{\vee}\otimes 1)\circ (J\otimes 1) \circ R_{\mathbb{E}_{\g}^E} = 0
\end{equation*}

Therefore, we get an element $\alpha_{\mathbb{E}_{\g}^E}$ of total degree 2.

\begin{prop}
\begin{enumerate}
    \item $d_{\CO(h)}\alpha_{\mathbb{E}_{\g}^E} = 0$, hence we get a cocycle in the Chevalley-Eilenberg complex of $(\h^{\perp}\otimes \underline{\End}(E), \mathbb{E}^{\h^{\perp}\otimes \underline{\End}(E)})$.
    \item The cohomology class $[\alpha_{\mathbb{E}_{\g}^E}]$ in the 
    $L_{\infty}$-algebroid cohomology $H^{\bullet}(\h, \h^{\perp}\otimes \underline{\End}(E))$ is independent of the extension $\mathbb{E}_{\g}^E$. We call $[\alpha_{\mathbb{E}_{\g}^E}]$ the {\it Atiyah class} of the $L_{\infty}$-pair $(\g, \h)$ with respect to $E$. 
    \item For the canonical $\h$-module $(\g/\h)$, there is a canonical Atiyah class
    \begin{equation*}
        [\alpha^{\g/\h}] \in H^2(\h, \h^{\perp}\otimes \underline{\End}(\g/\h)) = H^2(\h,  \underline{\Hom}(\g/\h \otimes \g/\h, \g/\h)).
    \end{equation*}
\end{enumerate}
\end{prop}
\begin{proof}

First we need a lemma
\begin{lem}
Let $x \in \CA(\g, \underline{\End}(E))$ satisfy $(\iota^{\vee}\otimes 1)(x) = 0$, then
\begin{equation}
    [D^{\h, E}, (J\otimes 1)(x)] = (J\otimes 1)[D^{\g, E}, x].
\end{equation}
\end{lem}
\begin{proof}
It suffices to prove for $x$ homogeneous. We have 
\begin{align*}
    (J\otimes 1)[D^{\g, E}, x] =& (J\otimes 1)(D^{\g, E}\circ x - (-1)^{|x|}x \circ D^{\g, E})\\
    =& (\iota^{\vee} \otimes 1)(D^{\g, E}) \circ (J\otimes 1)(x) + (J\otimes 1)(D^{\g, E}) \circ (\iota^{\vee} \otimes 1)(x)\\
    &-(-1)^{|x|}\big( (\iota^{\vee} \otimes 1)(x) \circ (J\otimes 1)(D^{\g, E}) + (J\otimes 1)(x) \circ (\iota^{\vee} \otimes 1)(D^{\g, E})\big)\\
    =&D^{\h, E} \circ (J\otimes 1)(x)  -(-1)^{|x|}(J\otimes 1)(x) \circ D^{\h, E}\\
    =&[D^{\h, E}, (J\otimes 1)(x)].
\end{align*}
\end{proof}
Now, let us prove the proposition. By the previous commutative diagram $(\iota^{\vee} \otimes 1) \circ R_{\mathbb{E}_{\g}^E} = R_{\mathbb{E}_{\h}^E} = 0$, hence $(J\otimes 1)(R_{\mathbb{E}_{\g}^E}) \in \CA(\h, \h^{\vee}\otimes \underline{\End}(E))$. We have
\begin{align*}
\mathbb{E}_{\h}^{\h^{\perp}\otimes \underline{\End}(E)}(\alpha_{\mathbb{E}_{\g}^E}) =& (d_{\CO(\h)} + D^{\h^{\perp}} + [D^{\h, E}, -])((J\otimes 1)(R_{\mathbb{E}_{\g}^E})) \\
=&(d_{\CO(\h)} + D^{\h^{\perp}})((J\otimes 1)(R_{\mathbb{E}_{\g}^E})) + [D^{\h, E}, (J\otimes 1)(R_{\mathbb{E}_{\g}^E})] \\
=& (J\otimes 1)(d_{\CO(\g)} R_{\mathbb{E}_{\g}^E} + [D^{\g, E}, (R_{\mathbb{E}_{\g}^E})])
\end{align*}
where the last step follows from the Bianchi identity.

Next, let us look at (2). Consider another $\Z$-connection $\mathbb{E}' = d_{\CO(\g)} + {D^{\g, E}}'$ lifts the flat $\Z$-connection $\mathbb{E}^{\h, E}_{\h}$. Let $\omega = \mathbb{E}^{\g, E}_{\g} - \mathbb{E}' = D^{\g, E} - {D^{\g, E}}' \in \CA(\g, \underline{\End}(E))^1$, and we have $(\iota^{\vee} \otimes 1)(\omega) = 0$ and $(J\otimes 1)(\omega) \in \CA(\h, \h^{\perp} \otimes \underline{\End}(E))$, which implies $(J\otimes 1)(\omega^2) = 0$.  Now we have
\begin{align*}
    \alpha_{\mathbb{E}_{\g}^E} - \alpha_{\mathbb{E}'} =& (J\otimes 1)(R_{\mathbb{E}_{\g}^E} - R_{\mathbb{E}'})\\
    =& (J\otimes 1)(d_{\CA(\g)}(D^{\g, E}) + (D^{\g, E})^2 - d_{\CA(\g)}({D^{\g, E}}') - ({D^{\g, E}}')^2)\\
    =& (J\otimes 1)( d_{\CA(\g)} \omega  + \omega^2 +[{D^{\g, E}}', \omega])\\
    =& (d_{\CA(\h)} + D^{\h, \h^{\perp}})((J\otimes 1)(\omega) + [D^{\g, E}, (J\otimes 1)(\omega)]\\
    =& (d_{\CA(\h)} + D^{\h, \h^{\perp}} + [D^{\g, E}, -])((J\otimes 1)(\omega)) = \mathbb{E}_{\h}^{\h^{\perp}\otimes \underline{\End}(E)}((J\otimes 1)(\omega))
\end{align*}
which implies that the cohomology classes of $\alpha_{\mathbb{E}_{\g}^E} $ and $ \alpha_{\mathbb{E}'}$ are the same.

Finally, (3) follows from the standard identification $\h^{\perp} \simeq (\g / \h)^{\vee}$.

\end{proof}
Next, we shall construct the Atiyah classes from another way. Again, let $(\g, \h)$ be an $L_{\infty}$-pair. On $\h^{\vee}$, there exists a coadjoint $\h$-module structure $\mathbb{E}_{\h}^{\h^{\vee}} = d_{\CA(\h)} + D^{\h, \h^{\vee}}$ which is dual to the adjoint $\infty$-representation of $\h$ on itself. There is a natural $\h$-module structure on $\g$, and similarly on $\g^{\vee}$. We have a short exact sequence of $\h$-modules
\begin{equation*}
    0 \to \h^{\perp} \to \g^{\vee} \stackrel{\iota^{\vee}}{\to} \h^{\vee}\to 0
\end{equation*}
Now consider $(E, \mathbb{E})$ an $\h$-module, then we have 
\begin{equation*}
    0 \to \h^{\perp}\otimes \underline{\End}(E) \to \g^{\vee}\otimes \underline{\End}(E) \stackrel{\iota^{\vee}\otimes 1}{\to} \h^{\vee}\otimes \underline{\End}(E)\to 0
\end{equation*}
which induces a short exact sequence of dga's
\begin{equation*}
    0 \to \CA(\h,\h^{\perp}\otimes \underline{\End}(E)) \to \CA(\h, \g^{\vee}\otimes \underline{\End}(E)) \stackrel{1\otimes \iota^{\vee}\otimes 1}{\to} \CA(\h, \h^{\vee}\otimes \underline{\End}(E))\to 0
\end{equation*}
Hence we have a long exact sequence of $L_{\infty}$-algebroid cohomology
\begin{align*}
    &\cdots H^1(\h,\h^{\perp}\otimes \underline{\End}(E)) \to H^1(\h, \g^{\vee}\otimes \underline{\End}(E)) \stackrel{1\otimes \iota^{\vee}\otimes 1}{\to} H^1(\h, \h^{\vee}\otimes \underline{\End}(E))\\
    &\stackrel{\delta}{\to }H^2(\h,\h^{\perp}\otimes \underline{\End}(E))\to H^2(\h, \g^{\vee}\otimes \underline{\End}(E)) \stackrel{1\otimes \iota^{\vee}\otimes 1}{\to} H^2(\h, \h^{\vee}\otimes \underline{\End}(E))\cdots
\end{align*}
\begin{lem}
    The element $(d^{\dr}_{\h} \otimes 1)(D^{\h, E}) \in \CA(\h, \h^{\vee}\otimes \underline{\End}(E))$ is a degree 1 cocycle.
\end{lem}

\begin{proof}
     Since $d^{\dr}_{\CA(\h)}$ is a derivation on $\CA(\h, \h^{\vee})$, $d^{\dr}_{\CA(\h)}\otimes 1$ is a derivation on $\CA(\h, \h^{\vee}\otimes \underline{\End}(E))$. Hence,
     \begin{align*}
         (d^{\dr}_{\CA(\h)}\otimes 1)({(D^{\h, E})}^2) =& (d^{\dr}_{\CA(\h)}\otimes 1)(D^{\h, E}) \circ D^{\h, E} + D^{\h, E}\circ (d^{\dr}_{\CA(\h)}\otimes 1)(D^{\h, E})\\
         =& [D^{\h, E}, (d^{\dr}_{\CA(\h)}\otimes 1)D^{\h, E}]
     \end{align*}
     It is easy to verify that
     $$(\mathbb{E}_{\h}^{\h^{\vee}} \otimes 1)(d^{\dr}_{\CA(\h)}\otimes 1)(D^{\h, E}) = (d^{\dr}_{\CA(\h)}\otimes 1)(d_{A(\h)}D^{\h, E})
     $$
     
     Now
     \begin{align*}
         \mathbb{E}_{\h}^{\h^{\vee}\otimes \underline{\End}(E)}((d^{\dr}_{\CA(\h)}\otimes 1)(D^{\h, E})) &= (\mathbb{E}_{\h}^{\h^{\vee}} + [D^{\h, E}, -])((d^{\dr}_{\CA(\h)}\otimes 1)(D^{\h, E}))\\
         &=(d^{\dr}_{\CA(\h)}\otimes 1)(d_{A(\h)}D^{\h, E}) + (d^{\dr}_{\CA(\h)}\otimes 1)({(D^{\h, E})}^2)  \\
         &= (d^{\dr}_{\CA(\h)}\otimes 1)(d_{A(\h)}D^{\h, E} + {(D^{\h, E})}^2) = 0
     \end{align*}
     since $d_{A(\h)}D^{\h, E} + {(D^{\h, E})}^2 = 0$ which is the Maurer-Cartan equation.
\end{proof}

\begin{prop}
    The cohomology class $\delta[(d^{\dr}_{\h} \otimes 1)(D^{\h, E})] \in H^2(\h,\h^{\perp}\otimes \underline{\End}(E))$ is the same as the Atiyah class $\alpha^E_{\g, \h}$.
\end{prop}
\begin{proof}
     First we choose an element $\beta \in \CA(\h, \g^{\vee}\otimes \underline{\End}(E))$ of degree 1 such that $(1\otimes \iota^{\vee} \otimes 1)(\beta) = (d^{\dr}_{\h} \otimes 1)(D^{\h, E})$. Then we get $\alpha \in \CA^2(\h, \h^{\perp} \otimes \underline{\End}(E))$ by
     \begin{equation*}
         \alpha = \mathbb{E}_{\h}^{\g^{\vee}\otimes \underline{\End}(E)}(\beta) = (d_{\CA(\h)} + D^{\g^{\vee}}+[D^{\h, E}, -])(\beta)
     \end{equation*}
     whose cohomology class is the one given by $\delta[(d^{\dr}_{\h} \otimes 1)(D^{\h, E})]$. We want to show that $[\alpha]$ actually agrees with the Atiyah class $[\alpha^E_{\g, \h}]$. We will do this by extending the $\Z$-connection $\mathbb{E}_{\h}^E$ given by the $\h$-module structure of $E$ and extend it to a $\Z$-connection $\mathbb{E}^E_{\g}$ over $\g$, and show that the resulting Atiyah cocycle $\alpha^E_{\mathbb{E}^E_{\g}}$ coincides with $\alpha$.
     
     First, we want to find an element $D^{\g, E} \in \CA(\g, \underline{\End}(E))$ with $(J\otimes 1)(D^{\g, E}) = \beta$ and $(\iota^{\vee} \otimes 1)(D^{\g, E}) = D^{\h, E}$. By surjectivity of $J$, we can find some   $K^{\g, E} \in \CA(\g, \underline{\End}(E))$ such that $(J\otimes 1)(K^{\g, E}) = \beta$. Now 
     \begin{align*}
         (d^{\dr}_{\CA(\h)}\otimes 1)(\iota^{\vee} \otimes 1)(K^{\g, E}) =& (1\otimes \iota^{\vee} \otimes 1)(J\otimes 1)(K^{\g, E})\\
         =&(1\otimes \iota^{\vee} \otimes 1)(\beta)\\
         =&(d^{\dr}_{\h} \otimes 1)(D^{\h, E})
     \end{align*}
e     Now $(\iota^{\vee} \otimes 1)(K^{\g, E}) - D^{\h, E} = \phi$ for some $\phi \in \Gamma(\underline{\End}(E))$ as $\ker d^{\dr}_{\CA(\h)} \simeq \cinf(M)$. Hence, we could let $D^{\g, E} = K^{\h, E} - \phi$.

Now 

\begin{align*}
    \alpha^E_{\mathbb{E}^E_{\g}} =& (J\otimes 1)(d_{\CA(\g)}D^{\g, E} +(D^{\g, E})^2)\\
    =& (\mathbb{E}_{\h}^{\g^{\vee}}\circ (J\otimes 1) (D^{\g, E}) + [(\iota^{\vee}\otimes 1)(D^{\g, E}),(J\otimes 1) (D^{\h, E})]\\
    =& (d_{\CA(\h)} + D^{\g^{\vee}})(J\otimes 1) (D^{\g, E}) +[D^{\h, E},(J\otimes 1) (D^{\h, E})]\\
    =& (d_{\CA(\h)} + D^{\g^{\vee}}+[D^{\h, E}, -])(\beta)=\alpha
\end{align*}
Hence, $[\alpha]$ agrees with $[\alpha^E_{\mathbb{E}^E_{\g}}]$.
\end{proof}

Let $\h$ be an $L_{\infty}$-algebroid and $(E, \mathbb{E})$ a $\h$-module, we have the canonical Abelian extension $( \g = \h \oplus E, \mathbb{E})$ of $\h$ along $E$, which induces an $L_{\infty}$-pair $(\g, \h)$. The Atiyah class $\alpha_{\mathbb{E}}$ is trivial in this case. Hence, we see that the Atiyah class measures the nontriviality of the extension of $\h$ to $\g$. 

Next, we are going to see some simple situation when does Atiyah classes vanish.

\begin{prop}
Let $(\g, \h)$ be an $L_{\infty}$-pair and $(E, \mathbb{E} = d_{\CA(h)} + D^{\h, E})$ be an $\h$-module, then the Atiyah class $[\alpha_{\mathbb{E}_{\g}^E}] \in H^{\bullet}(\h, \h^{\perp}\otimes \underline{\End}(E))$ vanishes if any of the following equivalent condition is met:
\begin{enumerate}

    \item There exists a $\Z$-connection $\mathbb{E}'$ over $\g$ on $E$ extending $\mathbb{E}$ such that the Atyiah cocycle $\alpha_{\mathbb{E}_{\g}^E}$ relative to $\alpha_{\mathbb{E'}}$ vanishes.
    
    \item There exists a degree 1 cocycle $\phi \in \CA(\h, \g^{\vee}\otimes \underline{\End}(E))$ such that $(1\otimes \iota^{\vee} \otimes 1)(\phi) = (d_{\h}^{\dr} \otimes 1)(D^{\h, E})$.
    
    \item There exists an $\h$-module morphism $\{\phi_k: \sym \h \otimes \g \to \underline{\End}(E)[1]\}_{k \ge 0}$ from $\g$ to $\underline{\End}(E)[1]$ extending the canonical $\h$ module morphism $\{\phi_k^{\h, E}\}_{k \ge 0}$ from $\h$ to $\underline{\End}(E)[1]$, i.e.
    \begin{equation*}
        (\phi_k)\circ (1\otimes j) = \phi_k^{\h, E} : \sym \h \otimes \g \to \underline{\End}(E)[1]
    \end{equation*}

\end{enumerate}
\end{prop}
\begin{proof}
     (1) $\Rightarrow$ (2): Clearly (3) implies that the Atiyah class of $E$ vanishes. Hence  $\delta[(d^{\dr}_{\h} \otimes 1)(D^{\h, E})] = 0$. We can find a degree 1 element $\tilde{\phi}\in \CA(\h, \g^{\vee}\otimes \underline{\End}(E))$ such that  $[(1\otimes \iota^{\vee} \otimes 1)(\tilde{\phi}] = [(d^{\dr}_{\h} \otimes 1)(D^{\h, E})]$. Then we can find an element $\beta \in \CA(\h, \h^{\vee}\otimes \underline{\End}(E))$ of degree 0 such that $\mathbb{E}^{\h^{\vee}\otimes \underline{\End}(E)}_{\h}(\beta) = (1\otimes \iota^{\vee} \otimes 1)(\tilde{\phi}) - (d^{\dr}_{\h} \otimes 1)(D^{\h, E})$, and therefore a $\gamma \in \CA^1(\h, \g^{\vee}\otimes \underline{\End}(E))$ maps to $\h$, i.e. $(1\otimes \iota^{\vee} \otimes 1)(\gamma) = \beta$.
     
     Now we let $\phi = \Tilde{\phi} - \beta$, by an easy calculation, we have 
     $(1\otimes \iota^{\vee} \otimes 1)(\phi) = (d_{\h}^{\dr} \otimes 1)(D^{\h, E})$.
     
     (2) $\Rightarrow$ (1): Given $\phi \in \CA^1(\h, \g^{\vee}\otimes \underline{\End}(E))$, we can find a $D^{\g, E}\in \CA(L, \underline{\End}(E)))$ such that $(J\otimes 1)(D^{\g, E}) = \beta$ and $(\iota^{\vee} \otimes 1)(D^{\g, E}) = D^{\h, E}$. Now $\mathbb{E}' = d_{\CA(\g)} + D^{\g, E}$ is a $\Z$-connection extending $\mathbb{E}$, and the associated Atiyah cocycle

     \begin{align*}
         \alpha^E_{\mathbb{E}'} =& (J\otimes 1)(R_\mathbb{E}')\\
         =& (J\otimes 1)(d_{\CA(\g)}D^{\g, E} +(D^{\g, E})^2 )\\
         =& \mathbb{E}_{\h}^{\g^{\vee}}( (J\otimes 1) (D^{\g, E})) + [(\iota^{\vee}\otimes 1)(D^{\g, E}),(J\otimes 1) (D^{\g, E})]\\
         =& \mathbb{E}_{\h}^{\g^{\vee}} (\phi) + [D^{\h, E},\phi] = 0
     \end{align*}

     (2) $\Leftrightarrow$ (3): Note that $\phi \in \CA^1(\h, \g^{\vee}\otimes \underline{\End}(E))$ consists of a family of map $\phi_k: \sym \h \otimes \g \to \underline{\End}(E)[1]$.

\end{proof}
\subsection{Scalar Atiyah classes and Todd classes}
Let $(E, \mathbb{E})$ be an $\h$-module and $(\g, \h)$ an $L_{\infty}$-pair. We define the {\it scalar Atiyah classes}\index{scalar Atiyah class} of the $L_{\infty}$-pair to be 
\begin{equation*}
    c_k(\g, \h) = \frac{1}{k!}\bigg(\frac{i}{2\pi}\bigg)^k \str(\alpha_{\g, \h})\in H^{k}(\h, \h^{\perp})
\end{equation*}

Let $\ber: \Gamma(\underline{\End}(E)) \to \cinf(M)$ be the {\it Berezinian map}\index{Berezinian map} (superdeterminant), then we define the {\it Todd class}\index{Todd class} of an $L_{\infty}$-pair $(\g, \h)$ to be 
\begin{equation*}
    \td_{\g, \h} = \ber\Bigg(\frac{\alpha_{(\g, \h)}}{1 - e^{-\alpha_{(\g, \h)}}} \Bigg) \in \bigoplus_{k\ge 0} H^{k}(\h, \h^{\perp})
\end{equation*}
\begin{example}
    Let $X$ be a compact Kahler manifold. Consider the $L_{\infty}$-pair $(T_{\C}X, T^{0,1}_X)$, then the natural map of sheaf cohomology $\bigoplus_k H^k(X, \Omega_X^k) \to \bigoplus_k H^{2k} (X, \C)$ sends the scalar Atiyah classes $c_k(T_X^{1,0})$ and the Todd class $\td_{T_X^{1,0}}$ of the $L_{\infty}$-pair $(T_{\C}X, T^{0,1}_X)$ to the $k$-th Chern characters $\ch_k(X)$ and the Todd class $\td_X$ of $X$ respectively.
\end{example}

\subsection{Infinitesimal ideal system of $L_{\infty}$-algebroids}

In this section, we will define infinitesimal ideal systems associated to an $L_{\infty}$-pair, and show that there is a natural infinitesimal ideal systems associated to an $L_{\infty}$-algebroid fibration. The infinitesimal ideal system structure is actually related to the vanishing to Atiyah classes.

\begin{defn}
Consider an $L_{\infty}$-pair $(\g, \h)$, we define an {\it infinitesimal ideal system}\index{infinitesimal ideal system} to be a triple $(\CF_M, \h, \mathbb{E})$ such that $\CF_M \subset T_M$ is an involutive locally free subsheaf of the tangent sheaf $T_M$, $\rho(\h) \subset \CF_M$, and $\mathbb{E}$ is a flat $\Z$-connection over $\CF_M$ on $\g/\h$ which satisfies
\begin{enumerate}
    \item If $g\in \g$ is $\mathbb{E}$-flat, then $[g, h_1,\cdots, h_{i-1}] \in \h$ for all $h_i \in h$ and all $i$-brackets for $i\ge 2$.
    \item If $g_1, \cdots, g_i \in \g$ are $\mathbb{E}$-flat, then $[g_1, \cdots, g_i]$ is also $\mathbb{E}$-flat.
    \item If $g\in \g$ is $\mathbb{E}$-flat, then $\rho(g)$ is $\nabla^{F_M}$-flat, where
    $\nabla^{F_M}$ is the Bott connection on $T_M/\CF_M$.
\end{enumerate}
\end{defn}
This is a direct generalization of infinitesimal ideal systems in Lie algebroids.
\begin{prop}
 For any $L_{\infty}$-algebroid fibration $\phi:\g \to \h$ over a $\cinf$-map $f:M \to N$, there exist an infinitesimal ideal system $(\CF, \h, \mathbb{E})$ associated to it.
\end{prop}


\begin{defn}
We define the Atiyah class of an infinitesimal ideal system $(\CF_M, \h. \mathbb{E})$ in an $L_{\infty}$-algebroid $\g$ to be the Atiyah class of the flat $\Z$-connection $\mathbb{E}$.
\end{defn}
\begin{prop}
Let $(\g, \h)$ be an $L_{\infty}$-pair on $M$. If there exists an infinitesimal ideal system $(\CF_M, \h. \mathbb{E})$ in $\g$, such that the quotient $(\g/\h)/\mathbb{E} \to M/\CF_M$ exists and is smooth, then the Atiyah class of the infinitesimal ideal system
vanishes. 
\end{prop}
\begin{example}[Simple foliations]
    For example, if $\CF$ is a simple foliation, i.e. the leaf space of $\CF$ is a manifold, then the Atiyah class associated to $L_{\infty}$-pair $(TM, \CF)$ vanishes.
\end{example}


\part{Singular foliations and $L_{\infty}$-algebroids}
\section{Singular foliation and their homotopy theory}
\subsection{Foliations}
A foliation is a partition of a manifold into immersed submanifolds.
\begin{defn}[\cite{MM03}]
    Let $M$ be a smooth manifold. A (regular) {\it foliation}\index{foliation} $\CF$ of codimension $q$ on $M$ can be described in the following equivalent data:
    \begin{enumerate}
        \item A {\it foliation atlas}\index{foliation atlas} $\{\phi_i: U_i\to \R^{n-q}\times \R^q\}$ of $M$ for which the change-of-coordinates diffeomorphisms $\phi_{ij}$'s are globally of the form
        $$
        \phi_{ij}(x,y) = \big(g_{ij})(x,y), h_{ij}(y) \big)
        $$
        with respect to the decomposition $\R^n = \R^{n-q} \times \R^q$, where $n = \dim M$.  Note that each leaf is partitioned into {\it plaques}\index{plaques}, which are connected components of the submanifolds $\phi^{-1}_i(\R^{n-1}\times \{y\}), y \in \R^q$. The plaques globally glue to {\it leaves}\index{leaf}, which are immersed submanifolds of $M$. We call the first $n-q$ directions in the decomposition the {\it leaf directions}, and the last $q$ directions the {\it transversal directions}.
        \item An open cover $\{U_i\}$ of $M$ with submersions $s_i : U_i \to \R^q$ such that there are diffeomorphisms 
        $$
        \gamma_{ij}: s_j(U_i\cap U_j) \to s_i(U_i\cap U_j)
        $$
        with $\gamma_{ij}\circ s_j|_{U_i\cap U_j} = s_i|_{U_i\cap U_j}$, which is necessarily unique. $\gamma_{ij}$'s satisfy the cocycle condition $\gamma_{ij}\circ \gamma_{jk} = \gamma_{ik}$, which is called the {\it Haefliger cocycle}\index{Haefliger cocycle} representing $\CF$.
        \item An integrable sub-bundle $F$ of $TM$ of rank $n-q$, i.e. for any $X, Y \in \Gamma(F)$, $[X, Y]\in \Gamma(F)$. We usually denote $\Gamma(F)$ by $\CF$ and $F$ by $T\CF$.
        \item A locally trivial {\it differential ideal}\index{differential ideal} $\mathcal{J} = \oplus_{k=1}^n \mathcal{J}^k$ of rank $q$ in the de Rham dga $\Omega^{\bt}(M)$. 
    \end{enumerate}
	When there is a foliation on a manifold $M$, we denote by $(M, \CF)$ a {\it foliated manifold}\index{foliated manifold}.
\end{defn}

The first two conditions are descriptions of $\CF$ by local charts, which tells us that locally a foliation is decomposed into two distinction directions: the leaf direction and the transversal direction. We denote $M/\CF$ the {\it leaf space} by quotienting equivalence relations such that $x\sim y$ if $x$ and $y$ lie on the same leaf.

\begin{defn}
	Let $M$ be a smooth manifold. We define a {\it complex foliation} to be an involutive sub-bundle of the complexified tangent bundle $T_{\C}M = TM\otimes_{\R} \C$.
\end{defn}
A complex foliation $\CF$ is called {\it real} if and  $\overline{T\CF} = T\CF$. In this case, we can define a real foliation $T\CF_{\R} = T\CF \cap TM$.
\begin{defn}
	A map between two foliated manifolds $f:(M_1, \CF_1)\to (M_2, \CF_2)$ is called {\it foliated} if it preserves the foliated structure, i.e. $f$ maps leaves of $(M_1, \CF_1)$ into leaves of $(M_2, \CF_2)$.
\end{defn}
We denote $\mfd^{\fol}$ the category of foliated manifolds where the morphisms are foliated maps, and $\mfd^{\fol}_{\C}$ the category of complex-foliated manifolds. 
\begin{example}
    Let $M$ be a smooth manifold. The two most simple foliation on $M$ is given by (1) foliation by whole manifold, i.e. $\CF = T_M$, and (2) foliation by points , i.e. $\CF = M\times \{0\}$. For the first case, the transversal direction is trivial, i.e. 0-dimensional. For the second case, the transversal direction is the whole manifold, which is again foliated by points. Hence, we see that the transversal direction always locally admits a foliation by points.
\end{example}

\begin{example}[Product foliations]
Given two foliations $(M_1, \CF_1)$ and $(M_2, \CF_2)$, we can form their product $(M_1\times M_2, \CF_1\times \CF_2)$, where $\CF_1\times \CF_2$ is given by $F_1\times F_2\subset TM_1\times TM_2 = T(M_1\times M_2)$.
\end{example}

We say $f$ is {\it transversal}\index{transeversal} to $\CF$ if $f$ is transverse to all leaves of $\CF$ in the image of $f$, i.e. for any $x\in N$,
$$
(df)_x(T_xN) + T_{f(x)}(\CF) = T_{f(x)} M
$$
\begin{example}[Pullback foliations]
Let $(M, \CF)$ be a manifold with foliation $\CF$. Consider a smooth map $f: N\to M$ transversal to $\CF$. It is not hard to show that the pull-back $f^*\CF$ is a foliation on $N$. This example will be important when we consider the 'derived' counterpart of the foliation.
\end{example}
\begin{example}[Flat bundles]
    Let $G$ be a group acting freely and properly discontinuesly on a connected manifold $\tilde{M}$ and $\Tilde{M}/G = M$. For example, we can take $\tilde{M}$  be the universal cover of $M$ and $G = \pi_1(M, x)$ for some $x\in M$. Here we let $G$ be a right action on $\tilde{M}$. Suppose $G$ also acts on the left on some manifold $F$, then we can form a quotient space $E = \tilde{M}\times_G F$ from the product space $\tilde{M}\times F$ by identifying $(yg,z)\sim (y, gz)$. It is easy to show that $E$ is a manifold, and we have the following commutative diagram
    \begin{center}
       \begin{tikzcd}
    \tilde{M}\times F  \arrow[r] \arrow[d, "\pr_1"]
    
    & E  \arrow[d, "\pi"] \\
     \tilde{M} \arrow[r]
    & M
\end{tikzcd}
    \end{center}
    The projection $\pr_1$ induces a submersion $\pi$, which gives $E$ a fiber bundle structure over $M$ with fiber $F$.
    
    This construction also produces a foliation $\CF(\pr_2)$ on $\tilde{M}\times F $, which is given by the submersion $\pr_2: \tilde{M}\times F \to p$. By the construction, $\CF(\pr_2)$ is $G$-invariant, hence we can form a quotient foliation $\CF = \CF(\pr_2)/G$ on $E$. 
    
    We can look into more detail about leaves of $\CF$. Let $x \in F$, and $G_x\subset G$ the isotropy group of the $G$-action at $x$, then the leaf of the foliation $F$ associated to $\tilde{M}\times \{z\}$ is diffeomorphic $\Tilde{M}/G_x$.

\end{example}
\subsection{Singular foliation}
We also often see foliated structures with singularities, i.e. the dimensions of leaves are not constant.

\begin{defn}[\cite{AC09}, \cite{LLS20}]
	Let $M$ be a smooth manifold. A {\it singular foliation}\index{singular foliaiton} $\CF$ on $M$ is a locally finitely generated subsheaf of $\CO_M$-modules of the tangent sheaf $T_M$ which is involutive, i.e. closed under Lie brackets.
\end{defn}
Equivalently, we can characterize a singular foliation $\CF$ as a locally finitely generated $O_M$-submodule of $\Gamma(TM)$. Clearly, regular foliations are singular foliations, since sub-bundles of $TM$ are finitely generated. By a result of Hermann \cite{Her60}, a singular foliation  on $M$ induces a partition of $M$ into leaves.

A {\it singular sub-foliation}\index{singular sub-foliaiton} $\CF'$ of a singular foliation $\CF$ is a singular foliation such that, for all open sets $U\subset M$, we have $\CF'(U) \subset \CF(U)$. 

\begin{defn}
    Let $(M, \CF)$ be a singular-foliated manifold and $x \in M$. The {\it tangent space of the leaf}\index{tangent space!of the leaf} at $x\in M$ is the image $F_x$ of $\CF$ in $T_xM$. The {\it fiber}\index{fiber!of a singular foliation} of $\CF$ at $x$ is $\CF_x = \CF / I_x \CF$, where$  I_x = \{f\in \CO_M:f(x) = 0\}$.
\end{defn}
Let $\ev_x:\CF \to T_x M$ be the evaluation map of $\CF$ at $x$. Clearly $\ev_x$ vanishes on $I_x \CF$, therefore it descends to a map $\tilde{\ev_x}: \CF_x \to F_x \subset T_x M$. $\ker {\ev_x}$ is a Lie subalgebra of $\CF$ and $I_x \CF$ is an ideal in this Lie algebra. It follows that $\ker \tilde{\ev_x} = \ker {\ev_x}/ I_x \CF$ is a Lie algebra, and we call this Lie algebra the {\it isotropy Lie algebra}\index{isotropy Lie algebra}.

Below are some basic results for fibers and tangent spaces of the leaves of singular foliations:
\begin{prop}[\cite{AS06}]
Let $(M, \CF)$ be a singular-foliated manifold, and $x\in M$. We have
\begin{enumerate}
    \item Let $X_1, \cdots, X_k \in \CF$ whose images in $\CF_x$ form a basis of $\CF_x$, then there is a neighborhood $U$ of $x$ such that $\CF_U$ is generated by $X_1, \cdots, X_k$.
    \item $\dim F_x$ is lower semi-continuous and $\dim \CF_x$ is upper semi-continuous.
    \item The set 
    $$
    U = \{ x\in M: \tilde{\ev}_x:\CF_x \to F_x\}\text{ is an isomorphism}
    $$
    is the set of continuity of $x \mapsto \dim F_x$. $U$ is open and dense. $F|_U\subset TM|_U$ is a sub-bundle, hence $\CF_U$ is a regular foliation.
\end{enumerate}
\end{prop}
\begin{proof}
    See \cite[Proposition 1.5]{AS06}.
\end{proof}

\subsection{Holonomy and monodromy}
Holonomy of a regular foliation is defined as germs of local diffeomorphisms of transversals along a path on a leaf. It turns out to be one of the most import concepts related to foliations, for example, we can construct the {\it holonomy groupoid}\index{holonomy groupoid} of a foliation. A related notion is the {\it monodromy}\index{monodromy}, which describes the leafwise homotopy classes of paths.

Let $\Diff_x(M)$ denote the group of diffeomorphisms of a manifold fixing $M$. 
\begin{defn}[\cite{MM03}]
    Let $(\M, \CF)$ be a foliated manifold. Let $x, y$ be some points on some leaf $L$ of $\CF$. Let $S, T$ be {\it transversal sections}\index{transeversal section} (or {\it transversal}\index{transversal}) at $x$ and $y$, then for any path $\alpha:x\to y$ we can associate a germ of a diffeomorphism 
    $$
    \hol^{S,T}(\alpha): (S,x) \to (T,y)
    $$
    which is called the {\it holonomy} of $\alpha$ with respect to the transversal sections $S$ and $T$. For details about the construction of $\hol^{S,T}(\alpha)$, see \cite[Section 2.1]{MM03}.
\end{defn}
Two easy but import properties of the holonomy are \begin{enumerate}
    \item Homotopic paths induces the same holonomy;
    \item Holonomy is independent of the choice of transversals by identifying the holonomy of different transversals along the constant path.
\end{enumerate}
Let $x\in L$ and $S$ be a transversal at $x$. By the above properties, we have a group homomorphism
$$
\hol: \pi_1(L, x) \to \Diff_0(\R^q)
$$
by the independence of the choice of $T$, we get the {\it holonomy homomorphism}\index{holonomy homomorphism} $\hol$
\begin{equation*}
    \hol: \pi_1(L, x) \to \Diff_x(T)\simeq \Diff_0(\R^q)
\end{equation*}
which is defined up to conjugations in $\Diff_0(\R^q)$, where $q =\codim \CF$. We call the image of $\hol$ the {\it holonomy group}\index{holonomy group} of $L$ at $x$, which is determined up to an inner automorphism of $\Diff_0(\R^q)$. As a direct consequence, we have a short exact sequence 
\begin{equation*}
    1\to K \hookrightarrow \pi_1(L,x) \stackrel{\hol}{\to} \hol(L,x) \to 1
\end{equation*}
We will look at this sequence again when we generalize to higher holonomies. Similar to the case of homotopy, we say two path $\alpha, \beta: x\to y$ lying in the same leaf $L$ are in the same {\it holonomy class} if $\hol(\alpha^{-1}\beta) = \id$.

\begin{defn}[\cite{MM03}]
    Taking the differential at 0 gives a homomorphism $d_0:  \Diff_0(\R^q) \to \GL(q, \R)$. We call the composition 
    \begin{equation*}
        d\hol = d_0\circ \hol: \pi_1(L,x) \to \GL(q,\R)
    \end{equation*}
    the {\it linear holonomy homomorphism}\index{linear holonomy homomorphism} of $L$ at x, and we call the image the {\it linear holonomy group}\index{linear holonomy group}.
\end{defn}
An important construction in foliations is associating various groupoids to a foliation. The most important two are {\it holonomy groupoids} and {\it monodromy groupoids}.
\begin{defn}[\cite{MM03}]
    Let $(\M, \CF)$ be a foliated manifold. Define the {\it monodromy groupoid}\index{monodromy groupoid} $\mon(\CF)$ of $\CF$ to be a groupoid over $M$ whose arrows are homotopy classes of paths along leaves of $\CF$. Similarly, define the {\it holonomy groupoid}\index{holonomy groupoid} $\hol(\CF)$ of $\CF$ to be a groupoid over $M$ whose arrows are holonomy classes of paths along leaves of $\CF$
\end{defn}
Both holonomy groupoids and monodromy groupoids are Lie groupoids, hence they are powerful tools in studying the geometry and topology of foliations.

Next, we move to singular foliations.
\begin{defn}
    Let $(M, \CF)$ be a singular-foliated manifold. A {\it slice}\index{slice} $\mathcal{T}$ at $x$ is an embedded submanifold $\mathcal{T}\subset M$ such that $x\in \mathcal{T}$ and $T_x \mathcal{T}\oplus F_x = T_x M$.
\end{defn}
This is similar to the definition of local transversals in regular foliations.

Consider a path $\gamma:[0,1]\to M$ from $x$ to $y$ which lie in a single leaf $L$ of $\CF$. Fixed two slices $\mathcal{T}_x$ and $\mathcal{T}_y$ at $x$ and $y$ respectively. For each time $t$, we lift $\dot{\gamma}(t)$ to a vector field $X_t$ lying in $\CF$, such that the flow of the time-dependent vector field $\{X_t\}$ maps $\mathcal{T}_x$ to $\mathcal{T}_y$. However, if $\CF$ is not a regular foliation, the map $f\in \Hom(\mathcal{T}_x, \mathcal{T}_y)$ will dependent on the extension. Hence, we want to modify this such that we are not affected by the choice of extensions.

First, we want to make some notations. Let $\aut_{\CF}(M)$ be the subgroup of local diffeomorphisms of $M$ preserving $\CF$. Let $\exp$ denote the space of time-one flows of time-dependent vector fields in $\CF$. Recall we have the following exact sequence
\begin{equation*}
    0\to \g_x \to \CF_x \stackrel{\ev_x}{\to }F_x \to 0
\end{equation*}
where $\g_x = \CF(x)/I_x \CF$ is the isotropy Lie algebra. We have that both $\exp(I_x\CF)$ and $\exp(\CF(x))$ are subgroups of  $\aut_{\CF}(M)$. Denote the restriction of $\CF$ to the slice $\mathcal{T}_x$ by $\CF_{\mathcal{T}_x}$, i.e.  $\CF_{\mathcal{T}_x} = \CF|_{\CF_{\mathcal{T}_x}}\cap T\mathcal{T}_x$. Finally, we denote $\germ \aut_{\CF}(\mathcal{T}_x, \mathcal{T}_y)$ the space of germs of local diffeomorphisms from $(\mathcal{T}_x, \CF_{\mathcal{T}_x}) $ to $(\mathcal{T}_y, \CF_{\mathcal{T}_x})$.
\begin{thm}[]
    The class $\Gamma(-, 1): \mathcal{T}_x\to \mathcal{T}_y$ in the quotient 
    \begin{equation}
        \frac{\germ \aut_{\CF}(\mathcal{T}_x, \mathcal{T}_y)}{\exp(\CF_{\mathcal{T}_x})}
    \end{equation}
    is independent of the choice of extension $\Gamma$.
\end{thm}
\begin{proof}
    See \cite[Proposition 2.3]{AZ12}.
\end{proof}

This seems to be a good candidate to define the holonomy transformation. However, as pointed out in \cite{AZ12}, $\exp(\CF(x))$ is too large to do linearization. Hence, we replace $\CF(x)$ by $I_x \CF$ and get the following definition:
\begin{defn}[\cite{AZ12}]
    Let $(M, \CF)$ be a singular-foliated manifold, and $x, y \in M$ lying in some leaf $L$. We fix slices $\mathcal{T}_x$ and $\mathcal{T}_y$ at $x$ and $y$ respectively. We define a {\it holonomy transformation}\index{holonomy!transformation} from $x$ to $y$ to be an element of 
    \begin{equation}
        \frac{\germ \aut_{\CF}(\mathcal{T}_x, \mathcal{T}_y)}{\exp(I_x \CF_{\mathcal{T}_x})}
    \end{equation}
\end{defn}
\begin{lem}
    Let $x$ be a point on a regular leaf, then $\exp(I_x \CF_{\mathcal{T}_x})$ is trivial.
\end{lem}
\begin{proof}
    In this case, $\mathcal{T}_x$ is equipped with a trivial singular foliation, i.e. foliation by points. Hence, the flow is trivial.
\end{proof}
Therefore, we see that for regular foliation, the holonomy transformations reduce to the ordinary holonomy transformations $\germ \Diff(\mathcal{T}_x, \mathcal{T}_y)$. It is easy to see that we can form a topological groupoid $\holtrans(\CF)$ over $M$ with morphisms being
\begin{equation*}
   \bigcup_{x,y} \frac{\germ \aut_{\CF}(\mathcal{T}_x, \mathcal{T}_y)}{\exp(I_x \CF_{\mathcal{T}_x})}
\end{equation*}
There is a natural map from the holonomy groupoid in the sense of Androulidakis and Skandalis in \cite{AS06}, which justifies the correctness of the definition of holonomy transformations. Let's first review basics about holonomy groupoids defined by Androulidakis and Skandalis.

\begin{defn}[\cite{AS06}]

Let $(M, \CF_M), (N, \CF_N)$ be two singular foliated manifolds. A {\it bisubmersion}\index{bisubmersion} is a manifold $P$ with two surjective submersions  $s: P\to M$ and $t: P\to N$ such that 
$$s^{-1} \CF_M = t^{-1} \CF_N = \Gamma(\ker(s)_*) + \Gamma(\ker(t)_*)$$

We have the following diagram
    \begin{center}
        \begin{tikzcd}
             & {(P, \CF)} \arrow[ld, "s"] \arrow[rd, "t"] &              \\
{(M, \CF_M)} &                                                    & {(N, \CF_N)}
\end{tikzcd}
    \end{center}
where $\CF = s^{-1} \CF_M = t^{-1} \CF_N$ is the pullback singular foliation on $P$.

\end{defn}
\begin{defn}[]
    A {\it morphism between bisubmersions}\index{morphism!bisubmersion} $(U, s_U, t_U), (V, s_V, t_V)$ is a smooth map $f:U \to V$ such that for all $u\in U$, we have $s_V(f(u)) = s_U(u)$ and 
    $t_V(f(u)) = t_U(u)$.
    
    A {\it local morphism between bisubmersions}\index{local morphism!bisubmersion} $(U, s_U, t_U), (V, s_V, t_V)$ is a smooth map $f:U' \to V$ for some $U'\subset U$ such that for all $u\in U'$, we have $s_V(f(u)) = s_U(u)$ and 
    $t_V(f(u)) = t_U(u)$.
\end{defn}

\begin{defn}
    Let $\mathcal{U} = (U_i, s_i, t_i)_{i\in I}$ be a family of bisubmersions.
    \begin{enumerate}
        \item A bisubmersion $(U,s,t)$ is said to be {\it adapted to $\mathcal{U}$ at $u\in U$} if there exists an open subset $U'\subset U$ containing $u$ and a morphism of bisubmersion $U'\to U$.
        \item A bisubmersion $(U,s,t)$ is said to be {\it adapted to $\mathcal{U}$} if $(U,s,t)$ is adapted to $\mathcal{U}$ at $u\in U$ for all $u\in U$.
        \item We call the family $\mathcal{U} = (U_i, s_i, t_i)_{i\in I}$ an {\it atlas}\index{atlas} if \begin{enumerate}
            \item $\cup_{i\in I}s_i(U_i) = M$.
            \item Any elements of $\mathcal{U}$ is still adapted to $\mathcal{U}$ under taking inverses and compositions.
        \end{enumerate}
    \end{enumerate}
\end{defn}

Recall that an atlas of a manifold allows us to reconstruct the manifold, similarly an atlas of a singular foliation $(M, \CF)$ allows us to reconstruct the foliated structure by a groupoid over $M$, and this is the first step in constructing the holonomy groupoid.
\begin{thm}[\cite{AS06}]
    Let $\mathcal{U} = (U_i, s_i, t_i)_{i\in I}$ be an atlas of a singular foliated manifold $(M, \CF)$.
    \begin{enumerate}
        \item Let $G = \coprod_{i\in I} U_i / \sim$ where $\sim$ is the equivalence relation generated by local morphisms, i.e. $u\in U_i $ is equivalent to $v\in U_j$ if there exists a local morphism from $U_i\to U_j$ which takes $u$ to $v$. There are maps $s,t: G \to M$ such that $s\circ q_i = s_i$ and $t\circ q_i = t_i$, where $Q= (q_i)_{i\in I}:\coprod_{i\in I}  \to G $ is the quotient map.
         \item For any $(U, s_U, t_U)$ adapted to $\mathcal{U}$, there exists a map $q_U: U\to G$ such that for every local morphism $f: U'\subset U \to U_i$ and every $u\in U'$, we have $q_U(u) = q_i(f(u))$.
         \item There exists a (topological) groupoid structure on $G$ over $M$ with source and target maps $s$ and $t$ defined before, and $q_i(u)q_j(v) = q_{U_i\circ U_j}(u,v)$ .
    \end{enumerate}
\end{thm}
Hence, given any atlas, we can construct a groupoid which encodes information about the singular foliation. Apparently, we can simply take all possible bisubmersions to be an atlas, which is called the {\it full holonomy atlas}\index{atlas!full holonomy}. However, this atlas is obviously too big which will make the arrow space of the groupoid to be incredibly large and nasty. The second choice is to take all leaf-preserving bisubmersions, i.e. bisubmersions $(U,s,t)$ that $s(u)$ and $t(u)$ lying in the same leaf for all $u\in U$. This is called the {\it leaf-preserving atlas}\index{atlas!leaf-preserving}, which is much smaller than the full holonomy atlas. \cite{AS06} constructs an atlas $\mathcal{W}$ which is as minimal as possible.  
\begin{prop}[\cite{AS06}]\label{pathhol}
    Let $x\in M$ and $X_1,\cdots, X_n\in \CF$ be vector fields whose images at $x$ form a basis of $\CF_x$ (the fiber of $\CF$ at $x$). For $y=(y_1,\cdots, y_n)\in \R^n$,  set $\phi_y = \exp(\sum y_i X_i)\in \exp \CF$, i.e.  the image of $y$ under the time-1 flow of the vector field $\sum y_i X_i$. Let $W_0 = \R^n \times M, s_0(y,x) = x, s_0(y,x) = \phi_y(x)$. Then
    \begin{itemize}
        \item There exists a neighborhood $W$ of $(0,x)$ in $W_0$ such that $(W, s|_W, t|_W)$ is a bisubmersion.
        \item For any bisubmersions  $(V,s_V, t_v)$ carries the identity of $M$ at some $v\in V$, then there exists a local morphism from  $(V,s_V, t_v)$ to some $(W, s_W, t_W)\in \mathcal{W}$ at $v$ which sends $v$ to $(0,x)$.
    \end{itemize}
\end{prop}

\begin{defn}[\cite{AS06}]
    We define the {\it path holonomy atlas}\index{atlas! path holonomy} of a singular foliated manifold $(M, \CF)$ to be the maximal atlas generated by a cover of $M$ by $s$-connected bisubmersions of the form in Proposition \ref{pathhol}.
\end{defn}
The corresponding groupoid of the path holonomy atlas is the smallest due to (2) in Proposition \ref{pathhol}, which implies that this atlas is adapted to any other atlas.

\begin{defn}
    Let $(M, \CF)$ be a singular-foliated manifold. We define the {\it holonomy groupoid}\index{holonomy groupoid!of singular foliations} of $\CF$ to be the (topological) groupoid associated to the path holonomy atlas of $\CF$. We denote the holonomy groupoid of $\CF$ by $\hol(\CF)$ or $\hol^{AS}(\CF)$.
\end{defn}
    The topology of holonomy groupoids in this definition is usually pretty bad. For example, \cite[Example 3.7]{AS06} consider the singular foliation on $\R^2$ generated by the action of $\SL(2,\R)$, whose holonomy groupoid has a highly non-Hausdorff arrow space, for example, for any $x$ outside the origin, the sequence $\{(x/n, x/n)\}$ will converge to all $(g,0)$ for stabilizers $g$ of $x$. Therefore, in general we don't expect holonomy groupoids to be Lie groupoids. However, we do have the following local smoothness results.
    \begin{thm}[\cite{Deb13}]
         Let $(M, \CF)$ be a singular-foliated manifold. The $s$-fibers of $\hol(\CF)$ are smooth manifolds.
    \end{thm}
    \begin{cor}[\cite{AZ11}]
    The transitive groupoid $\hol_L(\CF)$  is smooth and integrates the Lie algebroid $A_L = \cup_{x\in L} \CF_x$, where $\hol_L(\CF) = \hol(\CF)|_{s^{-1}L}=\hol(\CF)|_{t^{-1}L}$
    \end{cor}
    
    Hence a natural question is can we find a higher categorical geometric object such that the 1-truncation equals the holonomy groupoid?
    
    We will answer this question later, and let's return to the holonomy transformation first.
\begin{thm}[\cite{AZ12}]
    Let $(M, \CF)$ be a singular-foliated manifold, and $x, y \in M$ lying in some leaf $L$. We fix slices $\mathcal{T}_x$ and $\mathcal{T}_y$ at $x$ and $y$ respectively. There is a natural injective map
    \begin{equation*}
        \Phi^y_x: \hol^{AS}(\CF)^y_x \to \holtrans(\CF)^y_x=\frac{\germ \aut_{\CF}(\mathcal{T}_x, \mathcal{T}_y)}{\exp(I_x \CF_{\mathcal{T}_x})}
    \end{equation*}
    where $\hol^{AS}(\CF)$ denotes the holonomy groupoid in the sense of Androulidakis and Skandalis \cite{AS06}. Moreover, $\Phi^y_x$ assembles to a global groupoid morphism
    $$
    \Phi: \hol(\CF) \to \holtrans(\CF)
    $$
\end{thm}

\begin{defn}
    Let $f:\mathcal{T}_x \to \mathcal{T}_y$ be a holonomy transformation. Suppose $f$ is also an embedding, then we call $f$ a {\it holonomy embedding}\index{holonomy embedding}.
\end{defn}

\subsection{Hausdorff Morita equivalences}
There are various notions of two (singular) foliations to be equivalent. Garmendia and Zambon \cite{GZ19} proposed a notion called {\it Hausdorff Morita equivalence} of singular foliations, which is constructed to be compatible to Androulidakis and Skandalis's construction of holonomy groupoids in \cite{AS06}.
\begin{defn}\cite{AS06}

\end{defn}

\begin{defn}[\cite{GZ19}]
    Let $(M, \CF_M), (N, \CF_N)$ be two singular foliated manifolds. We say $(M, \CF_M)$ and $(N, \CF_N)$ are {\it Hausdorff Morita equivalent}\index{Hausdorff Morita equivalent} if there exists a manifold $P$ and two surjective submersions with connected fibers $\pi_M: P\to M$ and $\pi_N: P\to N$ such that $\pi_M^{-1} \CF_M = \pi_N^{-1} \CF_N$. We have the following diagram
    \begin{center}
        \begin{tikzcd}
             & {(P, \CF)} \arrow[ld, "\pi_M"] \arrow[rd, "\pi_N"] &              \\
{(M, \CF_M)} &                                                    & {(N, \CF_N)}
\end{tikzcd}
    \end{center}
where $\CF = \pi_M^{-1} \CF_M = \pi_N^{-1} \CF_N$ is the pullback singular foliation on $P$.

\end{defn}
This definition is almost the same as the bisubmersion except that we do not require $\pi_M^{-1} \CF_M = \pi_N^{-1} \CF_N=\Gamma(\ker(\pi_M)_*) + \Gamma(\ker(\pi_N)_*)$. Hence, Hausdorff Morita equivalences are weaker notions than bisubmersions. 

The following theorem gives basic properties of Hausdorff Morita equivalences.
\begin{thm}{\cite{GZ19}}
    Let two singular foliated manifolds $(M, \CF_M), (N, \CF_N)$ be Hausdorff Morita equivalent. We have
    \begin{enumerate}
        \item There is a homeomorphism between the leaf spaces of $(M, \CF_M)$ and the leaf space of $(N, \CF_N)$, which maps the leaf though some $x \in M$ to the leaf of $\CF_N$ containing $\pi_N(\pi^{-1}(x))$, and preserves the codimension of leaves and the property of being an embedded leaf.
        \item Consider $x\in M, y\in N$. Pick slices $S_x$ at $x$ and $S_y$ at $y$, then the singular foliated manifolds $(S_x, \iota_{S_x}^{-1})$ and $(S_y, \iota_{S_y}^{-1})$ are diffeomorphic.
        \item Consider $x\in M, y\in N$. The isotropy Lie algebras $\g_x^{\CF_M}$ and $\g_y^{\CF_N}$ are isomorphic.
    \end{enumerate} 
\end{thm}

The next theorem justifies Hausdorff Morita equivalence is the correct notion which preserves holonomy groupoids.
\begin{thm}[\cite{GZ19}]
    If two singular foliated manifolds $(M, \CF_M), (N, \CF_N)$ are Hausdorff Morita equivalent, then their holonomy groupoids (in the sense of Androulidakis and Skandalis in \cite{AS06}) are Morita equivalent as open topological groupoids.
\end{thm}

Next, let's look at the behavior of Hausdorff Morita equivalences under pullbacks.
\begin{prop}\label{he1}
    Pullbacks of Hausdorff Morita equivalences are Hausdorff Morita equivalences.
\end{prop}

\begin{proof}
 Let $f:(M, \CF_M)
    \to (N, \CF_N)$, $g: (L, \CF_L) \to (N, \CF_N)$ be morphisms in $ \mfd^{\sfol}$. Suppose the pullback singular foliated manifold $(M\times_N L,\CF_M\times_{\CF_N} \CF_L)$ exists and  $f$ is a Hausdorff Morita equivalence, then the induces map $f':(M\times_N L,\CF_M\times_{\CF_N} \CF_L) \to (L, \CF_L)$ is also a Hausdorff Morita equivalence.
\end{proof}

\begin{proof}
Let $(M, \CF_M) \stackrel{\pi_M}{\leftarrow} (P,\CF)\stackrel{\pi_N}{\rightarrow} (N, \CF)$ be a morphism representing the Hausdorff Morita equivalence. We have the following commutative diagram
\begin{center}
    \begin{tikzcd}
                       & {(P',\CF')} \arrow[r] \arrow[ld, "\pi_M'"'] \arrow[ldd, "\pi_N'"] & {(P,\CF)} \arrow[ld, "\pi_M"'] \arrow[ldd, "\pi_N"] \\
{(M\times_N L,\CF_M\times_{\CF_N} \CF_L) } \arrow[r, "g'"] \arrow[d, "f"] & {(M, \CF_M)} \arrow[d, "f"]                        &                           \\
{(L, \CF_L)} \arrow[r, "g"]           & {(N, \CF)}                                  &                          
\end{tikzcd}
\end{center}
Since $\pi_M, \pi_N$ are submersions, we see their pullbacks exist, and the left triangle in the diagram is actually a composition of pullbacks, which implies that all three squares are pullbacks. In particular.  $(\pi_M')^{-1}(\CF_M\times_{\CF_N} \CF_L) = (\pi_N')^{-1}(\CF_L)$ follows from the composition pullback. Clearly the fiber of $\pi_M'$ and $\pi_N'$ since $\pi_M$ and $\pi_N$'s are.
\end{proof}

\subsection{Homotopy theory of singular foliations}

\begin{defn}
    Let $f:(M, \CF_M)
    \to (N, \CF_N)$, $g: (L, \CF_L) \to (N, \CF_N)$ be morphisms in $ \mfd^{\sfol}$. We say $f$ is {\it foliated transverse}\index{foliated transverse} to  $g$ if the natural map $(d_f\times d_g) (\CF_M \times \CF_L) \to f^*\CF_N \times_N g^* \CF_N$ is surjective. Here we use $f^*\CF_N \times_N g^* \CF_N$  to denote the pullback of $\CF_N$ on $M\times_N L$.
\end{defn}

\begin{prop}
Let $f:(M, \CF_M)
    \to (N, \CF_N)$, $g: (L, \CF_L) \to (N, \CF_N)$ be morphisms in $ \mfd^{\sfol}$. Suppose the $f$ is foliated transverse to $g$. If the pullback $(M\times_N L)$ exists, then the pullback $(\CF_M\times_{\CF_N} \CF_L)$ is a singular foliation on $(M\times_N L)$.
    \begin{center}
	\begin{tikzcd}
	(M\times_N L,\CF_M\times_{\CF_N} \CF_L)  \arrow[r, "g'"]\arrow[dr, phantom, "\ulcorner", very near start] \arrow[d,"f'"] & (M, \CF_M) \arrow[d, "f"] \\
	(L, \CF_L) \arrow[r,"g"]           & (N,\CF_N)         
	\end{tikzcd}
\end{center}
\end{prop}
\begin{proof}
    First let's look at the involutivity. Denote $\CF_M\times_{\CF_N} \CF_L$ by $\CF$.  Obviously, $g'$ is foliated, i.e. $d_{g'}(\CF) \subset \CF_M$. Let $X, X' \in \CF_M$. Write $d_{g'}(X) = \sum f_i Y_i\circ g'$ and $d_{g'}(X') = \sum f_i' Y_i'\circ g'$, then
    \begin{equation*}
        d_{g'}([X, X']) = \sum f_if_j' [Y_i, Y_j']\circ g' + \sum X(f_j')Y_j' \circ g' - \sum X'(f_i)Y_i \circ g'
    \end{equation*}
    hence we see the pullback is closed under brackets.
    
    Next, we want to show $\CF$ is locally finitely generated. By restricting to sufficient small open subsets of $M$, $L$ and $N$, we can assume $\CF_M$, $\CF_L$, and $\CF_N$ are finitely generated, and tangent bundles of $M$, $L$, and $N$ are trivial.
    
    \begin{equation*}
        \CF = \CF_M\times_{\CF_N} \CF_L
        =\CF_M \tensor[_{f^{-1}(\CF_N)}]{\times}{_{g^{-1}(\CF_N)}}  \CF_L
    \end{equation*}
    Note that \begin{align*}
        f^{-1}(\CF_N) =& f^*(\CF_N)\times_{\Gamma(M, f^*(TN))} \Gamma(M, TM)\\
        g^{-1}(\CF_N) =& g^*(\CF_N)\times_{\Gamma(L, g^*(TN))} \Gamma(L, TL)
    \end{align*}
    By foliated transversality, we see there exists a section
    $$
    s: \CF_N \to f^{-1}(\CF_N) \times g^{-1}(\CF_N)\subset \CF_M \times \CF_L
    $$
    It follows that $\CF_M \tensor[_{f^{-1}(\CF_N)}]{\times}{_{g^{-1}(\CF_N)}}  \CF_L$ is finitely generated.
\end{proof}
Similarly, we can define foliated submersions.
\begin{defn}
    Let $f:(M, \CF_M)
    \to (N, \CF_N)$ be morphisms in $ \mfd^{\sfol}$. We say $f$ is a {\it foliated submersion}\index{foliated submersion} if the natural map $d_f \CF_M \to f^*\CF_N$ is surjective.
\end{defn}

By this definition, a foliated submersion is then foliated transversal to any  foliated maps. Hence, a direct corollary is:
\begin{cor}\label{pb3}
    Let $f:(M, \CF_M)
    \to (N, \CF_N)$, $g: (L, \CF_L) \to (N, \CF_N)$ be morphisms in $ \mfd^{\sfol}$. Suppose $f$ is a surjective foliated submersion. If the pullback $(M\times_N L)$ exists, then the pullback $(\CF_M\times_{\CF_N} \CF_L)$ is a singular foliation on $(M\times_N L)$. Moreover, the induced map $f': (M\times_N L, \CF_M\times_{\CF_N} \CF_L)\to (N, \CF_N)$ is a surjective foliated submersion.
\end{cor}

Next, we will construct a finite dimensional model for the path object of singular foliations. Recall that, in general, $\cinf$ path spaces for finite dimensional manifolds are infinite dimensional, even if we restrict to $C^1$ paths. One way to remedy this is to consider only those 'short paths'. We will follow the construction in \cite[Section 2.1] {BLX21}.

First, we can some connection $\nabla$ on $M$. Let $\exp^{\nabla}$ denote the exponential map with respect to $\nabla$. Let $I=(a,b) \supset [0,1]$.

\begin{prop}[\cite{BLX21}]\label{ps1}
There exists a manifold $P_g M$, which is called the {\it manifold of short geodesic paths}\index{manifold of short geodesic paths} in $M$, which parametrizes a family of geodesic paths, such that
\begin{center}
    \begin{tikzcd}
              & {M} \arrow[ld, "0"] \arrow[rd, "\Delta"] \arrow[dd, "const"] &               \\
{TM}  &                                     & {M\times M}  \\& {{P_gM} } \arrow[lu, "\gamma(0)\times \gamma'(0)"] \arrow[ru, "\gamma(0)\times \gamma(1)"']            & 
\end{tikzcd}
\end{center}
where the two lower diagonal maps are open embeddings. 
\end{prop}
\begin{proof}
    Consider $U\subset TM$ an open neighborhood of the zero section, where the $\exp^{\nabla}(tv)$ is defined for any $t\in I, v\in U$. $U$ parametrize a family of geodesic paths with domain $I$:
    \begin{align*}
        U\times I &\to M\\
        (x, v, t)\mapsto \gamma_{x,v}(t) &= \exp^{\nabla}_x(tv)
    \end{align*}
    We can restrict to a smaller open neighborhood $V\subset V$ such that the both maps $U\to TM$ and $V\to M\times M$ are open embeddings. Now take $P_gM = V$. 
\end{proof}
Hence, $P_g M$ is diffeomorphic to an open neighborhood of the zero section of $TM$. The evaluation map $\ev_0: P_g M \to M$ is just the restriction of the projection $TM \to M$, and $\ev_1:P_g M \to M$ is given by the exponential map
\begin{equation*}
    v_x \mapsto \exp^{\nabla}_x v_x
\end{equation*}
\begin{lem}\label{ps2}
Let $(M, \CF)$ be a singular-foliated manifold. Then there exists a associated foliation $\CF_g$ on $P_g M$ which projects to $\CF$ along the projection $P_gM \to M$.
\end{lem}
\begin{proof}
    It suffices to construct $(M, \CF)$ locally. The connection gives a splitting $TTM \simeq TM\oplus vTM$, where $vTM$ denotes the vertical tangent bundle of $TM$. Hence, locally, $P_gM \simeq U\times \R^n$ for some $U\subset M$ and $n= \dim M$. Pick $X_1,\cdots, X_k$ to be generators of $\CF|_U$, and $e_1,\cdots, e_n$ be the local coordinate sections of $vTM|_U$. Now we define $\CF_g|_U$ to be the module generated by $X_1,\cdots, X_k, e_1,\cdots, e_n$. The involutivity follows directly from our construction. Note the projection $P_gM \to M$ is clearly a foliated map which just kills $e_i$'s, which then maps $\CF_g$ to $\CF$. 
\end{proof}
\begin{lem}\label{ps3}
    Both $\ev_0$ and $\ev_1$ are foliated submersions.
\end{lem}
\begin{proof}
    Follows from the construction that \begin{equation*}
        (X_1,\cdots, X_k, e_1,\cdots, e_n)\mapsto (X_1,\cdots, X_k)
    \end{equation*}
    is surjective.
\end{proof}

\begin{prop}
      There exists an incomplete category of fibrant objects structure on $\mfd^{\sfol}$, where
      \begin{enumerate}
          \item Fibrations are surjective foliated submersions.
          \item Weak equivalences are Hausdorff Morita equivalences.
          \item Path objects are foliated short geodesic path space.
      \end{enumerate}
\end{prop}
\begin{proof}
    (2) follows from Corollary \ref{ps3}. (3) then follows from (4) and Proposition \ref{he1}. 
    
    For (4), clearly isomorphisms are Hausdorff Morita equivalences, let's prove 2-out-of-3 property. By construction, the composition of Hausdorff Morita equivalences is the fiber product which is again a Hausdorff Morita equivalences, this directly implies that, if $(M, \CF_M)\simeq (N, \CF_N)$ and $(L, \CF_L)\simeq (N, \CF_N)$, then we have $(M, \CF_M)\simeq (L, \CF_L)$. By symmetry of the Hausdorff Morita equivalences, this will generate all cases of 2-out-of-3.
    
    For (5), clearly isomorphisms are surjective foliated submersions, and composition of foliated submersions are again foliated submersions by definition.
    
    For (6), given a singular foliated manifold $(M, \CF)$, we construct its path object $(M, \CF)^{\Delta[1]}$ to be $(P_g M, \CF_g)$. By Theorem \ref{ps1} and Lemma \ref{ps2}, we have the following factorization
    \begin{equation*}
	(M, \CF) \stackrel{\iota}{\longrightarrow} (P_g M, \CF_g) \stackrel{(\gamma(0), \gamma(1))}{\longrightarrow} 	(M, \CF) \times 	(M, \CF)
	\end{equation*}
	which composes to the diagonal map. $(\gamma(0), \gamma(1))$ are fibrations by Lemma \ref{ps3}. To see $\iota$ is a Hausdorff Morita equivalence, notice that the $\iota$ is an embedding, and $p:P_gM \to \iota(M)$ is a submersion. By construction $p_{-1}(\iota(\CF))$ is exactly $\CF_g$. Hence, $(P_gM, \CF)$ itself gives a Hausdorff Morita equivalence by $(\id, p)$.
    
    Finally, the trivial map $(M, \CF) \to *$ is clearly a surjective foliated submersion, hence a fibration.
\end{proof}

We denote the $\infty$-category of singular foliated manifolds presented by this iCFO by $\mfdi^{\sfol} = \mfd^{\sfol}[W^{-1}]$. We will use this later in constructing algebraic $K$-theory of singular foliations.
\subsection{Algebraic $K$-theory of singular foliations}
In this section, we will construct the {\it algebraic $K$ theory sheaves} \index{algebraic $K$ theory sheaves} $\mathsf{K}$ on $\mfd^{\sfol}$ following \cite{Bun18} for regular foliations. Then we can calculate Algebraic $K$-theory of $(M,\CF)$ by taking homotopy groups of $\mathsf{K}(M,\CF)$:
\begin{equation*}
    K^{\bt}(M,\CF) = \pi_{-\bt}(
    \big( \mathsf{K}(M,\CF)\big)
\end{equation*}

We consider $\cat$ with its Cartesian symmetric monoidal structure. Let $W$ denote the class of categorical equivalence, then we get a symmetric monoidal category $\cat[W^{-1}]$. We denote the category of commutative algebras in  $\cat[W^{-1}]$ by $\calg(\cat[W^{-1}])$. 

\paragraph{Construction of algebraic $K$ theory sheaves}
\begin{enumerate}
    \item Let $M$ be a manifold, we denote $\vect(M)$ the category of vector bundles over $M$. Any map $f: M \to N$ induces a functor $f^*:\vect(N)\to \vect(M)$. Hence, we get a stack $\vect$ on the site $\mfd$ with open covering topology. Similarly, we get a stack $\vect^{{\C}\sfol}$ by pullback along the forgetful functor 
\begin{equation*}
    F: \mfd^{{\C}\sfol} \to \mfd
\end{equation*}
We will write $\vect$  for $\vect^{{\C}\sfol}$ for simplicity.

\item Similarly, we can consider the category of pairs $(V, \nabla)$ of a vector bundle $V\to M$. Denote the resulting symmetric monoidal stack (with the Cartesian symmetric monoidal structure) by  $\vect^{\nabla}$.

\item Let $(M,\CF)\in \mfd^{{\C}\sfol}$. Denote $\vect^{\Flat}(M, \CF)$ the category of pairs $(V, \nabla_{I})$ of a vector bundle $V\to M$ and a flat partial connection $\nabla_{I}$ on $(M,\CF)$. A foliated map $f:(M,\CF) \to (M', \CF')$ induces a functor $f' :\vect^{\Flat}(M', \CF')\to \vect^{\Flat}(M, \CF)$. We get a stack $\vect^{\Flat}$ on the site $\mfd^{{\C}\sfol}$.

\item 
Finally, let $(M,\CF)\in \mfd^{{\C}\sfol}$. Denote $\vect^{\Flat,\nabla}(M, \CF)$ the category of pairs $(V, \nabla_{\CF})$ of a vector bundle $V\to M$ and a flat $\CF$-connection $\nabla_{\CF}$ on $(M,\CF)$. A foliated map $f:(M,\CF) \to (M', \CF')$ induces a functor $f' :\vect^{\Flat, \nabla}(M', \CF')\to \vect^{\Flat,\nabla}(M, \CF)$. We get a symmetric monoidal stack $\vect^{\Flat,\nabla}$ on the site $\mfd^{{\C}\sfol}$.

We have the following commutative diagrams of stacks by forgetful maps
\begin{center}
    \begin{tikzcd}
              & {\vect^{\Flat,\nabla}} \arrow[ld] \arrow[rd] &               \\
{\vect^{\Flat}} \arrow[rd] &                          & {\vect^{\nabla}} \arrow[ld] \\
              & {\vect}                       &              
\end{tikzcd}
\end{center}
in $\sh_{\calg(\cat(W^{-1}))}(\mfd^{{\C}\sfol})$. 

\item Now we can apply the {\it  K-theory machine}\index{K-theory machine} developed in \cite{BNV13}, we get a commutative diagram of presheaves of spectra
\begin{center}
    \begin{tikzcd}
{\mathsf{K}(\vect^{\Flat,\nabla})} \arrow[r] \arrow[d] & {\mathsf{K}(\vect^{\Flat})} \arrow[d] \\
{\mathsf{K}(\vect^{\nabla})} \arrow[r] \arrow[d] & {\mathsf{K}(\vect)} \arrow[d] \\
{\widehat{\ku}^{\nabla} = s\big(\mathsf{K}(\vect^{\nabla})\big)} \arrow[r]           & {\widehat{\ku}^{} = s\big(\mathsf{K}(\vect^{})\big)}      \end{tikzcd}
\end{center}
 Here $s:\psh \to \sh$ denotes the sheafification functor.
\end{enumerate}
The K-theory machine developed in \cite{BNV13} is basically a composition,
\begin{align*}
    &\calg(\cat(W^{-1})) \to \calg(\gpd(W^{-1})) \to \common(s\set[W^{-1}) \\
    &\to \comgrp(s\set[W^{-1}) \simeq \mathsf{Sp}_{\ge 0} \to \mathsf{Sp}
\end{align*}
where we \begin{enumerate}
    \item First take the groupoid underlying $\cat$.
    \item Applying nerve to get a commutative monoid in the category of spaces $s\set[W^{-1}]$, i.e. an $E_{\infty}$-space.
    \item Then applying the group completion, we get a commutative group in the category of spaces $s\set[W^{-1}]$, i.e. a grouplike $E_{\infty}$-space.
    \item Finally, we apply the functor which maps a commutative group in spaces to the corresponding connective spectrum whose $\infty$-group is this group.
\end{enumerate}
For more details about the K-theory machine $\mathsf{K}$, see \cite[Definition 6.1, Remark 6.4]{BNV13}. 
Let $L$ and $\mathcal{H}^{\Flat}$ denote the sheafification and homotopification functors.
\begin{defn}
    We define the following sheaves of spectra 
    \begin{align*}
        \mathbf{K} =& \mathcal{H}^{\Flat}(L(\mathsf{K}(\vect^{\Flat}))) \in \sh_{\Sp}^h(\mfd^{{\C}\sfol})\\
        \mathbf{K}^{\nabla} =& L(\mathsf{K}(\vect^{\Flat}, \nabla)) \in \sh_{\Sp}(\mfd^{{\C}\sfol})
    \end{align*}
    and for $i\in \Z$, we define the algebraic $K$-theory of a singular foliation $(M, \CF)$ \index{algebraic $K$-theory! of a singular foliation} by
    \begin{equation*}
        K^i(M,\CF) = \pi_{-i}(\mathbf{K}(M, \CF))
    \end{equation*}
\end{defn}
Note that $\mathsf{K}(\vect^{\Flat})$ is a homotopy invariant, hence we expect that homotopification will preserve this invariance. Therefore, the homotopification might not be necessary.

\section{Higher groupoids arised in singular foliations}
\subsection{Leaf spaces of singular foliations and \v{C}ech $\infty$-groupoids} 
Recall that for a foliated manifold $(M, \CF)$, the leaf space of $\CF$ is a space $\CF$ which is a quotient of $M$ by identifying points within the same leaves. We want to construct a smooth model for the leaf space of a singular foliation. 
\begin{defn}
    We define a {\it transversal basis}\index{transversal basis} for $(M, \CF)$ as a family $\mathcal{U}$ of slices $U$ such that given any slice $V$ at $x$, we can find a $U$ at $y$ and $x, y$ lying in the same leaf, and there exists a holonomy embedding $h: V\hookrightarrow U$.   
\end{defn}

Given a point $x$, we can take a 

\begin{defn}
    Let $(M, \CF)$ be a singular foliated manifold and $\mathcal{U}$ a transversal basis of for $(M, \CF)$. We define the {\it \v{C}ech $\infty$-groupoid}\index{\v{C}ech $\infty$-groupoid} $\cech_{\bt}(\CF)$ whose  $k$-simplices are
    \begin{equation*}
        \cech(\CF)_k = \coprod_{U_0\stackrel{h_1}{\to}\cdots\stackrel{h_k}{\to} U_k} U_0
    \end{equation*}
    where $h_i: U_{i-1}\to U_i$ are holonomy embeddings and $U_i\in \mathcal{U}$ for all $i$.
    
    There structure maps $d_i$'s and $s_i$'s are defined as 
\end{defn}

\begin{defn}

\end{defn}
For a regular foliation, the standard model for the leaf space is the classifying space of the holonomy groupoid. If $\CF$ is regular, then our construction reduces to \cite{CM00}, and we have the following isomorphism
\begin{thm}[\cite{CM00}]
    For a regular foliated-manifold $(M, \CF)$. There is a natural isomorphism
    $$
    \check{H}^{\bt}_{\mathcal{U}}(M/\CF) \simeq H^{\bt}(B\hol(M\CF); \R)
    $$ between the \v{C}ech-de Rham cohomology and the cohomology of the classifying space of the holonomy groupoid. The left hand side is independent of the choice of $\mathcal{U}$.
\end{thm}
Hence, we can regard $\cech(\CF)_{\bt}$ as a model for the leaf space of a singular foliated manifold $(M, \CF)$.
\subsection{Holonomy $\infty$-groupoids}

In this section, we are going to construct a higher groupoid model which extends the holonomy groupoids in the sense of Androulidakis and Skandalis.

Recall that, when we construct the groupoid associated to an atlas of bisubmersions, we take the quotient of equivalence relations induced by local morphisms of bisubmersions. This is the reason we get a crappy arrow space for the holonomy groupoid. One natural idea is that, instead of brutally quotient the (local) equivalence, we keep the gluing data, and take a 'nerve' similar to the case when we construct classifying spaces.

Let $(M, \CF)$ be a singular-foliated manifold. In this section, we fix the atlas $\mathcal{U}= (U_i, s_i, t_i)_{i\in I}$ to be the path holonomy atlas of $(M, \CF)$.

\begin{defn}
    We say a morphism between bisubmersions $f: (U, s_U, t_U)\to (V, s_V, t_V)$ is an {\it equivalence}\index{equivalence! of bisumsersion morphisms} if there exists a morphism $g:(V, s_V, t_V)\to (U, s_U, t_U)$. We say $f$ is an {\it isomorphism} \index{!isomorphism of bisumsersion morphism} if $f\circ g = \id, g\circ f = id$.
    
    We say a morphism between bisubmersions $f: (U, s_U, t_U)\to (V, s_V, t_V)$ is a {\it local equivalence (isomorphism)}\index{local equivalence! of bisumsersion morphisms}\index{local isomorphism! of bisumsersion morphisms} if there exists a morphism $f': U'\to V'$ where $U'\subset U, V'\subset V$ which is an equivalence (isomorphism).
\end{defn}
First, let's recall the following lemma about the local morphisms.
\begin{lem}[\cite{AS06}]
    Let $(U, s_U, t_U), (V, s_V, t_V)$ be bisubmersions and let $u\in U, v\in V$ with $s_U(u) = s_V(v)$. Then:
    \begin{enumerate}
        \item If the identity local diffeomorphism is carried by $U$ at $u$ and by $V$ at $v$, then there exists an open neighborhood $U'$ of $u$ in $U$, and a morphism $f: U' \to V$ such that $f(u) = v$.
        \item If there exists a local diffeomorphism is carried by $U$ at $u$ and by $V$ at $v$, then there exists an open neighborhood $U'$ of $u$ in $U$, and a morphism $f: U' \to V$ such that $f(u) = v$.
        \item If there exists a morphism of bisubmersions $f:U \to V$ such that $f(u) = v$, then there exists an open neighborhood $V'$ of $v$ in $V$, and a morphism $g: V' \to U$ such that $g(v) = u$.
    \end{enumerate}
\end{lem}
\begin{lem}Let $(U, s_U, t_U), (V, s_V, t_V)$ be bisubmersions and let $u\in U, v\in V$ with $s_U(u) = s_V(v)$, then any local morphisms around $u$ are local equivalences.
\end{lem}
\begin{proof}
    Let $f': U'\to V$ denote the morphism induced by $f$. Applying the previous lemma, we can restrict to an open neighborhood  $V'\subset V$ with a morphism $g: V'\to U'$ such that $g(v) = u$, which realizes a local equivalence.
\end{proof}

\begin{cor}
Local morphisms are local equivalences.
\end{cor}
\begin{proof}
Let $f:(U,s_U, t_U) \to (V, s_V, t_V)$ be a local morphism. Pick some $u\in U'$ where $f:U'\to V$ is the morphism defined by the local morphism $f$, and let $v = f(u)$. Then clearly $s_U(u) = s_V(v)$ and we can apply the previous lemma. 
\end{proof}

\begin{prop}
    Let $f: (U, s_U, t_U) \to (V, s_V, t_V)$ be a local morphism of bisubmersions which send $u \to v$, then there exists a fiber product
\end{prop}

\begin{lem}
    Let $f: (U, s_U, t_U) \to (V, s_V, t_V)$ be a local morphism of bisubmersions. Then $(V, t_V, s_V)$ is locally isomorphic to $(U, s_U, t_U)$.
\end{lem}

\begin{defn}
    We say two bisubmersions $(U, s_U, t_U), (V, s_V, t_V)$ of an atlas $\mathcal{U}$ are {\it $s$-sufficiently close}\index{ sufficiently close! $s$}, if there exists some  $(W, s_W, t_W)\in \mathcal{U}$ such that $s_W(U\times_{s_U, t_W} W) \cap s_V(V)$ is not empty. Similarly, we can define the notion of {\it $t$-sufficiently close}\index{ sufficiently close! $t$}.
\end{defn}

\begin{lem}\label{equi}
    Let $\mathcal{U}= (U_i, s_i, t_i)_{i\in I}$ an atlas of $(M, \CF)$.  Let $f: (U, s_U, t_U) \to (V, s_V, t_V)$ be a local morphism of bisubmersions of elements in $\mathcal{U}$.
Then there exists a fiber product $(U\times_{s_U, t_Wc} W, s_W, t_U)$ locally equivalent to $(V, s_V, s_V)$ which consists of identity diffeomorphisms on $s_V(V)$ for some  $(W, s_W, t_W)\in \mathcal{U}$.
\end{lem}
\begin{proof}
    Pick $u\in U'$ in the domain of local morphism, and $v = f(u)$. We have the following commutative diagram
    \begin{center}
        \begin{tikzcd}
       & u \arrow[ld, "s_U"'] \arrow[rd, "t_U"] \arrow[d, "f"] &        \\
s_U(u) & v \arrow[l, "s_V"] \arrow[r, "t_V"']                                & t_U(u)
\end{tikzcd}
    \end{center}
    Since the inverse of $(V, s_V, t_V)$ is adapted to $\mathcal{U}$. Without loss of generality, we can let  $(V, s_V, t_V)^{-1} = (V, t_V, s_V)\in \mathcal{U}$.
    Now let $(W, s_W, t_W) = (V, t_V, s_V)$. Let's consider the fiber product
    \begin{equation*}
        (W\times_{{s_W, t_U}} U, s_U, t_W)
    \end{equation*}

    We want to show this fiber product is locally equivalent to $(U, \id, \id)$. First, let $V'$ denote the domain of the induced morphism  $g: (V, s_V, t_V)\to (U, s_U, t_U)$. It suffices to show there exists a local morphism from $(U,\id, \id)$ to $(W'\times_{{s_W, t_U}} U',s_U, t_W)$  around $u\in U'$. Define $\phi:  W\times_{{s_W, t_U}}U \to V $ by
    \begin{equation*}
        \phi\big((w,u)\big) = f(u)
    \end{equation*}
    
    Then 
    \begin{equation*}
        s_V\Big(\phi\big((w,u)\big)\Big)  = s_V(f(u)) = s_V(v) = s_U\big((w,u) \big)
    \end{equation*}
    \begin{equation*}
         s_V(v)=t_W(v) = t_W\big((w,u) \big)
    \end{equation*}
    Therefore, we see that there exists a local morphism from the fiber product $(W\times_{{s_W, t_U}} U, s_U, t_W)$ to $(V, s_V, s_V)$. Therefore, we have a local equivalence.
     
    Note that, by assumption, we only have $h(w) = v$ at $v = f(u)$. But if we restrict to a small enough neighborhood of $V$, since the exponential map is uniquely determined locally, if we have an identity 
    Now we pick a $(W, s_W, t_W)\in \mathcal{U}$ such that $t_W(w) = s_U(u)$ and $s_W(w) = s_V(v)$ for some $w\in W$. This is possible since $(U, s_U, t_U)$ and $(V, s_V, t_V)$ are sufficiently close. Now we have $s_U\times_{s_U, t_W} W((u,w)) = s_V(v)$. Hence, we can apply previous lemmas to get $(U\times_{s_U, t_W} W, s_W, t_U)$ are locally equivalent to $(V, s_V, t_V)$ by the local morphism $g: V\to W$ defined by
    \begin{equation*}
        g(v) = f(u)
    \end{equation*}
    Clearly we have $S_W(F((u,w))) = s_W(w) = S_V(v)$
    
    For general case, we can connect $s_U(U)$ and $s_V(V)$ by taking fiber products with a series of bisubmersions $W_1, \cdots, W_k$.
    \begin{equation*}
        Z = U\times_{s_U, t_{W_1}} W_1\cdots \times_{s_{W_{k-1}}, t_{W_k}} W_k
    \end{equation*}
    Denote $(Z, s_Z, t_Z)$ the resulting fiber product.
    Notices that all fiber products are adapted by the atlas. Let $z\in Z$ such that $s_Z(z) = s_V(v)$, then by adeptness, there exists a $(U_i, s_i, s_j) \in \mathcal{U}$ with a local morphism $Z\to U_i$ at $z$. By previous lemmas, we see $U_i$ is locally equivalent to $Z$ at $z$. By construction, $U_i$ 
\end{proof}

\begin{proof}
Pick $U'\subset U$ such that $f:U'\to V$ is a morphism, and pick $u\in U', v = f(u)$.
Let's take two bisubmersion $(W, s_W, t_w)$ and $(X, s_X, t_X)$ such that
\begin{align*}
    s
\end{align*}

    By definition $U\times_{s_U, s_V} V = \{(u,v)\in U\times V| s_u(u) = s_v(v) \} $. Since $f$ is a morphism between $(U, s_U, t_U)$ and $(V, s_V, t_V)$, we have 
    \begin{align*}
        s_V(f(u)) &= s_U(u) =  s_V(v)
        t_V(f(u)) = t_U(u)
    \end{align*}
    hence 
    $$
    s_W(u,v) = t_v(v)
    $$
    
    there exists a fiber product $(W = U\times_{s_U, s_V} V, s_W = t_V, t_W = t_U)$ of $(U, s_U, t_U)$ and  $(V, s_V, t_V)^{-1} $
\end{proof}

Now let's construct the {\it holonomy $\infty$-groupoids} which enhances the holonomy groupoids of singular foliations.

\begin{defn}
    Let $(M, \CF)$ be a singular foliation manifold with an atlas $\mathcal{U}= (U_i, s_i, t_i)_{i\in I}$. We define a simplicial manifold $\hol^{\infty}_{\bt}(\CF)$ by
    \begin{align*}
        \hol^{\infty}_k(\CF) &= \coprod_{i_1,\cdots, i_k \in I} U_{i_1}\times_{s_{U_{i_1}}, t_{U_{i_2}}} U_{i_2}\times\cdots \times U_{i_{k-1}}\times_{s_{U_{i_{k-1}}}, t_{U_{i_k}}} U_{i_k} &k\ge 1\\
        \hol^{\infty}_1(\CF) &= M
    \end{align*}
    with face maps generated by
    \begin{align*}
        d_l(u_1,u_2\cdots, u_k) &= (u_2\cdots, u_k) & l = 0\\
        d_l(u_1,u_2\cdots, u_k) &= (u_1,u_2\cdots,u_{l-1}, (u_l, u_{l+1}),\cdots, u_k) & 1\le l \le k-1\\
        d_l(u_1,u_2\cdots, u_k)  &=(u_1,u_2\cdots, u_{k-1}) & l = k\\
    \end{align*}
    where $(u_l, u_{l+1}) \in U_{i_l}\times_{{s_{U_{i_l}}}, t_{U_{i_{l+1}}}} U_{i_{l+1}}$. 
\end{defn}

\begin{prop}
    $\hol^{\infty}_{\bt}(\CF)$ is a Lie $\infty$-groupoid.
\end{prop}
\begin{proof}
    Note that given any length $k$ fiber product
    \begin{equation*}
        U_{i_1}\times\cdots \times U_{i_k}
    \end{equation*}
    we can get its inverse simply by
    \begin{equation*}
        U_{i_k}\times\cdots \times U_{i_1}
    \end{equation*}
    with all $s_{i_l}$ and  $t_{i_l}$ switched, i.e. we take the inverse of each $(U_i, s_i, s_j)$ and then take the fiber product in the reverse order. This implies that all $k$-simplices are invertible. In particular, given a horn
    \begin{equation*}
        \Lambda^k[n] \to \hol^{\infty}_{\bt}(\CF)
    \end{equation*}
    we get a map $f:u_1 \to u_k$ as an element of $U_{i_1}\times \cdots U_{i_k}$
    , a map $g:u_k \to u_n$ as an element of $U_{i_{k+1}}\times \cdots U_{i_n}$, and a map $h:u_n\to u_1$ as an element of $U_{i_n}\times U_{i_1}$. Now we can simply take $\phi= g^{-1}\circ h \circ f^{-1}$, which gives the desired horn filling. Hence, the induced map
    \begin{equation*}
        \hol^{\infty}_{n}(\CF) \to M_{\Lambda^k[n]}\hol^{\infty}_{\bt}(\CF)
    \end{equation*}
    is clearly a surjective submersion.
\end{proof}
We call $\hol^{\infty}_{\bt}(\CF)$ the {\it holonomy $\infty$-groupoid} of the singular foliation $\CF$. Next, we shall justify the correctness of our construction. 

\begin{prop}
    The 1-truncation $\tau^{\le 1} \hol^{\infty}_{\bt}(\CF)$ of the holonomy $\infty$-groupoid is equivalent to the holonomy groupoid $\hol(\CF)$ in the sense of Androulidakis and Skandalis.
\end{prop}

\begin{proof}
    Recall that the arrow space of the holonomy groupoid is the quotient of path holonomy atlas by the relation such that $u\in U_i $ is equivalent to $v\in U_j$ if there exists a local morphism from $U_i\to U_j$ which takes $u$ to $v$. Now suppose $u\sim v$ for $u\in U_i, v\in U_j$, and there exists a local morphism $f: U_i\to U_j$ send $u$ to $v$. Then by Lemma \ref{equi}, $(v,u) \in U_j^{-1}\times_{t_{U_j}, t_{U_i}} U_i=W$ and there exists a local morphism near $(v,u)$ which maps to identity diffeomorphisms near $s_{U_j}(v)$. In addition, at $(v,u)$,
    \begin{equation*}
        s_{W}\big( (v,u)\big) = t_{U_j}(v) = t_{U_i}(u)= t_{W}\big( (v,u)\big)
    \end{equation*}
\end{proof}
\section{Holomorphic singular foliations}
\subsection{$\li$-algebroids for holomorphic singular foliations}
Let $(X, \CF)$ be a holomorphic singular foliation, i.e. $\CF$ is a locally finitely generated involutive coherent $\CO_X$-module. Since $T_{X}$ has a natural $L_{\infty}$ structure (indeed Lie structure), $\CF$ is closed with this $L_{\infty}$ structure. Denote the Dolbeault dga $(A^{0, \bullet}, \delb, 0)$ by $A$. Recall that a sheaf $\CS$ of $\cinf(X)$-modules is called {\it $\delb$-analytic coherent} if $\CS$ locally admits a resolution by finitely generated free modules and equips with a flat $\delb$-connection.

By a simple application of Hilbert's Syzygy theorem, we have that $\CF$ admits a local resolution.

\begin{prop}[\cite{LLS20}]
	Any holomorphic singular foliations on a complex manifold $X$ of dimension $n$ locally admits a finite resolution by  finitely generated
	free $\CO_X$-modules of length $\le n$. 
\end{prop}
Equivalently, we can regard $X$ has local resolution by trivial vector bundles.

The $\delb$-connection is directly inherited the flat $\delb$-connection on the $T_X$. Therefore, $\CF$ is {\it $\delb$-analytic coherent} \index{coherent!$\delb$-analytic}. In (3), the bundle $F$ is formed by taking all tangent vectors tangent to leaves, whereas in (4) a $k$-forms is in $\mathcal{J}$ if it vanishes on any $k$-tuples of vectors tangent to leaves.
Laurent-Gengoux, Lavau, Strob in \cite{LLS20} considered the case when a singular foliation on a smooth manifold admits a (global) resolution by vector bundles, and they proved that we can always construct an $L_{\infty}$-algebroid associated to that singular foliation:
\begin{thm}[\cite{LLS20}]
Give a foliation $\CF$ which admits a resolution by vector bundles $\CF_{\bt}$, there exists a {\it {universal} $L_{\infty}$}-algebroid \index{{\it {universal} $L_{\infty}$}-algebroid} $\g \in \lalgd_{\cinf{M}}$ whose linear part is the given resolution $\CF_{\bt}$. Here universal means that $\g$ is the terminal object in the category of $\lalgd_{\cinf{M}}/\CF$ which consists of $L_{\infty}$-algebroids resolving $\CF$. 
\end{thm}
We would like to generalize this to holomorphic singular foliations. 
Though we don't have global resolutions for holomorphic singular foliations, we can glue the local resolutions by higher homotopical information \cite{TT76}\cite{TT78}\cite{Blo05}\cite{Wei16}. As natural to coherent sheaves on complex manifolds, we will use the Dolbeault enhancement introduced by Block in \cite{Blo05}:

\begin{thm}[\cite{Blo05}]
	Let $X$ be a complex manifold, and $\g = T^{0, 1} X$ be the Dolbeault Lie algebroid. The homotopy category of the dg-category $\cohmod_{\CE(\g)}=\repi_{\g,A}$ is equivalent to the bounded derived category of chain complexes of sheaves of $\CO_X$-modules with coherent cohomology on 
	$X$
	
\end{thm}
As a corollary, for any coherent analytic sheaf $\CF$, there exists a cohesive module $(E, \mathbb{E})\in \cohmod_{A}$, which is unique up to quasi-isomorphism, corresponds to $\CF$. In fact, take $\CF^{\infty} = \CF \otimes_{\CO_X} \cinf(X)$ which is the $\delb$-coherent sheaf associated with $\CF$. we have a projective resolution
\begin{equation*}
0\to  E_{-n} \stackrel{d_n}{\longrightarrow}  \cdots \stackrel{d_2}{\longrightarrow} E_{-1} \stackrel{d_1}{\longrightarrow} E_{0} \stackrel{\rho}{\longrightarrow} \CF^{\infty} \longrightarrow 0
\end{equation*}
Tensoring the above sequence with the Dolbeault dga gives a resolution $E_{\bullet}\otimes_{\cinf(X)}A\to \CF^{\infty}\otimes_{\cinf(X)}A$. Denote   $\CF^{\infty}\otimes_{\cinf(X)}A$ by $F$. We can equip the dg-$A$-module $\E_{\bullet}=E_{\bullet}\otimes_{\cinf(X)}A$ with a $\Z$-connection $\nabla$ and then get a cohesive module. Denote the associated map $E_0 \to 0$ by $d_0$.

\subsection{Lifting dg-$A$-Module structure to                                             $L_{\infty}$-algebroid structure}
We shall construct $L_{\infty}$-structure on $\E_{\bullet}$ following similar methods in \cite{BFLS97} and \cite{LLS20}.

For simplicity, sometimes we will denote the $n$-ary bracket $[-,\cdots,-]_n$ by $l_n$. We regard $F\in \dgmod_A$ concentrated in degree zero, then $\rho$ naturally extends to a chain map. Let $C_{\bullet}$ and $B_{\bullet}$ denote the $d$-cycles and $d$-boundaries, respectively. We have $F\simeq H_0(\E_{\bullet})$. The existence of contracting homotopy $s:\E_{\bullet}\to \E_{\bullet-1}$ specifies a homotopy inverse $\delta: F \to \E_{\bullet}$.

\subsubsection{Construction of $l_2$}
First, we want to construct $l_2$ on $\E_{0}$. 
\begin{lem}
	There exists a  skew-symmetric $A$-linear map $\tilde{l}_2: \E_{0} \otimes \E_{0} \to \E_{0}$ satisfying
	\begin{enumerate}
		\item $ \tilde{l}_2(c_1, b_1) = 0$
		\item $[c_1, a\cdot c_2]= a[c_1, c_2] + 
		\rho(c_1)(a) c_2$
		\item $\sum_{\sigma \in \unsh(2,1)} (-1)^{\sigma} \tilde{l}_2(\tilde{l}_2( c_{\sigma(1)},c_{\sigma(2)}), c_{\sigma(3)}) \in B_0$
	\end{enumerate}
where $c_i \in \E_0$, $b_i \in B_0$, $a\in A$.
\end{lem}

\begin{proof}
	Define $\tilde{l}_2 = \delta \circ [-,-] \circ (\rho\otimes \rho)$. It is clearly skew-symmetric. Property (1) is satisfied since $\rho(b_1) = 0$.
	We claim that $[-,-] = \rho \circ \tilde{l}_2 \circ (\delta \otimes \delta)$. In fact, 
	\begin{equation*}
		d_0 \circ \tilde{l}_2 \circ (\delta \otimes \delta)= d_0\circ \delta \circ [-,-] \circ (\rho\otimes \rho)\circ (\delta \otimes \delta) 
	\end{equation*}
	since $ d_0 \circ \delta = \id$ the result follows.

	For (2), we want to construct $\delta$ explicitly. In fact, pick an open neighborhood $U\subset X$ such that $F|_U$ is finitely generated by $\{f_i\otimes a_j\}$. Pick $\{e_i\otimes a_j\}$ such that $\rho{(e_i\otimes a_j)} = f_i\otimes a_j$. Then we can define $[e_i, e_j] = \sum_{l=1}^k c_{ij}^l e_k$ where $c_{ij}^l$ comes from $[f_i, f_j] = \sum_{l=1}^k c_{ij}^l f_k$. Extend the brackets to all $E_0$ by Leibniz rules, i.e.
	$$
	[e_i, a\cdot e_j] = a[e_i, e_j] + \rho(e_i)(a)e_j
	$$. Note that here we regard $\rho(e_i \otimes a_j )\in F = \CF \otimes_{\CO_X} A $ sits inside $T_X\otimes_{\CO_X} A \subset \der_k(A)$. Finally, we glue all the local brackets by partition of unity.
	
	Next, we want to show (3) holds. On $\CF$, the Jacobi identity implies that, for $f_i\in \CF$
	$$
	\sum_{\sigma \in \unsh(2,1)}^{}(-1)^{\sigma} [[f_{\sigma(1)}, f_{\sigma(2)}], f_{\sigma(3)}]=0
	$$
	Since $[-, -] = \rho \circ \tilde{l}_2 \circ (\delta \otimes \delta)$, the left hand side becomes
	\begin{align*}
	&\sum_{\sigma \in \unsh(2,1)}^{}(-1)^{\sigma} \big(\rho \circ \tilde{l}_2 \circ (\delta \otimes \delta)\big) \circ \big( \rho \circ (\tilde{l}_2\otimes \1 ) \circ (\delta \otimes \delta \otimes \1) \big)(f_{\sigma(1)}, f_{\sigma(2)}, f_{\sigma(3)})\\     
	=& \sum_{\sigma \in \unsh(2,1)}^{}(-1)^{\sigma} \big(\rho \circ \tilde{l}_2 \circ (\delta \otimes \delta)\big) \circ \big( \rho \circ (\tilde{l}_2(\delta f_{\sigma(1)} \otimes \delta f_{\sigma(2)}) )\otimes f_{\sigma(3)}\big)\\
	=& \sum_{\sigma \in \unsh(2,1)}^{}(-1)^{\sigma}  \rho \circ \tilde{l}_2  \Big( \delta \circ \rho \circ (\tilde{l}_2(\delta f_{\sigma(1)} \otimes \delta f_{\sigma(2)}) ), \delta f_{\sigma(3)}\Big)
	\end{align*}
	Note that $d_0\circ \delta = \1_{\CF}$ and there exists a chain homotopy $s:\E_{\bullet} \to \E_{\bullet-1}$ with $\delta \circ d - \1_{\E_{\bullet}}= s\circ d + d \circ s$. Hence, we get, 
	\begin{align*}
	&\sum_{\sigma \in \unsh(2,1)}^{}(-1)^{\sigma}  \rho \circ \tilde{l}_2  \Big( (\1_{\E_0}+  d_1\circ s ) \circ (\tilde{l}_2(\delta f_{\sigma(1)} \otimes \delta f_{\sigma(2)}) ), \delta f_{\sigma(3)}\Big)\\
	=&\sum_{\sigma \in \unsh(2,1)}^{}(-1)^{\sigma} \rho \circ \tilde{l}_2  \Big(  (\tilde{l}_2(\delta f_{\sigma(1)} \otimes \delta f_{\sigma(2)}) ), \delta f_{\sigma(3)}\Big)\\
	&+ \sum_{\sigma \in \unsh(2,1)}^{}(-1)^{\sigma}  \rho \circ \tilde{l}_2  \Big( d_1\circ s  \circ (\tilde{l}_2(\delta f_{\sigma(1)} \otimes \delta f_{\sigma(2)}) ), \delta f_{\sigma(3)}\Big)
	\end{align*}
	By property (1), the second term is 0. Therefore, we have
	$$\rho	\Big(\sum_{\sigma \in \unsh(2,1)}^{}(-1)^{\sigma}  \circ \tilde{l}_2  \Big(  (\tilde{l}_2(\delta f_{\sigma(1)} \otimes \delta f_{\sigma(2)}) ), \delta f_{\sigma(3)}\Big)\Big)=0
	$$
	    which implies the term inside belongs to $B_0$.
	\end{proof}
    Next, we extend $\tilde{l}_2$ to a chain map $l_2: E_{\bullet}\otimes E_{\bullet} \to E_{\bullet}$. Consider $e_1\otimes e_0 \in E_1\otimes E_0$, then $d(_1\otimes e_0) = de_1 \otimes e_0$. We define $l_2\big(d(x_1\otimes x_0)\big)= \tilde{l}_2\big(d(x_1\otimes x_0)\big) =\tilde{l}_2\big((dx_1\otimes x_0)\big)=0$. Pick some $x_1 \in B_1$ such that $dx_1=0$, then we define $l_2(e_1\otimes e_0) = x_1$. Extend this to $E_0\otimes E_1$ skew-symmetrically. Finally, by induction, we can extend $l_2$ to all $E_i\otimes E_j$.
    
    It turns out that $l_1 = d$ and  $l_2$ satisfy a higher homotopy identity by introducing a new map $l_3$.
    
    \begin{prop}
    	There exists an almost Lie algebroid\index{almost Lie algebroid} structure on $\E_{\bullet}$.
    \end{prop}

    \begin{lem}
    There exists a degree one skew-symmetric map $l_3: \bigotimes_3 E_{\bullet} \to E_{\bullet} $ such that $l_1l_3 + l_3l_1 + l_2l_2 = 0$.
    \end{lem}
    \begin{proof}First let $e= e_1\otimes e_2\otimes e_3$ in degree 0, then by previous lemma, $l_2l_2e = b \in B_0$, hence we can find a $z\in E_1$ with $dz = b$. Now we define $l_3(e) = -z$, then $l_1l_3 + l_3l_1 + l_2l_2=0$ since $l_3l_1(e)=0$.
    	
    Next we proceed by induction, suppose the $l_3$ is constructed up to degree $k$, then $l_2l_2 + l_3l_1$ is then defined for degree $k+1$ elements. Compose with $l_1$ we have $l_1(l_2l_2 + l_3 l_1)=l_2l_2 l_1 + l_1 l_3 l_1$ used the fact that $l_1$ and $l_2$ commutes. Using the equality on degree $k$, we get $l_2l_2 l_1 + l_1 l_3 l_1 = -l_3 l_1 l_1 = 0$, hence we see $(l_2l_2 + l_3 l_1)$ on degree $k + 1$ element is a boundary $b$, then we can define the image of  $l_3$  to be a preimage of $b$ under $l_1$ just as before.
    	\end{proof}

Next, we want to construct $l_i$ for all $i\ge 3$. The idea is still the same as before. We first consider the degree 0 elements. Suppose we have construct $l_i$ for $1 \le i <n$ which satisfy strong homotopy Jacobi identities. For simplicity, we denote $\sum_{i+j = n+1}^{} (-1)^{i(j-1)} l_jl_i$ which is already summed under appropriate unshuffles. We will use the following lemma:
\begin{lem}
	Let $\{l_i\}$ define an $L_{\infty}$ structure, then we have 
	$$l_1 \sum_{i+j = n+1, i,j >1}^{} (-1)^{i(j-1)} l_jl_i = (-1)^{(n-1) }\sum_{i+j = n+1, i,j, >1}^{} (-1)^{i(j-1)} l_jl_il_1$$
\end{lem}
Hence on $\bigotimes^n X_0$, both $l_1 l_1 l_n$ and $l_1 l_n l_1$ vanishes, therefore we know $l_1l_n$ is a cycle, so we can define $l_n$ by acyclicity of $E_{\bullet}$.

Next, suppose $l_n$ has been constructed on all degree $<k$ elements.  $$
l_1 \sum_{i, j > 1, i+j = n+1}^{} (-1)^{i(j-1)} l_jl_i = (-1)^{(n-1)}\sum_{i, j >1, i+j = n+1}^{} (-1)^{i(j-1)} l_jl_il_1=0
$$
Hence, $\sum_{i, j > 1, i+j = n+1}^{} (-1)^{i(j-1)} l_jl_i$ is a cycle $b$ in $E_{n-3}$. Let $z\in E_{n-2}$ such that $l_1z = b$, then with appropriate care of signs from unshuffles, we can define $l_n = z$.

\begin{thm}\label{holo}
	Given a holomorphic foliation $\CF$ on a compact complex manifold $\CF$ , there exist an $L_{\infty}$-algebroid $\g$ over $A$, where the linear part of $\g$ corresponds to the cohesive module $E^{\bullet}$ associated to $\CF^{\infty} = \CF\otimes_{\CO_X} \cinf{X}$
\end{thm}

\subsection{Cofibrant replacement}
Since $F \in \dgmod_A$, $(\E, \nabla)$ is essentially a cofibrant replacement in the model category $\dgmod_A$. Note that we have a Quillen adjunction $\lalgd_A\leftrightarrow \dgmod_A$, a natural question is whether we can lift $(E, \mathbb{E})$ to a cofibrant object in $\lalgd_A$.

Since every strict morphism between $L_{\infty}$-morphism is an $\infty$-morphism, we have a function $\iota: \lalgd_A \to \lalgdi_A$. 

\subsection{Perfect singular foliations}

Following the idea of our construction of $L_{\infty}$-algebroid, we refine the notion of singular foliation as follows.
\begin{defn}
    Let $M$ be a $\cinf$-manifold, we define a {\it perfect singular foliation}\index{singular foliation!perfect} $\CF$ to be a subsheaf $\CF$ of $\CO_M$-module of the tangent sheaf $T_M$ such that:
    \begin{enumerate}
        \item (Perfectness)$\CF$ is a (strict) perfect $\CO_M$-module, i.e. there exists a (global) local resolution by finite projective $\CO_M$-modules
        \begin{equation*}
            0\to E_{-d}\to E_{d-1}\to \cdots \to E_{-1}\to E_0 \to \CF \to 0  
        \end{equation*}
        \item (Involutivity) $\CF$ is closed under brackets.
    \end{enumerate}
\end{defn}
In this definition, we replace the local finite generativity by local finite presentivity, which allows us to do many operations in homological algebras and homotopical operations as we did in the case of holomorphic singular foliations previously. For perfect singular foliations, we can always endow $L_{\infty}$-algebroid structures by \cite{LLS20}. Note that what we defined here is using the existing of a global resolution, which we might weaken to exist local resolutions. In that case, we call $\CF$ a {\it weakly perfect singular foliation}\index{perfect singular foliation! weakly}. In the case of compact manifolds, these two definition agrees. For weakly perfect foliations, we cannot construct $L_{\infty}$-structure directly, but we can follow the similar idea of Theorem \ref{holo}, and use {\it twisted perfect complex}\index{twisted perfect complex} (c.f. \cite{Wei16}\cite{TT76}) resolving $\CF$, and then use similar freely generating method as Theorem \ref{holo} to construct an $L_{\infty}$-algebroid structure.

We can easily define {\it perfect complex singular foliation} to be a subsheaf $\CF$ of the complexified tangent sheaf $T^{\C}_M$ which satisfies perfectness and involutivity.

We can also generalize this to derived manifolds:
\begin{defn}
    Let $(X, \CO_X)$ be a derived manifold, we define a {\it perfect singular foliation} $\CF$ to be a complex of subsheaf of $\CF$ of $\CO_X$-module of the tangent sheaf $T_X$ complex such that:
    \begin{enumerate}
        \item (Perfectness)$\CF$ is a perfect $\CO_X$-module, i.e. there exists a resolution by a double complex of finite projective $\CO_X$-modules
        \begin{equation*}
            0\to E_{-d}\to E_{d-1}\to \cdots \to E_{-1}\to E_0 \to \CF \to 0  
        \end{equation*}
        \item (Involutivity) $\CF$ is closed under brackets.
    \end{enumerate}
\end{defn}
For example, in coisotropic reduction (c.f. Example \ref{brst}), if $\CF$ has constant dimension in each stratum, then  $\CF$ is a perfect singular foliation on the derived critical locus.

Following the ideas for perfect singular foliations, we can generalize the elliptic involutive structures
\begin{defn}
Let $M$ be a $\cinf$-manifold, and $\CF\subset T^{\C}_M$ a complex singular foliation. We say $\CF$ is an {\it elliptic singular foliation}\index{singular foliation!elliptic} if 
\begin{enumerate}
    \item $\CF$ is a complex perfect singular foliation.
    \item $\CF + \overline{\CF} = T^{\C}_M$.
\end{enumerate}
\end{defn}
Clear, elliptic involutive structures are elliptic singular foliations, which we can also call {\it elliptic regular foliations}\index{elliptic regular foliation}.

\section{Elliptic involutive structures}
\subsection{Elliptic involutive structures and foliations}
Let $M$ be a compact manifold.
\begin{defn}
    A complex Lie algebroid $A$ is {\it elliptic} if its associated dga $\sym A^{\vee}[-1]$ is an elliptic complex.
\end{defn}
Note that here the ellipticity is equivalent to require $\rho(A) + \overline{\rho(A)} = T_{\C} M$.
\begin{defn}
    Let $M$ be a smooth manifold. An {\it elliptic involutive structure}\index{elliptic involutive
    structure} (EIS) consists of the following data:
    \begin{enumerate}
        \item An involutive sub-bundle $V$ of the complexified tangent bundle $T_{\C}M$.
        \item $V$ is an elliptic Lie algebroid.
    \end{enumerate}
\end{defn}
\begin{example}
    The complexified tangent Lie algebroid is clearly a trivial EIS.
\end{example}
\begin{example}
    Take $V = T^{0,1} M$ the anti-holomorphic tangent Lie algebroid, then an $EIS$ on $V$ corresponds to a complex structure on $M$.
 \end{example}
\begin{thm}[Newlander-Nirenberg]
    Let $V$ be an elliptic involutive structure on $M$. Then, locally, there exist on $M$ real coordinates $(t_1, \cdots, t_d)$ and complex coordinates $(z_1,\cdots, z_n)$ such that
\begin{equation}
    V = \Span \big\{\frac{\del}{\del t_1}, \cdots,\frac{\del}{\del t_d}, \frac{\del}{\del \overline{z}_1}, \cdots, \frac{\del}{\del \overline{z}_n}  \big\} = \Span \big\{\frac{\del}{\del z_1}, \cdots, \frac{\del}{\del z_n}^{\perp} \big\}
\end{equation}
Thus, locally, V looks like the product distribution $\R^d \oplus T^{0,1}\C^n$, where $\C^n$ has its
standard complex structure.
\end{thm}
The real part of an EIS is always a foliation on $M$, hence we also call an EIS an elliptic (regular) foliation\index{elliptic regular foliation}. 

\begin{prop}
Let $\CO_V$ denote the structure sheaf of an $EIS$ on $M$, then there exists an equivalence of categories
\begin{equation}
    \{\text{locally free sheaves of $\CO_V$ modules}\} \simeq \{\text{ finitely generated projective $V$-modules} \}
\end{equation}
\end{prop}
 
 \begin{proof}
 First let $\mathcal{E}$ be a locally free sheaf of $\CO_V$-module. Note that this condition is equivalent to that there exists a trivializing cover $\{U_i\}$ such that the transition functions of $E$ corresponding to this cover take value in $\CO_V$, where $E = \Gamma(E)\otimes_{\CO_V}\cinf(M)$ is the vector bundle corresponding to $\mathcal{E}$. Hence, in order to construct a flat $V$-connection on $E$, we just need to let the frame on $U_i\times \C^r$ to be parallel.

 Let $E$ be finitely generated projective $V$-modules, i.e. a vector bundle with a flat $V$-connection. We want to show that, for any $x\in M$, there exists a parallel local frame on some neighborhood of $x$. Let $\nabla$ be the flat connection on $E$, and $\{e_i\}$ a local frame on some neighborhood $U$ of $x$, then $\nabla e_i = \omega^j_i e_j$ where $\omega$ is the connection 1-form. By Newlander-Nirenberg, we can let $U$ be small such that $V|_U = \Span_{\C}\{dz^1, \cdots, dz^m\}^{\perp}$, where $z^1, \cdots, z^m$ are some complex coordinates on $E|_U$.
 
 Now we have $E|_U \simeq U\times \C^r$. Let $u^1, \cdots, u^r$ be complex coordinates on $\C^r$. Consider a distribution 
 \begin{equation*}
     V'= \Span_{\C} \bigg\{\frac{\del}{\del \overline{u}^k}, v - \omega^j_i(v)u^i\frac{\del}{\del u^k}| v \in V \subset T_{\C}(U\times \C^r)\bigg\}
     \subset T_{\C}(U\times \C^r)
 \end{equation*}
 By flatness of $\nabla$, we can show that this distribution is involutive. By our construction, $V' + \overline{V'} = T_{\C} U\times \C^r$, therefore, we get an EIS on $U\times \C^r$.
 According to our construction, the space of 1-form annihilate this distribution is 
 \begin{align*}
     (V')^{\perp} &= V^{\perp} + \Span_{\C} \{du^j +\tilde{\omega}^j_i u_i: j = 1,\cdots, r \}\\
 \end{align*}
 Note that $V^{\perp} = \{dz^1, \cdots, dz^m\}$ by Newlander-Nirenberg, and by same reason we also have ${V'}^{\perp} = \{d\tilde{z}^1, \cdots, d\tilde{z}^m\}$ for some coordinates $\{\tilde{z}^i\}$. Hence, we have $d\tilde{z}^j = F^j_i dz^i + G^j_i(du^i + \tilde{\omega}^i_k u^k)$. With possibly rearranging indices, we have $G^j_i \in GL_r(\C)$ in some neighborhood of zero-section. Differentiating previous equation and set $u^k = 0$ for all $k$, we get
 
 \begin{equation*}
     0 = dF^j_i\wedge dz_i + dG^j_i \wedge du^i - G^j_i \tilde{\omega}^i_k \wedge du^k
 \end{equation*}
 on $U \times \{0\}$. Pulling back to $V^{\vee}\otimes (^
 (T^{1,0}\C^r)^{\vee}$, we get   \begin{equation*}
     d_VG^j_i  - G^j_k \tilde{\omega}^k_i = 0
 \end{equation*}
 
 Now let $\tilde{\sigma}_i =(G^i_j)^{-1} \sigma_j$ on $U\times \{0\}$. Next, we want to show that $\tilde{\sigma}_k$'s are all parallel.
 
 First, note that
 \begin{align*}
     \nabla\sigma_j =& \omega^k_j \otimes \sigma_k = d_VG^k_j \otimes \tilde{\sigma}_k + G^k_j \nabla \tilde{\sigma}_k\\
     =&G^k_i\omega^i_j \otimes \tilde{\sigma}_k + G^k_j \nabla \tilde{\sigma}_k\\
     =&\omega^i_j \otimes {\sigma}_i + G^k_j \nabla \tilde{\sigma}_k
 \end{align*}
 Hence $ G^k_j \nabla \tilde{\sigma}_k = 0$, which implies $ \nabla \tilde{\sigma}_k = 0$ for all $k$. Therefore, $\{\tilde{\sigma}_k\}$ is a parallel local frame for $E$.
 \end{proof}
 
\subsection{Sheaf of $\CO_V$-modules}
\begin{defn}
    Let $(X, V)$ be a compact manifold with an elliptic involutive structure $V$. An {\it $V$-analytic} sheaf \index{$V$-analytic sheaf}on $X$ is a sheaf of $\CO_V$-modules. An $V$-analytic sheaf $\CF$ is called a coherent $V$-analytic sheaf if each point of $X$ is contained in a neighborhood $U$ such that $\CF|_U$ is the cokernel of a morphism $\CO_V^{\oplus m} \to \CO_V^{\oplus n} $ between free finite rank $V$-analytic sheaves.
\end{defn}

\begin{lem}
    Let $k \in \N$. If for every open set $U\subset M$ and every positive integer $m$, every morphism $\CO_V^{\oplus m}|_U \to \CO_V^{}|_U$ has locally finitely generated kernel, then every morphism $\CO_V^{\oplus m}|_U \to \CO_V^{k}|_U$ also has locally finitely generated kernel.
\end{lem}

\begin{proof}
We proceed by induction on $k$. The case $k= 1$ is just the assumption. Now assume $k > 1$ and we have already proved the case for all $j < k$, i.e. for every $U\subset M$, every morphism $\CO_V^{\oplus m}|_U \to \CO_V^{j}|_U$ has finitely generated kernel. Now we fix a $U$, and consider a morphism $f: \CO_V^{\oplus m}|_U \to \CO_V^{k}|_U$. We can regard $f = (g, h)$, where $g: \CO_V^{\oplus m}|_U \to \CO_V^{k - 1}|_U$ and $h: \CO_V^{\oplus m}|_U \to \CO_V^{1}|_U$ is constructed by composing $f$ with projection on the first $k - 1$ entries and the last entry respectively. By induction hypothesis, $\ker g$ is locally finitely generated. Hence, for any $x\in U$, we can find a neighborhood $V \subset U$ such that there exists some $p \in \N$ and a morphism $\phi: \CO_V^{\oplus p}|_V \to \CO_V^{m}|_V$, and $g\circ \phi$ is an exact sequence. Note that $h \circ \phi: \CO_V^{\oplus p}|_V \to \CO_V|_V$ also has a locally finitely generated kernel, by shrinking $V$ is necessary, there exists some $q \in \N$ and a morphism $\psi: \CO_V^{\oplus q}|_V \to \CO_V^{p}|_V$ which surjects on $\ker (h \circ \phi)$. Now $\ker f = \ker g \cap \ker h = \phi(\ker h\circ \phi) = \im \phi\circ \psi$, thus $\ker f$ is also locally finitely generated.
\end{proof}

\begin{lem}
Let $(p_1, \cdots, p_m)$ be an $m$-tuple of monic polynomials in $\CO_{n-1}[z_n]_0$ and $d = \max_i \deg p_i$, and let $\rho:\CO_{n}[z_n]_0^{\oplus m}$ be the morphism $\rho(f_1,\cdots, f_m) = \sum_i p_if_i$. Let $\mathcal{K}_d \subset \ker \rho$ be the subspace generated by $m$-tuples of polynomials of degree at most $d$ in $\CO_{n-1}[z_n]_0$, then $\mathcal{K}_d$ generated $\ker \rho$ as a $(\mathcal{O}_n)_0$-module. 
\end{lem}

\begin{prop}[Oka]
Let $U \subset M$ be a trivializing open set. Then each $\CO_V|_U$-module morphism $f: \CO_V|_U^{\oplus m} \to \CO_V|_U^{\oplus k}$ has locally finitely generated kernel.
\end{prop}
\begin{proof}
By previous lemma, it suffices to prove the case for $k = 1$. Recall that by Newlander-Nirenberg, $V$ locally look like the distribution $\R^d \oplus T^{0,1}\C^n$, so $V$ is split to a foliated direction spanned by $\frac{\del}{\del t_i}$'s and a transverse direction spanned by  $\frac{\del}{\del \overline{z}_j}$'s. Hence, $\CO_V|_U$ consists of $\cinf$ functions which are constant along $t_i$'s directions and holomorphic along $z_j$'s directions. Let $d= \dim \Span\{\frac{\del}{\del t_i} \}$  and $n = \dim \Span\{\frac{\del}{\del  \overline{z}_j} \}$. If $n = 0$, then the $f$ simplifies to a linear map between finite dimensional vector space, hence the kernel is also finite dimensional. We will proceed by induct on $n$. Suppose the proposition is proved for all $n < k$. Without loss of generality, $U\simeq \R^{\dim l}$, where $l = \dim M =2(d + n)$, and $x\in U$ sits at the origin in $\R^l$. $f$ has the form $f(s_1,\cdots, s_m) = \sum f_i s_i$ for an $m$-tuple of functions $\{f_i\} \subset \CO_V|_U$. By modifying coordinates if necessary, we can assume that the germ at $0$ of each $f_i$ vanishes at finite order. Note that $f_i$ are constant along $t$'s direction, i.e. we can regard $f(t_1,\cdots, t_d, z_1, \cdots z_n) = f(z_1, \cdots z_n)$ near 0. By Weierstrass preparation theorem, we can then write $f_i = u_i p_i$, where $u_i$ is a unit and $p_i$ is a Weierstrass polynomial. We can replace these germs by their representatives in some neighborhood $U'$ of $0$. Shrinking $U'$ if necessary, we can assume that $u_i$'s are non-vanishing in $U'$. Now the map $(t_1,\cdots, t_d,z_1,\cdots, z_n) \to (t_1,\cdots, t_d,u_1z_1,\cdots, u_nz_n)$ is an automorphism of $\CO^{\oplus m}|_U$ which maps the kernel of $f$ to the kernel of the map determined by the $m$-tuple $(p_1, \cdots, p_m)$. Hence, without loss of generality, we may replace $f_i$'s by $p_i$'s, i.e. $f = \sum f_i g_i$.
We may assume $U'$ has the form $\R^{2d}\times U'' \times U'''$ for some $U''\subset \C^n, U''' \subset \C$.
Let $d = \max_i \deg p_i $, and let $\mathcal{K}_d$ denote the sheaf on $U''$ defined as follows: for each open subset $W \subset U''$, $\mathcal{K}_d(W) \subset \ker f$ is the subspace consisting of $m$-tuple of polynomials of degree less than or equal to $d$ in $\CO_{n - 1}(U'')[z_n]$, where $\CO_{n - 1}(U'')$ denotes the sheaf of holomorphic functions on $U''$. We need to show that $\mathcal{K}_d$ is locally finitely generated as a $\CO_{V}(U'')$-modules. For each neighborhood $W \subset U''$, the space of $m$-tuples $(q_1,\cdots, q_m)$ forms a free module of rank $(d+1)m$ over $\CO_{n - 1}(U'')$, where each $q_i$'s is a polynomial of degree less than or equal to $d$. Hence, $f$ gives an $\CO_n(U'')[z_n]$-module morphism from a rank $(d + 1)m$ free module to a rank $(2d + 1)$ free module which consists of polynomial in $z_n$ of degree at most $2d$. Note that $\mathcal{K}_d$ is exactly the kernel of this restricted morphism. By induction hypothesis, we get that $\mathcal{K}_d$ is locally finitely generated.

Finally, we need to show that $\mathcal{K}_d$ actually generate the whole $\ker f$ as a sheaf of $\CO_V$-modules. It suffices to show that the stalk $(\mathcal{K}_d)_0$ at the origin generated the stalk $(\ker f)_0$ over $(\CO_n)_0$. Note that here we regard $\CO_n$ locally constant in the real directions. The $p_i$'s may not be Weierstrass polynomials when taking the germs at the origin, but they are still monic polynomials in $z_n$. The proof then follows from the previous lemma.

\end{proof}

\begin{cor}
Let $U$ be a trivializing open neighborhood of $(X,V)$, and $\mathcal{F}$ is a locally finitely generated sheaf of submodules of $\CO_V^{\oplus k}$, then $\mathcal{F}$ is $V$-coherent. In particular, given a morphism of $V$-analytic sheaves $\phi: \CO_V^{\oplus m}(U) \to \CO_V^{\oplus k}(U)$, then both $\ker \phi$ and $\im \phi$ are $V$-coherent.
\end{cor}

\begin{proof}
Let $x\in U$ be arbitrary, then there exists $W\ni x$ such that $\CF$ is the image of a morphism $\phi: \CO_V^{\oplus m}(W) \to \CO_V^{\oplus k}(W)$ of sheaves of $\CO_V$-modules. By Oka's theorem, its kernel is also locally finitely generated. With shrinking $W$ if necessary, we can assume there exists a morphism $\psi: \CO_V^{\oplus p}(W) \to \CO_V^{\oplus m}(W)$ which surjects on $\ker \phi$ on $W$. Therefore, $\CF$ is $V$-coherent.

Next, by assumption, $\im \phi$ is already finitely generated on $U$, and the kernel is also locally finitely generated, which implies $\im \phi$ is also $V$-coherent.
\end{proof}
\begin{cor}
Let $U$ be a trivializing open neighborhood of $(X, V)$, and $\mathcal{F}$ and $\mathcal{G}$ are $V$-coherent sheaves of submodules of $\CO^{\oplus m}_V(U)$, then $\mathcal{F}\cap \mathcal{G}$ is also $V$-coherent.
\end{cor}

\begin{proof}
For every point $x\in U$, there exists an open neighborhood $W\ni x$ such that $\CF|_W$ and $\CG|_W$ are images of some morphisms $\phi: \CO_V^{\oplus p}|_W \to \CO_V^{\oplus m}|_W$ and $\psi: \CO_V^{\oplus q}|_W \to \CO_V^{\oplus m}|_W$ respectively. Consider the map $\theta: \CO_V^{\oplus (p+q)}|_W \to \CO_V^{\oplus m}|_W$ by $\theta(f\oplus g) = \phi(f) - \psi(g)$. By Oka's theorem, $\ker \theta$ is $V$-coherent. Note that $\CF \cap \CG$ is the image of $\ker \theta$ under $\phi$. With shrinking $W$ if necessary, we can choose finitely many generators for $\ker \theta$, whose image under $\phi$ will then generate $\CF\cap \CG$. Hence, $\CF\cap \CG$ is locally finitely generated over $\CO_V$. By similar reasoning in the previous corollary, $\CF\cap \CG$ is $V$-coherent.
\end{proof}

\begin{prop}
Let $(X, V)$ be a manifold equipped with an elliptic involutive structure $V$, and $\CF$ a coherent $V$-analytic sheaf, then locally $\CF$ admits finite resolution of length less than or equal to $n + 1$  by free sheaves of modules.
\end{prop}
Let $M$ be a module over $\CA^{\bullet}$, then it localizes to a sheaf of $\CA^{\bullet}_X$-modules by taking $M_X(U) = M \otimes_{\CA^{\bullet}} \CA^{\bullet}_X(U)$. Let $(E^{\bullet}, \nabla) \in \mathcal{P}_{\CA^{\bullet}}$, define a double complex of sheaves $\mathcal{E}^{p,q}$ by $\mathcal{E}^{p,q}(U) = E^p \otimes_{\CA^{\bullet}} \CA^{q}_X(U)$. Define $(\mathcal{E}_X^{\bullet}, \nabla) = (\sum_{p+q = \bullet} \mathcal{E}^{p,q}_X, \nabla)$. Note that $\mathcal{E}_X^{\bullet}$ is a complex of soft sheaves of $\CO_V$-modules.

Analogous to Pali's definition of $\overline{\del}$-coherent analytic sheaves, we define $\overline{\del}_V$-coherent analytic sheaves\index{coherent!$\overline{\del}_V$-coherent analytic sheaf} for elliptic involutive structures.
\begin{defn}
    Let $(X, V)$ be an elliptic involutive structure. We define a {\it $\overline{\del}$-coherent analytic sheaf} $\CF$ to be a sheaf of modules over the sheaf of $\cinf$-functions $\cinf_X$ with
    \begin{enumerate}
        \item Finiteness: $\CF$ has locally finite resolution by finitely generated free modules over $\cinf_X$.
        \item $V$-analytic: $\CF$ is equipped with a flat $\overline{\del}_V$-connection, i.e. an operator $\overline{\del}_V: \CF \to  \CF \otimes_{\cinf_X} A^1_X$ with $\overline{\del}_V^2 = 0$. 
    \end{enumerate}
\end{defn}

\begin{prop}
The functor $\alpha: \ho \mathcal{P}_{A^{\bullet}} \to \dperf(X, \CO_V) \simeq \dcoh(X, \CO_V)$ defined by
\begin{equation}
    \alpha:(E^{\bullet}, \nabla) \mapsto (\mathcal{E}_X^{\bullet}, \nabla)
\end{equation}
is fully faithful.
\end{prop}
\begin{proof}
    Let $U = \{(t_1,\cdots, t_d, z_1,\cdots, z_n)| \ |t_i| < r, |z_i| < r\}$ be a polydisc in $X$. We want to show that there exists a small polydisc $V$ such that there exists a gauge transformation $\phi:\mathcal{E}^{\bullet}|_V \to \mathcal{E}^{\bullet}|_V$ of degree 0 such that 
    $\phi \circ \nabla \circ \phi^{-1} = \tilde{\nabla}^0 + \overline{\del}_V$. Hence, $\mathcal{E}^{\bullet}_X$ is gauge equivalent to a complex of finitely generated projective $V$-modules, that is, $H^p\big((\mathcal{E}^{\bullet, 0}), \nabla^0 \big)$ is $\overline{\del}_V$-coherent with $\overline{\del}_V$-coherent connection $\nabla^1$ for each $p$. Note that $U = U_1\times U_2$ where $U_1$ is contractible and $U_2$ is Stein, so there is no higher cohomology with respect to $\nabla^1$, then we are left with $V$-analytic sections over $U$, which are then coherent.
    
    The construction of $\phi$ follows from the proof of integrability of holomorphic structures on vector bundles. As we are in a polydisc $U$, we can write the $\mathbb{Z}$-connection $\nabla$ as $\nabla = \nabla^0 + \overline{\del}_V + J$ where 
    $$ J: \mathcal{E}^{p, q}(U) \to \bigoplus_{i \le p} \mathcal{E}^{i, q+(p-i) +1}(U)
    $$ is a $\cinf_X(U)$-linear map. Decompose $J$ as $J = J'\wedge d\overline{z}_1 + J''$ with $\iota_{\frac{\del}{\del \overline{z}_1}} J' =\iota_{\frac{\del}{\del \overline{z}_1}} J'' = 0$. For simplicity, write $\overline{\del}_i = d\overline{z}_i \wedge \frac{\del}{\del \overline{z}_i}$. We want to find a $\phi_1$ with $\phi_1(\overline{\del}_1 + J'\wedge d\overline{z}_1) \phi_1^{-1} = \overline{\del}_1$. It suffices to solve $\phi_1^{-1}\overline{\del}_1(\phi_1) = J'\wedge d\overline{z}_1 $ and treat $t_1,\cdots, t_d, z_2, \cdots, z_n$ as variables. Now we set $\nabla_1 = \phi_1(\nabla^0 + \overline{\del}_V + J' + J'')\phi_1^{-1}$. We can write $\nabla_1 = \nabla^0_1 + \overline{\del}_1 + \overline{\del}_{\ge 2} + J_1$. We claim that:
    \begin{enumerate}
        \item $\nabla^0_1\circ \nabla^0_1 = 0.$
        \item $\nabla^0_1$ and $J_1$ are holomorphic in $z_1$.
        \item $\iota_{\frac{\del }{\del \overline{z}_1}} J_1 = 0$
    \end{enumerate}
    Notice that 
    \begin{align*}
        0 =& \iota_{\frac{\del }{\del \overline{z}_1}}(\nabla_1 \circ \nabla_1)\\
        =&\iota_{\frac{\del }{\del \overline{z}_1}}(\nabla^0_1\circ \overline{\del}_1 + \overline{\del}_1 \circ \nabla^0_1 + J_1\circ \overline{\del}_1 +  \overline{\del}_1 \circ J_1 )\\
        =& \iota_{\frac{\del }{\del \overline{z}_1}}(\overline{\del}_1(\nabla^0_1) + \overline{\del}_1(J_1))
    \end{align*}
    For degree reason in the $p$-direction, the two summand in the bracket must both be zero. Therefore, we have proved the claim.
    
    Next, we shall iterate this procedure. Write $J_1 = J_1' \wedge d\overline{z}_2 + J_1''$ with $\iota_{\frac{\del}{\del \overline{z}_1}} J_1'= \iota_{\frac{\del}{\del \overline{z}_2}} J_1' =\iota_{\frac{\del}{\del \overline{z}_1}} J_1'' = \iota_{\frac{\del}{\del \overline{z}_2}} J_1''=0$, and we want to find a $\phi_2$ with $\phi_2(\overline{\del}_2 + J_1'\wedge d\overline{z}_2) \phi_2^{-1} = \overline{\del}_2$ and $\phi_2(\overline{\del}_1)\phi_2^{-1}= \overline{\del}_1$. Then it suffices to solve $\phi_2^{-1}\overline{\del}_2(\phi_2) = J_1'\wedge d\overline{z}_2 $. Note that $\phi_2$ is holomorphic in $z_1$ since $J_1'$ is. Now set $\nabla_2 = \phi_2 \circ \nabla_1 \circ \phi_2^{-1}$, and we can write $\nabla_2$ as $\nabla_2 =  \nabla^0_2 + \overline{\del}_1 + \overline{\del}_{2}+ \overline{\del}_{\ge 3} + J_2$ with $\iota_{\frac{\del }{\del \overline{z}_1}} J_2 = \iota_{\frac{\del }{\del \overline{z}_2}} J_2 =0$
    
    Continuing this fashion, we will arrive at $\nabla_n = \nabla^0_n + \overline{\del} + J_n$, where $\iota_{\frac{\del }{\del \overline{z}_i}} J_n =0$ for all $i = 1, \cdots, n$. Hence, now it suffices to deal with the real directions. Again, write $d_i = dt_i \wedge \frac{d}{d {t}_i}$ and $J_n = J_n'\wedge dt_1 + J_n''$ with $\iota_{\frac{\del}{\del t_1}}J_n' = \iota_{\frac{\del}{\del t_1}}J_n'' = 0$. We want to find $\psi_1$ with $\psi_1(d_1 + J_n'\wedge dt_1) \psi_1^{-1} = d_1$ and $\psi_1(\overline{\del}_i)\psi_1^{-1}= \overline{\del}_i$, which can be done by solving $\psi_1^{-1}d_2(\psi_1) = J_n'' \wedge dt_1$. Now set $\nabla_{n+1} = \psi_1\circ \nabla_n \circ \psi_1^{-1}$. We can write $\nabla_{n+1} = \nabla_{n+1}^0 + \overline{\del}+d_1 + d_{\ge 2} + J_{n+1}$. Similar to the previous argument, we can easily show that 
    \begin{enumerate}
        \item $\nabla_{n+1}^0\circ \nabla_{n+1}^0 = 0.$
        \item $\nabla_{n+1}^0$ and $J_{n+1}$ are flat in $t_1$.
        \item $\iota_{\frac{\del }{\del t_1}} J_{n+1} = 0$
    \end{enumerate}

    Iterating this procedure, we will reach $\widetilde{\nabla} = \nabla_{n+d} + \delb + d = \nabla_{n+d} + \delb_V$.

\end{proof}
\begin{lem}
    On an elliptic involutive structure $(X, V)$, $\cinf_X$ is flat over $\CO_V$.
\end{lem}

\begin{thm}[Block\cite{Blo05}]\label{quasicohesive}
    Suppose $(A^{\bullet}, d,c)$ is a curved dga. Let $X = (X, \nabla)$ be a quasi-cohesive module over $A^{\bullet}$, then there is an object $E = (E^{\bullet}, \nabla')$ in $\mathcal{P}_{A^{\bullet}}$ such that $\tilde{h}_X$ is quasi-isomorphic to $h_E$, under either of the two following conditions:
    
    \begin{enumerate}
        \item X is a quasi-finite quasi-cohesive module.
        \item $A^{\bullet}$ is flat over $A^{0}$ and there exists a bounded complex $(E, \nabla'^0)$ of finitely generated projective right $A^0$-modules and an $A^0$-linear quasi-isomorphisms $e^0: (E, \nabla'^0) \to (X, \nabla^0)$.
    \end{enumerate}
\end{thm}

\begin{prop}
Let $(\mathcal{E}^{\bullet}_X, d)$ be a complex of sheaf of $\CO_V$-modules with coherent $V$-analytic cohomology, then there exists a cohesive $A^{\bullet}$-module $E = (E^{\bullet}, \nabla)$unique up to quasi-isomorphism, and $\alpha(E)$ is  quasi-isomorphic to $(\mathcal{E}^{\bullet}_X, d)$. In addition, for two such complexes $(\mathcal{E}^{\bullet}_1, d)$ and $(\mathcal{E}^{\bullet}_2, d)$, the corresponding cohesive modules $E_1$ and $E_2$ satisfies

\begin{equation}
    \ext_{\CO_V}^k(\mathcal{E}^{\bullet}_1, \mathcal{E}^{\bullet}_2) \simeq H^k(\Hom_{\mathcal{P}_{A^{\bullet}}}(E_1, E_2))
\end{equation}
\end{prop}
\begin{proof}
 Without loss of generality, we can assume $(\mathcal{E}^{\bullet}_X, d)$ is a perfect complex over $\CO_V$. Define $\mathcal{E}^{\bullet}_{\infty} =\mathcal{E}^{\bullet}_X\otimes_{\CO_V}\cinf_X$. By flatness of $\cinf_X$ over $\CO_V$, $(\mathcal{E}^{\bullet}_X, d)$ is a perfect complex of $A_X$-modules, and the map $(\mathcal{E}^{\bullet}_X, d)\to (\mathcal{E}^{\bullet}_X\otimes_{\CO_V}A_X^{\bullet}, d\otimes 1 + 1\otimes \overline{\del}_V)$ is a quasi-isomorphism. By proposition, there exists a (strict) perfect complex $(E^{\bullet}, \nabla)$ of $\cinf(X)$-modules and a quasi-isomorphism $e^0:(E^{\bullet}, \nabla) \to (\Gamma(X, \mathcal{E}^{\bullet}_{\infty}), d)$. $(\Gamma(X, \mathcal{E}^{\bullet}_{\infty}), d)$ defines a quasi-cohesive module over $A^{\bullet}$, hence the result follows from Theorem \ref{quasicohesive}.
\end{proof}
In summary, we have just proved that

\begin{prop}
Let $(X, V)$ be a compact manifold $X$ with an elliptic involutive structure $V$, then there exists an equivalence of categories between $\dcoh(X)$, the bounded derived category of
complexes of sheaves of $\CO_V$-modules with coherent $V$-analytic cohomology, and $\ho \mathcal{P}_{\CA^{\bullet}}$, the homotopy category of the dg-category of cohesive modules over $\CA^{\bullet} = \sym V^{\vee}[-1]$, i.e.   
\begin{equation*}
    \dcoh(X, \CO_V) \simeq \ho \mathcal{P}_{\CA^{\bullet}}
\end{equation*}
\end{prop}

We can deduct analogous result of Pali for coherent $V$-analytic sheaves and $\overline{\del}_V$-coherent analytic sheaves. 
\begin{cor}[] The category of $V$-coherent analytic sheaves on $X$ is equivalent to the category of coherent $V$-analytic sheaves.
    
\end{cor}

\subsection{Cauchy-Riemann structures}
Let $M$ be a $2n + 1$ dimensional smooth manifold. An {\it almost Cauchy-Riemann structure} on $M$ is a sub-bundle $L$ of the complexified tangent bundle $T_{\C} M$ such that $L \cap \overline{L} = 0$. We say $L$ is a {\it Cauchy-Riemann structure} on $M$ is $L$ is involutive.

\begin{defn}
    Let $(M, L)$ be a CR manifold. Let $f\in \cinf_{\C}(M)$, then we say $f$ is a {\it Cauchy-Riemann} or simply {\it CR} function, if for all $Z \in \overline{L}$, $Z(f) = 0$
\end{defn}

\begin{prop}
Let $\CO_L$ denote the structure sheaf of a $CR$-structure on $M$, then there exists an equivalence of categories
\begin{equation}
    \{\text{locally free sheaves of $\CO_L$ modules}\} \simeq \{\text{ finitely generated projective $L$-modules} \}
\end{equation}
\end{prop}

\part{Higher monodromy and holonomy}

\section{Monodromy}
\subsection{Higher monodromy}
Let $(M, \CF)$ be a regular foliation, we consider the collection $\CG$ of all homotopy classes of paths lying in the same leaf, which form a smooth manifold and gives the {\it monodromy groupoid}\index{monodromy groupoid} of $\CF$. 

\begin{defn}
	A map $p: X \to Y$ is called {\it semi-locally simply connected}\index{semi-locally simply connected} if given any $x\in X$, there is a basic open neighborhood $V$ of $p(x)$ and a basic open neighborhood of $x$ such that $p(U) \subset V$ and the following diagram commutes and the lift exists.
	\begin{center}
		\begin{tikzcd}
		\del I^2 \arrow[r] \arrow[d]            & U \arrow[r] &  V\times_Y X \arrow[d] \\
		I^2 \arrow[rr] \arrow[rru, dotted] &             & V          
		\end{tikzcd}
	\end{center}
\end{defn} 

\begin{prop}
	Let $(M, \CF)$ be a regular foliation, then it induces a semi-locally simply connected map $p: M \to [M/\CF]$.
\end{prop}
\begin{proof}
	Let $x\in X$, then we can pick a open neighborhood $V$ of $p(x)$ which is contained in a single foliation chart $\R^q \times \R^{n-q}$, where $q=\dim \CF$ and $n = \dim M$.  $V\times_Y X$ equals the union of all leaves lying in $[M/\mathcal{F}]$, i.e. $V\times_Y X = \coprod_{x\in [M/\mathcal{F}]} L_x$. Note that $V\times_Y X \to V$ is a submersion hence the lift always exists. 
	\end{proof}
Similarly, we can show any submersions are semi-locally simply connected.
\begin{prop}
	Let $\pi: X\to Y$ be a submersion, then $\pi$ is semi-locally simply connected.
\end{prop}

\begin{cor}
	A smooth Serre fibration is semi-locally simply connected.
\end{cor}
\begin{proof}
	This follows from the fact that all smooth Serre fibrations are submersions.
	\end{proof}

\subsection{The monodromy $\infty$-groupoid $\mon_{\infty}(\CF)$}
\begin{defn}
	The {\it Monodromy  $\infty$-groupoid $\mon_{\infty}(\CF)$ of a foliation $(M, \CF)$, also denoted by $\Pi_{\infty} (\CF)$}, is a simplicial space whose $n$-simplices are $$\map(\Delta^n, \mathcal{F})= \map_{vert}(\Delta^n, M) = \{f: \Delta^n \to M| f(\Delta^n) \text{ lies in a single leaf} \}$$
	
	where $\Delta^n$ denotes the {\it geometric $n$-simplex}.
\end{defn}

\begin{prop}
	$\mon_{\infty}(\CF)$ is a simplicial space.
\end{prop}
\begin{proof}
	The topology on $n$-simplices $\map(\Delta^n, \mathcal{F})$ is inherited from the compact-open topology on $\map(\Delta^n, M)$. Given  $x\in \map(\Delta^n, \mathcal{F})$ and $V\subset [M/\mathcal{F}]$, a basic open neighborhood of $x$ has the form 
	$$<x, U> = \{y \in \map_{vect}(\Delta^n, M)\big|\exists h: \Delta^n \times \Delta^1 \to \pi^{-1}(U), h(-, 0) = x, h(-, 1)=y \}$$. Since the degeneracy maps $s^k: \Delta^n \to \Delta^{n+1}$ and face maps $d^k:\Delta^{n+1} \to  \Delta^n $ are all continuous, the face maps and degeneracy maps in $\mon_{\infty}(\CF)$ are all continuous. Therefore, $\mon_{\infty}(\CF)$ is a simplicial topological space.
	
	\end{proof}

Next, we want to explore the $\infty$-groupoid structure of $\mon_{\infty}(\CF)$. First, we start with a definition of topological $\infty$-groupoids.

\begin{defn}
	Let $X_{\bullet}$ be a simplicial space, we say $X_{\bullet}$ is a {\it topological $\infty$-groupoid}\index{topological $\infty$-groupoid} if all its structure maps are continuous, and diagrams of the following form commutes and the lift exists for $0\le i \le k, 0\le k < \infty$.
	\begin{center}
		\begin{tikzcd}
		\Lambda^i[k] \arrow[r] \arrow[d]          & X_{\bullet}  \\
		\Delta[k]\arrow[ru, dotted] &       
		\end{tikzcd}
	\end{center}
Note that here all maps are continuous, and $\Delta[n]$ here denotes the {\it standard $n$-simplex} rather than the geometric $n$-simplex.
\end{defn} 
We denote the space of $i$-th $k$-horns $\Lambda^i[k] \to X_{\bullet} $ by $ X(\Lambda^i[k])$ and the $n$-simplices $\Delta [k] \to X_{\bullet} $ by $X_n$. Since $L_x$
\begin{prop}
	$\mon_{\infty}(\CF)$ is a topological $\infty$-groupoid.
\end{prop}

\begin{proof} 
	Let $\sigma:\Lambda^i[k]\to \mon_{\infty}(\CF)$ be a $k$-horn. Note that
	$$
	X(\Lambda^i[k]) = \map(\big|\Lambda^i[k]\big|, \CF)
	$$
	where $\big|\Lambda^i[k]\big|$ denote the $i$-th geometric $k$-horn. Hence, the image of the standard $k$-horn $\big|\Lambda^i[k]\big|$ lies completely within a single leaf $L_x$ for some $x\in M$. Therefore, 
	$\sigma $ can be regard as a map from $ \Lambda^i[k]$ to $ \Pi_{\infty}(L_x)$, where $\Pi_{\infty}(L_x)$ denotes the fundamental $\infty$-groupoid of the leaf $L_x$. Therefore, the lift exists.
	\end{proof}

\subsection{Smooth monodromy $\infty$-groupoid $\mathcal{P}_{\infty}(\CF)$}

Next, we consider a smooth refinement of the $\mon_{\infty}(\CF)$. 

First, we recall the definition of $A$-path.
\begin{defn}
	Let $\pi: A\to M$ be a Lie algebroid with an anchor map $\rho:A\to T_M$. A $C^1$ curve $a: \Delta^1 \to A$ is called an $A$-path if
	$$ \frac{d}{dt}(\pi\circ a(t)) = \rho\big(a(t) \big)$$.
	
	If, in addition, $a(t)$ satisfies the following boundary conditions
	$$ a(0)=a(1) = 0, \quad \dot{a}(0) = \dot{a}(1) = 0$$
	then we say $a$ is an $A_0$-path.
\end{defn}
We want to generalize these to higher dimensions. Recall given a foliation $(M, \CF)$, there is an associated Lie algebroid $\CF \stackrel{\rho}{\rightarrow} T_M$ where the $\rho$ is simply the inclusion.

\begin{defn}
	Let $\sigma: \Delta^n \to \CF$ be a differentiable ($C^1$ or $C^{\infty}$) map such that
	\begin{enumerate}
		\item $T_M \big|_{ \pi \circ \sigma (\Delta^n)}  \subset \CF\big|_{ \pi \circ \sigma (\Delta^n)}$.
		\item For any piecewise smooth path $\gamma: \Delta^1 \to \im (\sigma)$, we have
		$$
		\frac{d}{dt}(\pi\circ \gamma)(t) = \gamma(t)
		$$
	\end{enumerate}
\end{defn}
We call such a map $\sigma: \Delta^n \to \CF$ a {\it $C^1$($C^{\infty}$) foliated $n$-simplex}. 

\begin{prop}
	The space of $C^1$ foliated $n$-simplices $P^n_{C^1} \mathcal{F}$ is a Banach manifold. The space of $C^{\infty}$ foliated $n$-simplices $P^{n}_{C^{\infty}} \CF$ is a Frech{\'e}t manifold.
\end{prop}

\begin{proof} First, consider the $C^{1}$ case. Let $(\sigma: \Delta^n \to \mathcal{F} )\in P^n_{C^1} \CF$. Pick up any Riemannian metric on $\CF$. Let $T_{\epsilon} \subset \gamma^* T_{\CF}$ consist of tangent vector of length $\le \epsilon$. For $\epsilon$ small, we have the exponential map	
$\exp: T_{\epsilon} \to \CF, (x, v) \mapsto \exp_{\sigma(x)} v$. Denote the $C^1$ section of $T_{\epsilon}$ by $PT_{\epsilon}$. Note that $\gamma^*T_{\CF}\simeq \Delta^1 \times \R^{n-q}$ is a trivial $(n-q)$ bundle, hence any trivialization will give a map $PT_{\epsilon} \to \map_{C^1}(\Delta^1, \R^n)$. Since $\map_{C^1}(\Delta^1, \R^n)$ is a Banach space, $PT_{\epsilon}$ gives a chart for $P^n_{C^1} A$.

Similarly, we can show the space of $C^r$ foliated $n$-simplices $P^n_{C^r} \CF$ is also a Banach manifold for $1<r<\infty$. Hence, $P^{n}_{C^{\infty}}\CF$ is a  Frech{\'e}t manifold.
	\end{proof}

The chart constructed in above gives a deformation of any smooth path. Given a path $\gamma: \Delta^1 \to M$, $T_{\epsilon} \simeq I \times D_{\epsilon}$ where $D_{\epsilon}$ is the $\epsilon$-disk in $\R^{n-q}$. For any section $\tilde{\gamma}\in \Gamma T_{\epsilon}$ lying above $\gamma$, the  exponential map then yield a smooth map $\phi: I\times D_{\epsilon} \to M$ by $\phi(t,u) = \exp_{\sigma(t)}\tilde{\gamma}(t, u)$. We call $\phi$ the {\it universal deformation of $\gamma$}.

Now we can define the smooth analogue of monodromy $\infty$-groupoid.
\begin{defn}
	Given a foliation $(M, \CF)$, define a {simplicial Frech{\'e}t manifold} $\mathcal{P}_{\infty}(\CF)$ whose $n$-simplices is $P^n_{C^{\infty}}$. Similarly, we define a simplicial Banach manifold
	$\mathcal{P}_{\infty}^{C^1}(\CF)$ whose $n$-simplices is $P^n_{C^{1}}$.
\end{defn}

\begin{rem}
	The definition of  $\mathcal{P}_{\infty}(\CF)$ is similar to {\it path $\infty$-groupoid} in literature. However, we don't mod out the thin homotopy classes of path in each level, since we want to keep the manifold structure. For example, the space of morphisms of the path 1-groupoid of a manifold is not a manifold in general.
\end{rem}

\begin{prop}
	$\mathcal{P}_{\infty}(\CF)$ and  $\mathcal{P}_{\infty}^{C^1}(\CF)$ are Lie $\infty$-groupoids.
\end{prop}
\section{Integrating derived $L_{\infty}$-algebroids}

\subsection{Integrating Lie algebroids}
In this section, we consider the integration of Lie algebroids to a Lie $\infty$-groupoids.

\begin{defn}
    Let $A\stackrel{\rho}{\to} TM$ be a Lie algebroid. Define a simplicial manifold $\CG_{\bullet}$ with
    $$
    \CG_n = \Hom_{\algd} (T\Delta_n), A).
    $$
\end{defn}

\begin{prop}
$\mathcal{G}_{\bullet} $ is a Lie $\infty$-groupoid.
\end{prop}
Hence, $\CG_{\bullet}$ presents the $\infty$-stack which globalizes $A$. The 2-truncation of $\CG_{\bullet}$ corresponds to the Lie 2-groupoid integrates $A$, which is equivalent to a Weinstein groupoid or a stacky Lie groupoid.  However, $\CG_n$ are all infinite dimensional for $n \ge 1$.

Parametrize $n$-simplex $\Delta_n$ by $\{1 \ge t_1\ge t_2 \ge \cdots \ge t_n \ge 0\}.$
\begin{prop}
Let $\alpha_1^t, \cdots, \alpha_n^t$ be time-dependent sections of $A$, where $t= (t_1, \cdots, t_n) \in \Delta_n$. Suppose the following equation holds:
\begin{equation}
    [\alpha_i, \alpha_j] = \frac{d \alpha_i}{d t_j} - \frac{d \alpha_j}{d t_i}
\end{equation}
for all $1\le i, j \le n$. Then there exists a family of time-dependent vector fields $X_i(x, t) = \rho(\alpha^t_i(x)) + \del t_i$. For any $x_0 \in M$, there exists a $n$-simplex $\sigma: \Delta_n \to M$ with
\begin{align*}
    \frac{d \sigma}{d t_i} =& X^t_i(\sigma(t))\\
    \sigma(0) =& x_0.
\end{align*}
In addition, Let $a_i(t) = \alpha_i^t(\sigma(t))$, then $a = \sum_{i=1}^n a_i dt_i: T\Delta_n \to A$ defines a Lie algebroid morphism.

\end{prop}

\begin{proof}
       First, we want to show the existence of the simplex $\sigma$. Applying the anchor map, we have 
    \begin{equation*}
        [\rho(\alpha_i), \rho(\alpha_j)] = 
    \end{equation*}
Hence $X_i + \del t_i$ mutually commutes as vector fields on $M \times \Delta_n$. Since we assumed that  for any $i\in \{1, \cdots, n\}$, the flow of $X_i$ exists on $\{ t_i = t_{i-1} \}$ up to time 1.   Note that $X_i$ is a $(t_1, \cdots, \hat{t_i}, \cdots, t_n )$-family of $t_i$-time-dependent vector fields. Denote the time-dependent flow of $X_i$ by $\Phi^{X_i}_{s_i', s_i}$. Hence, by previous observation we have
\begin{equation*}
    \Phi^{X_i^{(t_1, \cdots, \hat{t_i}, \cdots, s_j', \cdots, t_n )}}_{s_i', s_i} \circ \Phi^{X_j^{(t_1, \cdots, s_i, \cdots, \hat{t_j}, \cdots, t_n )}}_{s_j', s_j} = \Phi^{X_j^{(t_1, \cdots, s_i', \cdots, \hat{t_j}, \cdots, t_n )}}_{s_j', s_j} \circ \Phi^{X_i^{(t_1, \cdots, \hat{t_i}, \cdots, s_j, \cdots, t_n )}}_{s_i', s_i}
\end{equation*}

Hence, we can construct $\sigma$ as 
\begin{equation*}
    \sigma(t_1, \cdots, t_n) = \Phi^{X_n^{(t_1, \cdots, t_{n-1}, \hat{t_n} )}}_{t_n, 0}\circ \cdots \Phi^{X_i^{(t_1, \cdots, \hat{t_i}, 0, \cdots, 0 )}}_{t_i, 0}\circ \cdots \Phi^{X_1^{(t_1,  0, \cdots, 0 )}}_{t_1, 0}
\end{equation*}

Next, we want to show $a$ defines a Lie algebroid morphism, i.e. the induced map $a^*: \CE(A)^{\bullet} \to \Omega^{\bullet}\Delta_n$ is a dga morphism. Note that it suffices to check on degree 0 and 1, where all higher degree terms follow from Leibniz rules.

Let $\{e_1, \cdots, e_m\}$ be a basis of sections of $A$ around an open neighborhood $U$ of $\sigma(t_0) = x_0$. Let $c^l_{p,q}$ be the structure constants of $A$ on $U$, i.e. $[e_i, e_j] = \sum_{l=1}^m c^l_{p,q} e_l$ for $p,q = 1, \cdots m$. 
Write $\alpha^t_i = \sum_{p=1}^m \alpha^t_{i,p} e_p$, then let $a_i = \alpha_i \circ \sigma$
\begin{equation*}
    a(t) = \sum_{i=1}^n\sum_{p = 1}^m a_{i,p}(t)e_p(t) \otimes dt_i
\end{equation*}

First let $f\in \CE(A)^0 = \cinf(M)$, we have 
\begin{align*}
    <a^*\circ d_A(f), \del t_i>(t) &= <df(\sigma_t), \rho\circ a(\del t_i)>\\
    &= <df(\sigma_t), \rho\circ \alpha^t_i>\\
    &= <df, X^t_i>(\sigma(t))\\
    &= \frac{d}{dt_i}(f\circ \sigma(t))\\
    &= <d\circ a^*(f), \del t_t>(t).
\end{align*}
Hence $d_{\sigma} \circ a^* = a^*\circ d_A$ for degree 0.
Now let's verify the case for degree 1. Note that it suffices to check the dual basis $e_i^*$.
\begin{align*}
    <h^*d_A(e_l), \del t_i \wedge \del t_j> &= <h^*\big(\sum_{p,q = 1}^m c^l_{p,q} e^*_p \wedge e^*_q \big), \del t_i \wedge \del t_j >\\
    &= <\big(\sum_{p,q = 1}^m c^l_{p,q} e^*_p \wedge e^*_q \big), h(\del t_i) \wedge h(\del t_j) >\\
    &= \frac{1}{2} \sum_{p,q = 1}^m c^l_{i,j}(a_{i,p}a_{j,q} -a_{i,q}a_{j,p})\\
    &= \frac{d a_{i,l}}{dt_j} - \frac{da_{j,l}}{dt_i}
\end{align*}
where the last equality is due to $[\alpha_i, \alpha_j] =  \frac{d \alpha_i}{d t_j} - \frac{d \alpha_j}{d t_i}$. Note that
\begin{align*}
    <h^*(e_l^*), \del t_i> =& <e_l^*, \sum_{p=1}^m a_{i,p}e_p>\\
    =& a_{i,l}
\end{align*}
we have 
\begin{align*}
    \frac{d a_{i,l}}{dt_j} - \frac{da_{j,l}}{dt_i} = <d_{\sigma}\circ h^*(e_l^*), \del t_i\wedge \del t_j>
\end{align*}
Therefore, we have $d_{\sigma}\circ h^* = h^*\circ d_A$.
\end{proof}

Next, we want to show that all Lie algebroid morphisms can be obtained in this way.
\begin{lem}\label{extension}
    Let $\alpha^0_1, \cdots, \alpha^0_n$ be a family of $t$-time-dependent sections of $A$ for $t\in \Delta_n$, which satisfy 
    \begin{equation*}
        [\alpha_i^0, \alpha_j^0] =  \frac{d \alpha_i^0}{d t_j} - \frac{d \alpha_j^0}{d t_i}
    \end{equation*}
    for $i,j = 1,\cdots, n$. Suppose we have another family of sections $\alpha_{n+1}$ which depends on $t'\in \Delta_{n+1}$ and satisfying
    \begin{align*}
        [\alpha_k, \alpha_{n+1}] =&  \frac{d \alpha_k}{d t_{n+1}} - \frac{d \alpha_{n+1}}{d t_k}\\
        \alpha_k|_{t_{n+1} = 0} =&  \alpha_k^0
    \end{align*}
    then $[\alpha_i, \alpha_j] =  \frac{d \alpha_i}{d t_j} - \frac{d \alpha_j}{d t_i}$ are satisfied for all $i,j = 1 \cdots, n+1$.
\end{lem}
\begin{proof}
    We only need to verify the case for $i , j= 1, \cdots, n$. Let $\phi_{i, j} = [\alpha_i, \alpha_j] -  \frac{d \alpha_i}{d t_j} - \frac{d \alpha_j}{d t_i}$. Differentiate $\phi_{i,j}$ in $t_{n+1}$ we get
    \begin{align*}
        \frac{d\phi_{i,j}}{dt_{n+1}} =& [\alpha_{n+1}, \phi_{i,j}]
    \end{align*}

\end{proof}

\begin{prop}
Let $a: \sum_{i=1}^n a_i dt_i: T\Delta_n \to A$ be a Lie algebroid morphism, then there exists a family of time-dependent sections $\alpha^t_1, \cdots, \alpha^t_n$ such that
\begin{align*}
    [\alpha_i, \alpha_j] &= \frac{d \alpha_i}{d t_j} - \frac{d \alpha_j}{d t_i}\\
    a_k &= \alpha_k \circ \sigma
\end{align*}
\end{prop}
\begin{proof}
The case for $n = 1$ is obvious. We shall use the previous lemma and prove by induction. Suppose we have shown the case for $n = k$. Let $a: T \Delta_{n+1} \to A$ be a Lie algebroid morphism.    

First, extend $a_1|_{\{0 = t_n = \cdots = t_2 \le t_1 \le 1\}}$ to a $t_1$-time-dependent section $\alpha_1$. Next, we extend $a_2|_{\{0 = t_n = \cdots = t_3 \le t_2 \le t_1 \le 1\} }$ to a $(t_1, t_2)$-time-dependent section $\alpha_2$. Then we construct $\alpha_1$ as solution to 
    $$
    [\alpha_2, \alpha_1] = \frac{d \alpha_2}{d t_1} - \frac{d \alpha_1}{d t_2}
    $$
    with initial condition $\alpha_1|_{\{0 = t_n = \cdots = t_3 \le t_2 \le t_1 \le 1\}}$ constructed as in the previous step. Continuing this fashion, we extend $\alpha_i$ to a $(t_1,\cdots, t_i)$-time-dependent section of $A$, which satisfied the equation
    $$
    [\alpha_k, \alpha_i] = \frac{d \alpha_k}{d t_i} - \frac{d \alpha_i}{d t_k}
    $$
    for $1\le k <i$, with initial conditions $\alpha_k|_{\{0 = \cdots = t_{i} \le \cdots \le t_1 \le 1\}} $. $a_i = \alpha_i \circ \sigma$ is obvious by construction.
\end{proof}

\subsection{Homotopy and monodromy}

\begin{defn}
    Let $a_0, a_1: T\Delta^n \to A$ be two $n$-simplices, we say $a_0$ and $a_1$ are {\it homotopic} if there exists a Lie algebroid morphism $h = \sum_{k=1}^{n+1} h_k dt_k: T\Delta^{n}\times T\Delta^1 \to A$.

    \begin{enumerate}
        \item $a_{\epsilon}= \sum_{k=1}^{n} h_k(t_1, \cdots, t_n, \epsilon) dt_k$ for $\epsilon = 0, 1;$
        \item $h_{n+1}$ vanishes on the boundary of $\Delta^n$.
    \end{enumerate}
    It is easy to show that homotopies define an equivalence relation on the space of $n$-simplices.
\end{defn}

\begin{lem}
    A map $T\Delta^n \to A$ which vanishes on $\del \Delta^n$ is homotopic to a map vanishing on $T\del \Delta^n$.
\end{lem}

\begin{proof}
    First choose a cut-off function $\tau \in \cinf(\R)$ with $\tau(0)=0, \tau(1) = 1$, and $\tau'(t) > 0$ for $t \in (0,1)$.
    Define $h: \Delta^n \times \Delta^1 \to \Delta^n$ by
    $$
    h(t_1,\cdots,t_n, t_{n+1}) = \Big(\big(1-\tau(t_{n+1}) \big)t_1 + \tau(t_{n+1})t_1, \cdots, \big(1-\tau(t_{n+1}) \big)t_n + \tau(t_{n+1})t_n \Big)
    $$
    write $\tau(t_1,\cdots, t_n) = \big( \tau(t_1), \cdots, \tau(t_n)\big).$ Then $ a\circ dh$ gives a homotopy between $a$ and $a\circ d\tau$. Clearly $a\circ d\tau$ vanishes on $T \del\Delta^n$
\end{proof}

Given two simplices $a^i: T \Delta^n \to A$, we want to define the concatenation of them. The idea is to concatenate in the $t_n$-direction, but we have to be careful since the naive concatenation might not be smooth. First, in order to concatenate two simplices, we assume $d_1 a^1 = d^2 a^2$. We define the concatenation $a^1 \odot_{t_1} a^2$ by 
$$
a^1 \odot_{t_1} a^2 = \begin{cases}
(dr_1^{\tau}* a^0) \circ dp_0 \quad t_1 \in [0, 1/2]\\
(dr_1^{\tau}* a^1) \circ dp_1
\end{cases}
$$
where $dr_n^{\tau}: T\Delta^n \to T \Delta^n$ is the tangent map to $r_1^{\tau} (t_1, \cdots, t_n) \to (\tau(t_1), \cdots, t_n)$, and $p_i$ are maps which reparametrize the first coordinates
\begin{align*}
    p_0(t_1, \cdots, t_n) &= \Big(\tilde{p}_0^{-1}(t_1), t_2 \cdots, t_n\Big)\\
    p_1(t_1, \cdots, t_n) &= \Big(\tilde{p}_1^{-1}(t_1), t_2\cdots, t_n\Big)
\end{align*}
where 
\begin{align*}
    \tilde{p}_0(t_1) =&  \frac{(1-t_1)t_1}{1-t_2} + \frac{(t_1 - t_2)(1+ t_2)}{2(1-t_2)}\\
    \tilde{p}_1(t_1) =&  \frac{(t_1 - t_2) t_1}{1 - t_2} + \frac{(1-t_1)(1 + t_2)}{2 (1-t_2)}
\end{align*}
Based on the property of $\tau$ which smoothen the boundary of the simplices, it is easy to get
\begin{lem}
    The concatenation map $a^1\circ_{t_1} a^2: T\Delta^n \coprod T\Delta^n \to A
    $ of $n$-simplices is smooth.
\end{lem}
Next, let's define the homotopy groups of Lie algebroids in terms of simplices.
Let $A$ be a Lie algebroid. Its monodromy $\infty$-groupoid is defined as the simplicial manifold $\mon A$ with $n$-simplices
\begin{equation*}
    \mon^{\infty}(A)_n = \Hom_{\algd}(T\Delta^n, A)
\end{equation*}
We define the {\it isotropy $n$-simplices} at $x$ or {\it $A$-spheres} based at $x$ to be
$$
\mon^{\infty}_x(A)_n = \{ g\in \mon_n(A): \pi (g|_{\del \Delta}) = x\}
$$
We define {\it trivial isotropy $n$-simplices} to be those isotropy simplices whose base simplices are contractible
$$
\mon^{\infty}_x(A)_n^0 = \{ g\in \mon^{\infty}_x(A)_n: \pi (g) =\sim  *\}
$$
Note that all isotropy simplices at $x$ lie in a single leaf. According to our construction, we have
\begin{lem}
    $\mon^{\infty}_x(A)_n^0$ is the connected component of the identity of $\mon^{\infty}_x(A)_n$. We have a short exact sequence
    \begin{equation*}
        1 \to \mon^{\infty}_x(A)_n^0 \to \mon^{\infty}_x(A)_n \to \pi_n(L_x) \to 1
    \end{equation*}
    where $L_x$ denote the leaf of the foliation that $x$ sits in.
\end{lem}
Recall that, when we restrict to a leaf, we have the following exact sequence of Lie algebroids
\begin{equation*}
    0 \to \g_{L_x} \to A_{L_x} \stackrel{\rho}{\to} T L_x \to 0
\end{equation*}
where $\g_{L_x}$ is the bundle of isotropy Lie algebra $\g_x$.
\begin{prop}
There exists a long exact sequence
\begin{equation}                                   
    \cdots \pi_{i+1}(L_x) \stackrel{\del}{\to} \mon^{\infty}_x(\g_{L_x})_i \to \mon^{\infty}_x(A_{L_x})_i\to \pi_i(L_x) \to \cdots
\end{equation}

\begin{proof}
    First, we need to construct the boundary map $\del: \pi_{i+1}(L_x) \to \mon^{\infty}_x(\g_{L_x})_i$. Let $[\sigma] \in \pi_{i+1}(L_x)$, i.e. $\sigma: \Delta^{n+1} \to L_x$ such that $\sigma(\del \Delta^n) = x$. Let $\sum_{i=1}^{n+1} a_i dt_i: T\Delta^n \to A_{L_x}$ be any Lie algebroid morphism which lifts $d\sigma: T\Delta^n \to TL_x$ such that $a_i|_{\Delta^n} = 0$ for $1\le i \le n$ and $a_{n+1}|_{\del \Delta^n} = 0$, where $\Delta^n = \{1 \ge t_1\ge \cdots \ge t_n \ge 0\}= d_{n+1} \Delta^n$. Let $\Lambda^{n+1}_{n+1}$ be the $(n+1)$-th horn of $\Delta^{n+1}$ as usual, then after simple reparametrization, $a(T\Lambda^{n+1}_{n+1})$ gives a map $\Delta_n \to \g_{L_x}$ since $\sigma$ is constant on the boundary. Therefore, we define $\del(\sigma) = [a(T\Lambda^{n+1}_{n+1})] \in \mon^{\infty}_x(\g_{L_x})_i$. Now it suffices to show, for two homotopic simplices $\sigma_1, \sigma_2$, their image under $\del$ are homotopic as $\g_{L_x}$-paths.
    
    Let $a^i: T\Delta^n \to A_{L_x}$ be some lifts of $\sigma^i$, we want to show $\del(a^i) = a^i(T\Lambda^{n+1}_{n+1})$ are homotopic. In order to do this, we will construct an explicit homotopy $h: T\Delta^n \times T\Delta^1 \to A_{L_x}$. By assumption, there exists a homotopy $\sigma^s(t) = \sigma(t,s):\Delta^{n+1} \times \Delta^1 \to L_x, s\in \Delta^1$ of $\sigma^1$ and $\sigma^2$. Choose a homotopy $a^s_{n+1}(t)=a_{n+1}(t,s): T\Delta^n \times T\Delta^1 \to A_{L_x}$ such that, 
    $$
    \rho(a_{n+1}(t,s)) = \frac{d \sigma^s(t) }{d t_{n+1}}, \quad a^{n+1}(t,s)|_{\del \Delta^n} = 0 
    $$
    Let $\alpha_{n+1}(t,s, \sigma(t,s))$ be the corresponding time-dependent sections of $A_{L_x}$. Consider the solutions of the following system
    \begin{equation*}
        \begin{cases}
        \frac{d \alpha_i}{dt_{n+1}} - \frac{d \alpha_{n+1}}{dt_{i}} = [\alpha_i, \alpha_{n+1}] \quad 1\le i \le n\\
        \frac{d \beta}{dt_{n+1}} - \frac{d \alpha_{n+1}}{ds} = [\beta, \alpha_{n+1}]\\
        \alpha_i|_{\del \Delta^{n+1}} = 0\\
        \beta|_{\del \Delta^{n+1}} = 0
        
        \end{cases}
    \end{equation*}
    Note that $\alpha^i(t,s, \sigma(t,s)) = a^i(t,s)$. By similar arguments as in Lemma \ref{extension}, we have 
    \begin{equation*}
        \frac{d \alpha_i}{ds} - \frac{d \beta}{dt_{i}} = [\alpha_i, \beta] \quad 1\le i \le n
    \end{equation*}
    Since $\beta|_{\del \Delta^n} = 0$, it gives a homotopy between $a^1$ and $a^2$.
    
    Next, we want to show the sequence is exact. It suffices to verify the exactness at $\mon^{\infty}_x(\g_{L_x})_i$ for all $i \ge 1$. By construction, the image of $\del$ consists of 
\end{proof}

\end{prop}
Now, if we glue everything in the above exact sequence by the leave $L_x$, and regard all objects as bundles of groups over $M$, we get
\begin{cor}
There exists a long exact sequence of bundle of groups
\begin{equation}                                   
    \cdots \pi_{i+1}(L_x) \stackrel{\del}{\to} \mon^{\infty}(\g_{L_x})_i \to \mon^{\infty}(A)_i\to \pi_i(L_x) \to \cdots
\end{equation}
\end{cor}
Now we have a relation between isotropy $n$-simplices, $A$-simplices, and simplices along the foliation, where each of them corresponds to $n$-simplices of $\infty$-groupoids. Once we add back the simplicial structures, we get

\begin{prop}
There exists a fiber sequence of Lie $\infty$-groupoids
\begin{equation}
\coprod_{L_x}\mon^{\infty}(\g_{L_x}) \to \mon^{\infty}(A)\to \Pi_{\infty}(\CF) 
\end{equation}

\end{prop}

 We call the boundary map $\del$ to be {\it monodromy morphism}\index{monodromy morphism}. As we know for $i=2$, it corresponds to the classical monodromy morphism and its image in $\mon^{\infty}(\g_{L_x})$ controls the integrability of Lie algebroids.

Recall $\mon^{\infty}(A)_{n}$ consists of all $C^1$ $n$-simplices $\tilde{\sigma}: \Delta^n \to A$ which sits above its projection $\sigma= \pi(\tilde{\sigma})$  in $M$, and satisfies
\begin{equation*}
    \rho(a(t)) = \frac{d(\pi\circ a(t))}{dt}
\end{equation*}

In order to help us study the smooth structures on $\mon^{\infty}(A)_n$, let's first look at the larger space $\tilde{P}_n(A)$ which consists of $C^1$ $n$-simplices $\tilde{\sigma}: \Delta^n \to A$ over some base $C^2$-simplices in $M$. It's easy to see that $\tilde{P}_n(A)$ is a Banach manifold, and $\mon^{\infty}(A)_{n}$ is a submanifold of $\tilde{P}_n(A)$. The tangent space of a simplex $\tilde{\sigma}(t)$ consists of all $C^0$ $A$-sections over $\sigma(t)$. Using a connection $\nabla$ on $A$, we can view an element in $T_{\tilde{\sigma}}\tilde{P}_n(A)$ as a pair $(u, \phi)$, where $u: \Delta^n \to A$ and $\phi: \Delta^n \to TM$ are both simplices over the base simplex $\sigma$.

\begin{lem}
$\delta \in T_a \mon^{\infty}(A)_n$ have decomposition $\tilde{\delta}(t) = \big(u(t), \phi(t) \big)$ such that 
\begin{equation*}
    \rho(u) = \overline{\nabla}_{a(t)}\phi(t)
\end{equation*}
for $t = (t_1, \cdots, t_n) \in \{1 \ge t_1 \ge \cdots \ge  t_n \ge 0\}$.

\end{lem}
\begin{proof}
Consider a smooth map $F: \tilde{P}_n(A) \to \tilde{P}_n(TM)$ defined by $F(\tilde{\sigma}(t)) = \rho(\sigma(t)) - D  \sigma(t)$. Here we use $D$ denote the gradient $D = \sum_i \frac{d}{dt_i}$ in order to distinguish it from the connection. Let $\tilde{P}_n^0(TM)$ denotes the submanifold of paths with zero value in the fiber in $\tilde{P}_n(TM)$. Then it suffices to show $F$ is a submersion onto $\tilde{P}_n^0(TM)$, and $F^{-1}\tilde{P}_n^0(TM)= \mon^{\infty}(A)_n$. 
Let's restrict to the differential of $F$ onto $\tilde{P}_n^0(TM)$
\begin{equation*}
    dF: T_{\tilde{\sigma}}\tilde{P}_n(A) \to T_{0_{\sigma}}\tilde{P}(TM)
\end{equation*}
Here $0_{\sigma}$ denotes the canonical lift of $\sigma :\Delta^n \to M$ to $\Delta \to TM$ with zero values in the fiber. 
Note that $T_{0_{\sigma}}\tilde{P}(TM)$ consists of all smooth sections of $TTM$ over $\sigma$. The image of $x \in \Delta^n$ of $dF$ is 
    $$T_{0_x}\tilde{P}(TM) \simeq \bigoplus_{i=1}^{n} T_{0_x} TM
    \simeq  \bigoplus_{i=1}^{n} T_xM\oplus T_xM
    $$
    by the canonical splitting of $T_{0_x} T_M$. We claim that, for any connection $\nabla$ which splits $\tilde{\sigma}$ as $(u, \phi)$, the vertical and horizontal components of the splitting $T_{0_x} TM
    \simeq  T_xM\oplus T_xM$ are $\rho(u) - \overline{\nabla}_a \phi$ and $\phi$ respectively.
    Let $m = \dim M, k = \dim A$. Let $x = \{x_1, \cdots, x_m\}$ be a local chart of $M$, and $\{ \frac{\del}{\del x_i} \}$ be a local basis of $TM$, then denote the horizontal and vertical basis of $T_{0_x}TM$ by  $\{ \frac{\del}{\del x_i} \}$ and $\{ \frac{\delta}{\delta x_i} \}$ respectively. Without loss of generality, we assume $\nabla$ is the standard flat connection.
    Now let $\tilde{\sigma}(t) = \sum_{i = 1}^n  \tilde{\sigma}_i(t) e_i$

\end{proof}

Next, we want to show that the homotopy of $n$-simplices actually induces a (infinite dimensional) foliation on $\mon^{\infty}(A)_n$
\begin{prop}
 There exists a foliation $\CF_n$ on $\CG_n$ with finite codimension.
\end{prop}

First, let $\sigma$ be an $n$-simplex in $M$, consider a subspace $\tilde{P}_n^{\sigma}(A)$ of $\tilde{P}_n(A)$ defined by
$$
\tilde{P}_n^{\sigma}(A) = \{\gamma \in \tilde{P}_n(A): \gamma(0) = 0, \gamma(t) \in A_{\sigma(t)} \}
$$
that is, $\tilde{P}_n^{\sigma}(A)$ consists of sections of $A$ over $\sigma$ with initial condition $\gamma(0) = 0$. Let $\nabla$ be a connection on $A$, and $\tilde{\sigma}$ be an $A$-simplex over $\sigma$.

\begin{prop}
The $n$-truncation of $\CG_{\bullet}$ is a Lie $n$-groupoid.
\end{prop}

\subsection{Integrating $L_{\infty}$-algebroids}
Let $\g$ be an $L_{\infty}$ algebroids over a smooth manifold $M$ which is positively graded. Define a simplicial manifold $\mon_{\bullet}$ whose $n$-simplices are 
\begin{equation*}
    \mon_n =\Hom_{\dgcalg}\big(\CE(\g), \Omega(\Delta_n)\big)= \Hom_{\lalgdi}(T\Delta_n), \g)
\end{equation*}
Note that when $\g$ is a Lie algebroid, then $\mon_{\bullet}$ coincide with the construction in previous section.

\begin{prop}\cite{SS19}
$\mon_{\bullet}$ is a Lie $\infty$-groupoid where each $\mon_n$ is a Frech{\'e}t manifold.
\end{prop}
In this section, we shall prove an enhancement to the above result, which gives an $n$-truncation of integration of $L_{\infty}$-algebroids.

\begin{prop}
Let $\g$ be an $L_{\infty}$ algebroid over a smooth manifold $M$, with the underlying dg Module being perfect and concentrated in degree $[-n, 0]$. Then $\g$ integrates to a Lie $n$-algebroid which is an $n$-truncation of $\CG_{\bullet}$.
\end{prop}

\subsection{Local holonomy $\infty$-groupoid}

In this section, we will study the higher holonomy defined by monodromy morphisms. First, we will study the local structures. Let $L$ be a singular leaf of $\CF$. 

Recall that a fibration $P: E_{\bullet} \to F_{\bullet}$ in the semi-model category $\lalgd_{A^{\bullet}}$ is a degree-wise surjection. In particular, if $F_{\bullet}$ is a Lie algebroid, the $P$ degenerates to a surjection $E_{-1} \to F_{-1}$. Now consider in a fibration in $\lalgd$, we define a fibration $P: E_{\bullet} \to F_{\bullet}$ to be a commutative diagram 

\begin{center}
	\begin{tikzcd}
	E_{\bullet} \arrow[r, "P"] \arrow[d, ""] & F_{\bullet}
	\arrow[d, ""] \\
	M \arrow[r, "p"]    &  N 
	\end{tikzcd}
\end{center}
such that $P$ is a degreewise $L_{\infty}$ surjection and $p$ is a surjective submersion. 

\begin{defn}
    Let $P: E_{\bullet} \to F_{\bullet}$ be a $L_{\infty}$-algebroid fibration. An {\it Ehresmann connection}\index{Ehresmann connection! for $L_{\infty}$-algebroid fibration} for $P$ is a graded vector sub-bundle $H_{\bullet} \subset E_{\bullet}$ such that $H_{\bullet} \oplus \ker(P) = E_{\bullet}$. 
\end{defn}
Given an Ehresmann connection, we can lift a section of $F_{\bullet}$ to a unique section of $E_{\bullet}$, which is called a {\it horizontal lift}. Moreover,  $\rho_{E_{\bullet}}(\sigma(a))$ is $p$-related to $\rho_{F_{\bullet}}(a)$.

\begin{example}
    Let $E_{\bullet}$ and $F_{\bullet}$ be the tangent Lie algebroids $TM$ and $TN$ respectively. Then we recover the usual definition of manifold fibrations (surjective submersion). On the other hand, let $E_{\bullet}$ and $F_{\bullet}$ be ordinary Lie algebras, we recover Lie algebra epimorphisms.
\end{example}

\begin{defn}
We say an Ehresmann connection is complete if the for any complete vector field $\rho(\alpha)$, $\rho(\sigma(\alpha))$ is complete, where $\sigma: \Gamma(F_{\bullet}) \to \Gamma(E_{\bullet})$ is a lift induced by the connection. 
\end{defn}

Let's look at the fiber of an $L_{\infty}$-fibration. By definition, we have a graded vector bundle $K_{\bullet} = \ker(P) \subset E_{\bullet}$ over $M$. We can then restrict the $k$-ary brackets on $E_{\bullet}$ to $K_{\bullet}$, i.e.

\begin{equation*}
    l_k^{K_{\bullet}}(e_1,\cdots,e_k) = \pi (l_k^{E_{\bullet}}(e_1,\cdots,e_k))
\end{equation*}
for $e_i \in \Gamma(K_{\bullet})$, and $\pi: E_{\bullet} \to K_{\bullet}$ is the projection map. 
\begin{prop}
Let   $P: E_{\bullet} \to F_{\bullet}$ be an $L_{\infty}$-algebroid fibration over $p: M\to N$, then $K_{\bullet} = \ker(P)|_{p^{-1}(x)}$ inherits an  $L_{\infty}$-algebroid structure over an $L_{\infty}$-algebroid fibration for any $x \in N$.
\end{prop}
\begin{proof}
Fix $x\in N$. Clearly e$K_{\bullet}$ is a graded vector bundle over $p^{-1}(x)$, hence it suffices to show that the brackets $l_k^{K_{\bullet}}$ is well-defined and satisfies the homotopy Jacobi identities. By an analogue of Frobenius theorem for $L_{\infty}$-algebroids, it suffices to show that 
$\ann (\ker P) \simeq \im: \CE(P)$, where $\CE(P): \CO(E_{\bullet}) \to \CO(F_{\bullet})$ is the induced map on Chevalley-Eilenberg algebras, is $d^{\CE}$-closed. This follows from the fact that
$$
d^{\CE(E_{\bullet})}\circ f (\CO(F_{\bullet})) = f \circ d^{\CE(F_{\bullet})}(\CO(F_{\bullet}))
$$
and the homotopy Jacobi identity follows directly from $E_{\bullet}$.
\end{proof}
Next, we will show that $\ker(P)|_{p^{-1}(l)}$ can be patched together when we have a complete Ehresmann connection. We will need the following lemma.
\begin{lem}[\cite{LR19}]
    Let $\g$ be an $L_{\infty}$-algebroid over $A = \cinf(M)$. 
    Let $X$ be a degree zero vector field on $\CO(\g)$, i.e. a degree zero element in the tangent complex $T_{\g}$, then
    \begin{enumerate}
        \item For all fixed $t\in \R$, $X$ admits a time-$t$ flow $\Phi^X_t: \CO(\g) \to \CO(\g)$ if and only if the induced vector field $\tilde{X}$ on $M$ admits a time-$t$ flow.
        \item Assume $X$ is $d_{\CO(\g)}$-closed, i.e. $[d_{\CO(\g}, X] = 0$. Then the flow $\Phi^X_t: \CO(\g) \to \CO(\g)$ is an $L_{\infty}$-morphism for any admissible $t$.
        \item Assume $X$ is exact, i.e. there exists a $Y$ such that $[d_{\CO(\g)}, Y] = X$, then there exists an $\li$-morphism $\Phi^Y: \CO(\g) \otimes \R \to \CO(\g)$ defined in a small neighborhood of $\CO(\g) \otimes \{0\}$, such that the restriction $\Phi^X_t: \CO(\g) \to \CO(\g)$ if the flow of $[d_{\CO(\g)}, Y]$ at time $t$ for all admissible time $t$. Also, we have that all $\Phi^X_t$'s are homotopic $\li$-morphisms.
    \end{enumerate}
\end{lem}
\begin{proof}
    See \cite{LR19} Lemma 1.6.
\end{proof}
Given an $\li$-algebroids fibration $P: E_{\bullet} \to F_{\bt}$, we can regard a section of $E_{0}$ as a degree -1 vector field on $\CO(E_{\bt})$ by left contraction (note the degree shift here in $\CO(E_{\bt})$). Let $X\in \Gamma(F_0)$, then the degree zero vector fields $[\iota_{\sigma(X)}, d_{\CO(E_{\bt})}]$ and        $[\iota_{X}, d_{\CO(F_{\bt})
}]$   are $P$-related for every Ehresmann connection $K_{\bt}$.
\begin{lem}
    Let $K_{\bt}$ be an Ehresmann connection for an $\li$-algebroid fibration $P: E_{\bullet} \to F_{\bt}$, then $K_{\bt}$ is complete if and only if for any $X\in \Gamma(F_0)$, the time-$t$ flow of $[\iota_{\sigma(X)}, d_{\CO(E_{\bt})}]$ is defined if and only if the time $t$ flow of 
    $[\iota_{X}, d_{\CO(F_{\bt})}]$ is defined.
\end{lem}
\begin{proof}
By (1) in the previous lemma, the flow of a degree 0 vector field exists if and only if its induced vector fields on the base manifold exists. Note that the induced vector fields of $[\iota_{\sigma(X)}, d_{\CO(E_{\bt})}]$ and $[\iota_{X}, d_{\CO(F_{\bt})}]$ are $\rho_{E_{\bt}}(\sigma(X))$ and $\rho_{F_{\bt}}(X)$ respectively. Therefore, the result follows directly from the definition.
\end{proof}
A complete Ehresmann connection allows us to identify different fibers.

\begin{lem}
    Let $K_{\bt}$ be an Ehresmann connection for an $\li$-algebroid fibration $P: E_{\bullet} \to F_{\bt}$. Suppose the anchor map of $F_{\bt}$ is surjective, then the fibers $\mathcal{T}_x$ and $\mathcal{T}_y$ for $x, y  \in N$ are isomorphic as $L_{\infty}$-algebroids.
\end{lem}
\begin{proof}
     By assumption, there exists a vector field $Z$ on $N$ whose time 1 flow globally, and maps $x$ to $y$. By surjectivity of $\rho_{F_{\bt}}$, we can lift $Z$ to a section $X$ of $F_0$. From previous lemma, we know $[\iota_{\sigma(X)}, d_{\CO(E_{\bt})}]$ is $p$-related to $Z$. Since the time-1 flow of $X$ is well-defined, we know the time-1 flow $\Phi_1$ of $[\iota_{\sigma(X)}, d_{\CO(E_{\bt})}]$ is well-defined and is an $\li$-isomorphism. Note that $Phi_1$ induces a diffeomorphism $\phi_1: M\to M$ which is over the time-1 flow of $X$, hence it maps the fiber $p^{-1}(x)$ to $p^{-1}(y)$. There, we see $\Phi_1$ restricts to an isomorphism from $\mathcal{T}_x$ to $\mathcal{T}_y$.
\end{proof}
    For general $\rho_{F_{\bt}}$, we have the following local result.
\begin{lem}
    Let $K_{\bt}$ be an Ehresmann connection for an $\li$-algebroid fibration $P: E_{\bullet} \to F_{\bt}$. Let $x,y$ lie in a single leaf in the singular foliation associated to $F_{\bt}$ on $N$, then the fibers $\mathcal{T}_x$ and $\mathcal{T}_y$ for $x, y  \in N$ are isomorphic as $L_{\infty}$-algebroids.
\end{lem}

\begin{proof}
    It suffices to consider the anchor map of $F_{\bt}$ is not surjective. Then result follows from the previous lemma by replacing $F_{\bullet}$ with $F_{\bt}|_L$, where $L$ is the leaf containing both $x$ and $y$.
\end{proof}

\begin{cor}
    $\ker(P)_{p^{-1}(x)}$ glues to an $L_{\infty}$-algebroid $\mathcal{T}$ over $M$.
\end{cor}
\begin{proof}
    Classical partition of unity type argument.
\end{proof}

By construction $dp \circ\rho_{K_{\bullet}}(\Gamma(K_0)) =\rho_{F_{\bullet}} \circ P(\Gamma(K_0)) = 0 \subset \Gamma(TN)$, so $K_{\bullet}$ restricts to the fiber of $p$, i.e. given $x\in N$, we have an $L_{\infty}$-algebroid $K_{\bullet}^l$ over $p^{-1}(x)$.

\begin{example}
    In \cite{BZ11}, Lie algebroid fibrations corresponding to $L_{\infty}$-fibration in our sense for Lie algebroids as $L_{\infty}$-algebroids with a complete Ehresmann connection.
\end{example}
\begin{prop}
Let $P: E_{\bullet} \to F_{\bullet}$ be an $L_{\infty}$-algebroid fibration over $p:M \to N$. Let $T_x$ denote the fiber of $P$ over $x\in N$. Suppose $P$ admits a complete Ehresmann connection. There exists a long exact sequence
\begin{equation}
    \cdots \mon^{\infty}_{i+1}(F_{\bullet})_x \stackrel{\del}{\to} \mon^{\infty}_i(T_x)_y \to \mon^{\infty}_i(E_{\bullet})_y\to \mon^{\infty}_i(F_{\bullet})_x \to \cdots
\end{equation}
\end{prop}

Note that the boundary map is exactly the monodromy homomorphism $$\del: \mon^{\infty}_{i+1}(F_{\bullet})_x \stackrel{}{\to} \Gamma(\mon^{\infty}_i(T_x))$$.

Now consider $E_{\bullet}$ to be the universal $L_{\infty}$-algebroid associated to $\CF$, and let $L$ be a locally closed leaf.

\begin{defn}
    Let $L$ be a locally closed singular leaf of $\CF$. An {\it Ehresmann $\CF$-connection}\index{Ehresmann connection!for singular foliations} consists of a triple $($ an Ehresmann connection $H$ of a projection $p:M_L \to L$ of a neighborhood $M_L \subset M$ of $L$. We say an Ehresmann $\CF$-connection $(M_L, p, H)$ is complete near $L$ is $H$ is complete.

\end{defn}
Recall that we have the following local splitting property of singular foliations

\begin{thm}[local splitting]
    Let $(M, \CF)$ be a singular foliation. Let $x \in M$ be arbitrary, $k = \dim(F_x)$, and $\hat{S}$ a slice at $x$, i.e. an embedded submanifold of $M$ such that $T_x \hat{S} \oplus F_x = T_x M$. Then there exists an open neighborhood $U$ of $x$ in $M$ and a foliated diffeomorphism $(U, \CF|_U) \simeq (I^k, TI^k) \times (S, \CF_S$, where $S = \hat{S} \cap U$, $\CF|_U$ is the restriction of $\CF$ to $U$, $I = (-1, 1)$, and $\CF_S = \CF|_U \cap \Gamma(TS)$.
\end{thm}

We have an analog for the Ehresmann $\CF$-connection
\begin{prop}[local splitting by Ehresmann $\CF$-connection]
Let $(M_L, p, H)$ be a complete Ehresmann $\CF$-connection for a locally closed leaf $L$. For every $x\in L$, there exists a neighborhood $U \subset L$ and a foliated diffeomorphism $p^{-1}(U) \simeq U \times p^{-1}(L)$ which intertwining $\CF|_{p^{-1}(U)}$ and the product foliation $\Gamma(TU) \times \CF_{p^{-1}(l)}$.
\end{prop}
\begin{proof}
Without loss of generality, $L\simeq I^k$ where $k = \dim L$. It suffices to show that given a complete Ehresmann $\CF$-connection $(M_L, p, H)$ for $L$, then there actually exists a flat complete Ehresmann $\CF$-connection. 

We shall proceed by induction. The case for $k = 1$ is trivial since any dimension 1 distribution is integrable. Suppose we have proved the result for some $k \in \N$. Let consider a complete Ehresmann $\CF$-connection over $L = I^{k+1}$. Let $(t_1, \cdots, t_k, t_{k+1})$ be a coordinate for $I^{k+1}$. Since $H$ is complete, the horizontal lift of $\frac{\del}{\del t_{k+1}}$ is complete, and its flow $\Psi_t: p^{-1}(I^k \times \{s\}) \to p^{-1}(I^k \times \{s + t\})$ preserves $\CF$, where $t, s, s+t \in I$. Hence $(\Psi_t)_*: \CF|_{p^{-1}(I^k \times \{0\})} \simeq \CF|_{p^{-1}(I^k \times \{t\})}$. By induction hypothesis, the projection $p^{-1}(I^k \times \{0\}) \to (I^k \times \{0\})$ admits a flat complete Ehresmann $\CF$-connection $(M_L, p, H')$. Now we can use $\Psi_t$ to transport $H'$ to get a flat complete Ehresmann $\CF$-connection for all $p^{-1}(I^k \times \{t\}) \to (I^k \times \{t\})$. Therefore, we have a new distribution
\begin{equation*}
    H'' = <H(\frac{\del}{\del t_{k + 1}})> \oplus H'
\end{equation*}
Then it's easy to verify that $H''$ is actually the flat complete Ehresmann $\CF$-connection we need.
\end{proof}

Consider the tangent Lie algebroid $TL$ of a locally closed singular leaf $L$. 
\begin{lem} Let $E_{\bullet}$ be the $L_{\infty}$-algebroid resolving $\CF$. Suppose $L$ admits an Ehresmann $\CF$ connection $(M_L, p, H)$, then the induced map $E_{\bullet}\to TL$ is an $L_{\infty}$-algebroid fibration.
\end{lem}
\begin{proof}
 $P: E_{\bullet} \to TL$ is the composition of the anchor map $\rho: E_{\bullet} \to TM$ and the surjection $dp: TM\to TL$, which is clearly an $L_{\infty}$-morphism. 
\end{proof}

\begin{lem}
    Suppose $L$ admits an Ehresmann $\CF$-connection $(M_L, p, H)$, then there exists an Ehresmann connection for $P:E_{\bullet} \to TL$ such that the only non-trivial term is $H_{0} \subset E_0$ and $\rho(H_{0}) = H$. Moreover, $H_{0}$ is complete if and only if $H$ is complete.
\end{lem}
\begin{proof}
    We will do the construction locally and then the result follows from the standard partition of unity argument. Let $x\in L$, by assumption there exist a neighborhood $U_x$ and $k$ vector fields $X_1, \cdots, X_k$ generate $H$, where $k = \dim L$. Let $e_i \in \Gamma(E_0)$ be a lift of $X_i$ for $1\le i \le k$, i.e. $dp \circ \rho(e_i) = X_i$. Then $\{\e_i\}$ generate the desired $H$ on $U_x$. 
\end{proof}

Suppose $L$ admits a complete Ehresmann  $\CF$-connection $(M_L, p, H)$, and denote the fiber over $x$ to be $\mathcal{T}_x$, then we have an exact sequence
$$
\mathcal{T}_x \to E_{\bullet} \to TL
$$
Applying the previous result, we have 
\begin{prop}
Let $L$ be a locally closed singular leaf of $\CF$ which admits a complete Ehresmann $\CF$-connection $(M_L, p, H)$, then $P: E_{\bullet} \to TL$ is an $L_{\infty}$-algebroid fibration over $p:M_L \to L$. Let $\mathcal{T}_x$ denote the fiber of $P$ over $x\in L$, which corresponds to $L_{\infty}$-algebroid of the transversal foliation at $x$. There exists a long exact sequence
\begin{equation}
    \cdots \pi_{i+1}(L,x) \stackrel{\del}{\to} \mon^{\infty}_i(\mathcal{T}_x)_y \to \mon^{\infty}_i(E_{\bullet})_y\to \pi_i(L,x) \to \cdots
\end{equation}
\end{prop}

Note that the boundary map is exactly the monodromy homomorphism $$\del: \pi_{i+1}(L,x) \stackrel{}{\to} \Gamma(\mon^{\infty}_i(\mathcal{T}_x)_y)$$.

Recall that, $\mon_0(E_{\bullet}) = M/E_{\bullet}$ which corresponds to the leaf space of $\CF$.
When $i = 0$, $\mon^{\infty}_0(\mathcal{T}_x) = p^{-1}(x)/\CF|_{\mathcal{T}_x}$, and we have an identification $\Gamma(\mon^{\infty}_0(\mathcal{T}_x) = \Diff(p^{-1}(x)/\CF|_{\mathcal{T}_x})$ which is the bijections of the leaf space induced by diffeomorphisms. Hence, the image of $\del_1$ is the holonomy group $\hol(\CF)_x$.

\begin{defn}
    Define the {\it $n$-th holonomy}\index{holonomy!$n$-th holonomy of singular foliations} of $L$ to be the image of the $n$-th monodromy morphism $\hol_n(\CF, L) = \del(\mon_{n+1}(TL))$.
\end{defn}

\begin{prop}
There is a natural simplicial structure on $\hol_n(\CF, L)$, which assembles to a Lie $\infty$-groupoid $\hol_{\bullet}(\CF, L)$. We call $\hol_{\bullet}(\CF, L)$ the holonomy $\infty$-groupoid of $\CF$ at $L$.
\end{prop}

\begin{example}[embedded submanifold]
    Let's consider $L$ to be a simply connected embedded submanifold of $M$. Consider $\CF$ to be a singular foliation generated by all the vector field tangent to $L$. Let $U$ be a tubular neighborhood of $L$ in $M$, and $NL$ the normal bundle of $NL$. Let $f: \CF|_U \to NL$ be a foliated diffeomorphism which send $L$ to the zero section of $L$. Then the Atiyah Lie algebroid $\At(NL)$ of $NL$ is a Lie algebroid of minimal rank of $\CF$. Recall the $\At(NL)$ consists of covariant differential operators on $\Gamma(NL)$. 
    
    Now let's look at the long exact sequence of holonomy. Take some $x \in L$ and $y \in p^{-1}(x)$. There are two cases:
    \begin{enumerate}
        \item First consider $y \not = 0$. The transverse foliation consists of a fiber $V \simeq \R^q$ of $NL$, where $q$ is the codimension of $L$, with a regular leaf $V -\{0\}$ and a singular leaf $\{0\}$. Hence, we have \begin{equation*}
        \mon_{n}(\mathcal{T}_x, y) = \pi_n(V-\{0\}, y) \simeq \pi_n (S^{q -1}, y)
    \end{equation*} 
    Therefore, the $i$-th monodromy morphism reduces to
    \begin{equation*}
        \del: \pi_i(L,x) \to \Gamma(\pi_n (S^{q -1}))
    \end{equation*}
    which corresponds to the exact sequence of the fibration
    \begin{equation*}
        (V-\{0\}) \to (NL - L)\to L.
    \end{equation*}
    \item Next, let's look at $y = 0$. Note that $\At(NL)$ restricts to $\gl(V)$ on $V$, hence by classical Lie integration theory, we know that 
    \begin{equation*}
        \mon_n(\mathcal{T}_x, y) = \begin{cases}
        \widetilde{\GL(V)} \quad &n = 1\\
        \pi_n(\widetilde{\GL(V)}) & n > 1
        \end{cases}
    \end{equation*}
    \end{enumerate}

\end{example}

\section{Higher foliations}
\subsection{Tangent $\infty$-stack}
Assigning a manifold $M$ its tangent bundle $TM$ gives a functor $T: \mfd \to \mfd$. Associated to $T$ there is a natural projection $\pi : T \to \id$ given by $TM \to M$. Precompose $T$ with the Yoneda embedding $\y:\mfd^{\op} \to \pshi(\mfd)$ gives a $\infty$-functor $T^*: \pshi(\mfd) \to \pshi(\mfd)$, i.e.

\begin{center}
	\begin{tikzcd}
	X_{\bullet} \arrow[r, "d"] \arrow[d, "d"] & X_{\bullet}
	\arrow[d, "d"] \\
	Y_{\bullet} \arrow[r, "d"] \arrow[r, "d"]      &  Y_{\bullet}                 
	\end{tikzcd}
\end{center}

\begin{lem} There exists an $\infty$-functor $T^*: \pshi(\mfd) \to \pshi(\mfd)$ canonically associated to $T$.
	\end{lem}

\begin{lem}
	$T^*$ restricts to an $\infty$-functor $T^*: \shi(\mfd) \to  \shi(\mfd)$.
\end{lem}
\begin{proof}
	If follows from $T$ preserves open covers and pullbacks of covers are covers. 
	\end{proof}

Recall that two $\infty$-functors $F:\CC \to \CD$ and $G:\CD \to \CC$ if there exist a unit $\infty$-transformation $\epsilon: \id_{\CD}\to F\circ G$ such that the composition
$$
\Hom_{\CC}(F(x), y) \stackrel{\hom_{\CC}(G, G)}{\longrightarrow }\Hom_{\CD}(G\circ F(x), G(y)) \stackrel{\hom_{\CC}(G\epsilon, \id)}{\longrightarrow} \Hom_{\CD}(x, G(y))
$$
is an equivalence of Kan complex.

\begin{prop}
	$T^*:\shi(\mfd) \to \shi(\mfd) $ admits a left $\infty$-adjoint $T: \shi(\mfd) \to \shi(\mfd)$. 
\end{prop}	

\begin{proof}The $\infty$-category of Kan complexes $\igpd$ is homotopically complete, hence we can form the left $\infty$-adjoint $T^{pre}:\pshi(\mfd) \to \pshi(\mfd)$ on $\infty$-presheaves. The inclusion $i:\shi(\mfd) \to \pshi(\mfd)$ admits the left adjoint of the stackification functor $p:\shi(\mfd) \to \shi(\mfd)$. Now we define $T = p\circ T^{pre} \circ i$.
	\end{proof}
We call $T$ the {\it tangent $\infty$-stack functor}. We also want to extend $\pi: T \to \id_{\mfd}$ to $\infty$-stacks. Precomposing Yoneda embedding with $\pi$ determines a natural transformation $\pi^*: \id_{\pshi(\mfd)} \to T^*$. Take the left adjoint of $\pi^*$ to be $\pi^{pre}$. Now define $\pi: T\to \id_{\shi(\mfd)}$ to be $$\Hom_{[\shi(\mfd),\shi(\mfd)]}(p,p)\circ \pi^{pre} \circ \Hom_{[\pshi(\mfd),\pshi(\mfd)]}(i,i)
$$.

\begin{prop}
	There exist an $\infty$-natural equivalence $\epsilon: T\circ \y 	\Rightarrow \y \circ T$. 
\end{prop}
\begin{proof}Follows from the construction.
	
	\end{proof}
We have a fully faithfully embedding $|-|:\ilgpd \to \shi(\mfd)$ from Lie $\infty$-groupoids to $\infty$-stacks. We can form the tangent $\infty$-groupoid functor $T^{gpd}: \ilgpd \to \ilgpd$ by taking degreewise tangent bundle along with differentials.
\begin{defn}
	Let $X_{\bullet}$ be a Lie $\infty$-groupoid. Define the tangent groupoid functor $T^{gpd}$ to be the unique functor which sends $X_{\bullet}$ to $TX_{\bullet}$, where $TX_{i}= T(X_i)$ and all structure maps are getting by taking differentials.
\end{defn}

\begin{prop}Let $X_{\bullet}, Y_{\bullet} \in \ilgpd$, and $|X_{\bullet}|, |XY_{\bullet}|$ denote their associated stacks.
	We have a commutative square
	
\begin{center}
	\begin{tikzcd}
	X_{\bullet} \arrow[r, "d"] \arrow[d, "d"] & X_{\bullet}
	 \arrow[d, "d"] \\
	Y_{\bullet} \arrow[r, "d"] \arrow[r, "d"]      &  Y_{\bullet}                 
	\end{tikzcd}
\end{center}
\end{prop}

\begin{proof}
	\end{proof}
\subsection{$\infty$-vector fields on $\infty$-stack}
\begin{defn}
Let $\mathfrak{X}\in \shi(\mfd)$, we define an ($\infty$-){\it vector field} on $X$ to be a pair (X, $\epsilon_X$), where $X$ is a morphism $X: \mathfrak{X} \to T\mathfrak{X}$, and $\epsilon_X: \pi_{\mathfrak{X}}\circ X\Rightarrow \id_{\mathfrak{X}}$ is an equivalence.
\end{defn}
Similarly, we can define vector fields on a Lie $\infty$-groupoid.

\begin{defn}
	Let $G_{\bullet}\in \ilgpd$. An ($\infty$-)vector field on $G_{\bullet}$ is a morphism $X: G_{\bullet} \to TG_{\bullet}$ such that $\pi_{G_{\bullet}}\circ X = \id_{{G_{\bullet}}}$.
\end{defn}
\begin{prop}
    Given a Lie $\infty$-groupoid $G_{\bullet}$, and denote $G$ the associated $\infty$-stack, then we have an equivalence of category $\vect(G_{\bullet}) \simeq \vect(G)$
\end{prop}
\subsection{Higher foliations}

\subsubsection{Foliations on stacks}

Let $M$ be a smooth manifold. A regular foliation is defined as an involutive sub-bundle of the tangent bundle $TM$. This easily generalized to Lie $\infty$-groupoids.
\begin{defn}
	Let $X_{\bullet}$ be a Lie $\infty$-groupoid. We define a   {\it $\infty$-foliation}\index{foliation!higher} $\CF$ on $X_{\bullet}$ to be a sub-Lie $\infty$-groupoid $A_{\bullet}$ of the tangent $\infty$-groupoid $TX_{\bullet}$, where at each level, $A_i \subset TX_i$ is an involutive sub-bundle of $TX_i$.
\end{defn}
Recall that a singular foliation $\CF$ on a smooth manifold $M$ is defined as a subsheaf of the tangent $T_M$ which is involutive and locally finitely generated as a $\cinf(X)$-module. Replace the $\infty$-foliation degreewise by a singular foliation, we get a higher notion of singular foliation.
\begin{defn}
	Let $X_{\bullet}$ be a Lie $\infty$-groupoid. We define a singular  ($\infty$-)foliation $\CF$ on $X_{\bullet}$ to be a simplicial set $\CF_{\bullet}$, where at each level, $\CF_i$ is a subsheaf of the tangent $T_M$ which is involutive and locally finitely generated as a $\cinf(X)$-module.
\end{defn}
Here the simplicial set $\CF$ is actually a simplicial sheaf of $\cinf(X)$-modules, by applying the forgetful functor from sheaf of $\cinf(X)$-modules to sheaf of sets, we can regard $\CF$ as an element of $s\sh(\mfd)$.

Similarly, we can consider foliations on $\infty$-stacks.

\subsubsection{Foliation on 1-stack}

\part{Higher Riemann-Hilbert correspondence for foliations}

In this chapter, we study more in depth about the foliation dga (algebroid). Recall that, for a smooth manifold, we have the {\it de Rham theorem}\index{de Rham theorem}: given a manifold $M$, the singular cohomology groups $H^{\bt}(M, \R)$ and the de Rham cohomology groupoids $H^{\bt}_{\dr}(M, \R)$ are isomorphic, i.e.
\begin{equation*}
    H^{\bt}(M, \R)\simeq H^{\bt}_{\dr}(M, \R)
\end{equation*}
In other words, the singular cochain dga $C^{\bt}(M, \R)$ and de Rham dga $\mathcal{A}^{\bt}(M, \R)$ are {\it quasi-isomorphic}.

However, this quasi-isomorphism is not an dga quasi-isomorphism, since the product structure is not preserved. However, Guggenheim \cite{Gug77} proved that this quasi-isomorphism lifts to an $A_{\infty}$-quasi-isomorphism, where the product structure is preserved up to a higher homotopy coherence. We first study foliated dga's and prove an $A_{\infty}$ de Rham theorem for foliations. 

On the other hand, the similar method can be applied to modules over foliated dga's (algebroids). Recall that the classical {\it Riemann-Hilbert correspondence}\index{Riemann-Hilbert correspondence} (for manifolds) established the following equivalences:
\begin{enumerate}
    \item Local systems over $M$.
    \item Vector bundles with flat connections over $M$.
    \item Representations of the fundamental group of $M$.
\end{enumerate}

Following Chen's iterated integrals \cite{Che77}\cite{Gug77} and Igusa's integration of superconnections \cite{Igu09}, Block-Smith \cite{BS14} proves a higher Riemann-Hilbert correspondence for compact manifolds: the dg category of cohesive modules over the de Rham dga is $A_{\infty}$-quasi-equivalent to the dg-category of $\infty$-local systems over $M$:
\begin{equation*}
    \Mod^{\coh}_{\mathcal{A}} \simeq_{A_{\infty}} \Loc^{\dg}_{\ch_k}(M)
\end{equation*}
where the left-hand side is equivalent to the dg category of $\infty$-representations of the tangent Lie algebroid $TM$, and the right-hand side is equivalent to the dg category of the $\infty$-representations of the fundamental $\infty$-groupoid $\Pi^{\infty}(M)$. Notice that $\Pi^{\infty}(M)$ is equivalent to the integration of $TM$ by the Lie integration functor we mentioned before. Thus, we have the following homotopy-commutative square

\begin{center}
    \begin{tikzcd}
{T_M} \arrow[d, "\rep^{\infty}"] \arrow[r, "\int"] & {\Pi^{\infty}(M)} \arrow[d, "\rep^{\infty}"] \\
{\Mod^{\coh}_{\mathcal{A}}} \arrow[r]           & {\Loc^{\dg}_{\ch_k}(M)}   
\end{tikzcd}
\end{center}
Hence we can really understand the Riemann-Hilbert Correspondence as an equivalence between $\infty$-representations of $L_{\infty}$-algebroids and $\infty$-representations of the integration of $L_{\infty}$-algebroids, i.e. Lie $\infty$-groupoids. 
We apply this idea to the case of foliations and prove a higher Riemann-Hilbert correspondence for foliation, and construct the integration functor from the $\infty$-representations of $L_{\infty}$-algebroids and $\infty$-representations Lie $\infty$-groupoids.

\section{Algebras and modules of foliations}
\subsection{$\mathcal{D}$-module and foliations}
Let $X$ be a smooth manifold and $\CF$ a regular foliation on $X$. We consider a $\D_X$-module\index{$\D_X$-module} associated $\D_{\CF} = \D_X / \D_X\cdot \CF $, which stands for linear differential operator normal to $\CF$. 

First, consider $\CF = T_X$, then $\D_{\CF} / \D_X\cdot T_X \simeq \CO_X$.

Next, consider a general regular foliation $\CF$. Let $\A^{\bullet}= \A^{\bullet}(M, \widehat{\sym}(\CF^{\perp}))= \Omega^{\bullet}(X)\otimes \widehat{\sym}(\CF^{\perp})$. Let $q$ be the codimension of $\CF$. Let $\{x_i\}_i, 1\le i \le n$ coordinates for $U_i$ and $\hat{x}_i = dx_i, 1\le i \le q$ be a basis of $\CF^{\perp}$.

\begin{lem}
	$\nabla = d - \sum_{i=1}^q dx_i\wedge \frac{\del}{\del \hat{x}_i} $ is a flat connection on $\A^{\bullet}= \A^{\bullet}(M, \widehat{\sym}(\CF^{\perp}))$.
\end{lem}

\begin{proof}First, we show $\nabla^2 = 0$. 	Let $\sum_{\alpha}^{}f_{\alpha} \hat{x}^{\alpha} \in  \widehat{\sym}(\CF^{\perp})$, where $\alpha$'s are multi-indices, and $f_{\alpha} \in \cinf(M)$
	\begin{align*}
	\nabla\sum_{\alpha}^{}f_{\alpha}\hat{x}^{\alpha} ) =& ( d - \sum_{j=1}^q dx_j\wedge \frac{\del}{\del \hat{x}_j}) (\sum_{\alpha}^{}f_{\alpha}\hat{x}^{\alpha} ) \\
	=&
	\sum_{\alpha}\bigg(\Big(\sum_{i=1}^{n}dx_i \wedge  \frac{\del f_{\alpha}}{ \del x_i} \hat{x}^{\alpha}\Big) - \sum_{i=1}^q \Big( dx_i \wedge f_{\alpha}\frac{\del \hat{x}^{\alpha}}{\del \hat{x}_i}\Big)\bigg)
	\end{align*}
		Let's look at these two terms in the summand independently. Note that $d$ on the first term is just 0, hence
	\begin{align*}
	\nabla\Big(\sum_{i=1}^{n}dx_i \wedge  \frac{\del f_{\alpha}}{ \del x_i} \hat{x}^{\alpha}\Big) =& -\sum_{j=1}^q\Big(dx_j \wedge\frac{\del}{\del\hat{x_j}}\Big)\Big(\sum_{i=1}^{n}dx_i \wedge  \frac{\del f_{\alpha}}{ \del x_i} \hat{x}^{\alpha}\Big)\\
	=& \sum_{i=1}^{n} \sum_{j=1}^{q}dx_i\wedge dx_j \wedge \frac{\del f_{\alpha}}{ \del x_i} \frac{\del \hat{x}^{\alpha}}{\del \hat{x}_j}
	\end{align*}
	Next, for the second term, 
	\begin{align*}
	\nabla \Big(- \sum_{j=1}^q dx_i \wedge f_{\alpha}\frac{\del \hat{x}^{\alpha}}{\del \hat{x}_i}\Big)=& \sum_{j=1}^q dx_i \wedge \Big( \sum_{k=1}^ndx_k\wedge \frac{\del f_{\alpha}}{\del x_k}  \Big) \frac{\del \hat{x}^{\alpha}}{\del \hat{x}_i} - \sum_{i=1}^q dx_i \wedge \Big(\sum_{l=1}^q dx_l \wedge f_{\alpha} \frac{\del}{\del \hat{x_l}}  \frac{\del \hat{x}^{\alpha}}{\del \hat{x}_i}\Big)\\
	=& \sum_{j=1}^q \sum_{k=1}^n  dx_i \wedge  dx_k\wedge \frac{\del f_{\alpha}}{\del x_k}   \frac{\del \hat{x}^{\alpha}}{\del \hat{x}_i} - \sum_{i=1}^q\sum_{l=1}^q dx_i \wedge dx_l\wedge  f_{\alpha} \frac{\del}{\del \hat{x_l}}  \frac{\del \hat{x}^{\alpha}}{\del \hat{x}_i}
	\end{align*}
	The first term cancels the term in the precious equation and the second is clearly vanished. Hence, we see $\nabla$ is flat.
	
	Next, we want to show $\nabla$ is well-defined. Let $U$ and $V$ be two foliated neighborhoods with nonempty intersection. Let $\{x_i\}_{i=1}^n$ and $\{y_i\}_{i=1}^n$ be coordinates on $U$ and $V$ respectively. consider $\phi: U_i \to U_j$ be a transition function between foliated neighborhoods. Note that $\phi: U\simeq \R^q\times \R^{n-q} \to  \R^q\times \R^{n-q}\simeq V$ has the following form
	\begin{enumerate}
		\item $\phi_i(x)=\phi_i(x_1,\cdots, x_q)$ for $1\le i\le q$.
	    \item $\phi_i(x)=\phi_i(x_1,\cdots, x_n)$ for $q+1 \le i\le n$.
	\end{enumerate} 
    By our assumption, $\phi_i(x) = y_i$.
    \begin{align*}
    d - \sum_{i=1}^q dy_i\wedge \frac{\del}{\del \hat{y}_i} =& \sum_{i=1}^n dy_i\wedge \frac{\del}{\del y_i}- \sum_{i=1}^q dy_i\wedge \frac{\del}{\del \hat{y}_i}\\
    &= \sum_{i=1}^n \Big(\sum_{j=1}^n dx_j \wedge\frac{\partial \phi_i}{ \del x_j}\Big)\wedge \frac{\del}{\del y_i}- \sum_{i=1}^q \Big(\sum_{l=1}^n dx_j \wedge\frac{\partial \phi_i}{ \del x_j}\Big)\wedge \Big(\sum_{k=1}^{q}\frac{\del \hat{x}_k}{\del \hat{y}_i}\wedge \frac{\del}{\del \hat{x}_k}\Big) 
    \end{align*}
    
    Note that $\frac{\del \phi_i}{\del x_j} = 0$ for $1\le i \le q$ and $q+1 \le j \le n$. Hence, the above equation becomes
    
    \begin{align*}
    \sum_{i=1}^n \Big(\sum_{j=1}^q dx_j \wedge\frac{\partial \phi_i}{ \del x_j}\Big)\wedge \frac{\del}{\del y_i}- \sum_{i=1}^q \Big(\sum_{l=1}^q dx_l \wedge\frac{\partial \phi_i}{ \del x_l}\Big)\wedge \Big(\sum_{k=1}^{q}\frac{\del \hat{x}_k}{\del \hat{y}_i}\wedge \frac{\del}{\del \hat{x}_k}\Big) \\
    =\sum_{j=1}^q dx_j \wedge \Big( \sum_{i=1}^n \frac{\partial \phi_i}{ \del x_j}\wedge \frac{\del}{\del y_i}\Big)-  \sum_{l=1}^q\sum_{k=1}^{q} dx_l \wedge \frac{\del}{\del \hat{x}_k} \wedge \Big(\sum_{i=1}^q \frac{\partial \phi_i}{ \del x_l}\wedge \frac{\del \hat{x}_k}{\del \hat{y}_i}\Big) 
    \end{align*}
    Note that ${\del \hat{x}_k}/{\del \hat{y}_i}= {\del (dx_k)}/{\del (dy_i)}= \del x_k/\del y_i$, and then $\Big(\sum_{i=1}^q \frac{\partial \phi_i}{ \del x_l}\wedge \frac{\del \hat{x}_k}{\del \hat{y}_i}\Big) $ equals a diagonal matrix which restricts to $I_q$ on the top left $q\times q$ submatrix and all the other entries are 0. Hence, we get
    \begin{align*}
    \sum_{j=1}^q dx_j \wedge \Big( \sum_{i=1}^n \frac{\partial \phi_i}{ \del x_j}\wedge \frac{\del}{\del y_i}\Big)-  \sum_{l=1}^q\sum_{k=1}^{q} \delta_{lk} dx_l \wedge \frac{\del}{\del \hat{x}_k} 
    =  \sum_{j=1}^q dx_j \wedge  \frac{\del}{\del x_j}-  \sum_{m=1}^q dx_m \wedge \frac{\del}{\del \hat{x}_m} 
    \end{align*}

	\end{proof}

   Next, let us look at the cohomology of the dga $\A^{\bullet}= \A^{\bullet}(M, \widehat{\sym}(\CF^{\perp}))$
   
   \begin{lem}The 0-th cohomology of $\A^{\bullet}$ equal $\cinf$ functions on $M$ which is constant on leaves, i.e. 
   	$$H^{0} (\A^{\bullet}(M, \widehat{\sym}(\CF^{\perp})) \simeq \CO_{\CF}$$
   \end{lem}

\begin{proof}
	Let $\sum_{\alpha}^{}f_{\alpha} \hat{x}^{\alpha} \in  \widehat{\sym}(\CF^{\perp})$, where $\alpha$'s are multi-indices and $f_{\alpha} \in \cinf(M)$, then from the previous proof
	\begin{align*}
	\nabla(\sum_{\alpha}^{}f_{\alpha}\hat{x}^{\alpha} ) 
	=&
	\sum_{\alpha}\bigg(\Big(\sum_{i=1}^{n}dx_i \wedge  \frac{\del f_{\alpha}}{ \del x_i} \hat{x}^{\alpha}\Big) - \sum_{i=1}^q \Big( dx_i \wedge f_{\alpha}\frac{\del \hat{x}^{\alpha}}{\del \hat{x}_i}\Big)\bigg)\\
	=&
	\sum_{\alpha}\bigg(\Big(\sum_{i=1}^{n}dx_i \wedge  \frac{\del f_{\alpha}}{ \del x_i} \hat{x}^{\alpha} \Big)- \sum_{i=1}^q \Big( dx_i \wedge f_{\alpha} \alpha_i \hat{x}^{\alpha -1_i}\bigg)
	\end{align*}

	Here $\alpha = (\alpha_1, \cdots, \alpha_i, \cdots, \alpha_n)$, and $1_i$ denotes the multi-index with 1 at the $i$-th entry and 0's elsewhere. Now we can group the coefficients of $dx_i \wedge \hat{x}^{\alpha}$, so we get
	\begin{align*}
	\sum_{i=1}^{q} dx_i \wedge   	\sum_{\alpha} \Big(\frac{\del f_{\alpha}}{ \del x_i} - (\alpha_i+1)f_{\alpha + 1_i} \hat{x}^{\alpha} \Big) - \sum_{i= q+1}^n dx_i \wedge \Big(\sum_{\alpha}^{} \frac{\del f_{\alpha}}{ \del x_i} \hat{x}^{\alpha}\Big)
	\end{align*}
	Hence $\ker\nabla$ consists of sections $\sum_{\alpha}^{}f_{\alpha} \hat{x}^{\alpha}$ of $\widehat{\sym}(\CF^{\perp})$ where $f_{\alpha}$ satisfies $\frac{\del f_{\alpha}}{ \del x_i} - (\alpha_i+1)f_{\alpha + 1_i} = 0$ for $1\le i \le q$ and $\frac{\del f_{\alpha}}{ \del x_i}=0$ for $q+1\le i \le n$. The first condition implies $\sum_{\alpha}^{}f_{\alpha} \hat{x}^{\alpha}$ is holonomic and the second condition implies that $f_{\alpha}$'s are constant along leaves.
	
	\end{proof}
\begin{cor}
	The cohomology of  $\A^{\bullet}$ is isomorphic to $\Omega^{\bullet}_M\otimes_{\cinf(M)} \CO_{\CF}$.
\end{cor}

\subsection{Sheaf of constant functions along leaves}
We denote the sheaf of (smooth) functions on $M$ which are constant along leaves of $\CF$ by $\underline{\R}_{\CF}$. Regard $(M, \underline{\R}_{\CF})$ as a ringed space, then the sheaf of $\cinf$-functions $\cinf_M$ on $M$ is a sheave of $\underline{\R}_{\CF}$-module. We have the following conjecture:

	\paragraph{Conjecture} Given a foliation $(M, \CF)$, $\cinf_M$ is flat over $\underline{\R}_{\CF}$.
	
This is a proposition encoding differential geometric properties into a simple algebraic form, which will be useful in proving many results later. We won't prove it in this paper, and we shall use other method ($\cinf$-rings or topological algebras) to get rid of our issues. In this chapter, we will prove a partial result on this conjecture.

 The problem is local. It suffices to show $\cinf_{M,x}$ is flat over $\underline{\R}_{\CF, x}$ for all $x \in M$. Picking a foliation chart and a foliated neighborhood $ U \simeq \R^q \times \R^{n-q}$, then $\cinf_{M}(U)\simeq \cinf(\R^n)$ and $\underline{\R}_{\CF}\simeq \cinf(\R^{n-q})$. Hence, it suffices to show the following lemma
	\begin{lem}$ \cinf(\R^n)$ is a flat $ \cinf(\R^{n-q})$ module for $q\ge 0$.
	\end{lem}
The module structure is induced by the projection $p: \R^n \to \R^{n-q}$. Given $a_1, \cdots, a_k\in  \cinf(\R^{n-q})$ and  $b_1, \cdots, b_k\in  \cinf(\R^{n})$ such that $\sum_i a_ib_i = 0$, we want to show that there exist functions $G_1, \cdots, G_r \in \cinf_{\R^n}$ and $c_{ij}\in \R^{n-q}$, such that $\sum_{j=1}c_{ij}G_j$ for all $i$ and $\sum_i a_i c_{ij} = 0$ for all $j$.
	
	 Let's first consider the simple case $n=2, q=1$. We start by the following lemma, which is a special case of flatness when $k = 1$.
	 
	 \begin{lem}\label{lem1} Let $h\in \cinf(\R)$ to be strictly positive for $x<0$ and 0 for $x \ge 0$, and $g\in \cinf(\R^2)$. Let $\cinf(\R^2)$ as a $\cinf(\R)$-module induced by the projection $\R^2 \to \R$. Suppose $hg =0$, then $g(x,y) = c(x)G(x, y)$, where $c(x) \in \cinf(\R)$ vanishes on $x\le 0$ and $G(x,y) \in \cinf(\R^2)$ .
	 	\end{lem}
 	
 	\begin{proof}
 		Consider two sets $I = \{f\in \cinf(\R)|\text{$f=0$ when $x\le 0$ and $f> 0$ when $x> 0$} \}$, $J = \{f\in \map(\R^+, \R)| \text{$f(x)/x^n \to 0$ as $x\to 0$ for all $n>0$}\}$.
 		
 		\begin{lem} For any $f\in J$, there exists a $g\in I$ such that $f/g\to 0$ as $x\to 0$.
 			\end{lem}
 		\begin{proof}
 			Consider a bump function $r\in \cinf(\R)$ such that $r= 1$ for $x\le 0$ and $r=0$ for $x \ge 1$. Let $\{a_k\}$ be a monotonically decreasing sequence such that $\sum_i a_k < \infty$ and $a_k>0$ for all $k$. Define $\theta(x) = \sum_{i=1}^{\infty} r(x/a_i)$. Let 
 			\begin{equation}
 			g(x) = \begin{cases}
 			x^{\theta(x)} \quad &x>0\\
 			0 &x\le 0
 			\end{cases}.
 			\end{equation}
 			We claim that $g(x)$ satisfied the requirement in lemma. Outside any open neighborhood of $0$, there will be only finitely many non-zero summands in $\theta$, hence $g(x)$ is smooth and bounded outside any open neighborhood of $0$. Clearly $g(x)/x^n\to 0$ as $x\to 0$. We just need to check all derivative of $g$. Let $\theta_n(x) = \sum_{i=1}^n a_k$ and $g_n = x^{\theta_n}$, we have 
 			$$
 			\frac{d^l}{dx^l}g_n =\frac{d^l}{dx^l} x^{\theta_n} = \theta_n(x)\cdots(\theta_n(x)-l+1) x^{\theta_n(x) - l}
 			$$
 			for $x>0$. For each $n$, there exists an $\epsilon_n$ such that for $0<x<\epsilon_n$, $\frac{d^l}{dx^l} x^{\theta_{n+1}}<\frac{d^l}{dx^l} x^{\theta_n}$. For each $l$, if we pick $n$ large enough, e.g. $l >n$, then $\frac{d^l}{dx^l} x^{\theta_n} = 0$. Let $\epsilon>0$ be arbitrary, we want to show that there exists  $x_0>0$ and $N\in \N$ such that for all $n>N, 0<x<x_0$, $|\frac{d^l}{dx^l} x^{\theta_n(x)}-\frac{d^l}{dx^l} x^{\theta(x)}| < \epsilon$, that is $x^{\theta_n}\to x^{\theta}$ uniformly on $[0, x_0]$. Consider $n>l$ to be sufficiently large, then there exist $x_1> 0$ such that $\frac{d^l}{dx^l} x^{\theta_n} < \epsilon$ for $0\le x \le x_1$. Now
 			\begin{align*}
 				\frac{d^l}{dx^l} x^{\theta_{n+1}}& = \theta_{n+1}(x)\cdots(\theta_{n+1}(x)-l+1) x^{\theta_{n+1}(x) - l}\\
 				& =\frac{\theta_{n+1}(x)\cdots(\theta_{n+1}(x)-l+1)}{\theta_n(x)\cdots(\theta_n(x)-l+1)}\theta_n(x)\cdots(\theta_n(x)-l+1) x^{\theta_n(x) - l}x^{r(x/a_{n+1})} \\
 				& \le x^{\theta_{n+1}(x) - l}  \Bigg(\frac{\theta_{n+1}(x)-l+1}{\theta_{n}(x)-l+1} \Bigg)^l \frac{d^l}{dx^l} x^{\theta_{n}}   
 			\end{align*}
 		    Note that $\frac{\theta_{n+1}(x)-l+1}{\theta_{n}(x)-l+1}$ is monotonically decreasing as $n$ increases and as $x$ decreases. Let $\delta >0$ such that $(\frac{\theta_{n+1}(x)-l+1}{\theta_{n}(x)-l+1})^l < 1+\delta$ for all $0\le 0 \le x_1$, then we can find an $x_2<x_1$ such that $x^{r(x/a_{n+1})} <1/(1+\delta)$. Therefore, $\frac{d^l}{dx^l} x^{\theta_{n+1}} < \epsilon$ on $[0, x_2]$. By induction, we get  $\frac{d^l}{dx^l} x^{\theta_{k}} < \epsilon$ for all $k\ge n$ and $0\le x \le x_2$. 
 		    
 		    Now picking a sequence $\{x_m\}$ such that all $x_m \le x_2$ and $x_m \to 0$ monotonically, then $\lim_{m\to 0} g(x_m) = 0$ by continuity, and 
 		    $$
 		    \lim_{n\to \infty}\lim_{m\to 0}\frac{d^l}{dx^l}g_n(x_m) = \lim_{n\to \infty}\frac{d^l}{dx^l}g_n(0) = 0
 		    $$,
 		    then by the uniform convergence, we can change the order of limits and get
 			$$
 			\lim_{m\to 0}\lim_{n\to \infty}\frac{d^l}{dx^l}g_n(x_m) = \lim_{m\to 0}\frac{d^l}{dx^l}g(x_m) = 0
 			$$
 			
 			Now we have shown that $g\in J$. Since $f=0$ for all $x\le 0$, all derivatives of $f$ must vanish at infinite order at $0$, hence for each $n>0$, there exists a decreasing sequence $\{\epsilon_n\}$ and $\epsilon_k <1$ for all $k$ such that $|f|\le x^{n}$ for all $x\in [0, \epsilon_n]$. Now we just need to pick $a_k < \epsilon_{n+1}$ for all $k\ge n$ which ensure that $g(x)/x^{n-1}>1$ for $x \in (\epsilon_{n+1},\epsilon_n)$. Hence, $|f(x)/g(x)| \le x $  on $(\epsilon_{n+1},\epsilon_n)$ for all $n$, which implies $f/g \to 0$ as $x \to 0$.

%
%
%
 			\end{proof}
 		\begin{cor}
 			Given a sequence $\{f_i\}$ in $J$, there exists a $g\in I$ such that $f_k/g \to 0$ as $x \to 0$ for all $k$.
 		\end{cor}
 	
 		\begin{cor}
 		Given a sequence $\{f_i\}$ in $J$, there exists a $g\in I$ such that $f_k/g^n \to 0$ as $x \to 0$ for all $k, n \in \N$.
 	\end{cor}
 
 Now consider $f_{ijk}(x) = \sup_{|y| \le k}\big(\frac{d^{i+j}}{dx^i dy^j}g(x,y)\big)$. Then by previous lemma, we can find an $a(x)$ such that $f_{ijk}/a^n \to 0$ as $x \to 0^+$ for all $i, j, k, n$. Now, take $$G(x, y)= \begin{cases}
                   g(x,y)/a(x)\quad & x >0\\
                   0 & x\le 0
            \end{cases}
 $$
 It suffices to verify $G^{(n)}(x,y)\to 0$ as $x\to 0^+$. Expand by quotient rules, we have
 $$
  			\bigg(\frac{g(x,y)}{a(x)}\bigg)^{(n)} =   \frac{1}{a} \left( g^{(n)} -\sum_{j=1}^{n} \binom{n}{j} \bigg(\frac{g}{a}\bigg)^{(n-1-j)}a^{(j)} \right).
 			$$
which goes to $0$ by induction. Hence, $g=aG$ is the desired factorization.                                                                                                                                                                                                                                                                                                                                                                                                                                                                                                                                                                                                          
 
\end{proof}
 		    
 \begin{prop}
 	Let $p: \R^2 \to \R$ be a submersion. Let $h \in \cinf(\R)$ and $g\in \cinf{(\R^2)}$ such that $hg=0$, then there exist an $a\in \cinf(\R)$ and $G\in \cinf(\R^2)$ such that $h = aG$ and $ah = 0$.
 \end{prop}
\begin{proof}It suffices to consider the case $p$ is the projection. Without loss of generality, we restrict the domain on an open neighborhood $U$ of the origin on $\R$. Let $V = h^{-1}(0) \subset U$. If $V$ is nowhere dense, then $g$ must vanish on $p^{-1}(U\setminus V)\simeq (U\setminus V) \times \R$, which have to vanish on $\overline{U\setminus V}\times \R = \overline{U}\times \R$. In this case, $g$ vanishes on $\overline{U}\times \R$, then we just take $a$ be $a(x) = 0$ and $G$ arbitrary. Now suppose $x_0 \in U$ is contained in a closed interval $V'$. Let $U'\supset V'$ be a open neighborhood of $V'$ such that $g$ vanishes on $\overline{U'\setminus V'}$. Note that $h$ can still vanish on a nowhere dense set on $U'\setminus V'$.  Then $g$ vanishes at infinite order at $\del V'$. By previous lemma, there exist an $a\in \cinf(\R)$ vanishes on $x\in \overline{U'\setminus V'}$ and $g = aG$ for some $G\in \cinf(\R^2)$.
	\end{proof}

\begin{lem}
	Let $h_1,\cdots, h_k \in \cinf(\R)$ to be strictly positive for $x<0$ and 0 for $x \ge 0$, and $g_1, \cdots, g_k \in \cinf(\R^2)$ . Let $\cinf(\R^2)$ as a $\cinf(\R)$-module induced by the projection $\R^2 \to \R$. Suppose $\sum_i h_i g_i =0$, then $g_i(x,y) = c_{ij}(x)G_j(x, y)$, where $c_{ij}(x) \in \cinf(\R)$ vanishes on $x\le 0$ and $G_j(x, y)\in \cinf(\R^2)$.
\end{lem}
\section{ $A_{\infty}$ de Rham theorem for foliations}
\subsection{de Rham theorem for foliations}
\begin{thm}[de Rham theorem for foliations\index{de Rham theorem!for foliations}]\label{dr}
	Given a foliation $(M, \CF)$, there exists a isomorphism
	\begin{equation}
	H^{\bullet}(M, \extp^{\bullet}\CF^{\vee}) \simeq H^{\bullet}(M, C^{\bullet}(\CF))
	\end{equation} 
\end{thm}

Consider the codimension $q$ product foliation $\R^{n-q} \times \R^q$ on $\R^n$, we can build two new product foliations $\R^{n-q+1}\times \R^q$ and $\R^{n-q}\times \R^{q+1}$ out of it. Let $(x_1,\cdots, x_{n-q}, x_{n-q+1}, \cdots, x_n)$ be the canonical coordinates on $\R^n$, then $\CF = \{\del_{x_1}, \cdots, \del_{{n-q}} \}$. Therefore we have $\extp^{\bullet}\CF^{\vee} \simeq (\extp^{\bullet} T^{\vee}\R^{n-q})\times \R^q$, which yields $H^{\bullet}(\R^n, \CF) \simeq \Omega^{\bullet}(\R^{n-q})$. Therefore, $H^{\bullet}(\R^{n+1}, \R^{n-q}\times \R^{q + 1})\simeq H^{\bullet}(\R^{n+1}, \R^{n-q+1}\times \R^q)$. On the other hand, we also have $H^{\bullet}(\R^{n+1}, \R^{n-q+1}\times \R^{q })\simeq H^{\bullet}(\R^{n}, \R^{n-q}\times \R^q)$ by Poincare lemma for $\R^n$.
\begin{lem}[Poincare lemma for foliations\index{Poincare lemma!for foliations}]
	Consider the codimension $q$ product foliation $(\R^n, \CF) = (\R^n, \R^{n-q} \times \R^q)$, then 
	\begin{equation}
	H^{i}(\R^n,\extp^{\bullet}\CF) ) \simeq
	\begin{cases}
	\cinf(\R^q) &\quad i = 0  \\
	0 & i \not= 0 
	\end{cases}	
	\end{equation}.
\end{lem}
\begin{proof}
	 By previous observation,$H^{i}(\R^n,\extp^{\bullet}\CF) )\simeq H^i(\Omega^{\bullet}(\R^{n-q}))$, then the result follows from Poincare lemma for $\R^n$.  
	\end{proof}

\begin{defn}
	We define the (smooth) $\CF$-foliated singular $n$-chain $C_n(\CF, G)$ of a foliation $(M, \CF)$ to be the free Abelian group generated by (smooth) foliated $n$-simplices $\sigma: \Delta^n \to \CF$ with coefficient in some Abelian group $G$. Define the differential $d_n: C_n(\CF, G) \to C_{n-1}(\CF, G)$ by $d_n = \sum_{i=0}^{n} (-1)^i \delta_i$, where $\delta_i$ is the $i$-th face map. We call $\big(C_*(\CF, G), d\big)$ the foliated singular chain complex.
\end{defn}

\begin{defn}
	We define the (smooth) foliated singular cochains $C^{\bullet}(\CF)$ to be $\cinf$ function on the monodromy $\infty$-groupoid $\mon_{\infty} \CF$ associated to $\CF$, i.e. $C^n(\CF) = \cinf(\mon_n\CF, \R)$. 
\end{defn}

\begin{lem}{\label{poincare cochain}}
	Consider the codimension $q$ product foliation $(\R^n, \CF) = (\R^n, \R^{n-q} \times \R^q)$, then 
\begin{equation}
H^{i}(\R^n,C^{\bullet}\CF) ) \simeq
\begin{cases}
\cinf(\R^q) &\quad i = 0  \\
0 & i \not= 0 
\end{cases}	
\end{equation}.
\end{lem}
\begin{proof}
	Given a $k$-simplex  $\sigma: \Delta^n \to \CF$ ($0\le k \le n-q-1$), define $K: C_{k} \to C_{k+1}$ by $K\sigma(\sum_{j=0}^{k+1}t_jx_j) = (1 - t_{q+1})\sigma(\sum_{j=0}^q\frac{t_j}{1-t_{q+1}} x_j)$ which sends a foliated $k$-simplex to $k+1$ simplex, then by standard calculation we have $\del K- K \del = (-1)^{q+1}$. Let $L$ be the adjoint of $K$, then $(-1)^{k+1}(dL-Ld)=1$, which gives the result.
	\end{proof}

 On the other hand, $C^n\CF$ are soft since  $C^n\CF$ are sheaves of  $C^0{\CF} \simeq \cinf(M)$-modules.

\begin{proof}[Proof of de Rham theorem]
	By Poincare lemma, we have $0\to \R_{\CF} \to \Gamma(\extp^{\bullet}\CF)$ which is a resolution of $\underline{\R}_{\CF}$ by fine sheaves. Note that $C^n(\CF)$'s are sheaves of $C^0(\CF)\simeq \cinf(M)$-modules, which are soft since $\cinf(M)$ is. By lemma , $C^{\bullet}(\CF)$ is a soft resolution of $\underline{\R}_{\CF}$. Then integration over chains gives the desired quasi-isomorphism.
	\end{proof}

Next, we are going to show the quasi-isomorphism between the dga of $\CF$-foliated forms and the dga of smooth singular $\CF$-cochains actually lifts to an $A_{\infty}$-quasi-isomorphism
\begin{equation*}
    \phi:(\extp^{\bullet}\CF^{\vee}, -d, \wedge) \to  (C^{\bullet}(\CF), \delta, \cup)
\end{equation*}
The $\phi$ is defined as a composition of two maps
\begin{equation*}
    \mathsf{B} \big((\extp^{\bullet}\CF^{\vee})[1]\big)\stackrel{F}{\to} \Omega^{\bt}(P\CF) \stackrel{G}{\to} C^{\bt}(\CF)[1]
\end{equation*}
here $\mathsf{B}$ is the bar construction. The first map is similar to Chen's iterated integral map, and the second map is similar to Igusa's construction in \cite{Igu09}.

\subsection{Riemann-Hilbert correspondence}

\begin{thm}[Riemann-Hilbert correspondance for foliation]
	Let $(M, \CF)$ be a manifold with foliation $\CF$, then the following categories are equivalent
	\begin{enumerate}
	    \item The category of foliated local systems $\Loc(\CF)$.
	    \item The category of vector bundles with flat $\CF$-connection.
	    \item The category of the representations of the fundamental groupoid.
	\end{enumerate}
\end{thm}

Let $P\CF$ denote the Frech{\'e}t manifold $P^{1}_{C^{\infty}} \CF$ which consists of smooth path along leaves. We parametrize geometric $k$-simplex $\Delta^n$ by $t = (1 \ge t_1\ge t_2 \cdots \ge t_k \ge 0)$. First we have a map of evaluation on a path
$$
\ev_k: P\CF \times \Delta^k \to M^k: (\gamma, (t_1, \cdots, t_k)) \mapsto \big(\gamma(t_1), \cdots, \gamma(t_k)\big)
$$
The image of $\ev_k$ fixing $\gamma$ lies in a single leaf. Along with the natural inclusion $P\CF \subset PM$, the following diagrams commutes

\begin{center}
\begin{tikzcd}
{ P\CF \times \Delta^k} \arrow[r] \arrow[rd] \arrow[d, hook] & {\coprod_{x\in M} L_x}           \\
{ PM \times \Delta^k} \arrow[r]                            & {M^k} \arrow[u]
\end{tikzcd}
\end{center}

\begin{defn}
	We define  $T_{\CF}P\CF$ to be a vector bundle whose fiber at $\gamma \in P\CF$ is the vector space of all $\cinf$-sections $I\to \CF$ along $\gamma$. We define the dual bundle $T^{\vee}_{\CF} P \CF$ of  $T_{\CF}P\CF$ to be the vector bundle whose fiber at $\gamma$ is the space of all bounded linear functionals, i.e. $T^{\vee}_{\CF, \gamma} P\CF = \Hom(T^{}_{\CF, \gamma} P\CF, \R)$.  
\end{defn}
We denote the $\cinf$-section of $T^{\vee}_{\CF, \gamma} P\CF$ by $\Omega_{\CF}^1 P\CF$, and the exterior algebra of $\Omega_{\CF}^1 P\CF$ by $\Omega_{\CF}^{\bullet} P\CF$.

\begin{lem}
	Let $f\in \cinf(PM)$ and $\gamma_0 \in M$, there exists a unique section $D_{} f \in \Omega_{\CF}^1 P\CF.$
\end{lem}
\begin{proof}Let $\eta \in T_{\CF, \gamma_0}P\CF$. Take an one-parameter deformation $\gamma_s$ of $\gamma_0$ such that $\frac{\del}{\del s} \gamma_s = \eta$, then we can define $Df\big|_{\gamma_0} (\eta)= \frac{\del}{\del s}\big|_{s=0}(f\circ \gamma_s)$. We want to show this gives a unique bounded linear functional on $T_{\CF, \gamma_0} P\CF$. The boundedness and linearity is obvious.  
	\end{proof}
\begin{cor}
	For any smooth deformation $\gamma_s$ of $\gamma_0$, we have the following chain role
	$$
	Df\big|_{\gamma_0} \big(\frac{\del}{\del s}\bigg|_{s = 0} \gamma_s\big) = \frac{\del}{\del s}\bigg|_{s = 0} (f\circ\gamma_s)
	$$.
\end{cor}

Next, we want to define higher differentials on $\Omega_{\CF}^{\bullet} P\CF$. A key observation is that $T_{\CF}P\CF$ is involutive. Given two elements $\eta, \zeta \in \Gamma(T_{\CF, \gamma_0}P\CF)$, we can regard them as sections $I \to \CF$ along $\gamma_0$. Then, by involutivity of $\CF$, $[\gamma, \eta]$ is still a section $I\to \CF$. Using this fact, we can define all higher differential on $\Omega_{\CF}^{\bullet} P\CF$ simply by Chevalley-Eilenberg formula.

$\ev_1$ induces a smooth map $T\ev_1\big|_{\gamma,t}:T_{\CF, \gamma}P\CF \to \CF_{\gamma(t)}$. Given a vector bundle $V$ on $M$, we can get a pullback bundle $W_t$ along $\ev_1$ at time $t$, i.e. $W_t = \ev_{1}^* V_{\gamma(t)}$. Hence, $W$ is a vector bundle on $P\CF \times \Delta_1$.
\begin{lem}
	$\ev_1^* \Gamma (\CF^{\vee})$ lies in $\Gamma(T^{\vee}_{\CF} P\CF)$.
\end{lem}
\begin{proof}
	\end{proof}
\subsection{Chen's iterated integral}
Let $\pi: \CF \times \Delta^k \to \CF$ be the projection on the first factor.
   Define the push forward map
   \begin{equation*}
       \pi_*: \extp^{\bt}\big((\CF \times \Delta^k)^{\vee}\big) \to \extp^{\bt}(\CF^{\vee})
   \end{equation*}
   by
   \begin{equation*}
       \pi_*\big(f(x,t)dt_{i_1}\cdots dt_{i_k}dx_{j_1}\cdots dx_{j_s} \big) = \Big(\int_{\Delta_k} f(x,t)dt_{i_1}\cdots dt_{i_k} \Big)dx_{j_1}\cdots dx_{j_s}
   \end{equation*}
   Note that here $\CF \times \Delta^k$ is a foliation on $M\times \Delta^k$ which extends $\CF$ trivial along the $\Delta^k$ direction, i.e. $\CF\times \Delta^k := \CF \times T\Delta^k$.
   
   If $M$ is compact, we have
   $$
   \int_M \pi_*(\alpha) = \int_{M\times \Delta^k} \alpha
   $$
   for all $\alpha \in \extp^{\bt}(\CF{\vee}\otimes \Delta^k)$.
   
   \begin{lem}\label{a6}
       $\pi_*$ is a morphism of left $\extp^{\bt}\CF^{\vee}$-modules of degree $-k$, i.e for every $alpha \in \extp^{\bt}(\CF^{\vee})$ and $\beta \in \extp^{\bt}\big((\CF \times \Delta^k)^{\vee}\big)$, we have 
       \begin{equation}
           \pi_*(\pi^* \alpha \wedge \beta) = (-1)^{|\alpha| k} \alpha \wedge \pi_* \beta
       \end{equation}
       
       In addition, let $\del \pi$ be the composition
       \begin{equation*}
           \CF \otimes \del \Delta^k \stackrel{\id \otimes \iota}{\longrightarrow} \CF \otimes \Delta^k \stackrel{\pi}{\longrightarrow} M
       \end{equation*}
       Then we have 
       \begin{equation*}
           \pi_* \circ d - (-1)^k d\circ \pi_* = (\del \pi)_* \circ (\id \times \iota)^*
       \end{equation*}
   \end{lem}
    \begin{proof}. 
        Similar to \cite[Lemma 3.5]{AS12}. Note that we just need to restrict to integration along leaves.
    \end{proof}
	Next, we shall construct Chen's iterated integral map. Let $a_1[1]\otimes \cdots \otimes a_n[1]$ be an element of $\mathsf{B} \big((\extp^{\bullet}\CF^{\vee})[1]\big)$. Given a path $\gamma: I \to \CF \in P\CF$, we define a differential form on $P\CF$ by
	\begin{enumerate}
	    \item Pull back each $a_i$ to $M^k$ via the $i$-th projection map $p_i: M^k \to M$, then we get a wedge product $p_1^* a_1\wedge \cdots p_k^* a_k$.
	    \item Pullback $p_1^* a_1\wedge \cdots p_k^* a_k$ to a form on $P\CF \times \Delta^k$ via $\ev_k$.
	    \item Push forward through $\pi$ to get a form on $P\CF$.
	    \item Finally, correct the sign by multiplying $\spadesuit = \sum_{1 \le i < k} (T(a_i) - 1)(k-i)$ where $T(a_i)$ denotes the total degree of $a_i$.
	\end{enumerate}
	In summary,
	\begin{defn}[Chen's iterated integrals on foliated manifold\index{iterated integrals}]
	    Let $(M, \CF)$ be a foliated manifold, define Chen's iterated integral map from the bar complex of the suspension of foliation algebra to the foliated path space by
	    \begin{equation}
	        \mathsf{C}(a_1[1]\otimes \cdots \otimes a_k[1]) = (-1)^{\spadesuit}\pi_*\big( \ev_k^*(p_1^* a_1\wedge \cdots p_k^* a_k)\big)
	    \end{equation}
	\end{defn}
	\begin{rem}
	    Note that if any of the $a_i$'s is of degree 0, then the iterated integral vanishes. This follows from the observation that the form $\ev_k^*(p_1^* a_1\wedge \cdots p_k^* a_k) \in \Omega^{\bt}(\CF \otimes \Delta^k)$ is annihilated by vector fields $\frac{\del}{\del t_i}$'s, $1 \le i \le k$, which forces the push forward along $\pi: P\CF \times \Delta_k \to \CF$ vanishing.
	\end{rem}
	
	\begin{lem}
	    $\mathsf{C}$ is natural, i.e. for any foliated map $f:(M, \CF_1) \to (N, \CF_2)$, the diagram
	    \begin{center}
        \begin{tikzcd}
        \mathsf{B} \big((\extp^{\bullet}\CF_1^{\vee})[1]\big) \arrow[r, "\mathsf{C}"]  & \Omega^{\bt}(P\CF_1)          \\
        \mathsf{B}\big((\extp^{\bullet}\CF_2^{\vee})[1]\big) \arrow[r, "\mathsf{C}"] \arrow[u, "\mathsf{B}f"]                           & \Omega^{\bt}(P\CF_2)   \arrow[u, "(Pf)^*"]
        \end{tikzcd}
        \end{center}
	\end{lem}
	\begin{proof}
	Since $f$ is foliated,
	\begin{align*}
	(Pf)^* \ev_k^*(p_1^* a_1\wedge \cdots p_k^* a_k)
	     &= ((f\otimes \id)\circ \ev_k)^*(p_1^* a_1\wedge \cdots p_k^* a_k)\\
	     &=\ev_k^* f^*(p_1^* a_1\wedge \cdots p_k^* a_k)\\
	    &= \ev_k^* ((p_1\circ f)^* a_1\wedge \cdots (p_k\circ f)^* a_k)
	\end{align*}
	\end{proof}

\begin{lem}\label{a4}
    Let $a_1[1]\otimes \cdots \otimes a_k[1]\in \mathsf{B} \big((\extp^{\bullet}\CF^{\vee})[1]\big) $ be an element of the bar complex, then we have 
    \begin{align}
        d(\mathsf{C}(a_1[1]\otimes \cdots \otimes a_k[1])) =& \mathsf{C}\big(\overline{D}(a_1[1]\otimes \cdots \otimes a_k])\big) + \ev_1^*(a_1)\wedge \mathsf{C}(a_2[1]\otimes \cdots \otimes a_k[1]) \\
        &- (-1)^{|a_1|+\cdots |a_{k-1}|}\mathsf{C}(a_1[1]\otimes \cdots \otimes a_{k-1}[1])\wedge\ev_0^*(a_n)
    \end{align}
    here $\overline{D}$ is the differential of the foliation dga $(\extp^{\bullet}\CF^{\vee}, -d, \wedge)$.
\end{lem}
\begin{proof}
    Note that by lemma, 
    \begin{align*}
        d(\mathsf{C}(a_1[1]\otimes \cdots \otimes a_k[1])) =& (-1)^{\spadesuit}\Big( (-1)^k(\pi_* d\big(\ev_k^* (p_1^* a_1\wedge \cdots p_k^* a_k) \big) \\
        &+ (-1)^{k+1}((\del\pi)_* (\id \otimes \iota)^*\big(\ev_k^* (p_1^* a_1\wedge \cdots p_k^* a_k)\Big)\\
        =&\sum_{i=1}^k(-1)^{|a_1|+\cdots |a_{i-1}|}\mathsf{C}(a_1[1]\otimes \cdots\otimes (-da_i)[1]\otimes \cdots \otimes a_k[1]) \\
        &+ \Big(\sum_{i=1}^{k - 1}(-1)^{|a_1|+\cdots |a_{i}|}\mathsf{C}(a_1[1]\otimes \cdots\otimes (a_i\wedge a_{i+1})[1]\otimes \cdots \otimes a_k[1])\\
        &+\ev_1^*(a_1)\wedge \mathsf{C}(a_2[1]\otimes \cdots \otimes a_k[1]) \\
        &- (-1)^{|a_1|+\cdots |a_{k-1}|}\mathsf{C}(a_1[1]\otimes \cdots \otimes a_{k-1}[1])\wedge\ev_0^*(a_n)\Big)
    \end{align*}
\end{proof}

Let $\cinf_{+,\del I}(I)$ be the space of differentiable maps from $I\to I$ which are monotonically increasing and fixing the boundary $\del I$.

\begin{defn}
    We call a differential form $\alpha \in \Omega^{\bt}(P\CF)$ is {\it reparametrization invariant} if $\alpha$ is invariant under any reparametrization $\phi \in \cinf_{+,\del I}(I)$, i.e.
    $$
    \phi^* \alpha = \alpha
    $$
    Denote the subcomplex of invariant forms by $\Omega^{\bt}_{\inv}(P\CF)$
\end{defn}

\begin{lem}
    The image's of Chen's map on foliation
    \begin{equation*}
         \mathsf{C}:\mathsf{B} \big((\extp^{\bullet}\CF_1^{\vee})[1]\big) \to   \Omega^{\bt}(P\CF_1)    
    \end{equation*}
    lies in $\Omega^{\bt}_{\inv}(P\CF)$
\end{lem}
\begin{proof}

\end{proof}
\subsection{Cube's to simplices}

In this section, we shall construct a map 
\begin{equation}
    \mathsf{S}: \Omega^{\bt}(P\CF) \to C^{\bt}(\CF)[1]
\end{equation}
which is based on Igusa's construction from cubes to simplices \cite{Igu09}.
Recall that in this chapter, we parametrize the $k$-simplex by
$$
\Delta^k = \{(t_1, \cdots, t_k)\in \R^k|1 \ge t_1\ge t_2 \cdots \ge t_k \ge 0\} \subset \R^k
$$
The coface maps $\del_i: \Delta^k \to \Delta^{k+1}$ are given by
\begin{equation}
    (t_1, \cdots, t_k) \mapsto
    \begin{cases}
        (1, t_1, \cdots, t_k) &\text{for $i = 0$}\\
        (t_1, \cdots, t_{i-1}, t_i, t_i, t_{i+1},\cdots,t_k) &\text{for $0<i < k+1$}\\
        (t_1, \cdots, t_k,0) &\text{for $i = k+1$}
    \end{cases}
\end{equation}
The codegeneracy maps $\epsilon_i: \Delta^k \to \Delta^{k-1}$ are given by
\begin{equation}
    (t_1, \cdots, t_k) \mapsto (t_1, \cdots,\hat{t_i},\cdots, t_k)
\end{equation}
The $i$-th vertex of $\Delta^k$ is the point
$$
(\underbrace{1, \cdots, 1}_\text{ $i$-times}, \underbrace{0, \cdots, 0}_\text{ $k-i$-times})
$$

Recall the smooth singular $\CF$-chains $C_{\bt}(\CF)$ is given by $C_{k}(\CF) = \cinf(\Delta^k, M)$. With structure map $d_i = \del_i^*, s_i = \epsilon_i^*$, we equip $C_{\bt}(\CF)$ a simplicial set structure, which is equivalent to the monodromy $\infty$-groupoids $\mon^{\infty} (\CF)$ of $\CF$

We  define  maps $P_i$ and $Q_i$ which send element of $\mon^{\infty} (\CF)$ to its back-face and front-face respectively, i.t. $P_i$ and $Q_i$ are pullbacks of 
\begin{align*}
    U_i:& \Delta^i \to \Delta^k, (t_1, \cdots, t_i) \mapsto (1, \cdots, 1, t_1, \cdots, t_i)\\
     V_i:& \Delta^i \to \Delta^k, (t_1, \cdots, t_i) \mapsto ( t_1, \cdots, t_i,0, \cdots, 0)
\end{align*}
respectively.
\begin{defn}
    Let $(M, \CF)$ be a foliated manifold, we define the dga of (smooth) singular $\CF$-cochains $(C^{\bt}(M), \delta, \cup)$ consisting of the following data:
    \begin{enumerate}
        \item The grade vector space $C^{\bt}(M)$ of linear functional on the vector space generated by $\mon^{\infty}(\CF)$.
        \item The differential $\delta$ is given by
        \begin{equation*}
            (\delta \phi)(\sigma) = \sum_{i=0}^k(d^*_i\phi)(\sigma)=\sum_{i=0}^k(\phi)(\del^*_i\sigma)
        \end{equation*}
        \item The product $\cup$ is given by the usual cup product
        \begin{equation*}
            (\phi \cup \psi)(\sigma) = \phi(V_i^*\sigma)\psi(U_j^*\sigma)
        \end{equation*}
    \end{enumerate}
\end{defn}
Define $\pi_k: I^k \to \Delta_k$ by the order preserving retraction, i.e $\pi_k(x_1,\cdots, x_k) = (t_1, \cdots, t_k)$ with $t_i = \max\{x_i, \cdots, x_k\}$ for each $k$.

Consider an element $\lambda_w: I\to I^{k}$ of $PI^k$ which is parametrized by a $w\in I^k$. In detail, if $w = (w_1, \cdots, w_{k-1})$, then $\lambda_k$ travels backwards through the $k+1$ points
\begin{equation*}
    0 \leftarrow w_1x_1 \leftarrow w_1 x_1 + w_2 x_2 \leftarrow \cdots \leftarrow \sum_{i = 1}^k w_i e_i
\end{equation*}
For more details, see \cite[Proposition 4.6]{Igu09}. 
Set $\lambda_{(k-1)}: I^{k-1} \to PI^k$ by sending $w$ to $\lambda_w$.

Finally, we define $\theta_{(k)}$ to be the composition
\begin{equation*}
    \theta_{(k)} = P\pi_k \circ \lambda _{(k-1)}: I^{k-1} \to P\Delta^k
\end{equation*}
We denote the adjoint of $\theta_{(k)}$ to be $\theta_{k}: I^k \to \Delta^k$.
\begin{rem}
    By construction, $\theta_{(k)}$ are piecewise linear but not smooth. We can correct it by reparametrization, for example, let the derivative vanish near the vertices. Since the image of Chen's map $\mathsf{C}$ is invariant, our construction for $\theta_{(k)}$ is well-defined.
\end{rem}
\begin{rem}
    It is easy to verify that 
    \begin{equation*}
        \int_{I^k} \theta^*_k \alpha = (-1)^k\int_{\Delta^k} \alpha
    \end{equation*}
    for any form $\alpha\in \Omega^{\bt}(\Delta^k)$.
\end{rem}

For each $i$, define  $\widehat{\del_i^-}$ to be the map which inserts a 0 between the $(i-1)$-th and $i$-th coordinates. Note that the $i$-th negative face operator is given by $\del_i^-(\theta_{(k)}) = \theta_{(k)}\circ \widehat{\del_i^-}$
\begin{lem}[\cite{Igu09}]\label{a1}
    For each $1 \le i \le k - 1$, we have the following commutative diagram
    \begin{center}
        \begin{tikzcd}
I^{k-2} \arrow[r, "\widehat{\del_i^-}"] \arrow[d, "\theta_{(k-1)}"] & I^{k-1} \arrow[r, "\theta_{(k)}"] & P(\Delta^k, v_k, v_0)   \\P(\Delta^{k-1}, v_{k-1}, v_0)     \arrow[rr, "\omega_i"]                         &             & P(\Delta^{k-1}, v_{k-1}, v_0)  \arrow[u, "P\del_i"]
\end{tikzcd}
    \end{center}
    that is,
$$
\del_i^-(\theta_{(k)}) = \theta_{(k)}\circ \widehat{\del_i^-} =P\widehat{\del_i^-} \circ \omega_i \circ \theta_{(k - 1)}
$$
Here $w_i$ is given by the following reparametrization:
for each $\gamma \in P(\Delta^{k-1}, v_{k-1}, v_0)$, $w_i(\gamma)$ is defined by
\begin{equation*}
    \omega_i(\gamma)(t)  =
    \begin{cases}
        \gamma(\frac{kt}{k - 1}) &\text{if $t\le \frac{j-1}{k}$}\\
        \gamma(\frac{j-1}{k - 1}) &\text{if $ \frac{j-1}{k}\le k \le \frac{k}{k}$} \\
        \gamma(\frac{kt - 1}{k - 1}) &\text{if $t\ge \frac{j}{k}$}
    \end{cases}
\end{equation*}
\end{lem}
\begin{proof}
    See \cite[Lemma 4.7]{Igu09}.
\end{proof}
Set $\widehat{\del_i^+}: I^{k-1} \to I^{k-1}$ to be the map which inserts 1 between the $(i-1)$-th and $i$-th places.
\begin{lem}[\cite{Igu09}]\label{a2}
    For each $1 \le i \le k - 1$, we have the following commutative diagram
    \begin{center}
        \begin{tikzcd}
I^{k-2} \arrow[r, "\widehat{\del_i^+}"] \arrow[d, "\simeq"] & I^{k-1} \arrow[r, "\theta_{(k)}"] & P(\Delta^k, v_k, v_0)   \\I^{i-1} \times I^{k-i-1}   \arrow[rr, "\theta_{(i)}\times \theta_{(k-i)}"]                         &             & P(\Delta^{i}, v_{i}, v_0) \times P(\Delta^{k-i}, v_{k-i}, v_0)  \arrow[u, "\mu_{i}"]
\end{tikzcd}
    \end{center}
    that is,
$$
\del_i^+(\theta_{(k)}) = \theta_{(k)}\circ \widehat{\del_i^+} =\mu_{i} \circ (\theta_{(i)}\times \theta_{(j)})
$$
where $\mu_{i,j}$ is the path composition map

\begin{equation*}
    \mu_{i}(\alpha, \beta)(t)  =
    \begin{cases}
        U_{k-i}\big(\beta(\frac{kt}{k - i})\big) &\text{if $ t\le \frac{k-i}{k}$}\\
        V_i\Big(\alpha\big(\frac{k}{i}(t-\frac{k-i}{k})\big)\Big) &\text{if $  k \ge \frac{k-i}{k}$}
    \end{cases}
\end{equation*}
\end{lem}
\begin{proof}
    See \cite[Lemma 4.8]{Igu09}.
\end{proof}

\begin{lem}[\cite{AS12}]\label{a3}
    Let $a_1,\cdots, a_n$ be forms on $\Delta^k$, then we have the following factorization
    \begin{align*}
        &\int_{P(\Delta^{i}, v_{i}, v_0) \times P(\Delta^{k-i}, v_{k-i}, v_0) } (\mu_i)^*\mathsf{C}(a_1[1]\otimes \cdots \otimes a_n) \\
        &= \sum_{l = 1}^n \Big(\int_{P(\Delta^{i}, v_{i}, v_0) } \mathsf{C}(V_i^* a_1[1]\otimes \cdots \otimes V_i^* a_l[1]) \Big)\\
        &\times \Big(\int_{P(\Delta^{k-i}, v_{k-i}, v_0) } \mathsf{C}(U_{k-i}^* a_{l+1}1[1]\otimes \cdots \otimes U_{k-i}^* a_n[1]) \Big)
    \end{align*}
\end{lem}
\begin{proof}
    See \cite[Lemma 3.19]{AS12}.
\end{proof}
We define the map $\mathsf{S}:\Omega^{\bt}(P\CF)\to  C^{\bt}(\CF)[1]$ to be
\begin{equation*}
    \mathsf{S}(\alpha) = \int_{I^{k-1}} (\theta_{(k)})^*P\sigma^*\alpha
\end{equation*}
for $\alpha \in \Omega^{\bt}(P\CF)$.

\subsection{$A_{\infty}$ de Rham theorem for foliation}
Next, we will prove the $A_{\infty}$-enhancement of the de Rham theorem for foliations.
\begin{thm}[$A_{\infty}$ de Rham theorem for foliation\index{de Rham theorem!$A_{\infty}$ for foliation }]
	Let $(M, \CF)$ be a foliated manifold, there exists an $A_{\infty}$-quasi-isomorphism between $\big(\Omega^{\bullet}(\CF), -d, \wedge \big)$ and $\big(C^{\bullet}(\CF), \delta, \cup \big)$
\end{thm}
We have already constructed the map
$$
\mathsf{S}\circ \mathsf{C}:\mathsf{B} \big((\extp^{\bullet}\CF_1^{\vee})[1]\big) \to   \Omega^{\bt}(P\CF_1)\to  C^{\bt}(\CF)[1]
$$

\begin{lem}\label{a5}
    Let $a_1, \cdots, a_n$ be $\CF$-foliated forms, then we have the following identity
    \begin{align*}
        \mathsf{S}(d(\mathsf{C}(a_1[1]\otimes \cdots \otimes a_n[1]))) =& \delta'(\mathsf{S}(\mathsf{C}(a_1[1]\otimes \cdots \otimes a_n[1]))) +\\
        &\sum_{l=1}^{n-1}\mathsf{S}((\mathsf{C}(a_1[1]\otimes \cdots \otimes a_l[1])))\cup' \mathsf{S}((\mathsf{C}(a_{l+1}[1]\otimes \cdots \otimes a_n[1])))
    \end{align*}
    Here $\delta'$ and $\cup'$ are differential and product of the dga of singular $\CF$-cochains at the level of suspensions.
\end{lem}
\begin{proof}
    We follow \cite[Proposition 3.22]{AS12}. Consider $\alpha = \mathsf{C}(a_1[1]\otimes\cdots\otimes a_n[1]) \in \Omega^{\bt}(P\CF_1)$, and $\sigma\in \mon^{\infty}(\CF)^k$ a simplex. 
    We want to compute 
    \begin{equation*}
        \int_{I^{k-1}}d (\theta_{(k)})^*P\sigma^*\alpha = \int_{\del I^{k-1}} \iota^* (\theta_{(k)})^*P\sigma^*\alpha
    \end{equation*}
    Recall $\widehat{\del^{\pm}_i}$ are the canonical embeddings of $I^{k-2}$ into $I^{k-1}$ as top and bottom faces. Then the right-hand side of the above equation breaks to 
    \begin{equation*}
        \sum_{i=1}^{k-1} (-1)^i\int_{I^{k-2}} (\del_i^-)^* (\theta_{(k)})^*P\sigma^*\alpha -\sum_{i=1}^{k-1} (-1)^i\int_{I^{k-2}} (\del_i^+
        )^* (\theta_{(k)})^*P\sigma^*\alpha
    \end{equation*}
    By Lemma \ref{a1} and properties of Chen's map, we have 
    \begin{align*}
        \int_{I^{k-2}} (\del_i^-)^* (\theta_{(k)})^*P\sigma^*\alpha =& \int_{I^{k-2}}  (\theta_{(k-1)})^*(P\del_i^*\sigma)^*\alpha
    \end{align*}
    On the other hand, by Lemma \ref{a2} and Lemma \ref{a3},
    \begin{align*}
        \int_{I^{k-2}}& (\del_i^+
        )^* (\theta_{(k)})^*P\sigma^*\alpha = \\
        &\sum_{l=0}^{n}\mathsf{S}((\mathsf{C}(a_1[1]\otimes \cdots \otimes a_l[1])))(V_i^*\sigma) \mathsf{S}((\mathsf{C}(a_{l+1}[1]\otimes \cdots \otimes a_n[1])))(U^*_{k-i \sigma})
    \end{align*}
    
    Summing up all the items yields the desired result.
\end{proof}

Now we describe our proposed $A_{\infty}$-map. Let $(M, \CF)$ be a foliated manifold, we define a series of maps $\phi_n: (\Omega^{\bt}(\CF)[1])^{\otimes n} \to C^{\bt}(\CF)[1]$ by
\begin{enumerate}
    \item For $n = 1$,
    \begin{equation*}
        (\phi_1(a[1])(\sigma) = (-1)^k \int_{\Delta^k}\sigma^* \alpha
    \end{equation*}
    \item For $n>1$,
    $$\phi_n(a_1[1]\otimes \cdots a_n[1]) = (\mathsf{S}\circ \mathsf{C})(a_1[1]\otimes \cdots a_n[1])
    $$
\end{enumerate}
Next we shall prove that $\phi_n$'s form an $A_{\infty}$-morphism. The case for $\CF = TM$ is proved by Guggenheim in \cite{Gug77}. We will follow the proof in \cite[Theorem 3.25]{AS12}.
\begin{prop}
$\phi_n$'s form an $A_{\infty}$-morphism from $\Omega^{\bt}(\CF)$ to $C^{\bt}(\CF)$ which induces a quasi-isomorphism. Moreover, this map is natural with respect to pullbacks along $\cinf$-maps.
\end{prop}
\begin{proof}
    Let $a_1[1]\otimes \cdots a_n[1] \in \mathsf{B} \big((\extp^{\bullet}\CF_1^{\vee})[1]\big)$.
    
    First consider the case $n\not = 2$. By lemma \ref{a4},
    \begin{align*}
        d(\mathsf{C}(a_1[1]\otimes \cdots \otimes a_k[1])) =& \mathsf{C}\big(\overline{D}(a_1[1]\otimes \cdots \otimes a_k])\big) + \ev_1^*(a_1)\wedge \mathsf{C}(a_2[1]\otimes \cdots \otimes a_k[1]) \\
        &- (-1)^{|a_1|+\cdots |a_{k-1}|}\mathsf{C}(a_1[1]\otimes \cdots \otimes a_{k-1}[1])\wedge\ev_0^*(a_n)
    \end{align*}
    By lemma \ref{a5},
    \begin{align*}
        \mathsf{S}(d(\mathsf{C}(a_1[1]\otimes \cdots \otimes a_n[1]))) =& \delta'(\mathsf{S}(\mathsf{C}(a_1[1]\otimes \cdots \otimes a_n[1]))) +\\
        &\sum_{l=1}^{n-1}\mathsf{S}((\mathsf{C}(a_1[1]\otimes \cdots \otimes a_l[1])))\cup' \mathsf{S}((\mathsf{C}(a_{l+1}[1]\otimes \cdots \otimes a_n[1])))
    \end{align*}
    Combining these two equations gives
    \begin{align*}
       (\mathsf{S}\circ \mathsf{C})\big(\overline{D}(a_1[1]\otimes \cdots \otimes a_n])\big)=& \delta'(\mathsf{S}(\mathsf{C}(a_1[1]\otimes \cdots \otimes a_n[1]))) \\
        &+\sum_{l=1}^{n-1}\mathsf{S}((\mathsf{C}(a_1[1]\otimes \cdots \otimes a_l[1])))\cup' \mathsf{S}((\mathsf{C}(a_{l+1}[1]\otimes \cdots \otimes a_n[1])))\\
        &-\mathsf{S}\big(\ev_1^*(a_1)\wedge \mathsf{C}(a_2[1]\otimes \cdots \otimes a_n[1])\big) \\
        &+ (-1)^{|a_1|+\cdots |a_{n-1}|}\mathsf{S}\big(\mathsf{C}(a_1[1]\otimes \cdots \otimes a_{n-1}[1])\wedge\ev_0^*(a_n)\big)
    \end{align*}
    The third term
    $$
    \mathsf{S}\big(\ev_1^*(a_1)\wedge \mathsf{C}(a_2[1]\otimes \cdots \otimes a_n[1])\big) = - \phi_1(a_1[1])\cup' (\mathsf{S}\circ \mathsf{C})(a_2[1]\otimes \cdots \otimes a_n[1])
    $$
    for $|a_1| = 0$.
    The fourth term
    \begin{align*}
        \mathsf{S}&\big(\mathsf{C}(a_1[1]\otimes \cdots \otimes a_{n-1}[1])\wedge\ev_0^*(a_n)\big) = \\
        & (-1)^{|a_1|+\cdots |a_{k-1}|}(\mathsf{S}\circ \mathsf{C})(a_1[1]\otimes \cdots \otimes a_{n-1}[1])\cup' \phi_1(a_n[1])
    \end{align*}
    for $|a_n| = 0$. These two terms vanish for $|a_1| > 0$ and $|a_n| >0$ respectively.
    
    Therefore, putting everything together, we have
    \begin{align*}
        (\mathsf{S}\circ \mathsf{C})\big(\overline{D}(a_1[1]\otimes \cdots \otimes a_n])\big)=& \delta'(\phi_n(a_1[1]\otimes \cdots \otimes a_n[1])) \\
        &+\sum_{l=1}^{n-1}\phi_l(a_1[1]\otimes \cdots \otimes a_l[1]))\cup' \phi_{n-l}(a_{l+1}[1]\otimes \cdots \otimes a_n[1])
    \end{align*}
    On the other hand, by definition
    \begin{align*}
        (\mathsf{S}&\circ \mathsf{C})\big(\overline{D}(a_1[1]\otimes \cdots \otimes a_n])\big)=\\
        &\sum_{i=1}^n (-1)^{|a_1|+\cdots |a_{i-1}|}\phi_n(a_1[1]\otimes\cdots\otimes a_{i-1}[1] \otimes (-da_i)[1]\otimes a_{i+1}[1]\otimes\cdots\otimes a_n[1])\\
        &+\sum_{i=1}^{n-1} (-1)^{|a_1|+\cdots |a_{i}|}\phi_{n-1}(a_1[1]\otimes\cdots\otimes a_{i-1}[1] \otimes (a_i\wedge a_{i+1})[1]\otimes a_{i+2}[1]\otimes\cdots\otimes a_n[1])
    \end{align*}
    Combining these two equations yields the desired $A_{\infty}$-structure maps.
    
    For $n=2$ and $|a_1| = |a_2| = 0$. Just noted that for two foliated functions, for $\phi$ to be an $A_{\infty}$-map, we only need to check $(a_1a_2)(x) = a_1(x)a_2(x)$.
    
    The quasi-isomorphism follows from the ordinary de Rham theorem for foliations (Theorem \ref{dr}). The naturality follows from the naturality of the maps $\mathsf{S}$ and $\mathsf{C}$.
    
\end{proof}

\begin{rem}
    It is easy to verified that, according to the construction,
    $\phi_1(f[1]) = f[1]$ for any $|f| = 0$, and $\phi_n(a_1[1]\otimes \cdots \otimes a_n[1])$ vanishes if any of the $a_i[1]$ in the argument is of degree 0.
\end{rem}

\begin{lem}
    The image of $\phi_n$'s lies in the dga of normalized $\CF$-cochains.                  
\end{lem}
\begin{proof}
Follows from \cite[Proposition 3.26]{AS12}. Note that by our construction of $\phi_n$, we just need to restricted to leaves.
\end{proof}

\section{Riemann-Hilbert correspondence for $\infty$-foliated local systems}
\subsection{Iterated integrals on vector bundles}
In this section, we shall generalize iterated integrals\index{iterated integrals! for vector bundles} in the previous section to the case of graded vector bundles (or dg modules).
Let $V$ be a graded vector bundle on $M$, denote 
\begin{align*}
   \iota: \big(\
\Gamma(\End(V)\otimes \bigwedge\nolimits^{\bullet}\CF^{\vee}) \big)^{\otimes_{\R}^k} &\to \Gamma \big(\End(V)^{\boxtimes^k}\otimes 
(\extp^{\bullet}\CF^{\vee})^{\boxtimes^k} \big) \\
a_1\otimes \cdots \otimes a_n &\mapsto a_1\boxtimes \cdots \boxtimes a_n
\end{align*}
 the canonical embedding.

 The pull back of $\ev_k$ induces 
$$\ev_k^*: \Gamma \big(\End(V)^{\boxtimes^k}\otimes 
(\extp^{\bullet}\CF^{\vee})^{\boxtimes^k} \big) \to \Gamma\big(\boxtimes_i \ev_1^*\End(V)_{t_i} \otimes (\extp^{\bullet}T^{\vee}_{\CF}P\CF)\times \Delta^k \big)
$$
Let $\mu$ denote the multiplication map on $\ev_1^*\End(V)_{t_i}$, i.e. 
$$
\mu: \ev_k^* \Gamma \big(\End(V)^{\boxtimes^k}\otimes 
(\extp^{\bullet}\CF^{\vee})^{\boxtimes^k} \big) \to \Gamma\big( p_0^*\End(V) \otimes (\extp^{\bullet}T^{\vee}_{\CF}P\CF)\times \Delta^k \big)
$$
where $p_0: P\CF \to M$ is the evaluation map at $t=0$.  Denote $\pi$ the projection map $\pi: P\CF \times \Delta^k \to P\CF$.
\begin{defn}
	We define the iterative integral  $$\int : \big(\
	\Gamma(\End(V)\otimes \bigwedge\nolimits^{\bullet}\CF^{\vee}) \big)^{\otimes_{\R}^k} \to \Gamma\big( p_0^*\End(V) \otimes (\extp^{\bullet}T^{\vee}_{\CF}P\CF) \big)$$
	on graded vector bundles $V$ over a foliation $\CF$ to be the composition 
	\begin{equation}
	\int a_1\otimes a_2 \otimes\cdots\otimes a_k = (-1)^{\spadesuit}\pi_{\star}\circ \mu\circ \ev_l^* \circ \iota (a_1\otimes a_2 \otimes\cdots\otimes a_k)
	\end{equation}
	with $\spadesuit = \sum_{1 \le i < k} (T(a_i) - 1)(k-i)$ where $T(a_i)$ denotes the total degree of $a_i$.
\end{defn}

\begin{lem}
    On $\Gamma(\End(V)\otimes \bigwedge\nolimits^{\bullet}\CF^{\vee}) \big)$, we have
     \begin{equation*}
           \pi_* \circ d - (-1)^k d\circ \pi_* = (\del \pi)_* \circ (\id \times \iota)^*
       \end{equation*}
       Let $\alpha \in\Gamma(\End(V)\otimes \extp^{\bullet}T^{\vee}_{\CF}P\CF) \big)$, $\beta\in \Gamma\big( \End(V) \otimes (\extp^{\bullet}T^{\vee}_{\CF}P\CF)\times \Delta^k \big)$,
       \begin{align*}
           \pi_*(\pi^* \alpha \circ \beta) &= (-1)^{k T(\alpha) } \alpha \circ \pi_* \beta\\
           \pi_*( \circ \beta\pi^* \alpha ) &=  \pi_* \beta \circ \alpha
       \end{align*}
\end{lem}
\begin{proof}
    Similar to Lemma \ref{a6}.
\end{proof}

\begin{lem}[Stoke's theorem\index{Stoke's theorem}]
\begin{align}
	d\int \omega_1\cdots \omega_r =& \sum_{i=1}^{r} (-1)^i T \omega_1 \cdots d \omega_i \omega_{i+1}\cdots \omega_r + \sum_{i=1}^{r-1}(-1)^i T \omega_1\cdots (T\omega_i \circ \omega_{i+1})\cdots \omega_r\\
	&+ p_1^* \omega_1 \circ \int \omega_2\cdots \omega_r - T ( \int \omega_1\cdots \omega_{r-1} )\circ p_0^* \omega_r
\end{align}
\end{lem}
\begin{proof}
    Similar to Lemma \ref{a4}. See also \cite[Proposition 3.3]{BS14}.
\end{proof}
\subsection{$\infty$-holonomy of $\Z$-connection over $\CF$
}

Let $V$ be a $\Z$-graded vector bundle with a $\Z$-connection $\nabla$ over $\A^{\bullet} = \extp^{\bullet}\CF^{\vee}$. Locally, $\nabla = d - \sum_{i=0}^m A_i$, where $A_i \in \End^{1-i}(V) \otimes_{\A_0}\A^i$. Let $\omega = \sum_{i=0}^m A_i$. We define $p$-th holonomy of $\nabla$ to be the iterative integral 
\begin{equation}
\Psi_p = \int \omega^{\otimes p} \in \Gamma \big((\extp^{k}T^{\vee}_{\CF}P\CF)\otimes \End^{-k}(V)\big)
\end{equation}
and $\Psi_0 = \id$ for $p=0$.
\begin{defn}
	Define the {\it $\infty$-holonomy}\index{holonomy!$\infty$-holonomy} associated to $\nabla$ to be $\Psi = \sum_{p=0}^{\infty} \Psi_p$.
\end{defn}
Since $\omega$ has total degree 1, 

\begin{align*}
	\Psi_p =& \sum_{i=1, j = p-i}^{p}(-1)^{i+1}\int(\omega^{\otimes i})d\omega (\omega^{\otimes j} ) + \sum_{i=1, j = p-i-1}^{p-1}(-1)^i\int (\omega^{\otimes i})(\omega_i \circ \omega_{i+1})(\omega^{\otimes j})\\
	&+ p_1^*\omega\circ \int \omega^{\otimes (r-1)} - \big(\int \omega^{\otimes r-1}\big)\circ p_0^*\omega^{}
	\end{align*}

Summing in $p$, we get 
\begin{align*}
	d \Psi = \bigg(\int \kappa + \big(\int \kappa \omega +- \int \omega \kappa\big) + \cdots + \sum_{i+j =p-1}(-1)^i\int \omega^i \kappa \omega^j + \cdots \bigg) + p_1^*\omega\circ \Psi - \Psi\circ p_0^*\omega^{}
	\end{align*}

If $\nabla$ is flat, then locally $\nabla^2 = (d -\omega)^2 = -d\omega - T\omega\circ \omega = -d\omega + \omega\circ \omega$.

Let $\sigma: \Delta^k \to \CF$ be a foliated simplex. We can regard it as a $k-1$-family of paths into $\CF$. We can break this into two parts. First we have a map $\theta_{(k-1)}: I_{k-1} \to P\Delta^k_{(v_k, v_0)}$, then there is a canonical map $P\sigma: P\Delta^k_{(v_k, v_0)} \to P\CF_{(x_k, x_0)}$. We define a series of map $\psi_k \in \End^{1-k}(V)$ by
$$\psi_k(\sigma) =
\begin{cases} \int_{I^{k-1}} (-1)^{(k-1)(K\Psi)} \theta^*_{(k-1)} (P\sigma)^* \Psi &\quad k \ge 1\\
(V_x, \nabla^0_x) &\quad k =0 
\end{cases}
$$
which is essentially the integral of $I^{k-1}$ of the pullback holonomy of the $\Z$-connection $\nabla$.

Now we define the Riemann-Hilbert functor $\RH: \mathcal{P}_{\A} \to \repi(\mon_{\infty}\CF)$. On objects we define $\RH_0 : \ob(\mathcal{P}_{\A}) \to \ob(\repi(\mon_{\infty}\CF))$ by $\RH_0\big((E^{\bullet}, \nabla) \big)(\sigma_k) = \psi_k(\sigma_k)$. We claim that the image of this functor are $\infty$-local systems. Note that, by our construction
\begin{align*}
    \RH_0\big((E^{\bullet}, \nabla) \big)_x =& E_x\\
    \RH_0\big((E^{\bullet}, \nabla) \big)(x) =& \mathbb{E}^0_x\\
    \RH_0\big((E^{\bullet}, \nabla) \big)(\sigma_{k>0}) =& \int_{I^{k-1}} (-1)^{(k-1)(K\Psi)} \theta^*_{(k-1)} (P\sigma)^* \Psi
\end{align*}

Write $F$ the image of $\RH_0\big((E^{\bullet})), \nabla) \big)$ for simplicity, i.e. $F(\sigma_k) = \RH_0\big((E^{\bullet}))(\sigma_k)$. Since $\mathbb{E}$ is flat, we have 
$$
d\Psi = - p_0^* A^0 \circ \Psi + \Psi \circ p_1^*A^0
$$
Integrate the left side and apply the Stoke's formula we get
\begin{equation*}
    -\hat{\delta}F - \sum_{i=1}^{k-1}(-1)^i F(\sigma_{0\cdots i})F(\sigma_{i\cdots k})
\end{equation*}
Plug in the integration of right side, we get
$$
\mathbb{E}^0\circ F(\sigma_k) - (-1)^k F(\sigma_k) - \sum_{i=1}^{k-1}(-1)^iF(\sigma_{(0\cdots \hat{i}\cdots k)})+\sum_{i=1}^{k-1}(-1)^i F(\sigma_{0\cdots i})F(\sigma_{i\cdots k})=0
$$
which is the $k$-th level of the Maurer-Cartan equations for $\infty$-local system condition. Therefore, $\RH_0$ is a well-defined map on objects.
\begin{thm}[\cite{Igu09}]
	The  image of an object under the functor $\RH$ is an $\infty$-representations of $\mon_{\infty}(\CF)$ if and only if $\nabla$ is flat.
\end{thm}

\begin{proof}
	By Theorem 4.10 in \cite{Igu09}, we have
	$$
	\psi_0(x_0)\phi_k(\sigma) + (-1)^k\psi_k(\sigma)\psi(x_k) = \sum_{i=1}^{k-1}(-1)^i\big(\psi_{k-1}(\sigma_{(0\cdots \hat{i}\cdots k)}) - \psi_i(\sigma_{(0\dots i)})\psi_{k-i}(\sigma_{(i\cdots k)})\big)
	$$
	which is equivalent to 
	$$
	 \sum_{i=1}^{k-1}(-1)^i \psi_{k-1}(\sigma_{(0\cdots \hat{i}\cdots k)}) - \sum_{i=0}^k  (-1)^i \psi_i(\sigma_{(0\dots i)})\psi_{k-i}(\sigma_{(i\cdots k)})
	$$
	i.e. $\hat{\delta}\psi + \psi\cup\psi = 0$.  
	
	For the other direction, we just go back from the definition of $\psi$, and found that $A_0 \Psi_k - \Psi_k A_0 = d\Psi_{k-1}$ must be equal for all $k$, which is equivalent to the flatness of $\nabla$.
	\end{proof}
	
	Now we proceed to $\RH$ on higher simplices 
$$
\RH_n: \mathcal{P}_{\A}(E_{n-1}^{\bullet}, E_n^{\bullet}) \otimes \cdots \otimes \mathcal{P}_{\A}(E_{0}^{\bullet}, E_1^{\bullet}) \to \repi(\mon_{\infty}\CF)(\RH_0(E_0^{\bullet}), \RH_0(E_n^{\bullet}))[1-n]
$$
by
\begin{equation}
\RH_n (\phi_n \otimes \cdots \otimes \phi_1)(\sigma_k) = \RH_0\big(C_{(\phi_n \otimes \cdots \otimes \phi_1)}(\sigma_k) \big)_{n+1,1}
\end{equation}
Next, we will need $\infty$-holonomy with respect to the pre-triangulated structure of $\mathcal{P}_{\CA}$. We follow the calculation in \cite[Section 3.5]{BS14}
of the following 
\begin{itemize}
    \item $\infty$-holonomy with respect to the shift.
    \item $\infty$-holonomy with the cone
\end{itemize}

\begin{prop}
$\RH$ is an $A_{\infty}$-functor.
\end{prop}
\begin{proof}
    We follow \cite[Theorem 4.2]{BS14}. Let $\phi = \phi_n\otimes \cdots \otimes \phi_1 \in \mathcal{P}_{\A}(E_{n-1}, E_n)\otimes\cdots \otimes \mathcal{P}_{\A}(E_{0}, E_1)$ be a tuple of morphisms, denote the holonomy of the associated to the generalized homological cone $C_{\phi}$ by $\Psi^{\phi_n\otimes \cdots \otimes \phi_1}$. Locally, we can write $D^{\phi} = d - \omega$. By ..., we have that on $P\CF(x_0, x_1)$, where $x_0, x_1$ lie in some leaf. By the $\infty$-holonomy for cones, we have
    
    \begin{align*}
        -& d\Psi^{\phi_n\otimes \cdots \otimes \phi_1}_{n+1, 1} - p_0^*\omega^0_{n+1, 1} \circ  \Psi^{\phi_n\otimes \cdots \otimes \phi_1}_{n+1, 1}  + \Psi^{\phi_n\otimes \cdots \otimes \phi_1}_{n+1, 1} \circ p_1^*\omega^0_{n+1, 1} =\\
        &\sum_{k=1}^{n-1}(-1)^{n - k - 1 - |\phi_n\otimes \cdots \otimes \phi_{k+2}|} \Psi^{\phi_n\otimes \cdots\otimes \phi_{k+1}\circ \phi_k\otimes \cdots \otimes \phi_1}_{n+1, 1} + \\
        &\sum_{k=1}^{n}(-1)^{n - k  - |\phi_n\otimes \cdots \otimes \phi_{k+1}|} \Psi^{\phi_n\otimes \cdots\otimes d\phi_k\otimes \cdots \otimes \phi_1}_{n+1, 1}
    \end{align*}
    Now applying $\int (-1)^{K(\Psi J(-)} \theta^*(P[-])^*(\Psi)$ to both sides of the equation. For simplicity, we denote $\phi_k\otimes \cdots \otimes \phi_l$ by $\phi_{k,l}$. We have
    \begin{align*}
        [& RH_0 (C_{\phi}) \cup RH_0(C_{\phi}) + \hat{\delta} \RH_0(C_{\phi})]_{n+1, 1} =\\
        &\sum_{k=1}^{n-1}(-1)^{n - k - 1 - |\phi_{n,k+2}|}\RH_{n-1}(\phi_n\otimes \cdots\otimes \phi_{k+1}\circ \phi_k\otimes \cdots \otimes \phi_1)\\
        &+ \sum_{k=1}^{n-1}(-1)^{n - k - |\phi_{n,k+1}|}\RH_{n}(\phi_n\otimes \cdots\otimes  d\phi_k\otimes \cdots \otimes \phi_1)
    \end{align*}
    By the matrix decomposition formulas in calculating the holonomy of cones in \cite[Section 3.5.2]{AS14}, we get 
    \begin{align*}
        [ RH_0 (C_{\phi}) \cup RH_0(C_{\phi})]_{n+1, 1} =& \sum_{i+j = n}\RH_0(C_{\phi_{n, i+1}})_{j+1, 1}\cup \RH_0(C_{\phi_{i, 1}})_{i+1, 1}[(j - \sum_{k = i+1}^n p_k)]\\
        &+ \RH_0(E_n)\cup\RH_0(C_{\phi})_{n+1, 1} + \RH_0(C_{\phi})_{n+1, 1}\cup \RH_0(E_0[n-|\phi|])\\
        =&\sum_{i+j = n}(-1)^{(j - \sum_{k = i+1}^n p_k)}\RH_j(C_{\phi_{n, i+1}})\cup \RH_i(C_{\phi_{i, 1}})\\
        &+ \RH_0(E_n)\cup\RH_0(C_{\phi})_{n+1, 1} + (-1)^{n-|\phi|}\RH_0(C_{\phi})_{n+1, 1}\cup \RH_0(E_0)\\
    \end{align*}
    By our construction,
    \begin{align*}
        D_{\Loc^{\dg}_{\CC}(K_{\bt})}\big(\RH_n(\phi)\big) =& \hat{\delta} \big(\RH_n(\phi)\big) + \RH_0(E_n)\cup \RH_n(\phi) \\
        &+(-1)^{n-|\phi|}\RH_n(\phi)\cup \RH_0(E_0)
    \end{align*}
    Put the last two equations into the one above, we get
    \begin{align*}
        & \Big(\sum_{i+j = n}(-1)^{(j - \sum_{k = i+1}^n p_k)}\RH_j(C_{\phi_{n, i+1}})\cup \RH_i(C_{\phi_{i, 1}})  \Big) + D\big(\RH_n(\phi_{n,1})\big)\\
       &=\sum_{k=1}^{n-1}(-1)^{n - k - 1 - |\phi_{n,k+2}|}\RH_{n-1}(\phi_n\otimes \cdots\otimes \phi_{k+1}\circ \phi_k\otimes \cdots \otimes \phi_1)\\
        & \quad + \sum_{k=1}^{n-1}(-1)^{n - k - |\phi_{n,k+1}|}\RH_{n}(\phi_n\otimes \cdots\otimes  d\phi_k\otimes \cdots \otimes \phi_1)
    \end{align*}
    which is the $A_{\infty}$-relation for an $A_{\infty}$-functor between two dg-categories. Therefore, $\RH$ is an $A_{\infty}$-functor.
    \end{proof}

\subsection{Riemann-Hilbert correspondence}
\begin{thm}
	The functor $\RH$ is an $A_{\infty}$-quasi-equivalence.
\end{thm}
First, we want to show $\RH$ is {\it $A_{\infty}$-quasi-fully-faithful}. Consider two objects $(E_1{\bt}, \mathbb{E}_1)$, $(E_1{\bt}, \mathbb{E}_1) \in  \mathcal{P}_{\A}$. The chain map
\begin{equation*}
    \RH_1:  \mathcal{P}_{\A}(E_1, E_2) \to \Loc^{\infty}_{\Ch_k}(\CF)(\RH_0(E_1), \RH_0(E_2))
\end{equation*}
induces a map on the spectral sequence.  In $E_1$-page, on $\mathcal{P}_{\A}$ side, $H^{\bt}\big((E_i, \mathbb{E}_i^0)\big)$ are vector bundles with flat connections, while on the other side $H^{\bt}\big((\RH(E_i), \mathbb{E}_i^0)\big)$ are $\CF$-local systems. In $E_2$-page, the map is 
\begin{align*}
    H^{\bt}\bigg(M, \Hom\Big( H^{\bt}\big((E_1, \mathbb{E}_1^0)\big), H^{\bt}\big((E_2, \mathbb{E}_2^0)\big)\Big)\bigg) \to\\
    H^{\bt}\bigg(M, \Hom\Big( H^{\bt}\big((\RH(E_1), \mathbb{E}_1^0)\big), H^{\bt}\big((\RH(E_2), \mathbb{E}_2^0)\big)\Big)\bigg)
\end{align*}
which is an isomorphism by the de Rham theorem for foliated local systems.
Next, we shall prove that  $\RH$ is {\it $A_{\infty}$-essentially surjective}. 

Let $F\in \Loc^{\dg}_{\CC}(\CF)$, we want to construct an object $(E^{\bullet}, \nabla) \in \mathcal{P}_{\A}$ whose image under $\RH_0$ is quasi-isomorphic to $F$. First notice that $\underline{\R}_{\CF}$ defines a representation of $\mon_{\infty}(\CF)$ by previous section, which can be viewed as $\infty$-local system over $\CF$. Regard $(M, \underline{\R}_{\CF})$ as a ringed space. Construct a complex of sheaves $(\underline{C}_{F}^{\bullet}, D)$ by
$$
\big(\underline{C}_{F}(U)^{\bullet}, D(U)\big) = \Loc^{\infty}_{\Ch_k}(\CF)\big(\underline{\R}_{\CF}|_U, F
|_U\big)^{\bullet}
$$
We claim that $\underline{C}_{F}^{\bullet}$ is soft. First notice that for $i>0$, $\underline{C}_{F}^{i}$ is a $\underline{C}_{F}^{0}$-module by cup products on open sets. By definition, $\underline{C}_{F}^{0}=\{\phi: (\Loc^{\infty}_{\Ch_k}(\CF))_0 \to \Ch_{\R}^0| \phi(x)\in \Ch_{\R}^0 \big(\underline{\R}_{F}(x), F (x)\big)  \}$ which is a sheaf of discontinuous sections, hence soft. Therefore, all $\underline{C}_{F}^{i}$'s are soft. Recall that for two $\infty$-representation of a Lie $\infty$-groupoid, the $E_1$ term of the spectral sequence is an ordinary representation. Hence, $\underline{C}_{F}^{\bullet}$ is a perfect complex of sheaves. Let $\underline{\A}^{\bullet}$ be the sheaf of  $\cinf$ sections of $\A^{\bullet} = \extp^{\bullet} \CF^{\vee}$.  $\underline{\A}^0 = {\cinf(M)}$ is flat over $\underline{\R}_{\CF}$ as $\cinf$-rings since locally the module of smooth functions on $M$ are $\cinf(\R^n)$ and the foliated functions are $\cinf(\R^{n-q})$ where $q = \codim \CF$, and  $\cinf(\R^{n-q}) \otimes_{\infty} \cinf(\R^{q}) \simeq \cinf(\R^{n})$ where $\otimes_{\infty}$ is the tensor product for $\cinf$-rings. Therefore, $ \underline{C}_{F}^{\infty} = \underline{C}_{F}^{\bullet} \otimes_{\underline{\R}_{\CF}} \underline{\A}^0$ is a sheaf of perfect $\A^0$-modules.

Again by the flatness of $\underline{\A}^0$ is flat over $\underline{\R}_{\CF}$. We have a quasi-isomorphism $(\underline{C}_{F}^{\bullet}, D) \simeq (\underline{C}_{F}^{\infty}\otimes_{\underline{\A}^0} \underline{\A}^{\bullet}, D\otimes 1 + 1 \otimes d)$. We need the following proposition from Proposition 2.3.2, Expos\`{e} II, SGA6, \cite{Ber06}. 
\begin{rem}
    One of the core tool in previous proof is the flatness of $\underline{\A}^0$  over $\underline{\R}_{\CF}$, or, in other words, 
    the flatness of $\underline{\A}^0$  over $H^0(\underline{\A})$. We then expect a natural extension of our results to arbitrary $L_{\infty}$-algebroids $\g$ with associated foliation dga $A$, and $A^0$ is flat over $H^0(A)$. A natural question will be, given any dga $A$ (which presents some geometric object), when is $A^0$ flat over $H^0$? Or we can consider an even more generalization, given a map of sheaves of algebras
    $$d:A^0 \to A^1$$,
    when $A^0$ is flat over $H^0(A) = \ker(d)$?
    We believe that this question is related to a more general phenomenon in noncommutative geometry. 
\end{rem}
\begin{prop}
Let $(X, \underline{\mathcal{S}}_X)$ be a ringed space, where $X$ is compact and $\underline{\mathcal{S}}_X$ is a soft sheaf of rings. Then
\begin{enumerate}
    \item The global section functor
    $$ \Gamma: \Mod_{\underline{\mathcal{S}}_X} \to \Mod_{\underline{\mathcal{S}}_X(X)}
    $$ is exact and establishes an equivalence of categories between the category of sheaves of right $\underline{\mathcal{S}}_X$-modules and the category of right modules over the global sections $\underline{\mathcal{S}}_X(X)$ of $\underline{\mathcal{S}}_X$.
    \item If $\underline{F} \in \Mod_{\underline{\mathcal{S}}_X}$ locally has finite resolutions by finitely generated free $\underline{\mathcal{S}}_X$-modules, then $\Gamma(X, \underline{F})$ has a finite resolution by finitely generated projective modules.
    \item The derived category of perfect complexes of sheaves $\dperf(\Mod_{\underline{\mathcal{S}}_X})$ is equivalent to the derived category of perfect complexes of modules $\dperf(\Mod_{\underline{\mathcal{S}}_X(X)})$.
    
\end{enumerate}
\end{prop}

By this theorem, there is a (strict) perfect complex of $\A^0$-modules $(E, \mathbb{E}^0)$ and a quasi-isomorphism $e^0:(E^{\bt}, \mathbb{E}^0) \to (F^{\bt}, \mathbb{F}^0) = (\Gamma(M, \underline{C}_{F}^{\infty}), D)$. We shall follow the argument of Theorem 3.2.7 of \cite{Blo05} to construct the higher components $\mathbb{E}^i$ of $\Z$-connection along with the higher components of a morphism $e^i$.

On $F^{\bt}$, we have a $\Z$-connection $$\mathbb{F} = D\otimes1 + 1\otimes d: F^{\bt} \to F^{\bt}\otimes_{\A^0} \A^{\bt}
$$. The idea is to transfer this $\Z$-connection to $E^{\bt}$ which is compatible with the quasi-isomorphism on $H^0$'s. Note that we have an induced connection 
$$
\mathbb{H}^k: H^k(F^{\bt}, \mathbb{F}^0) \to H^k(F^{\bt}, \mathbb{F}^0)\otimes_{\A^0} \A^{1}
$$
for each $k$. First we will transfer this connection to a connection on $H^k(E^{\bt}, \mathbb{E}^0)$, and we have the following commutative diagram
\begin{center}
	\begin{tikzcd}
	H^k(E^{\bt}, \mathbb{E}^0) \arrow[r, "\mathbb{H}^k"] \arrow[d, "e^0"] & H^k(E^{\bt}, \mathbb{E}^0) \otimes_{\A^0} \A^{1}
	\arrow[d, "e^0\otimes 1"] \\
	H^k(F^{\bt}, \mathbb{F}^0)\arrow[r, "\mathbb{H}^k"]    &  H^k(F^{\bt}, \mathbb{F}^0) \otimes_{\A^0} \A^{1}               
	\end{tikzcd}
\end{center}
Note that $e^0\otimes 1$ is a quasi-isomorphism since $\A^{\bt}$ is flat over $\A^0$. We need the following lemma.
\begin{lem}
    Given a bounded complex of finitely generated projective $\A^0$-modules $(E^{\bt}, \mathbb{E}^0)$ with connections $\mathbb{H}^k: H^k(F^{\bt}, \mathbb{F}^0) \to H^k(F^{\bt}, \mathbb{F}^0)\otimes_{\A^0} \A^{1}$ for each $k$, there exists connections
    $$
        \tilde{\mathbb{H}}^k: E^k \to E^k\otimes_{\A^0} \A^{1}
    $$
    lifting $\mathbb{H}^k$, i.e. 
    $$
    \tilde{\mathbb{H}}^k\mathbb{E}^0 = (\mathbb{E}^0\otimes 1)\tilde{\mathbb{H}}^k
    $$
    for each $k$
    and the connection induced on the cohomology is $\mathbb{H}^k$. 
\end{lem}
\begin{proof}
This is Lemma 3.2.8 in \cite{Blo05} and Lemma 4.6 in \cite{BS14}. Since $E^{\bt}$ is bounded, let $[N, M]$ be its magnitude. Pick some arbitrary connection $\nabla$ on $E^M$. Consider the following diagram whose rows are exact

\begin{center}
\begin{tikzcd}
E^M \arrow[r, "j"] \arrow[rd, "\theta"] \arrow[d, "\nabla"] & H^M(E^{\bt}, \mathbb{E}^0) \arrow[d, "\mathbb{H}^M"] \arrow[r] & 0          \\
E^M\otimes_{\A^0} \A^{1} \arrow[r, "j\otimes 1"]                            & H^M(E^{\bt}, \mathbb{E}^0)\otimes_{\A^0} \A^{1} \arrow[r] & 0
\end{tikzcd}
\end{center}
Here $\theta = \mathbb{H}^M \circ j - (j\otimes 1)\circ \nabla$ is $\A^0$-linear and $j\otimes 1$ is surjective. Now, by the projectivity of $E^M$, we can lift $\theta$ to a $\tilde{\theta}: E^M \to E^M \otimes_{\A^0} \A^{1}$ such that $(j\otimes 1)\tilde{\theta} = \theta$. Let $\Tilde{\mathbb{H}}^k = \nabla + \Tilde{\theta}$. Now replace $\nabla$ by $\Tilde{\mathbb{H}}^k$ and the above diagram still commutes.

Now pick some arbitrary connection $\nabla_{M-1}$ on $E^{M-1}$. Note that $\mathbb{E}^0\nabla_{M-1} = {\mathbb{H}}^{M-1}\mathbb{E}^0 = 0$ does not necessarily hold. Set $\mu = {\mathbb{H}}^{M-1}\mathbb{E}^0  - \mathbb{E}^0\nabla_{M-1}$, then $\mu$ is $\A^0$-linear. $\im \mu \subset \im \mathbb{E}^0 \otimes 1$ since $\im {\mathbb{H}}^{M-1}\mathbb{E}^0 \subset  \im \mathbb{E}^0 \otimes 1$ as ${\mathbb{H}}$ is a lift of $\mathbb{H}$. Now by the projectivity we can lift this map to a $\tilde{\mu}: E^{M-1} \to E^{M-1} \otimes_{\A^0} \A^{1}$ such that $(\mathbb{E}^0 \otimes 1)\circ \tilde{\mu} = \mu$. Now set $\Tilde{\mathbb{H}}^{M-1} = \nabla_{M-1} + \Tilde{\mu}$, then $(\mathbb{E}^0\otimes 1)\Tilde{\mathbb{H}}^{M-1} = \Tilde{\mathbb{H}}^{M-1}\mathbb{E}^0$. We have the following diagram
\begin{center}
    \begin{tikzcd}
E^N \arrow[r, "\mathbb{E}^0"]                                & E^{N+1} \arrow[r, "\mathbb{E}^0"]                                & \cdots \arrow[r, "\mathbb{E}^0"]          & E^{M-1} \arrow[r, "\mathbb{E}^0"] \arrow[d, "\nabla_{M-1}"] \arrow[rd, "\mu"] & E^M \arrow[d, "\tilde{\mathbb{H}}_M"] \\
E^N \otimes_{\A^0} \A^{1} \arrow[r, "\mathbb{E}^0\otimes 1"] & E^{N+1} \otimes_{\A^0} \A^{1} \arrow[r, "\mathbb{E}^0\otimes 1"] & \cdots \arrow[r, "\mathbb{E}^0\otimes 1"] & E^{M-1} \otimes_{\A^0} \A^{1} \arrow[r, "\mathbb{E}^0\otimes 1"]              & E^M \otimes_{\A^0} \A^{1}            
\end{tikzcd}
\end{center}
Now we continue in the same fashion and construct all $\Tilde{\mathbb{H}}^{k}$ with $(\mathbb{E}^0\otimes 1)\Tilde{\mathbb{H}}^{k} = \Tilde{\mathbb{H}}^{k}\mathbb{E}^0$ for all $k$.
\end{proof}

Now let's continue the proof of the main theorem. Set $\tilde{\mathbb{E}}^1 = (-1)^k \tilde{\mathbb{H}}_k$ on $E^k$ for each $k$. By our construction
\begin{equation*}
    \mathbb{E}^0 \tilde{\mathbb{E}}^1 + \tilde{\mathbb{E}}^1 \mathbb{E}^0 = 0
\end{equation*}
but $e^0\tilde{\mathbb{E}}^1 = \mathbb{F}^1 e^0$ might not hold. We will correct this by modifying $\tilde{\mathbb{E}}^1$. Consider the map $\psi = e^0\tilde{\mathbb{E}}^1 - \mathbb{F}^1 e^0: E^{\bt} \to F^{\bt}\otimes_{\A^0} \A^{1}$. It is easy to verify that $\psi$ is $\A^0$-linear and a map of chain complexes. Now we have the following diagram

\begin{center}
    \begin{tikzcd}
                                                     & (E^{\bt}\otimes_{\A^0} \A^{1}, \mathbb{E}^0 \otimes 1) \arrow[d, "e^0\otimes 1"] \\
E^{\bt} \arrow[r, "\psi"] \arrow[ru, "\tilde{\psi}"] & (F^{\bt}\otimes_{\A^0} \A^{1}, \mathbb{F}^0 \otimes 1)                          
\end{tikzcd}
\end{center}
here $e^0\otimes 1$ is a quasi-isomorphism since $e^0$ is a homotopy equivalence. $\tilde{\psi}$ is a lift of $\psi$ and there exists a homotopy $e^1:E^{\bt} \to  F^{\bt - 1}\otimes_{\A^0} \A^{1}$
$$
\phi - (e^0\otimes 1)\tilde{\psi} = (e^1 \mathbb{E}^0 + \mathbb{F}^0 e^1)
$$
Now let $\mathbb{E}^1 = \tilde{\mathbb{E}}$. We have
\begin{equation*}
    \mathbb{E}^0\mathbb{E}^1 + \mathbb{E}^1\mathbb{E}^0 =0
\end{equation*}
and 
\begin{equation*}
    e^1 \mathbb{E}^0 + \mathbb{F}^0 e^1 = e^0 \mathbb{E}^1 + \mathbb{F}^1 e^0
\end{equation*}
Now we have constructed the first two components $\mathbb{E}^0$ and $\mathbb{E}^1$ of the $\Z$-connection, and the first two components $e^0$ and $e^1$ of the quasi-isomorphism $E^{\bt}\otimes_{\A^0} \A^{\bt} \to F^{\bt}\otimes_{\A^0} \A^{\bt}$.

Now let's proceed to construct the rest components. Consider the mapping cone $C_{e^0}^{\bt}$ of $e^0$, i.e. $C_{e^0}^{\bt} = E[1]^{\bt} \oplus F^{\bt}$.  Now let $\mathbb{L}^0$ be defined as the matrix
$$
\begin{pmatrix}
\mathbb{E}^0[1] & 0 \\
e^0[1] & \mathbb{F}^0 
\end{pmatrix}
$$
Define $\mathbb{L}^1$ as the matrix 
$\begin{pmatrix}
\mathbb{E}^1[1] & 0 \\
e^1[1] & \mathbb{F}^1
\end{pmatrix}
$
Now $\mathbb{L}^0\mathbb{L}^0 = 0$ and $[\mathbb{L}^0, \mathbb{L}^1] = 0$ by construction. Let
\begin{equation*}
    D = \mathbb{L}^1\mathbb{L}^1 +
    \begin{pmatrix}
0 & 0 \\
\mathbb{F}^2e^0 & [\mathbb{F}^0, \mathbb{F}^2]
\end{pmatrix}
\end{equation*}

It is easy to check that $D$ is $\A$-linear, 
\begin{enumerate}
    \item $[\mathbb{L}^0, D] = 0$,\\
    \item $D|_{0\oplus F^{\bt}}$.
\end{enumerate}
Note that $(C_{e^0}^{\bt}, \mathbb{L}^0)$ is acyclic since it is a mapping cone of a quasi-isomorphism. By flatness of $\A^{\bt}$ over $\A^0$, $(C_{e^0}^{\bt}\otimes_{\A^0} \A^2, \mathbb{L}^0\otimes 1)$ is also acyclic. In addition,
\begin{equation*}
    \Hom^{\bt}_{\A^0}\big( (E^{\bt}, \mathbb{E}^0), (C_{e^0}^{\bt}\otimes_{\A^0} \A^2, \mathbb{L}^0\otimes 1) \big)
\end{equation*}
is a subcomplex of 
\begin{equation*}
    \Hom^{\bt}_{\A^0}\big(C_{e^0}^{\bt}, (C_{e^0}^{\bt}\otimes_{\A^0} \A^2, [\mathbb{L}^0, -]) \big) 
\end{equation*}
How $D\in \Hom^{\bt}_{\A^0}\big( (E^{\bt}, \mathbb{E}^0), (C_{e^0}^{\bt}\otimes_{\A^0} \A^2, \mathbb{L}^0\otimes 1) \big)$ is a cycle, so there exists some $\Tilde{\mathbb{L}}^2 \in \Hom^{\bt}_{\A^0}\big( (E^{\bt}, \mathbb{E}^0), (C_{e^0}^{\bt}\otimes_{\A^0} \A^2, \mathbb{L}^0\otimes 1) \big)$ such that $-D = [\mathbb{L}^0,\Tilde{\mathbb{L}}^2]$. Define $\mathbb{L}^2$ by
$$
\mathbb{L}^2 = \Tilde{\mathbb{L}}^2 + \begin{pmatrix}
0 & 0 \\
0 & \mathbb{F}^2
\end{pmatrix}
$$
We have 
\begin{align*}
    [\mathbb{L}^0,\mathbb{L}^2 ] =& \Big[\mathbb{L}^0,  \Tilde{\mathbb{L}}^2 + \begin{pmatrix}
0 & 0 \\
0 & \mathbb{F}^2
\end{pmatrix}\Big]\\
=& -D + \Big[\mathbb{L}^0,   \begin{pmatrix}
0 & 0 \\
0 & \mathbb{F}^2
\end{pmatrix}\Big]\\
=& -\mathbb{L}^1\mathbb{L}^1
\end{align*}
Therefore we get
\begin{align*}
    \mathbb{L}^0\mathbb{L}^2 + \mathbb{L}^1\mathbb{L}^1 + \mathbb{L}^2\mathbb{L}^0 = 0
\end{align*}
Following this pattern, we continue by setting
\begin{equation*}
    D = \mathbb{L}^1\mathbb{L}^2 + \mathbb{L}^2\mathbb{L}^1 + \begin{pmatrix}
0 & 0 \\
\mathbb{F}^3e^0 & [\mathbb{F}^0, \mathbb{F}^3]
\end{pmatrix}
\end{equation*}
Again it is easy to verify that $D$ is $\A^0$-linear, and 
\begin{enumerate}
    \item $[\mathbb{L}^0, D] = 0$,\\
    \item $D|_{0\oplus F^{\bt}}$.
\end{enumerate}
By the same reasoning as before, so there exists some $\Tilde{\mathbb{L}}^3 \in \Hom^{\bt}_{\A^0}\big( (E^{\bt}, \mathbb{E}^0), (C_{e^0}^{\bt}\otimes_{\A^0} \A^3, \mathbb{L}^0\otimes 1) \big)$ such that $-D = [\mathbb{L}^0,\Tilde{\mathbb{L}}^3]$. Define $\mathbb{L}^3$ by
$$
\mathbb{L}^3 = \Tilde{\mathbb{L}}^3 + \begin{pmatrix}
0 & 0 \\
0 & \mathbb{F}^3
\end{pmatrix}
$$
By easy verification we get $\sum_{i=0}^3\mathbb{L}^i \mathbb{L}^{3-i} = 0$.

Now suppose we have constructed $\mathbb{L}^0, \cdots, \mathbb{L}^n$ which satisfy
\begin{equation*}
    \sum_{i = 0}^k \mathbb{L}^i \mathbb{L}^{k-i} = 0
\end{equation*}
for $k = 0,\cdots, n$.
Then we define 
\begin{equation*}
    D = \sum_{i = 1}^n  \mathbb{L}^i \mathbb{L}^{n + 1 -i} + \begin{pmatrix}
0 & 0 \\
\mathbb{F}^{n+1}e^0 & [\mathbb{F}^0, \mathbb{F}^{n+1}]
\end{pmatrix}
\end{equation*}
Again we have $D$ is $\A^0$-linear, and 
\begin{enumerate}
    \item $[\mathbb{L}^0, D] = 0$,\\
    \item $D|_{0\oplus F^{\bt}}$.
\end{enumerate}
We can continue the inductive construction of $\mathbb{L}$ to get a $\Z$-connection satisfying $\mathbb{L} \mathbb{L} = 0$. Then we have constructed both components of the $\Z$-connection and the morphism from $(E^{\bt}, \mathbb{E})$ to $(F^{\bt}, \mathbb{F})$.

Now we have shown that $\RH$ is $A_{\infty}$-essentially surjective. Therefore, $\RH$ is an $A_{\infty}$-quasi-equivalence.

\begin{cor}
The $\infty$-category $\loci^{\infty}_{\Ch_k}\CF$ is equivalent to the $\infty$-category $\Modi_{\A}^{\coh}$, for $\A = \CE(\CF)$.
\end{cor}
\subsection{Integrate $\infty$-representations of $L_{\infty}$-algebroids}

$\RH$ is a functor from cohesive modules over the foliation dga $A$, which can also be regarded as cohesive modules over the foliation Lie algebroid $T\CF$. 
It is not hard to generalize the $\RH$ as a functor from cohesive modules over any $\li$-algebroids, where we only need to refine the iterated integrals to the corresponding vector bundles over the foliations, i.e. we only integrate along leaves of the (singular) foliations generated by  $\li$-algebroids. On the other hand, the monodromy $\infty$-groupoid of a perfect singular foliation $\CF$ is the truncation of the integration of the $\li$-algebroid $\g$ associated $\CF$. Therefore, given a perfect singular foliation $\CF$ with its associated $\li$-algebroid $\g$, we get the following commutative diagram

\begin{center}
    \begin{equation*}
        \begin{tikzcd}[ampersand replacement=\&]
{\g} \arrow[r, "\int"] \arrow[d, "\rep"] \& {\mon(\g)} \arrow[d, "\rep"]\arrow[r, "\tau"]\& {\mon(\CF)} \arrow[d, "\rep"] \\
{\Mod^{\coh}_{\g}} \arrow[r, "\RH"]                \& {\Loc^{\infty}(\g)} \arrow[r, "\tau"]     \& {\Loc^{\infty}(\CF)}     
\end{tikzcd}
    \end{equation*}
    \end{center}
    
where $\tau$ denotes the truncation functor.
A natural question to ask is when $\RH$ will be an $A_{\infty}$-quasi-equivalence, or induce an $\infty$-equivalence at the $\infty$-category level. This will be studied in a future paper.

\printindex

\nocite{*}
\printbibliography

@article{AC09,
  title={Representations up to homotopy and Bott's spectral sequence for Lie groupoids},
  author={C. A. Abad and M. Crainic},
  journal={Advances in Mathematics},
  year={2009},
  volume={248},
  pages={416-452}
}

@misc{AC11,
      title={Representations up to homotopy of Lie algebroids}, 
      author={Camilo Arias Abad and Marius Crainic},
      year={2011},
      eprint={0901.0319},
      archivePrefix={arXiv},
      primaryClass={math.DG}
}

@article{AHS78,
    author = "Atiyah, M. F. and Hitchin, Nigel J. and Singer, I. M.",
    title = "{Self duality in Four-Dimensional Riemannian Geometry}",
    doi = "10.1098/rspa.1978.0143",
    journal = "Proc. Roy. Soc. Lond. A",
    volume = "362",
    pages = "425--461",
    year = "1978"
}

@book{AD11,
  title={Introduction to Banach Spaces and Algebras},
  author={Allan, G.R. and Dales, H.G.},
  isbn={9780199206537},
  lccn={2011292228},
  series={Introduction to Banach Spaces and Algebras},
  year={2011},
  publisher={Oxford University Press}
}

@book{AG07,
  title={Pseudo-differential Operators and the Nash-Moser Theorem},
  author={Alinhac, S. and G{\'e}rard, P. },
  isbn={9780821834541},
  lccn={2006047985},
  series={Graduate studies in mathematics},
  year={2007},
  publisher={American Mathematical Society}
}

@article{AKSZ97,
  title={The Geometry of the Master Equation and Topological Quantum Field Theory},
  author={M. Alexandrov and M. Kontsevich and A. Schwarz and O. Zaboronsky},
  journal={International Journal of Modern Physics A},
  year={1997},
  volume={12},
  pages={1405-1429}
}

@misc{AMV21,
      title={Chern-Weil theory for $\infty$-local systems}, 
      author={Camilo Arias Abad and Santiago Pineda Montoya and Alexander Quintero Velez},
      year={2021},
      eprint={2105.00461},
      archivePrefix={arXiv},
      primaryClass={math.AT}
}

@book{Arm10,
  title={Groups and Symmetry},
  author={Armstrong, M.A.},
  isbn={9781441930859},
  lccn={87037677},
  series={Undergraduate Texts in Mathematics},
  year={2010},
  publisher={Springer New York}
}

@book{Arn92,
  title={Ordinary Differential Equations},
  author={Arnold, V.I.},
  isbn={9783540548133},
  lccn={91044188},
  series={Springer Textbook},
  year={1992},
  publisher={Springer Berlin Heidelberg}
}

@book{Arn97,
  title={Mathematical Methods of Classical Mechanics},
  author={Arnol'd, V.I.},
  isbn={9780387968902},
  lccn={97221119},
  series={Graduate Texts in Mathematics},
  year={1997},
  publisher={Springer New York}
}

@book{Arn13,
  title={Lectures on Partial Differential Equations},
  author={Cooke, R. and Arnold, V.I.},
  isbn={9783662054413},
  series={Universitext},
  year={2013},
  publisher={Springer Berlin Heidelberg}
}

@article{AS06,
author = {Androulidakis, Iakovos and Skandalis, Georges},
year = {2006},
month = {12},
pages = {},
title = {The holonomy groupoid of a singular foliation},
volume = {2009},
journal = {Journal für die reine und angewandte Mathematik (Crelles Journal)},
doi = {10.1515/CRELLE.2009.001}
}

@article{AS14,
  title={Higher holonomies: Comparing two constructions},
  author={C. A. Abad and F. Schaetz},
  journal={Differential Geometry and Its Applications},
  year={2014},
  volume={40},
  pages={14-42}
}

@article{AS12,
author = {Abad, Camilo and Schatz, F.},
year = {2012},
month = {07},
pages = {},
title = {The $A_{\infty}$ de Rham Theorem and Integration of Representations up to Homotopy},
volume = {2013},
journal = {International Mathematics Research Notices},
doi = {10.1093/imrn/rns166}
}

@inproceedings{Ati13,
  title={Geometry of Yang-Mills fields},
  author={F. Atiyah},
  year={2013}
}

@article{AZ11,
  title={Smoothness of holonomy covers for singular foliations and essential isotropy},
  author={Iakovos Androulidakis and M. Zambon},
  journal={Mathematische Zeitschrift},
  year={2011},
  volume={275},
  pages={921-951}
}

@article{AZ12,
  title={Holonomy transformations for singular foliations},
  author={Iakovos Androulidakis and M. Zambon},
  journal={Advances in Mathematics},
  year={2012},
  volume={256},
  pages={348-397}
}

@article{Bat79,
 ISSN = {00029947},
 author = {Majorie Batchelor},
 journal = {Transactions of the American Mathematical Society},
 pages = {329--338},
 publisher = {American Mathematical Society},
 title = {The Structure of Supermanifolds},
 volume = {253},
 year = {1979}
}

@article{BB72,
  title={Singularities of holomorphic foliations},
  author={Baum, Paul and Bott, Raoul},
  journal={Journal of differential geometry},
  volume={7},
  number={3-4},
  pages={279--342},
  year={1972},
  publisher={Lehigh University}
}

@book{BB12,
  title={Topology and Analysis: The Atiyah-Singer Index Formula and Gauge-Theoretic Physics},
  author={Bleecker, D.D. and Booss, B.},
  isbn={9781468406276},
  lccn={84026684},
  series={Universitext},
  year={2012},
  publisher={Springer New York}
}

@book{BCH14,
  title={An Introduction to Involutive Structures},
  author={Berhanu, S. and Cordaro, P.D. and Hounie, J.},
  isbn={9780511388149},
  series={New Mathematical Monographs},
  year={2014},
  publisher={Cambridge University Press}
}

@article{BD10,
title = {Mukai duality for gerbes with connection},
volume = {2010},
author = {Jonathan Block and Calder Daenzer},
address = {Berlin :},
issn = {0075-4102},
journal = {Journal für die reine und angewandte Mathematik.},
lccn = {00227020},
number = {639},
publisher = {Walter de Gruyter},
year = {2010},
}

@book{BD13,
  title={Representations of Compact Lie Groups},
  author={Br{\"o}cker, T. and Dieck, T.},
  isbn={9783662129180},
  series={Graduate Texts in Mathematics},
  year={2013},
  publisher={Springer Berlin Heidelberg}
}

@article{BBLE20,
  title={Deformation spaces and normal forms around transversals},
  author={Francis Bischoff and H. Bursztyn and H. Lima and E. Meinrenken},
  journal={Compositio Mathematica},
  year={2020},
  volume={156},
  pages={697 - 732}
}

@book{Ber06,
  title={Th{\'e}orie des Intersections et Th{\'e}or{\`e}me de Riemann-Roch: S{\'e}minaire de G{\'e}om{\'e}trie Alg{\'e}brique du Bois Marie 1966 /67 (SGA 6)},
  author={Berthelot, P. and Ferrand, D. and Grothendieck, A. and Jouanolou, J.P. and Jussila, O. and Illusie, L. and Kleiman, S. and Raynaud, M. and Serre, J.P.},
  isbn={9783540369363},
  series={Lecture Notes in Mathematics},
  year={2006},
  publisher={Springer Berlin Heidelberg}
}

@book{Bes,
  title={Einstein Manifolds},
  author={Besse, A.L.},
  isbn={9783540741206},
  lccn={2007938035},
  series={Classics in Mathematics},
  year={2007},
  publisher={Springer Berlin Heidelberg}
}

@article{BFLS97,
author = {Barnich, G. and Fulp, R. and Lada, Tom and Stasheff, Jim},
year = {1997},
month = {02},
pages = {},
title = {The sh Lie Structure of Poisson Brackets in Field Theory},
volume = {191},
journal = {Communications in Mathematical Physics},
doi = {10.1007/s002200050278}
}

@Inbook{BG17,
author="Behrend, Kai
and Getzler, Ezra",
editor="Bost, Jean-Beno{\^i}t
and Hofer, Helmut
and Labourie, Fran{\c{c}}ois
and Le Jan, Yves
and Ma, Xiaonan
and Zhang, Weiping",
title="Geometric Higher Groupoids and Categories",
bookTitle="Geometry, Analysis and Probability: In Honor of Jean-Michel Bismut",
year="2017",
publisher="Springer International Publishing",
address="Cham",
pages="1--45",
isbn="978-3-319-49638-2",
}

@inproceedings{BGV92,
  title={Heat Kernels and Dirac Operators},
  author={Nicole Berline and E. Getzler and M. Vergne},
  year={1992}
}

@misc{BK17,
      title={Quasi-coherent sheaves in differential geometry}, 
      author={Dennis Borisov and Kobi Kremnizer},
      year={2017},
      eprint={1707.01145},
      archivePrefix={arXiv},
      primaryClass={math.DG}
}

@misc{Blo05,
	title={Duality and equivalence of module categories in noncommutative geometry I},
	author={Jonathan Block},
	year={2005},
	eprint={math/0509284},
	archivePrefix={arXiv},
	primaryClass={math.QA}
}

@misc{BLX21,
      title={Derived Differentiable Manifolds}, 
      author={Kai Behrend and Hsuan-Yi Liao and Ping Xu},
      year={2021},
      eprint={2006.01376},
      archivePrefix={arXiv},
      primaryClass={math.DG}
}

@article{BMS87,
author = {A. Beilinson and R. MacPherson and V. Schechtman},
title = {{Notes on motivic cohomology}},
volume = {54},
journal = {Duke Mathematical Journal},
number = {2},
publisher = {Duke University Press},
pages = {679 -- 710},
year = {1987},
doi = {10.1215/S0012-7094-87-05430-5}
}

@misc{BN11,
      title={Simplicial approach to derived differential manifolds}, 
      author={Dennis Borisov and Justin Noel},
      year={2011},
      eprint={1112.0033},
      archivePrefix={arXiv},
      primaryClass={math.DG}
}

@article{BNV13,
author = {Bunke, Ulrich and Nikolaus, Thomas and Völkl, Michael},
year = {2013},
month = {11},
pages = {},
title = {Differential cohomology theories as sheaves of spectra},
volume = {11},
journal = {Journal of Homotopy and Related Structures},
doi = {10.1007/s40062-014-0092-5}
}

@inproceedings{Bos08,
  title={Lecture notes on groupoid cohomology},
  author={R. Bos},
  year={2008}
}

@article{BP13,
  title={On the category of Lie n-algebroids},
  author={G. Bonavolont{\`a} and N. Poncin},
  journal={Journal of Geometry and Physics},
  year={2013},
  volume={73},
  pages={70-90}
}

@misc{BP15,
      title={Smooth one-dimensional topological field theories are vector bundles with connection}, 
      author={Daniel Berwick-Evans and Dmitri Pavlov},
      year={2015},
      eprint={1501.00967},
      archivePrefix={arXiv},
      primaryClass={math.AT}
}

@book{Bre12,
  title={Sheaf Theory},
  author={Bredon, G.E.},
  isbn={9781461206477},
  lccn={96044232},
  series={Graduate Texts in Mathematics},
  year={2012},
  publisher={Springer New York}
}

@article{Bro73,
  title={Abstract homotopy theory and generalized sheaf cohomology},
  author={K. S. Brown},
  journal={Transactions of the American Mathematical Society},
  year={1973},
  volume={186},
  pages={419-458}
}

@article{BRS75,
     author = {Becchi, C. and Rouet, A. and Stora, R.},
     title = {Renormalization of Gauge Theories},
     journal = {Les rencontres physiciens-math\'ematiciens de Strasbourg -RCP25},
     note = {talk:10},
     publisher = {Institut de Recherche Math\'ematique Avanc\'ee - Universit\'e Louis Pasteur},
     volume = {22},
     year = {1975},
     language = {en},
  
}

@article{BS10,
author = {Bunke, Ulrich and Schick, Thomas},
year = {2010},
month = {11},
pages = {},
title = {Differential K-Theory: A Survey},
volume = {17},
isbn = {978-3-642-22841-4},
journal = {Springer Proceedings in Mathematics},
doi = {10.1007/978-3-642-22842-1_11}
}

@article{BS14,
  title={The higher Riemann–Hilbert correspondence},
  author={Jonathan Block and Aaron M. Smith},
  journal={Advances in Mathematics},
  year={2014},
  volume={252},
  pages={382-405}
}

@article{BSS76,
title = {On the de Rham theory of certain classifying spaces},
journal = {Advances in Mathematics},
volume = {20},
number = {1},
pages = {43-56},
year = {1976},
issn = {0001-8708},
doi = {https://doi.org/10.1016/0001-8708(76)90169-9},
author = {R Bott and H Shulman and J Stasheff}
}

@misc{BSW21,
      title={Coherent sheaves, superconnections, and RRG}, 
      author={Jean-Michel Bismut and Shu Shen and Zhaoting Wei},
      year={2021},
      eprint={2102.08129},
      archivePrefix={arXiv},
      primaryClass={math.AG}
}

@book{BT95,
  title={Differential Forms in Algebraic Topology},
  author={Bott, R. and Tu, L.W.},
  isbn={9780387906133},
  series={Graduate Texts in Mathematics},
  year={1995},
  publisher={Springer New York}
}

@article{
Bun18,
title = {Foliated manifolds, algebraic K-theory, and a secondary invariant},
journal = {Münster Journal of Mathematics},
author = {Bunke, Ulrich},
publisher = {Mathematisches Institut (Universität Münster)},
year = {2018},
edition = {[Electronic ed.]}
}

@article{BX06,
  title={Differentiable stacks and gerbes},
  author={K. Behrend and P. Xu},
  journal={Journal of Symplectic Geometry},
  year={2006},
  volume={9},
  pages={285-341}
}

@article{BZ11,
title = "Lie algebroid fibrations",
journal = "Advances in Mathematics",
volume = "226",
number = "4",
pages = "3105 - 3135",
year = "2011",
issn = "0001-8708",
doi = "https://doi.org/10.1016/j.aim.2010.10.006",
author = "Olivier Brahic and Chenchang Zhu",
}

@unpublished{BZ21,
author = "Jonathan Block and Qingyun Zeng",
title = "Singular foliation and characteristic classes",
note = "in preparation",
}

@inproceedings{Car11,
  title={Categorical properties of topological and differentiable stacks},
  author={D. Carchedi},
  year={2011}
}

@article{Car15,
author = {Carchedi, David},
year = {2015},
month = {04},
pages = {},
title = {On The Homotopy Type of Higher Orbifolds and Haefliger Classifying Spaces},
volume = {294},
journal = {Advances in Mathematics},
doi = {10.1016/j.aim.2016.03.007}
}

@article{CD02,
  title={Yang–Mills Algebra},
  author={A. Connes and M. Dubois-Violette},
  journal={Letters in Mathematical Physics},
  year={2002},
  volume={61},
  pages={149-158}
}

@book{Car92,
  title={Riemannian Geometry},
  author={do Carmo, M.P.},
  isbn={9780817634902},
  lccn={91037377},
  series={Mathematics (Birkh{\"a}user) theory},
  year={1992},
  publisher={Birkh{\"a}user Boston}
}

@book{Car12,
  title={Differential Forms and Applications},
  author={do Carmo, M.P.},
  isbn={9783642579516},
  lccn={94021965},
  series={Universitext},
  year={2012},
  publisher={Springer Berlin Heidelberg}
}

@misc{CCN21,
      title={Lie algebroids are curved Lie algebras}, 
      author={Damien Calaque and Ricardo Campos and Joost Nuiten},
      year={2021},
      eprint={2103.10728},
      archivePrefix={arXiv},
      primaryClass={math.AT}
}

@article{CF03,
	title = "Integrability of Lie brackets",
	author = "Marius Crainic and Fernandes, Rui Loja",
	year = "2003",
	month = "3",
	doi = "10.4007/annals.2003.157.575",
	language = "English (US)",
	volume = "157",
	pages = "575--620",
	journal = "Annals of Mathematics",
	issn = "0003-486X",
	publisher = "Princeton University Press",
	number = "2",
}

@inbook{CF11,
title = "Lectures on integrability of Lie brackets",
author = "Marius Crainic and Fernandes, {Rui Loja}",
year = "2011",
language = "English (US)",
volume = "17",
series = "Geom. Topol. Monogr.",
publisher = "Geom. Topol. Publ., Coventry",
pages = "1--107",
booktitle = "Lectures on Poisson geometry",
}

@book{Che13,
  title={Complex Manifolds without Potential Theory: with an appendix on the geometry of characteristic classes},
  author={Chern, S.},
  isbn={9781468493443},
  series={Universitext},
  year={2013},
  publisher={Springer New York}
}

@book{Che,
  title={Comparison Theorems in Riemannian Geometry},
  isbn={9780444107640},
  lccn={74083725},
  series={ISSN},
  author={ Cheeger, D. G. Ebin},
  year={1975},
  publisher={Elsevier Science}
}

@article{Che77,
  title={Iterated path integrals},
  author={Kuo-Tsai Chen},
  journal={Bulletin of the American Mathematical Society},
  year={1977},
  volume={83},
  pages={831-879}
}

@article{CLX19,
author = {Chen, Zhuo and Lang, Honglei and Xiang, Maosong},
title = {Atiyah Classes of Strongly Homotopy Lie Pairs},
journal = {Algebra Colloquium},
volume = {26},
number = {02},
pages = {195-230},
year = {2019},
doi = {10.1142/S1005386719000178},
}

@article{CL19,
author = {Corrêa, Maurício and Lourenço, Fernando},
year = {2019},
month = {07},
pages = {},
title = {Determination of Baum-Bott residues of higher codimensional foliations},
volume = {23},
journal = {Asian Journal of Mathematics},
doi = {10.4310/AJM.2019.v23.n3.a8}
}

@article{CM00,
  title={{\v C}ech-De Rham theory for leaf spaces of foliations},
  author={M. Crainic and I. Moerdijk},
  journal={Mathematische Annalen},
  year={2000},
  volume={328},
  pages={59-85}
}

@article{CM04,
  title={Deformations of Lie brackets: cohomological aspects},
  author={M. Crainic and I. Moerdijk},
  journal={Journal of the European Mathematical Society},
  year={2004},
  volume={10},
  pages={1037-1059}
}

@article{CR12,
author = {Carchedi, David and Roytenberg, Dmitry},
year = {2012},
month = {11},
pages = {},
title = {On Theories of Superalgebras of Differentiable Functions},
volume = {28},
journal = {Theory and Applications of Categories}
}

@article{CS17,
  title={A tour about existence and uniqueness of dg enhancements and lifts},
  author={A. Canonaco and P. Stellari},
  journal={Journal of Geometry and Physics},
  year={2017},
  volume={122},
  pages={28-52}
}

@misc{Dat16,
      title={D-modules and complex foliations}, 
      author={Hamidou Dathe},
      year={2016},
      eprint={1606.09060},
      archivePrefix={arXiv},
      primaryClass={math.DG}
}

@article{Deb00,
author = {Debord, Claire},
year = {2000},
month = {01},
pages = {},
title = {Local integration of Lie algebroids},
volume = {54},
journal = {Banach Center Publications},
doi = {10.4064/bc54-0-2}
}

@article{Deb01,
author = {Claire Debord},
title = {{Holonomy Groupoids of Singular Foliations}},
volume = {58},
journal = {Journal of Differential Geometry},
number = {3},
publisher = {Lehigh University},
pages = {467 -- 500},
year = {2001},
doi = {10.4310/jdg/1090348356},
URL = {https://doi.org/10.4310/jdg/1090348356}
}

@book{Deb05,
  title={Complex Tori and Abelian Varieties},
  author={Debarre, O.},
  isbn={9780821831656},
  lccn={2005049857},
  series={Collection SMF},
  year={2005},
  publisher={American Mathematical Society}
}

@article{Deb13,
  title={Longitudinal smoothness of the holonomy groupoid},
  author={C. Debord},
  journal={Comptes Rendus Mathematique},
  year={2013},
  volume={351},
  pages={613-616}
}

@article{DHI04, title={Hypercovers and simplicial presheaves}, volume={136}, DOI={10.1017/S0305004103007175}, number={1}, journal={Mathematical Proceedings of the Cambridge Philosophical Society}, publisher={Cambridge University Press}, author={DUGGER, DANIEL and HOLLANDER, SHARON and ISAKSEN, DANIEL C.}, year={2004}, pages={9–51}}

@book{Din07,
  title={Supersymmetry and String Theory: Beyond the Standard Model},
  author={Dine, M.},
  isbn={9781139462440},
  year={2007},
  publisher={Cambridge University Press}
}

@book{Dou88,
  title={Banach Algebra Techniques in Operator Theory},
  author={Douglas, R.G.},
  series={Pure and applied mathematics},
  year={1988},
  publisher={Acad. Press}
}

@book{Dun07,
  title={Motivic Homotopy Theory: Lectures at a Summer School in Nordfjordeid, Norway, August 2002},
  author={Dundas, B.I. and Levine, M. and Ostvaer, P.A. and Jahren, B. and {\O}stv{\ae}r, P.A. and R{\"o}ndigs, O. and Voevodsky, V.},
  isbn={9783540458951},
  lccn={2006933719},
  series={Universitext (Berlin. Print)},
  year={2007},
  publisher={Springer}
}

@article{Eas13,
author = {Eastwood, Michael},
year = {2013},
month = {01},
pages = {},
title = {A double fibration transform for complex projective space},
isbn = {9780821887387},
doi = {10.1090/conm/598/11999}
}

@book{Eis13,
  title={Commutative Algebra: with a View Toward Algebraic Geometry},
  author={Eisenbud, D.},
  isbn={9781461253501},
  series={Graduate Texts in Mathematics},
  year={2013},
  publisher={Springer New York}
}

@article{EE67,
author = {Clifford J. Earle and James Eells Jr.},
title = {{Foliations and Fibrations}},
volume = {1},
journal = {Journal of Differential Geometry},
number = {1-2},
publisher = {Lehigh University},
pages = {33 -- 41},
year = {1967},
doi = {10.4310/jdg/1214427879}
}

@book{Ell12,
  title={Calabi-Yau Manifolds and Related Geometries: Lectures at a Summer School in Nordfjordeid, Norway, June 2001},
  author={Ellingsrud, G. and Gross, M. and Huybrechts, D. and Olson, L. and Joyce, D. and Ranestad, K. and Stromme, S.A.},
  isbn={9783642190049},
  lccn={2002034350},
  series={Universitext},
  year={2012},
  publisher={Springer Berlin Heidelberg}
}

@article{ES93,
author = {Michael G. Eastwood and Michael A. Singer},
title = {{The Fröhlicher spectral sequence on a twistor space}},
volume = {38},
journal = {Journal of Differential Geometry},
number = {3},
publisher = {Lehigh University},
pages = {653 -- 669},
year = {1993},
doi = {10.4310/jdg/1214454485}
}

@book{Eva10,
  title={Partial Differential Equations},
  author={Evans, L.C.},
  isbn={9780821849743},
  lccn={2009044716},
  series={Graduate studies in mathematics},
  year={2010},
  publisher={American Mathematical Society}
}

@misc{Eyn14,
      title={The Frobenius theorem for Banach distributions on infinite-dimensional manifolds and applications in infinite-dimensional Lie theory}, 
      author={Jan Milan Eyni},
      year={2014},
      eprint={1407.3166},
      archivePrefix={arXiv},
      primaryClass={math.GR}
}

@book{Fol95,
  title={Introduction to Partial Differential Equations},
  author={Folland, G.B.},
  isbn={9780691043616},
  lccn={95032308},
  series={Mathematical Notes},
  year={1995},
  publisher={Princeton University Press}
}

@book{FGI05,
  title={Fundamental Algebraic Geometry: Grothendieck's FGA Explained},
  author={Fantechi, B. and Gottsche, L. and Illusie, L.},
  isbn={9780821842454},
  lccn={2005053614},
  series={Mathematical surveys and monographs},
  year={2005},
  publisher={American Mathematical Society}
}

@book{FH91,
  title={Representation Theory: A First Course},
  author={Harris, W.F.J. and Fulton, W. and Harris, J.},
  isbn={9780387974958},
  lccn={96165313},
  series={Graduate Texts in Mathematics},
  year={1991},
  publisher={Springer New York}
}

@book{FHT12,
  title={Rational Homotopy Theory},
  author={Felix, Y. and Halperin, S. and Thomas, J.C.},
  isbn={9781461301059},
  lccn={00041913},
  series={Graduate Texts in Mathematics},
  year={2012},
  publisher={Springer New York}
}

@article{Fio04,
  title={An introduction to the Batalin-Vilkovisky formalism},
  author={D. Fiorenza},
  journal={arXiv: Quantum Algebra},
  year={2004}
}

@article{FSS19,
author = {Fiorenza, Domenico and Schreiber, Urs and Stasheff, Jim},
year = {2010},
month = {11},
pages = {},
title = {Čech cocycles for differential characteristic classes: An $\infty$-Lie theoretic construction},
volume = {16},
journal = {Advances in Theoretical and Mathematical Physics},
doi = {10.4310/ATMP.2012.v16.n1.a5}
}

@book{GH14,
  title={Principles of Algebraic Geometry},
  author={Griffiths, P. and Harris, J.},
  isbn={9781118626320},
  series={Wiley Classics Library},
  year={2014},
  publisher={Wiley}
}

@book{GHL04,
  title={Riemannian Geometry},
  author={Gallot, S. and Hulin, D. and Lafontaine, J.},
  isbn={9783540204930},
  lccn={90022646},
  series={Universitext},
  year={2004},
  publisher={Springer Berlin Heidelberg}
}

@article{Get09,
	ISSN = {0003486X},
	author = {Ezra Getzler},
	journal = {Annals of Mathematics},
	number = {1},
	pages = {271--301},
	publisher = {Annals of Mathematics},
	title = {Lie Theory for Nilpotent $L_{\infty}$-Algebras},
	volume = {170},
	year = {2009}
}

@inproceedings{GJ09,
  title={Simplicial Homotopy Theory},
  author={P. Goerss and J. Jardine},
  booktitle={Modern Birkh{\"a}user Classics},
  year={2009}
}

@book{GR19,
  title={A Study in Derived Algebraic Geometry: Volume I: Correspondences and Duality},
  author={Gaitsgory, D. and Rozenblyum, N.},
  isbn={9781470452841},
  lccn={2016054809},
  series={Mathematical Surveys and Monographs},
  year={2019},
  publisher={American Mathematical Society}
}

@book{GS99,
  title={4-Manifolds and Kirby Calculus},
  author={Gompf, R.E. and Stipsicz, A.I. },
  isbn={9780821809945},
  lccn={99029942},
  series={Graduate studies in mathematics},
  year={1999},
  publisher={American Mathematical Society}
}

@book{GS12,
  title={Mathematical Concepts of Quantum Mechanics},
  author={Gustafson, S.J. and Sigal, I.M.},
  isbn={9783642557293},
  lccn={2003052641},
  series={Universitext},
  year={2012},
  publisher={Springer Berlin Heidelberg}
}

@book{GT01,
  title={Elliptic Partial Differential Equations of Second Order},
  author={Gilbarg, D. and Trudinger, N.S.},
  isbn={9783540411604},
  lccn={00052272},
  series={Classics in Mathematics},
  year={2001},
  publisher={Springer Berlin Heidelberg}
}

@article{Gug77,
author = {V. K. A. M. Gugenheim},
title = {{On Chen's iterated integrals}},
volume = {21},
journal = {Illinois Journal of Mathematics},
number = {3},
publisher = {Duke University Press},
pages = {703 -- 715},
year = {1977},
doi = {10.1215/ijm/1256049021}
}

@article{GZ19,
author = {Garmendia, Alfonso and Zambon, Marco},
year = {2019},
month = {02},
pages = {},
title = {Hausdorff Morita Equivalence of singular foliations},
volume = {55},
journal = {Annals of Global Analysis and Geometry},
doi = {10.1007/s10455-018-9620-6}
}

@book{Hal15,
  title={Lie Groups, Lie Algebras, and Representations: An Elementary Introduction},
  author={Hall, B.},
  isbn={9783319134673},
  series={Graduate Texts in Mathematics},
  year={2015},
  publisher={Springer International Publishing}
}

@book{Hat02,
  title={Algebraic Topology},
  author={Hatcher, A.},
  isbn={9780521795401},
  lccn={00065166},
  series={Algebraic Topology},
  year={2002},
  publisher={Cambridge University Press}
}

@book{Har77,
publisher = {Springer New York :},
series = {Algebraic Geometry},
title = {Algebraic Geometry /},
volume = {52},
year = {1977},
author = {Hartshorne, Robin},
address = {New York, NY :},
edition = {1st ed. 1977.},
isbn = {9781441928078},
issn = {0072-5285},
}

@article{Hen08,
author = {Henriques, André},
year = {2008},
month = {07},
pages = {},
title = {Integrating $L^{\infty}$-algebras},
volume = {144},
journal = {Compositio Mathematica - COMPOS MATH},
doi = {10.1112/S0010437X07003405}
}

@article{Hen18,
  title={Tangent Lie algebra of derived Artin stacks},
  author={B. Hennion},
  journal={Crelle's Journal},
  year={2018},
  volume={2018},
  pages={1-45}
}

@article{Hep09,
author = {Hepworth, Richard},
journal = {Theory and Applications of Categories [electronic only]},
language = {eng},
pages = {542-587},
publisher = {Mount Allison University, Department of Mathematics and Computer Science, Sackville},
title = {Vector fields and flows on differentiable stacks.},
volume = {22},
year = {2009},
}

@article{Her60,
  title={On the Differential Geometry of Foliations},
  author={R. Hermann},
  journal={Annals of Mathematics},
  year={1960},
  volume={72},
  pages={445}
}

@book{HL11,
  title={Elliptic Partial Differential Equations},
  author={Han, Q. and Lin, F.},
  isbn={9780821853139},
  lccn={2010051489},
  series={Courant lecture notes in mathematics},
  year={2011},
  publisher={Courant Institute of Mathematical Sciences, New York University}
}

@book{Hir12,
  title={Differential Topology},
  author={Hirsch, M.W.},
  isbn={9781468494495},
  series={Graduate Texts in Mathematics},
  year={2012},
  publisher={Springer New York}
}

@article{Hit80,
    author = "Hitchin, Nigel J.",
    title = "{LINEAR FIELD EQUATIONS ON SELFDUAL SPACES}",
    doi = "10.1098/rspa.1980.0028",
    journal = "Proc. Roy. Soc. Lond. A",
    volume = "370",
    pages = "173--191",
    year = "1980"
}

@article{HL02,
author = {Heitsch, James and Lazarov, Connor},
year = {2002},
month = {01},
pages = {437-468},
title = {Riemann-Roch-Grothendieck and torsion for foliations},
volume = {12},
journal = {Journal of Geometric Analysis},
doi = {10.1007/BF02922049}
}

@book{Hor03,
  title={The Analysis of Linear Partial Differential Operators I: Distribution Theory and Fourier Analysis},
  author={H{\"o}rmander, L.},
  isbn={9783540006626},
  lccn={2003050516},
  series={Classics in Mathematics},
  year={2003},
  publisher={Springer Berlin Heidelberg}
}

@misc{Hor12,
      title={Parallel Transport on Higher Loop Spaces}, 
      author={Ivan Horozov},
      year={2012},
      eprint={1206.5784},
      archivePrefix={arXiv},
      primaryClass={math.AT}
}

@article{Hov98,
  title={Monoidal model categories},
  author={Mark Hovey},
  eprint={9803002},
      archivePrefix={arXiv},
      primaryClass={math.AT},
  year={1998}
}

@book{Hov07,
  title={Model Categories},
  author={Hovey, M.},
  isbn={9780821843611},
  lccn={98345399},
  series={Mathematical surveys and monographs},
  year={2007},
  publisher={American Mathematical Society}
}

@article{Hue98,
  title={Lie-Rinehart algebras, Gerstenhaber algebras and Batalin-Vilkovisky algebras},
  author={J. Huebschmann},
  journal={Annales de l'Institut Fourier},
  year={1998},
  volume={48},
  pages={425-440}
}

@article{Hue03,
author = {Huebschmann, Johannes},
year = {2003},
month = {04},
pages = {},
title = {Lie-Rinehart algebras, descent, and quantization},
volume = {43},
journal = {Fields Inst. Commun.}
}

@book{Huy06,
  title={Complex Geometry: An Introduction},
  author={Huybrechts, D.},
  isbn={9783540266877},
  lccn={2004108312},
  series={Universitext},
  year={2006},
  publisher={Springer Berlin Heidelberg}
}

@book{Huy10,
  title={The Geometry of Moduli Spaces of Sheaves},
  author={Huybrechts, D. and Lehn, M.},
  isbn={9781139485821},
  series={Cambridge Mathematical Library},
  year={2010},
  publisher={Cambridge University Press}
}

@misc{Igu09,
	title={Iterated integrals of superconnections},
	author={Kiyoshi Igusa},
	year={2009},
	eprint={0912.0249},
	archivePrefix={arXiv},
	primaryClass={math.AT}
}

@article{Ito83,
  title={On the Moduli Space of Anti-Self-Dual Yang-Mills Connections on Kahler Surfaces},
  author={M. Itoh},
  journal={Publications of The Research Institute for Mathematical Sciences},
  year={1983},
  volume={19},
  pages={15-32}
}

@article{Ito85,
author = {Mitsuhiro Itoh},
title = {{The moduli space of Yang-Mills connections over a Kähler surface is a complex manifold}},
volume = {22},
journal = {Osaka Journal of Mathematics},
number = {4},
publisher = {Osaka University and Osaka City University, Departments of Mathematics},
pages = {845 -- 862},
year = {1985},
doi = {ojm/1200778770},
}

@article{Jac00,
  title={Transversely holomorphic foliations and CR structures},
  author={Jacobowitz, Howard},
  journal={Mat. Contemp},
  volume={18},
  pages={175--194},
  year={2000}
}

@book{Jar05,
  title={Local Homotopy Theory},
  author={Jardine, J.F.},
  isbn={9781493923007},
  lccn={2015933625},
  series={Springer Monographs in Mathematics},
  year={2015},
  publisher={Springer New York}
}

@article{Jar87,
  title={Simplicial presheaves},
  author={John F. Jardine},
  journal={Journal of Pure and Applied Algebra},
  volume={47},
  pages={35--87},
  year={1987}
}

@book{JLY10,
  title={Cohomology of Groups and Algebraic K-theory},
  author={Ji, L. and Liu, K. and Yau, S.T.},
  isbn={9781571461445},
  series={Advanced Lectures in Mathematics - International Press},
  year={2010},
  publisher={International Press}
}

@book{Joh91,
  title={Partial Differential Equations},
  author={John, F.},
  isbn={9780387906096},
  lccn={81016636},
  series={Applied Mathematical Sciences},
  year={1991},
  publisher={Springer New York}
}

@book{Jos13,
  title={Compact Riemann Surfaces: An Introduction to Contemporary Mathematics},
  author={Jost, J.},
  isbn={9783662034460},
  lccn={96037590},
  series={Universitext},
  year={2013},
  publisher={Springer Berlin Heidelberg}
}

@book{Jos17,
  title={Riemannian Geometry and Geometric Analysis},
  author={Jost, J.},
  isbn={9783319618609},
  series={Universitext},
  year={2017},
  publisher={Springer International Publishing}
}

@article{Joy11,
author = {Joyce, Dominic},
year = {2011},
month = {04},
pages = {},
title = {An introduction to $\cinf$ schemes and $\cinf$ algebraic geometry},
volume = {17},
journal = {Surveys in Differential Geometry},
doi = {10.4310/SDG.2012.v17.n1.a7}
}

@misc{Joy12,
      title={An introduction to d-manifolds and derived differential geometry}, 
      author={Dominic Joyce},
      year={2012},
      eprint={1206.4207},
      archivePrefix={arXiv},
      primaryClass={math.DG}
}

@book{Joy19,
  title={Algebraic Geometry over $\cinf$-Rings},
  author={Joyce, D.},
  isbn={9781470436452},
  lccn={2019033051},
  series={Memoirs of the American Mathematical Society},
  year={2019},
  publisher={American Mathematical Society}
}

@book{Kod12,
  title={Complex Manifolds and Deformation of Complex Structures},
  author={Akao, K. and Kodaira, K.},
  isbn={9781461385905},
  lccn={85009825},
  series={Grundlehren der mathematischen Wissenschaften},
  year={2012},
  publisher={Springer New York}
}

@book{Kar09,
  title={K-Theory: An Introduction},
  author={Karoubi, M.},
  isbn={9783540798903},
  series={Classics in Mathematics},
  year={2009},
  publisher={Springer Berlin Heidelberg}
}

@article{Kas07,
author = {Kasparov, G},
year = {2007},
month = {10},
pages = {513},
title = {The operator K-functor and extensions of C*-algebras},
volume = {16},
journal = {Mathematics of the USSR-Izvestiya},
doi = {10.1070/IM1981v016n03ABEH001320}
}

@inproceedings{Kel06,
  title={On differential graded categories},
  author={B. Keller},
  year={2006}
}

@inproceedings{Kor04,
  title={Elliptic Involutive Structures and Generalized Higgs Algebroids},
  author={Eric O. Korman},
  year={2014}
}

@article{KSV95,
author = {Takashi Kimura and Jim Stasheff and Alexander A. Voronov},
title = {{On operad structures of moduli spaces and string theory}},
volume = {171},
journal = {Communications in Mathematical Physics},
number = {1},
publisher = {Springer},
pages = {1 -- 25},
year = {1995},
doi = {cmp/1104273401},
URL = {https://doi.org/}
}

@article{KT74,
  title={Characteristic invariants of foliated bundles},
  author={F. Kamber and P. Tondeur},
  journal={manuscripta mathematica},
  year={1974},
  volume={11},
  pages={51-89}
}

@inproceedings{Kor14,
  title={Elliptic involutive structures and generalized Higgs algebroids},
  author={Eric O. Korman},
  year={2014}
}

@book{Lan05,
  title={Algebra},
  author={Lang, S.},
  isbn={9780387953854},
  lccn={01054916},
  series={Graduate Texts in Mathematics},
  year={2005},
  publisher={Springer New York}
}

@book{Lan95,
  title={Differential and Riemannian Manifolds},
  author={Lang, S.},
  isbn={9780387943381},
  lccn={94020828},
  series={Graduate texts in mathematics},
  year={1995},
  publisher={Springer-Verlag}
}

@misc{Lav18,
      title={Lie $\infty$-algebroids and singular foliations}, 
      author={Sylvain Lavau},
      year={2018},
      eprint={1703.07404},
      archivePrefix={arXiv},
      primaryClass={math.DG}
}

@inproceedings{Lee21,
  title={Mapping Spaces, Signatures, and Data},
  author={D. Lee},
  year={2021}
}

@misc{Lea19a,
      title={Obstructions to representations up to homotopy and ideals}, 
      author={Madeleine Jotz Lean},
      year={2019},
      eprint={1905.10237},
      archivePrefix={arXiv},
      primaryClass={math.DG}
}

@misc{Lea19b,
      title={Infinitesimal ideal systems and the Atiyah class}, 
      author={Madeleine Jotz Lean},
      year={2019},
      eprint={1910.04492},
      archivePrefix={arXiv},
      primaryClass={math.DG}
}

@book{Lee06,
  title={Riemannian Manifolds: An Introduction to Curvature},
  author={Lee, J.M.},
  isbn={9780387227269},
  lccn={97014537},
  series={Graduate Texts in Mathematics},
  year={2006},
  publisher={Springer New York}
}

@book{Lee10,
  title={Introduction to Topological Manifolds},
  author={Lee, J.M.},
  isbn={9781441979407},
  lccn={2011287064},
  series={Graduate Texts in Mathematics},
  year={2010},
  publisher={Springer New York}
}

@book{Lee13,
  title={Introduction to Smooth Manifolds},
  author={Lee, J.M.},
  isbn={9780387217529},
  lccn={2002070454},
  series={Graduate Texts in Mathematics},
  year={2013},
  publisher={Springer New York}
}

@misc{Li15,
      title={Higher Groupoid Actions, Bibundles, and Differentiation}, 
      author={Du Li},
      year={2015},
      eprint={1512.04209},
      archivePrefix={arXiv},
      primaryClass={math.DG}
}

@article{LO14,
          volume = {25},
          number = {5},
           month = {7},
          author = {M. Jotz Lean and C. Ortiz},
           title = {Foliated groupoids and infinitesimal ideal systems},
            year = {2014},
         journal = {Indagationes Mathematicae},
           pages = {1019--1053},
}

@article{LLS20,
  title={The universal Lie infinity-algebroid of a singular foliation},
  author={C. Laurent-Gengoux and Sylvain Lavau and T. Strobl},
  journal={Doc. Math.},
  year={2020},
  volume = "25",
  pages = "1571-1652"
}

@book{LM16,
  title={Spin Geometry (PMS-38), Volume 38},
  author={Lawson, H.B. and Michelsohn, M.L.},
  isbn={9781400883912},
  lccn={89032544},
  series={Princeton Mathematical Series},
  year={2016},
  publisher={Princeton University Press}
}

@misc{LMP20,
      title={Modules and representations up to homotopy of Lie $n$-algebroids}, 
      author={Madeleine Jotz Lean and Rajan Amit Mehta and Theocharis Papantonis},
      year={2020},
      eprint={2001.01101},
      archivePrefix={arXiv},
      primaryClass={math.DG}
}

@article{LR21,
     author = {Camille Laurent-Gengoux and Leonid Ryvkin},
     title = {The neighborhood of a singular leaf},
     journal = {Journal de l{\textquoteright}\'Ecole polytechnique {\textemdash} Math\'ematiques},
     pages = {1037--1064},
     publisher = {\'Ecole polytechnique},
     volume = {8},
     year = {2021},
     doi = {10.5802/jep.165},
     language = {en}
}

@misc{LR19,
      title={The holonomy of a singular leaf}, 
      author={Camille Laurent-Gengoux and Leonid Ryvkin},
      year={2019},
      eprint={1912.05286},
      archivePrefix={arXiv},
      primaryClass={math.DG}
}

@misc{Lur06,
      title={Stable $\infty$-Categories}, 
      author={Jacob Lurie},
      year={2006},
      eprint={0608228},
      archivePrefix={arXiv},
      primaryClass={math.CT}
}

@misc{Lur09b,
      title={Derived Algebraic Geometry V: Structured Spaces}, 
      author={Jacob Lurie},
      year={2009},
      eprint={0905.0459},
      archivePrefix={arXiv},
      primaryClass={math.CT}
}

@book{Lur09a,
  title={Higher Topos Theory (AM-170)},
  author={Jacob Lurie},
  isbn={9780691140490},
  lccn={2008038170},
  series={Annals of Mathematics Studies},
  year={2009},
  publisher={Princeton University Press}
}

@book{Lur17, 
 url={https://people.math.harvard.edu/~lurie/papers/HA.pdf}, title={Higher Algebra}, author={Jacob Lurie}, year = {2017}}

@book{Lur18, 
 url={https://www.math.ias.edu/~lurie/papers/SAG-rootfile.pdf}, title={Spectral Algebraic Geometry}, author={Jacob Lurie}, year = {2018}}

@book{Mac87,
  title={Lie Groupoids and Lie Algebroids in Differential Geometry},
  author={Mackenzie, K.},
  isbn={9780521348829},
  lccn={87010287},
  series={Lecture note series / London mathematical society},
  year={1987},
  publisher={Cambridge University Press}
}

@book{Mac05,
  title={General Theory of Lie Groupoids and Lie Algebroids},
  author={Mackenzie, K.},
  isbn={9780521499286},
  lccn={2005296814},
  series={London Mathematical Society Lecture Note Series},
  year={2005},
  publisher={Cambridge University Press}
}

@article{Mac03,
     author = {Mackaay, Marco},
     title = {A note on the holonomy of connections in twisted bundles},
     journal = {Cahiers de Topologie et G\'eom\'etrie Diff\'erentielle Cat\'egoriques},
     pages = {39--62},
     publisher = {Dunod \'editeur, publi\'e avec le concours du CNRS},
     volume = {44},
     number = {1},
     year = {2003},
     zbl = {1067.58003},
     mrnumber = {1961525}
}

@misc{Mac19,
      title={A characteristic map for the holonomy groupoid of a foliation}, 
      author={Lachlan MacDonald},
      year={2019},
      eprint={1910.02167},
      archivePrefix={arXiv},
      primaryClass={math.DG}
}

@book{Mal66,
  title={Ideals of Differentiable Functions},
  author={Malgrange, B.},
  series={Tata Institute of Fundamental Research: Studies in Mathematics},
  year={1966},
  publisher={University Press}
}

@article{Mas02,
author = {Xosé M. Masa},
title = {{Alexander-Spanier cohomology of foliated manifolds}},
volume = {46},
journal = {Illinois Journal of Mathematics},
number = {4},
pages = {979 -- 998},
year = {2002},
doi = {10.1215/ijm/1258138462}
}

@misc{Men19,
      title={Thin homotopy and the holonomy approach to gauge theories}, 
      author={Claudio Meneses},
      year={2019},
      eprint={1904.10822},
      archivePrefix={arXiv},
      primaryClass={math-ph}
}

@misc{Met03,
      title={Topological and Smooth Stacks}, 
      author={David Metzler},
      year={2003},
      eprint={0306176},
      archivePrefix={arXiv},
      primaryClass={math-DG}
}

@book{MM03,
  title={Introduction to Foliations and Lie Groupoids},
  author={Moerdijk, I. and Mrcun, J.},
  isbn={9781139438988},
  year={2003},
  series={Cambridge Studies in Advanced Mathematics},
  publisher={Cambridge University Press}
}

@article{Mos79,
  title={The differentiable space structures of Milnor classifying spaces, simplicial complexes, and geometric realizations},
  author={M. Mostow},
  journal={Journal of Differential Geometry},
  year={1979},
  volume={14},
  pages={255-293}
}

@book{Mor96,
  title={The Seiberg-Witten Equations and Applications to the Topology of Smooth Four-manifolds},
  author={Morgan, J.W.},
  isbn={9780691025971},
  lccn={lc95043748},
  series={Mathematical Notes - Princeton University Press},
  year={1996},
  publisher={Princeton University Press}
}

@misc{Mov08,
      title={A Note on Self-Dual Yang-Mills Theory}, 
      author={M. V. Movshev},
      year={2008},
      eprint={0812.0224},
      archivePrefix={arXiv},
      primaryClass={math-ph}
}

@article{MP97,
author = {Moerdijk, I. and Pronk, Dorette},
year = {1997},
month = {07},
pages = {3-21},
title = {Orbifolds, Sheaves and Groupoids},
volume = {12},
journal = {K-Theory},
doi = {10.1023/A:1007767628271}
}

@article{MQ86,
    author = "Mathai, Varghese and Quillen, Daniel G.",
    title = "{Superconnections, Thom classes and equivariant differential forms}",
    doi = "10.1016/0040-9383(86)90007-8",
    journal = "Topology",
    volume = "25",
    pages = "85--110",
    year = "1986"
}

@book{MR13,
  title={Models for Smooth Infinitesimal Analysis},
  author={Moerdijk, I. and Reyes, G.E.},
  isbn={9781475741438},
  series={SpringerLink : B{\"u}cher},
  year={2013},
  publisher={Springer New York}
}

@article{MR06,
  title={Cohomology of singular Riemannian foliations},
  author={Xos{\'e} M. Masa and A. Rodr{\'i}guez-Fern{\'a}ndez},
  journal={Comptes Rendus Mathematique},
  year={2006},
  volume={342},
  pages={601-604}
}

@book{MS74,
  title={Characteristic Classes},
  author={Milnor, J.W. and Stasheff, J.D.},
  isbn={9780691081229},
  lccn={lc72004050},
  series={Annals of mathematics studies},
  year={1974},
  publisher={Princeton University Press}
}

@book{MS12,
  title={J-holomorphic Curves and Symplectic Topology},
  author={McDuff, D. and Salamon, D.},
  isbn={9780821887462},
  lccn={12016161},
  series={American Mathematical Society colloquium publications},
  year={2012},
  publisher={American Mathematical Society}
}

@book{Mum13,
publisher = {Springer Berlin Heidelberg :},
title = {The Red Book of Varieties and Schemes},
year = {2013},
author = {Mumford, David. author.},
address = {Berlin, Heidelberg :},
isbn = {9783540504979},
}

@inproceedings{Mur06,
  title={Twistor Theory},
  author={M. K. Murray},
  year={2006}
}

@misc{MZ19,
      title={Deformation of singular foliations, 1: Local deformation cohomology}, 
      author={Philippe Monnier and Nguyen Tien Zung},
      year={2019},
      eprint={1904.06888},
      archivePrefix={arXiv},
      primaryClass={math.DG}
}

@article{Nui19,
  title={Homotopical Algebra for Lie Algebroids},
  author={J. Nuiten},
  journal={Applied Categorical Structures},
  year={2019},
  pages={1-42}
}

@inproceedings{Nui18,
  title={Lie algebroids in derived differential topology},
  author={J. Nuiten},
  year={2018}
}

@article{NS11,
title = {Equivariance in higher geometry},
journal = {Advances in Mathematics},
volume = {226},
number = {4},
pages = {3367-3408},
year = {2011},
issn = {0001-8708},
doi = {https://doi.org/10.1016/j.aim.2010.10.016},
author = {Thomas Nikolaus and Christoph Schweigert},
}

@article{NSS12,
author = {Nikolaus, Thomas and Schreiber, Urs and Stevenson, Danny},
year = {2012},
month = {07},
pages = {},
title = {Principal $\infty$-bundles - Presentations},
volume = {10},
journal = {Journal of Homotopy and Related Structures},
doi = {10.1007/s40062-014-0077-4}
}

@misc{Pal03,
      title={$\Bar{\del}$-coherent sheaves over complex manifolds}, 
      author={Nefon Pali},
      year={2003},
      eprint={0305422},
      archivePrefix={arXiv},
      primaryClass={math.AG}
}

@article{Pri13,
title = "Presenting higher stacks as simplicial schemes",
journal = "Advances in Mathematics",
volume = "238",
pages = "184 - 245",
year = "2013",
issn = "0001-8708",
author = "J.P. Pridham"
}

@article{Pri17,
author = {Pridham, J P},
address = {Oxford, UK :},
issn = {1753-8416},
journal = {Journal of topology.},
lccn = {2008252964},
number = {1},
publisher = {Oxford University Press,},
title = {Shifted Poisson and symplectic structures on derived N-stacks},
volume = {10},
year = {2017},
}

@misc{Pri20c,
      title={Non-commutative derived moduli prestacks}, 
      author={J. P. Pridham},
      year={2020},
      eprint={2008.11684},
      archivePrefix={arXiv},
      primaryClass={math.AG}
}

@article{Pri20b,
author = {Pridham, J.P.},
year = {2020},
month = {01},
pages = {106922},
title = {A differential graded model for derived analytic geometry},
volume = {360},
journal = {Advances in Mathematics},
doi = {10.1016/j.aim.2019.106922}
}

@misc{Pri20a,
      title={An outline of shifted Poisson structures and deformation quantisation in derived differential geometry}, 
      author={J. P. Pridham},
      year={2020},
      eprint={1804.07622},
      archivePrefix={arXiv},
      primaryClass={math.DG}
}

@book{PP05,
  title={Quadratic Algebras},
  author={Polishchuk, A. and Positselski, L.},
  isbn={9780821838341},
  lccn={2005048198},
  series={University lecture series},
  url={https://books.google.is/books?id=wcAFCAAAQBAJ},
  year={2005},
  publisher={American Mathematical Society}
}

@article{PW05,
author = {Pietro Polesello, Ingo Waschkies},
journal = {Homology, Homotopy and Applications},
language = {eng},
number = {1},
pages = {109-150},
publisher = {International Press, Somerville},
title = {Higher monodromy.},
volume = {7},
year = {2005},
}

@misc{Qia16,
      title={On the Bott-Chern characteristic classes for coherent sheaves}, 
      author={Hua Qiang},
      year={2016},
      eprint={1611.04238},
      archivePrefix={arXiv},
      primaryClass={math.DG}
}

@book{Qui67,
  title={Homotopical Algebra},
  author={Quillen, D.G.},
  series={Lecture notes in mathematics},
  year={1967},
  publisher={Springer-Verl}
}

@inproceedings{Qui73,
  title={Higher algebraic K-theory: I},
  author={D. Quillen},
  year={1973}
}

@article{Qui85,
title = {Superconnections and the Chern character},
journal = {Topology},
volume = {24},
number = {1},
pages = {89-95},
year = {1985},
issn = {0040-9383},
doi = {https://doi.org/10.1016/0040-9383(85)90047-3},
author = {Daniel Quillen}
}

@misc{Ost08,
      title={Homotopy theory of C*-algebras}, 
      author={Paul Arne Østvær},
      year={2008},
      eprint={0812.0154},
      archivePrefix={arXiv},
      primaryClass={math.AT}
}

@article{Pra85,
  title={How to define the differentiable graph of a singular foliation},
  author={J. Pradines},
  journal={Cahiers de Topologie et G{\'e}om{\'e}trie Diff{\'e}rentielle Cat{\'e}goriques},
  year={1985},
  volume={26},
  pages={339-380}
}

@book{RLL00,
  title={An Introduction to K-Theory for C*-Algebras},
  author={R{\o}rdam, M. and Larsen, F. and Laustsen, N.},
  isbn={9780521789448},
  lccn={00023597},
  series={An Introduction to K-theory for C*-algebras},
  year={2000},
  publisher={Cambridge University Press}
}

@book{Roe96,
  title={Index Theory, Coarse Geometry, and Topology of Manifolds},
  author={John Roe},
  isbn={9780821804131},
  lccn={lc96022058},
  series={CBMS Regional Conference Series},
  year={1996},
  publisher={American Mathematical Society}
}

@book{Roe98,
  title={Elliptic Operators, Topology, and Asymptotic Methods, Second Edition},
  author={Roe, J.},
  isbn={9780582325029},
  lccn={gb98062456},
  series={Chapman \& Hall/CRC Research Notes in Mathematics Series},
  year={1999},
  publisher={Taylor \& Francis}
}

@book{Ros12,
  title={Algebraic K-Theory and Its Applications},
  author={Rosenberg, J.},
  isbn={9781461243144},
  lccn={94008077},
  series={Graduate Texts in Mathematics},
  year={2012},
  publisher={Springer New York}
}

@book{RR06,
  title={An Introduction to Partial Differential Equations},
  author={Renardy, M. and Rogers, R.C.},
  isbn={9780387216874},
  lccn={2003042471},
  series={Texts in Applied Mathematics},
  year={2006},
  publisher={Springer New York}
}

@book{Rud78,
  title={Real and Complex Analysis},
  author={Rudin, W.},
  isbn={9780070995574},
  lccn={73015743},
  series={McGraw-Hill series in higher mathematics},
  year={1978},
  publisher={McGraw-Hill}
}

@book{Rud91,
  title={Functional Analysis},
  author={Rudin, W.},
  isbn={9780070542365},
  lccn={lc90005677},
  series={International series in pure and applied mathematics},
  year={1991},
  publisher={McGraw-Hill}
}

@misc{RV18,
      title={Smooth loop stacks of differentiable stacks and gerbes}, 
      author={David Michael Roberts and Raymond F. Vozzo},
      year={2018},
      eprint={1602.07973},
      archivePrefix={arXiv},
      primaryClass={math.CT}
}

@article{RZ20,
author = "Rogers, Christopher L and Zhu, Chenchang",
fjournal = "Algebraic & Geometric Topology",
journal = "Algebr. Geom. Topol.",
number = "3",
pages = "1127--1219",
publisher = "MSP",
title = "On the homotopy theory for Lie $\infty$-groupoids, with an application to integrating $L_{\infty}$–algebras",
volume = "20",
year = "2020"
}

@article{SS19,
    author = {Severa, Pavel and Siran, Michal},
    title = "{Integration of Differential Graded Manifolds}",
    journal = {International Mathematics Research Notices},
    year = {2019},
    month = {02},
    issn = {1073-7928},
    doi = {10.1093/imrn/rnz004},
    url = {https://doi.org/10.1093/imrn/rnz004},
    note = {rnz004},
    eprint = {https://academic.oup.com/imrn/advance-article-pdf/doi/10.1093/imrn/rnz004/27815769/rnz004.pdf},
}

@article{SSJ09,
author = {Sati, Hisham and Schreiber, Urs and Stasheff, Jim},
year = {2009},
month = {10},
pages = {},
title = {Twisted Differential String and Fivebrane Structures},
volume = {315},
journal = {Communications in Mathematical Physics},
doi = {10.1007/s00220-012-1510-3}
}

@article{Ste74,
  title={Accessibility and foliations with singularities},
  author={P. Stefan},
  journal={Bulletin of the American Mathematical Society},
  year={1974},
  volume={80},
  pages={1142-1145}
}

@article{Sus73,
  title={Orbits of families of vector fields and integrability of distributions},
  author={H. Sussmann},
  journal={Transactions of the American Mathematical Society},
  year={1973},
  volume={180},
  pages={171-188}
}

@misc{SV10,
      title={De Rham Theorem for $L^\infty$ forms and homology on singular spaces}, 
      author={L. Shartser and G. Valette},
      year={2010},
      eprint={1002.4143},
      archivePrefix={arXiv},
      primaryClass={math.MG}
}

@inproceedings{Smi11,
  title={The Higher Riemann-Hilbert Correspondence and Multiholomorphic Mappings},
  author={Aaron M. Smith},
  year={2011}
}

@book{Sco05,
  title={The Wild World of 4-Manifolds},
  author={Scorpan, A.},
  isbn={9780821837498},
  lccn={2004062768},
  year={2005},
  publisher={American Mathematical Society}
}

@article{Sch84,
author = {Claude Schochet},
title = {{Topological methods for $C^{\ast}$-algebras. III. Axiomatic homology.}},
volume = {114},
journal = {Pacific Journal of Mathematics},
number = {2},
publisher = {Pacific Journal of Mathematics, A Non-profit Corporation},
pages = {399 -- 445},
year = {1984},
doi = {pjm/1102708717},
}

@article{Seg72,
  title={Categories and cohomology theories},
  author={G. Segal},
  journal={Topology},
  year={1974},
  volume={13},
  pages={293-312}
}

@article{Seg78,
  title={Classifying spaces related to foliations},
  author={G. Segal},
  journal={Topology},
  year={1978},
  volume={17},
  pages={367-382}
}

@book{SR13,
  title={Basic Algebraic Geometry 1: Varieties in Projective Space},
  author={Shafarevich, I.R. and Reid, M.},
  isbn={9783642379550},
  year={2013},
  publisher={Springer Berlin Heidelberg}
}

@article{SL91,
 ISSN = {0003486X},
 author = {Reyer Sjamaar and Eugene Lerman},
 journal = {Annals of Mathematics},
 number = {2},
 pages = {375--422},
 publisher = {Annals of Mathematics},
 title = {Stratified Symplectic Spaces and Reduction},
 volume = {134},
 year = {1991}
}

@article{Spi08,
  title={Derived smooth manifolds},
  author={David I. Spivak},
  journal={Duke Mathematical Journal},
  year={2008},
  volume={153},
  pages={55-128}
}

@article{Sta98,
author = {Stasheff, Jim},
year = {1998},
month = {01},
pages = {},
title = {The (secret?) homological algebra of the Batalin-Vilkovisky approach},
isbn = {9780821808283},
doi = {10.1090/conm/219/03076}
}

@book{Str07,
  title={Partial Differential Equations: An Introduction},
  author={Strauss, W.A.},
  isbn={9780470054567},
  lccn={2010287546},
  year={2007},
  publisher={Wiley}
}

@article{Sul77,
  title={Infinitesimal computations in topology},
  author={D. Sullivan},
  journal={Publications Math{\'e}matiques de l'Institut des Hautes {\'E}tudes Scientifiques},
  year={1977},
  volume={47},
  pages={269-331}
}

@article{Suw88,
  title={$\mathcal D$-modules associated to complex analytic singular foliations},
  author={T. Suwa},
  journal={Journal of the Faculty of Science. University of Tokyo. Section IA. Mathematics},
  year={1988},
  volume={37},
  pages={297-320}
}

@book{Tao06,
  title={Nonlinear Dispersive Equations: Local and Global Analysis},
  author={Tao, T.},
  isbn={9780821841433},
  lccn={2006042820},
  series={Conference Board of the Mathematical Sciences. Regional conference series in mathematics},
  year={2006},
  publisher={American Mathematical Society}
}

@book{Tao10,
  title={An Epsilon of Room, I: Real Analysis},
  author={Tao, T.},
  isbn={9780821852781},
  lccn={2010036469},
  series={An Epsilon of Room},
  year={2010},
  publisher={American Mathematical Society}
}

@book{Tao11,
  title={An Introduction to Measure Theory},
  author={Tao, T.},
  isbn={9780821869192},
  lccn={2011018926},
  series={Graduate studies in mathematics},
  year={2011},
  publisher={American Mathematical Society}
}

@book{Tao14,
  title={Hilbert's Fifth Problem and Related Topics},
  author={Tao, T.},
  isbn={9781470415648},
  lccn={2014009022},
  series={Graduate Studies in Mathematics},
  year={2014},
  publisher={American Mathematical Society}
}

@book{Tho12,
  title={Elementary Topics in Differential Geometry},
  author={Thorpe, J.A.},
  isbn={9781461261537},
  series={Undergraduate Texts in Mathematics},
  year={2012},
  publisher={Springer New York}
}

@book{Tre04,
  title={Hypo-Analytic Structures (PMS-40), Volume 40: Local Theory (PMS-40)},
  author={Treves, F.},
  isbn={9781400862887},
  series={Princeton Mathematical Series},
  year={2014},
  publisher={Princeton University Press}
}

@article{Tre09,
  title={Exit paths and constructible stacks},
  author={David Treumann},
  journal={Compositio Mathematica},
  year={2009},
  volume={145},
  pages={1504 - 1532}
}

@article{TT76,
  title={A parametrix for $\Bar{\del}$ and Riemann-Roch in \v{C}ech theory},
  author={D. Toledo and Y. L. Tong},
  journal={Topology},
  year={1976},
  volume={15},
  pages={273-301}
}

@article{TT78,
  title={Duality and intersection theory in complex manifolds. I},
  author={D. Toledo and Y. L. Tong},
  journal={Mathematische Annalen},
  year={1978},
  volume={237},
  pages={41-77}
}

@article{TV02,
  title={Homotopical algebraic geometry. I. Topos theory.},
  author={B. Toen and G. Vezzosi},
  journal={Advances in Mathematics},
  year={2002},
  volume={193},
  pages={257-372}
}

@book{TV08,
title = {Homotopical algebraic geometry II : geometric stacks and applications /},
volume = {902},
author = {Toën, Bertrand, 1973-},
address = {Providence, R.I. :},
isbn = {9780821840993 (alk. paper)},
lccn = {2008060003},
publisher = {American Mathematical Society,},
series = {Memoirs of the American Mathematical Society},
year = {2008},
}

@misc{TV20,
      title={Algebraic foliations and derived geometry: the Riemann-Hilbert correspondence}, 
      author={Bertrand Toën and Gabriele Vezzosi},
      year={2020},
      eprint={2001.05450},
      archivePrefix={arXiv},
      primaryClass={math.AG}
}

@MISC{Tyu75,
    author = {I. V. Tyutin},
    title = {Gauge invariance in field theory and statistical physics in operator formalism,” preprint LEBEDEV-75-39},
    year = {1975}
}

@article{UUY,
  title={Homotopical algebra of $C^*$-algebras},
  author={Otgonbayar Uuye},
  journal={Journal of Noncommutative Geometry},
  year={2013},
  volume={7},
  pages={981-1006}
}

@article{Vit15,
author = {Vitagliano, Luca},
year = {2015},
month = {04},
pages = {1450026 [34 pages]},
title = {On the Strong Homotopy Associative Algebra of a Foliation},
volume = {17},
journal = {Communications in Contemporary Mathematics},
doi = {10.1142/S0219199714500266}
}

@book{Voi02,
  title={Hodge Theory and Complex Algebraic Geometry I: Volume 1},
  author={Voisin, C. and Schneps, L.},
  isbn={9781139437691},
  series={Cambridge Studies in Advanced Mathematics},
  year={2002},
  publisher={Cambridge University Press}
}

@article{Vor05,
  title={Higher derived brackets and homotopy algebras},
  author={T. Voronov},
  journal={Journal of Pure and Applied Algebra},
  year={2005},
  volume={202},
  pages={133-153}
}

@article{Vor10,
    author = "Voronov, Theodore Th.",
    editor = "Kielanowski, Piotr and Buchstaber, Victor and Odzijewicz, Anatol and Schlichenmaier, Martin and Voronov, Theodore",
    title = "{$Q$-manifolds and Higher Analogs of Lie Algebroids}",
    doi = "10.1063/1.3527417",
    journal = "AIP Conf. Proc.",
    volume = "1307",
    number = "1",
    pages = "191--202",
    year = "2010"
}

@book{War86,
  title={Foundations of Differentiable Manifolds and Lie Groups},
  author={Warner, F.W.},
  isbn={9780387908946},
  lccn={83012395},
  series={Graduate Texts in Mathematics},
  year={1983},
  publisher={Springer}
}

@misc{Wei18,
      title={The descent of twisted perfect complexes on a space with soft structure sheaf}, 
      author={Zhaoting Wei},
      year={2018},
      eprint={1605.07111},
      archivePrefix={arXiv},
      primaryClass={math.AG}
}

@article{Wei16,
author = {Wei, Zhaoting},
year = {2016},
month = {09},
pages = {},
title = {Twisted complexes on a ringed space as a dg-enhancement of perfect complexes},
volume = {2},
journal = {European Journal of Mathematics},
doi = {10.1007/s40879-016-0102-8}
}

@inbook{Wein07,
	address = {Boston, MA},
	author = {Weinstein, Alan},
	booktitle = {From Geometry to Quantum Mechanics: In Honor of Hideki Omori},
	doi = {10.1007/978-0-8176-4530-4_7},
	editor = {Maeda, Yoshiaki and Ochiai, Takushiro and Michor, Peter and Yoshioka, Akira},
	isbn = {978-0-8176-4530-4},
	pages = {93--109},
	publisher = {Birkh{\"a}user Boston},
	title = {The Integration Problem for Complex Lie Algebroids},
	url = {https://doi.org/10.1007/978-0-8176-4530-4_7},
	year = {2007}}

@book{Wel13,
  title={Differential Analysis on Complex Manifolds},
  author={Wells, R.O.},
  isbn={9781475739466},
  series={Graduate Texts in Mathematics},
  year={2013},
  publisher={Springer New York}
}

@article{Wol89,
     author = {Wolak, Robert},
     title = {Foliated and associated geometric structures on foliated manifolds},
     journal = {Annales de la Facult\'e des sciences de Toulouse : Math\'ematiques},
     pages = {337--360},
     publisher = {Universit\'e Paul Sabatier},
     address = {Toulouse},
     volume = {Ser. 5, 10},
     number = {3},
     year = {1989},
     zbl = {0698.57007},
     mrnumber = {1425491},
     language = {en}
}

@article{WY18,
  title={Bousfield Localization and Algebras over Colored Operads},
  author={David White and Donald Yau},
  journal={Applied Categorical Structures},
  year={2018},
  volume={26},
  pages={153-203}
}

@article{WY19,
author = {White, David and Yau, Donald},
year = {2019},
month = {08},
pages = {},
title = {Homotopical Adjoint Lifting Theorem},
volume = {27},
journal = {Applied Categorical Structures},
doi = {10.1007/s10485-019-09560-2}
}

@article{Yu12,
author = {Yu, Shilin},
year = {2012},
month = {11},
pages = {},
title = {The Dolbeault dga of the formal neighborhood of a diagonal},
volume = {9},
journal = {Journal of Noncommutative Geometry},
doi = {10.4171/JNCG/190}
}

@article{Zhu09,
  title={n-Groupoids and Stacky Groupoids},
  author={Chenchang Zhu},
  journal={International Mathematics Research Notices},
  year={2009},
  volume={2009},
  pages={4087-4141}
}

\end{document}